\documentclass[bezier,11pt]{article}

\usepackage[french]{babel}
\usepackage{inputenc} 
\usepackage[T1]{fontenc}

\usepackage{amsfonts}
\usepackage{amssymb}
\usepackage{amsmath}
\usepackage{mathrsfs}

\newcommand{\dcb}{\begin{array}{lll}}
\newcommand{\dce}{\end{array}}
\newcommand{\ebe}{\begin{enumerate}\setlength{\baselineskip}{13pt}\setlength{\parskip}{0pt}}
\newcommand{\dbe}{\end{enumerate}}

\newcommand{\ibegin}{\begin{itemize}\setlength{\baselineskip}{19pt}\setlength{\parskip}{7pt}}
\newcommand{\iend}{\end{itemize}}
\newcommand{\ok}{\rule{4pt}{6pt}}%

\newtheorem{Theorem}{Theorem}[section]
\newtheorem {Cor}[Theorem]{Corollary}
\newtheorem {definition}[Theorem]{Definition}
\newtheorem {pro}[Theorem]{Proposition}
\newtheorem {Lemma}[Theorem]{Lemma}
\newtheorem {rem}[Theorem]{Remark}
\newtheorem {assumption}[Theorem]{Assumption}
\newtheorem {exercise}[Theorem]{Exercise}

\newcommand {\bd}{\begin{definition}}
\newcommand {\ed}{\end{definition}}
\newcommand {\bl}{\begin{Lemma}}
\newcommand {\el}{\end{Lemma}}
\newcommand {\bcor}{\begin{Cor}}
\newcommand {\ecor}{\end{Cor}}
\newcommand {\brem}{\begin{rem} \rm }
\newcommand {\erem}{\end{rem}}
\newcommand{\bethe}{\begin{Theorem}}
\newcommand{\ethe}{\end{Theorem}}
\newcommand{\bpro}{\begin{pro}}
	\newcommand{\epro}{\end{pro}}
\newcommand {\bassumption}{\begin{assumption}}
\newcommand {\eassumption}{\end{assumption}}
\newcommand {\bex}{\begin{exercise} \rm}
\newcommand {\eex}{\end{exercise}}

\def \ind{1\!\!1\!}
\def\cro#1{\langle #1\rangle}

\newcommand{\stocint}{\centerdot}
\newcommand{\proof}{\textbf{Proof.}}

\newcommand{\transp}{{^\top\!}}

\begin{document}

\begin{center}
\textbf{\LARGE An introduction of the enlargement of filtration}

Lecture\footnote{a working version} at the African Mathematic School in Marrakech

October 19-23, 2015
\end{center}

\

One knows that the problem of filtration enlargement was originated from questions of Ito, Meyer and Williams. The theory of enlargement of filtration started from 1977 with the papers Barlow \cite{barlow} and Jeulin-Yor \cite{JY2} on the honest times and the paper Ito \cite{ito78} on the extension of stochastic integral. The theory of enlargement of filtration had experienced a strong development in the 80s. There were the publications of the books \cite{Jeulin80, JYexample} and the papers Jacod \cite{jacod2} and Yoeurp \cite{yoeurp}. In the recent years, this theory has seen a vigorous revival, because of its applications in mathematic finance, notably in the study of insider problem, of credit risk, of market arbitrage, etc.

Excellent references exist on the theory of enlargement of filtration. See for example \cite{DM3, JYC, jeanblanc, Jeulin80, MY, N, protter}. We can notice an evolution in style of the presentations of the theory. After the monograph \cite{Jeulin80}, the textbooks on the theory favor clearly two particular models (initial and progressive enlargement) and their formulas, leaving aside general results such as the drift operator or the connexion with Girsanov theorem. In compensation, to gain a more global insight into the matter,  the present introduction takes volontier a style in the spirit of \cite[Chapitre II]{Jeulin80}, concentrating on the general aspect of the theory.

By the word "general", we mean propositions which do not found on the specific ways that the enlargements are constructed. Hypothesis $(H)$ or Hypothesis $(H')$ are general propositions, while key lemma or invariance time (cf.\cite{BJR, CS}) are not. There is a general proposition which is not found in the main stream literature, namely the local solution method developed in Song \cite{SongThesis}. The local solution method was originated from the questions how general is the relationship between the enlargement of filtration and Girsanov formula, whether exists a unified proof of the formulas in \cite{jacod2} or in \cite{barlow, JY2, yoeurp}, also of those in the monograph \cite{Jeulin80}. The method is therefore formulated in a setting which covers all models involved in this purpose. As we will show below, this method provides actually a unified proof of the various enlargement filtration formulas. In this sense, the local solution method constitutes a general rule which restructures the theory of enlargement of filtration around a common idea, and provides a canonical way to view and to do the filtration enlargements. We devote the Part II for a detailed presentation of this method.

As mentioned before, two particular models are very well expounded in the literature, i.e. the initial enlargement model under Jacod criterion and the progressive enlargement model with honest times. These two models correspond to the two questions \cite{barlow, ito78, JY2} at the origin of the theory. These two models have simple measurability structures. The explicit knowledge of the measurability structure makes apparent the conditioning relationship between the original and the enlarged filtrations, transforming the filtrations computations into processes computations, leading to explicit formulas. We will not give specific presentations of these two models in this lecture. Instead, we will show how to obtain the formulas with the local solution method in these two models.

The lecture is organized into parts and sections.

\ebe
\item[Part I.]
Problem of enlargement of filtration
\ebe
\item
Preliminary
\ebe
\item
Basic setting
\item
Bichteler-Dellacherie theorem
\item
Bracket processes of semimartingales
\item
Stochastic integrals in different filtrations under different probability measures
\item
Isomorphism of stochastic basis
\item
Filtration transformation under a map
\item
Generalized Girsanov formula
\item
Regular conditional distribution
\item
Prediction process
\dbe

\item
Enlargement of filtration in general
\ebe
\item
Hypothesis $(H)$
\item
Hypothesis $(H')$ and drift operator
\item
The problem of "faux-amis"
\item
Invariance under isomorphism
\item
A panorama of specific models and formulas
\dbe

\dbe
\item[Part II.] Local solution method
\ebe
\item
Local solution method: Introduction
\item
Local solution method: Pieces of semimartingales and their aggregation to form a global semimartingale
\ebe
\item
Assumption and result
\item
Elementary integral
\item
Properties of processes satisfying Assumption \ref{partialcovering}
\dbe

\item
Local solution method: A method to find local solutions
\ebe
\item
The funcitons $h^u, u\geq 0$
\item
Computation of local solutions
\dbe

\item
Classical results as consequence of the local solution method
\ebe
\item
Jacod's criterion
\item
Progressive enlargement of filtration
\item
Honest time
\item
Girsanov theorem and formulas of enlargement of filtration
\dbe

\item
A different proof of \cite[Théorème(3.26)]{Jeulin80}
\item
A last passage time in \cite{JYC}
\ebe
\item
The problem setting
\item
Computing $\tau^u,\alpha^u$ for $0\leq u<1$
\item
The brackets computations
\item
The local solutions before time 1 with their drifts 
\item
The random measure generated by the drifts of the local solutions and its distribution function
\item
Ending remark
\dbe

\item
Expansion with the future infimum process
\ebe
\item
The setting
\item
Computation of $\alpha^u$
\item
Local solutions
\item
Aggregation of local solutions
\item
Computation of $V^+$
\dbe
\dbe

\item[Part III.]
Drift operator and the viability
\ebe
\item
Introduction
\item
Notations and vocabulary
\item
Three fundamental concepts
\ebe
\item
Enlargement of filtrations and hypothesis $(H')$
\item
The martingale representation property
\item
Structure condition
\dbe

\item
Main results
\ebe
\item
Drift multiplier assumption and full viability
\item
The theorems
\dbe
\item
Structure condition decomposed under the martingale representation property
\item
Solution of the continuous structure condition
\item
Solution of the accessible structure condition
\ebe
\item
Equations at the stopping times $T_{n}$
\item
Conditional expectation at predictable stopping times $T_{n}$
\item
An intermediate result
\item
The integrability and conclusion
\dbe

\item
Solution of the totally inaccessible structure condition
\ebe
\item
Equations at the stopping times $S_{n}$
\item
Conditional expectation at stopping times $S_{n}$
\item
Consequence of the totally inaccessible structure condition
\item
Final conclusions
\dbe

\dbe
\dbe

\pagebreak

\begin{center}
\Large \textbf{Part I. 
Problem of enlargement of filtration}
\end{center}
\pagebreak

\section{Preliminary}

The theory of enlargement of filtration is a specific branch of the stochastic calculus. All aspects of the stochastic calculus are involved in this theory. In fact, we can quote Dellacherie-Meyer \cite[Chapitre IV Section 3]{DM}, 
\begin{quote}
Contrairement à ce qu'on pourrait penser, cette théorie ne concerne pas les processus, mais consiste à étudier de manière quelque peu approfondie la structure déterminée, sur un espace mesurable $(\Omega,\mathcal{F})$ ou un espace probabilisé $(\Omega,\mathcal{F},\mathbb{P})$, par la donnée d'une filtration $(\mathcal{F}_t)_{t\geq 0}$.
\end{quote}

That said, the problem of filtration enlargement involves the stochastic calculus in a specific way. A series of results will be presented in this section, which form basic connexions between the stochastic calculus and the problem of filtration enlargement. They are selected notably because of their fundamental role in \cite{Jeulin80,SongThesis}.

\subsection{Basic setting}

We work on probability spaces equipped with filtrations (satisfying the usual conditions), that we call below stochastic basis. We employ the vocabulary of semimartingale calculus as defined in \cite{CS, DM, HWY, jacod, protter}.

In particular, a process $X=(X_t)_{t\geq 0}$  is said integrable, if each $X_t$ is integrable. It is said càdlàg, if for almost all $\omega$, the function $t\rightarrow X_t(\omega)$ is right continuous and has left limit at every point. It is said of class$(D)$, if the family of $X_T$, where $T$ runs over the family of all finite stopping times, is a uniformly integrable family (cf. \cite[p.106]{protter} and also \cite{DM2, HWY}).

The semimartingale space $\mathsf{H}^r, r\geq 1,$ is defined in \cite{emery, Jeulin80, meyer2}. We define $\mathsf{H}^r_{loc}$ as in \cite{jacod, HWY}. A sequence of semimartingales $(X_n)_{n\geq 1}$ is said to converge in $\mathsf{H}^1_{loc}$, if there exists a sequence $(T_m)_{m\geq 1}$ of stopping times tending to the infinity, such that, for every fixed $T_m$, the sequence of stopped processes $(X^{T_m})_{n\geq 1}$ converges in $\mathsf{H}^1$. The subspace of $\mathsf{H}^r$ consisting of local martingales is denoted by $\mathsf{M}^r$.

When $X$ is a special semimartingale, $X$ can be decomposed into $X=M+A$, where $M$ is a local martingale and $A$ is a predictable process with finite variation null at the origin. We will call $A$ the drift of $X$.

The semimartingale calculus depend on the reference probability measure and on the reference filtration. When different probability measures or filtrations are involved in a same computation, expressions such as $(\mathbb{P}, \mathbb{F})$, $\bullet^{\mathbb{P}\cdot\mathbb{F}}$, will be used to indicate the reference probability or filtration. This indication may however be ignored, if no confusion exists. 

The optional (resp. predictable) projection of a process $X$ is denoted by $^{o}\!X$ or $^{\mathbb{P}\cdot\mathbb{F}-o}\!X$ (resp. ${}^{p}\!X$ or $^{\mathbb{P}\cdot\mathbb{F}-p}\!X$). The optional (resp. predictable) dual projection of a process $A$ with finite variation is denoted by $X^{o}$ or $X^{\mathbb{P}\cdot\mathbb{F}-o}$ (resp. $X^{p}$ or $X^{\mathbb{P}\cdot\mathbb{F}-p}$). 

For a function $f$ and a $\sigma$-algebra $\mathcal{T}$, the expression $f\in\mathcal{T}$ will mean that $f$ is $\mathcal{T}$ measurable.  (In)equalities between random variables are almost sure equalities, unless otherwise stated.

\

\subsection{Bichteler-Dellacherie theorem}

To fully appreciate the literature of the filtration enlargement theory, especially the search of enlargement of filtration formulas, we should recall what is a semimartingale. A semimartingale is a "good integrator", as stated in the Bichteler-Dellacherie theorem (cf. \cite{Bich, Dell}). 

\bethe
An adapted cadlag process is a semimartingale, if and only if its corresponding integral, defined on the space of simple predictable processes endowed with the uniform topology, constitutes a continuous map for the convergence in probability.
\ethe

Note that the above continuity in probability for a process $X$ means$$
\lim_{n\rightarrow\infty} H_{n}{{\centerdot}} X = 0,\ \mbox{ in probability,} 
$$
for any sequence $H_{n}, n\in\mathbb{N}$, of simple predictable processes converging uniformly to zero. This property of the semimartingales is at the basis of the local solution method. See section \ref{lsmintro}. Right now, let us illustrate the application of this theorem under a filtration change by the {Stricker theorem}.

\bethe
If a semimartingale is adapted to a sub-filtration, it is semimartingale also in the smaller filtration.
\ethe

\textbf{Proof.}
This last theorem is true because the continuity of the Bichteler-Dellacherie theorem in a bigger filtration implies the same continuity in a smaller filtration. See \cite[Theorem 12.36]{HWY}. \ \ok

\subsection{Bracket processes of semimartingales}

We work on a stochastic basis $(\Omega, \mathcal{A}, \mathbb{F}, \mathbb{P})$ (with $\mathbb{F}$ a filtration, $\mathbb{P}$ a probability measure). The bracket process of a semimartingale, is the limit in probability of quadratic variations.

\bl
For a semimartignale $X$, for $t\geq 0$, $[X,X]_{t}$ is the limit in probability of the quadratic variations $X_{0}^2+\sum_{i=1}^k(X_{t_{i+1}} - X_{t_{i}})^2$, when the mesh of the partition $\tau = (t_{i}, 1\leq i\leq k)$ of $[0,t]$ tends to zero.
\el

(See \cite[Remark of Theorem 9.33]{HWY} for a proof.) Hence, the bracket process is an invariant under filtration change. This invariance plays an important role in the analysis of the problem of filtration enlargement.

The next results concerning bracket processes are borrowed from \cite[Chapitre I]{Jeulin80}.

\bl
Let $Z$ be an optional thin process such that $(\sum_{s\geq 0}Z^2_{s})^{1/2}$ is integrable. For $r\geq 1$, there exists a universal constant $\gamma_{r}$ such that$$
\|(\sum_{s\geq 0}({^{p}}\!Z)^2_{s})^{1/2}\|_{r}
\leq 
\gamma_{r}
\|(\sum_{s\geq 0}Z^2_{s})^{1/2}\|_{r}.
$$
\el

\textbf{Proof.} Note that $^{p}\! Z$ is a thin process too. Consider firstly the case $r=1$. With the convention $Z_{\infty}=0$, $^{p}\! Z_{\infty}=0$, we have$$
\mathbb{E}[(\sum_{s\geq 0}({^{p}}\!Z)^2_{s})^{1/2}]
=
\sup \mathbb{E}[(\sum_{n=0}^N({^{p}}\!Z)^2_{T_{n}})^{1/2}],
$$
where the $\sup$ is taken on the finite increasing families of predictable stopping times. We write$$
\mathbb{E}[(\sum_{n=0}^N({^{p}}\!Z)^2_{T_{n}})^{1/2}]
=
\sup_{H}\mathbb{E}[\sum_{n=0}^N {^{p}}\!Z_{T_{n}} H_{n}],
$$
where the $\sup$ is taken on the random vectors $H$ such that $\sum_{n=0}^N H_{n}^2\leq 1$. Putting $h_{n}=\mathbb{E}[H_{n}|\mathcal{F}_{T_{n}-}]$, $f_{t}=\sum_{n=0}^N Z_{T_{n}}\ind_{\{T_{n}\leq t\}}$ and $g_{t}=\sum_{n=0}^N h_{n}\ind_{\{T_{n}\leq t\}}$, we compute
$$
\dcb
&&\mathbb{E}[\sum_{n=0}^N {^{p}}\!Z_{T_{n}} H_{n}]
=
\mathbb{E}[\sum_{n=0}^N Z_{T_{n}} h_{n}]
=
\mathbb{E}[[f,g]_{\infty}]\\

&\leq&
\mathbb{E}[\sum_{s\geq 0}\frac{|\Delta_{s}f|}{[f,f]_{s}^{1/4}}\ [f,f]_{s}^{1/4} |\Delta_{s}g|]

\leq 
\mathbb{E}[\left(\int_{0}^\infty \frac{1}{[f,f]_{s}^{1/2}}d[f,f]_{s}\right)^{1/2}\ \left(\int_{0}^\infty [f,f]_{s}^{1/2}d[g,g]_{s}\right)^{1/2}]\\
&&\ \mbox{with Kunita-Watanabe inequality,} \\
&\leq&
2\left(\mathbb{E}[ [f,f]_{\infty}^{1/2}]\right)^{1/2}\ \left(\mathbb{E}[\int_{0}^\infty ([g,g]_{\infty}-[g,g]_{s-}) d[f,f]_{s}^{1/2}]\right)^{1/2}\\

&\leq&
2\left(\mathbb{E}[ [f,f]_{\infty}^{1/2}]\right)^{1/2}\ \left(\mathbb{E}[\int_{0}^\infty {^o}\!([g,g]_{\infty}-[g,g]_{-})_{s} d[f,f]_{s}^{1/2}]\right)^{1/2}.
\dce
$$
We now compute the optional projection $ {^o}\!([g,g]_{\infty}-[g,g]_{s-})$. Let $T$ be any stopping time. 
$$
\dcb
&&\mathbb{E}[([g,g]_{\infty}-[g,g]_{T-})\ind_{\{T<\infty\}}]
=
\mathbb{E}[\int_{T}^\infty d[g,g]_{s} + (\Delta_{T}g)^2\ind_{\{T<\infty\}}]\\
&=&
\mathbb{E}[\sum_{n=0}^N h_{n}^2\ind_{\{T<T_{n}\}}]
+
\mathbb{E}[\sum_{n=0}^N h_{n}^2\ind_{\{T=T_{n}\}}]

\leq
\mathbb{E}[\sum_{n=0}^N H_{n}^2]
+
\mathbb{E}[\sum_{n=0}^N H_{n}^2] \leq 2.
\dce
$$
This inequality valid for all stopping time implies that $ {^o}\!([g,g]_{\infty}-[g,g]_{-})\leq 2$, and consequently, $$
\mathbb{E}[\sum_{n=0}^N {^{p}}\!Z_{T_{n}} H_{n}]
\leq
2 \mathbb{E}[(\sum_{n=0}^N Z^2_{T_{n}})^{1/2}].
$$
The lemma is proved for $r=1$.

Next, consider $r>1$. Suppose that $(\sum_{s\leq \cdot}Z^2_{s})^{1/2}\in\mathsf{L}^r$. For any finite increasing predictable stopping times $T_{0}<T_{1}<\ldots<T_{N}$, the process $U=\sum_{n=0}^N (Z_{T_{n}}- {^{p}}\!Z_{T_{n}} )\ind_{[T_{n},\infty)}$ defines a martingale. For any bounded martingale $M$,  the bracket $[U,M]$ satisfies $$
\dcb
&&\mathbb{E}[ [U,M]_{\infty} ]
=
\mathbb{E}[ \sum_{n=0}^N (Z_{T_{n}}- {^{p}}\!Z_{T_{n}} )\Delta_{T_{n}}M ]
=
\mathbb{E}[ \sum_{n=0}^N Z_{T_{n}} \Delta_{T_{n}}M ]
\leq 
\mathbb{E}[ \left(\sum_{n=0}^N Z_{T_{n}}^2\right)^{1/2} [M,M]^{1/2}_{\infty} ]\\

&\leq& 
\| \left(\sum_{n=0}^N Z_{T_{n}}^2\right)^{1/2}\|_{r} \ \|[M,M]^{1/2}_{\infty} \|_{r'} \ \mbox{ with $\frac{1}{r} + \frac{1}{r'} =1$},\\

&\leq&
C \| \left(\sum_{n=0}^N Z_{T_{n}}^2\right)^{1/2}\|_{r} \ \| M_{\infty} \|_{r'}, \ \mbox{ by Burkholder-Davis-Gundy's inequality.}
\dce
$$
By the duality, the above inequality implies $$
\| [U,U]_{\infty}^{1/2} \|_{r}
\leq
C\| U_{\infty} \|_{r}
\leq
C\| (\sum_{n=0}^N Z_{T_{n}}^2)^{1/2}\|_{r}
\leq
C \| (\sum_{s\geq 0} Z_{s}^2)^{1/2}\|_{r}.
$$
(The constant $C$ changes from line to line.) From this last inequality, the lemma follows. \ \ok

\

\bpro\label{bracketcontrol}
Let $X$ be a special semimartingale with canonical decomposition $X=M+V$ (with $M$ a local martingale and $V$ a predictable process with finite variation). For $r\geq 1$, there exists a universal constant $\gamma_{r}$ such that$$
\dcb
\|[V,V]^{1/2}_{\infty}\|_{r}
\leq 
\gamma_{r}
\|[X,X]^{1/2}_{\infty}\|_{r},\\

\|[M,M]^{1/2}_{\infty}\|_{r}
\leq 
(1+\gamma_{r})
\|[X,X]^{1/2}_{\infty}\|_{r}.
\dce
$$
\epro

\textbf{Proof.}
We have ${^{\mathbb{F}\cdot p}}\!(\Delta X)=\Delta V$. The first inequality is, therefore, a direct consequence of the above lemma and of $\sum_{s}(\Delta X_{s})^2\leq [X,X]_{\infty}$. The second inequality results from $$
[M,M]^{1/2}_{\infty}
\leq [X,X]^{1/2}_{\infty}+[V,V]^{1/2}_{\infty}. \ \ok
$$

\brem
A typical application of Proposition \ref{bracketcontrol} in an enlargement $\mathbb{F}\subset \mathbb{G}$ setting is that, if a sequence $(X_{i})_{i\in\mathbb{N}}$ converges in $\mathsf{H}^r(\mathbb{F})$, if the $(X_{i})_{i\in\mathbb{N}}$ are $\mathbb{G}$ semimartingales, then the martingale parts of $(X_{i})_{i\in\mathbb{N}}$ in $\mathbb{G}$ converges in $\mathsf{H}^r(\mathbb{G})$. This means that, if the "$\mathbb{G}$ martingale part" is defined for a generating system in $\mathsf{H}^r(\mathbb{F})$, it will be defined on the whole $\mathsf{H}^r(\mathbb{F})$. Compared with the discussion in \cite{JYfa} on the faux-amis, we see that the $\mathbb{G}$ martingale part of an $\mathbb{F}$ local martingale is not the most important point in the study of the problem of filtration enlargement. Instead, the $\mathbb{G}$ drift part of an $\mathbb{F}$ local martingale will be the center point in that study. See section \ref{noncovering} for more information.
\erem

\

\subsection{Stochastic integrals in different filtrations under different probability measures}

The stochastic integral is another invariant under filtration change. It is an important feature, because, when we deal with the problem of filtration enlargement, a stochastic integral defined in a filtration under some probability measure may be considered in another filtration under another probability measure. More precisely, we have the following result, borrowed from \cite{JS3} (cf. also \cite[Theorem 12.37]{HWY} or \cite[Chapter IV, Theorems 25 and 33]{protter}).

Given a measurable space $(\Omega,\mathcal{A})$, let ${\mathbb{Q}}$ and
${\mathbb{P}}$ be two probability measures on $\mathcal{A}$. Let ${\mathbb{F}}=({\mathcal{F}}_t)_{t\geq
	0}$ and ${\mathbb{G}}=({\mathcal{G}}_t)_{t\geq 0}$ be two
right-continuous filtrations in $\mathcal{A}$. 
Let $X$ be a multi-dimensional càdlàg process and  $0\leq S\leq T$ be two given random variables. Consider the following assumption.

\bassumption  \label{a1} 
\ebe 

\item 
$({\mathbb{Q}},{\mathbb{F}})$ (resp.
$({\mathbb{P}},{\mathbb{G}})$) satisfies the usual conditions.

\item 
$S,T$
are ${\mathbb{F}}$-stopping times and
${\mathbb{G}}$-stopping times.

\item
$X$ is
a $({\mathbb{Q}},{\mathbb{F}})$ semimartingale and a
$({\mathbb{P}},{\mathbb{G}})$ semimartingale.

\item   
The probability ${\mathbb{P}}$ is equivalent to
${\mathbb{Q}}$ on
${\mathcal{F}}_\infty\vee{\mathcal{G}}_\infty$.

\item   
for any
${\mathbb{F}}$-predictable process $J$, $J\ind_{(S,T]}$ is a
${\mathbb{G}}$-predictable process. 
\dbe 
\eassumption

We denote by $\mathcal{I}({\mathbb{Q}},{\mathbb{F}},X)$ the family of all $\mathbb{F}$ predictable processes which are integrable with respect to $X$ in $({\mathbb{Q}},{\mathbb{F}})$. We define similarly the family $\mathcal{I}({\mathbb{P}},{\mathbb{G}},X)$.

\bpro\label{equalintegral}
Suppose Assumption \ref{a1}. Let $J$ be a multi-dimensional ${\mathbb{F}}$-predictable process in $ \mathcal{I}({\mathbb{Q}},{\mathbb{F}},X)\cap \mathcal{I}({\mathbb{P}},{\mathbb{G}},X)$. 
Then, the stochastic integral $J\ind_{(S,T]}\stocint X$ defined in the two senses gives the same process.
\epro

\proof
The lemma is true for elementary $\mathbb{F}$ predictable process $J$. Apply \cite[Theorem 1.4]{HWY} (monotone class theorem), \cite[Remark(ii) of Definition 4.8]{CS} and \cite[Lemma 4.12]{CS}, we see that the lemma is true for bounded $\mathbb{F}$ predictable process. Again by \cite[Remark(ii) of Definition 4.8]{CS} and \cite[Lemma 4.11]{CS}, the proposition is proved. \ok

\

\subsection{Isomorphism of stochastic basis}

Semimartingale property is an invariant under stochastic basis isomorphism. This property is an essential element in the conception of the local solution method. In the sense of this invariance, we say that the problem of filtration enlargement is a problem "in law".

This subsection is based on the paper \cite{BEKSY}.

\bd
With respect to two (complete) probability spaces $(\Omega,\mathcal{A},\mathbb{P})$ and $(\Omega',\mathcal{A}',\mathbb{P}')$, we call almost sure morphism from $(\Omega,\mathcal{A},\mathbb{P})$ to $(\Omega',\mathcal{A}',\mathbb{P}')$, any map from $\mathsf{L}^0(\Omega,\mathcal{A},\mathbb{P})$ to $\mathsf{L}^0(\Omega',\mathcal{A}',\mathbb{P}')$ such that, for any $n\in\mathbb{N}$, for any $U_{1}, \ldots, U_{n}$ in $\mathsf{L}^0(\Omega,\mathcal{A},\mathbb{P})$, for any real Borel functions $f$ on $\mathbb{R}^n$, we have $$
\Psi(f(U_{1}, \ldots, U_{n})) = f(\Psi(U_{1}), \ldots, \Psi(U_{n})),\ \mbox{ almost surely.}
$$ 
\ed

\bpro
The almost sure morphisms $\Psi$ are linear maps, increassing and continuous for the convergence in probability. They preserve the constants. For any $X\in \mathsf{L}^0(\Omega,\mathcal{A},\mathbb{P})$, the law of $\Psi(X)$ is absolutely continuous with respect to that of $X$. In particular, if $X$ is an indicator, it is so for $\Psi(X)$. The composition of two almost sure morphisms is again an almost sure morphism.
\epro

\textbf{Proof.} Apply the morphism property on linear functions to obtain the linearity. Apply then the morphism property on the expression $1+U$ once with linearity, once with the function $f(x)=1+x$. We obtain the preservation of the constants. Based on that, for any Borel binary relation $\mathfrak{R}$, the equality $\ind_{\mathfrak{R}}(X,Y)=1$ implies $\ind_{\mathfrak{R}}(\Psi(X),\Psi(Y))=1$, i.e., the preservation of $\mathfrak{R}$. In particular, $\Psi$ is increasing.

Consider the continuity. We only consider a sequence $(U_{n})_{n\in\mathbb{N}}$ such that $\sum_{n=0}^\infty \mathbb{E}[1\wedge |U_{n}|]<\infty$. This menas that $S=\sum_{n=0}^\infty 1\wedge |U_{n}| \in \mathsf{L}^0$ and $\sum_{n=0}^N 1\wedge |U_{n}|\leq S$, for any integer $N>0$. In other words, $\sum_{n=0}^N 1\wedge |\Psi(U_{n})|\leq \Psi(S)$, for any integer $N>0$, i.e., $(\Psi(U_{n}))_{n\in\mathbb{N}}$ converges in probability $\mathbb{P}'$.

If $\mathbb{P}[U\in B] = 0$ for a Borel set $B\subset \mathbb{R}$, we have $\ind_{B}(U)=0$ so that $\ind_{B}(\Psi(U))=0$, i.e. the law of $\Psi(U)$ is absolutely continuous with respect to that of $X$.\ \ok

\bd
We call imbedding of  $(\Omega,\mathcal{A},\mathbb{P})$ into $(\Omega',\mathcal{A}',\mathbb{P}')$, any almost sure morphism from $(\Omega,\mathcal{A},\mathbb{P})$ into $(\Omega',\mathcal{A}',\mathbb{P}')$, which preserves the probability laws, i.e., for any $X\in \mathsf{L}^0(\Omega,\mathcal{A},\mathbb{P})$, $\Psi(X)$ has the same law as that of $X$. If $\Psi$ is an one-to-one map  from $\mathsf{L}^0(\Omega,\mathcal{A},\mathbb{P})$ to $\mathsf{L}^0(\Omega',\mathcal{A}',\mathbb{P}')$, we call it an isomorphism. 
\ed

\bpro
For any an imbedding $\Psi$ of  $(\Omega,\mathcal{A},\mathbb{P})$ into $(\Omega',\mathcal{A}',\mathbb{P}')$, its restriction on any $\mathsf{L}^p(\Omega,\mathcal{A},\mathbb{P}), p\geq 1$, is an isometry. For any finite or denumerable family of random variables $(U_{i})_{i\in \mathtt{D}}$ on $\Omega$, $(U_{i})_{i\in \mathtt{D}}$ and $(\Psi(U_{i}))_{i\in \mathtt{D}}$ have the same law. For any sub-$\sigma$-algebra $\mathcal{B}$ in $\mathcal{A}$, there exists a unique sub-$\sigma$-algebra $\mathcal{B}'$ in $\mathcal{A}'$ such that $\Psi$ defines a bijection between $\mathsf{L}^0(\Omega,\mathcal{B},\mathbb{P})$ and $\mathsf{L}^0(\Omega',\mathcal{B}',\mathbb{P}')$, and, for two random variables $Y,Z$ on $\Omega$, the equality $Y=\mathbb{E}[Z|\mathcal{B}]$ holds, if and only if $\Psi(Y)=\mathbb{E}'[\Psi(Z)|\mathcal{B}']$. We denote $\mathcal{B}'$ by $\Psi(\mathcal{B})$.
\epro

\textbf{Proof.} Consider the statement about $(U_{i})_{i\in \mathtt{D}}$. If $\mathtt{D}$ is finite, the result follows from the definition. With the $\lambda$-$\pi$ lemma (cf. \cite[Theorem 1.2]{HWY}), the result extends to denumerable set $\mathtt{D}$. 

Consider $\mathcal{B}$. Let $\mathcal{B}'$ be the $\sigma$-algebra generated by the representatives of $\Psi(U)$ for $U\in \mathcal{B}$. Because of the law preservation, the map $\Psi$ restricted on $\mathsf{L}^0(\mathcal{B})$ is an isometric into its image $\Psi(\mathsf{L}^0(\mathcal{B}))$. The two spaces $\mathsf{L}^0(\mathcal{B})$ and $\Psi(\mathsf{L}^0(\mathcal{B}))$ are $F$-spaces. Hence, $$
\{U'\in \mathsf{L}^0(\mathcal{B}'): \exists U, U'=\Psi(U)\}
$$
is a functional montone class in the sense of \cite{RW}. By the monotone class theorem, $\Psi(\mathsf{L}^0(\mathcal{B}))\supset \mathsf{L}^0(\mathcal{B}')$. But the inverse inclusion is clearly true. We have actually an equality $\Psi(\mathsf{L}^0(\mathcal{B}))= \mathsf{L}^0(\mathcal{B}')$. The $\sigma$-algebra satisfying this equality mush be unique. \ \ok

Consider two stochastic basis $(\Omega,\mathbb{F},\mathbb{P})$ and $(\Omega',\mathbb{F}',\mathbb{P}')$, where $\mathbb{F}=(\mathcal{F}_{t}, t\geq 0)$ and $\mathbb{F}'=(\mathcal{F}'_{t}, t\geq 0)$ are filtrations satisfying respectively the usual conditions with respect to themself. 

\bd
We call an isomorphism of the stochastic basis $(\Omega,\mathbb{F},\mathbb{P})$  onto the stochastic basis $(\Omega',\mathbb{F}',\mathbb{P}')$, any imbedding $\Psi$ from $(\Omega,\mathcal{F}_{\infty},\mathbb{P})$ into $(\Omega',\mathcal{F}'_{\infty},\mathbb{P}')$, which satisfies $\Psi(\mathcal{F}_{t})=\mathcal{F}'_{t}, t\in\mathbb{R}_{+}$.
\ed

\bpro\label{isota}
Let $\Psi$ be an isomorphism of the stochastic basis $(\Omega,\mathbb{F},\mathbb{P})$ onto the stochastic basis $(\Omega',\mathbb{F}',\mathbb{P}')$. A random variable $T$ on $\Omega$ is an $\mathbb{F}$ (predictable, totally inaccessible) stopping time, if and only if $\Psi(T)$ is an $\mathbb{F}'$ (predictable, totally inaccessible) stopping time.
\epro

\textbf{Proof.}
To prove the principal result, we note that a random variable is a stopping time, if and only if it is the decreasing limit of a sequence of discrete stopping times. For a discrete random variable $T$, it is an $\mathbb{F}$ stopping time, if and only if $
\ind_{\{T=a\}} = \mathbb{E}[\ind_{\{T=a\}} | \mathcal{F}_{a}],
$
for all $a\geq 0$, if and only if $
\ind_{\{\Psi(T)=a\}} = \mathbb{E}[\ind_{\{\Psi(T)=a\}} | \mathcal{F}'_{a}],
$
so that $\Psi(T)$ is an $\mathbb{F}'$ stopping time.

To prove the predictable stopping time property, we apply \cite[Theorem 4.34]{HWY} of foretellable times. The property of totally inaccessible time can be directly deduced from the definition. \ \ok

\bpro\label{sbi2}
Let $\Psi$ be an isomorphism of the stochastic basis $(\Omega,\mathbb{F},\mathbb{P})$ onto the stochastic basis $(\Omega',\mathbb{F}',\mathbb{P}')$. Let $X=(X_{t},t\geq 0)$ be a process on $\Omega$. Then, $X$ has a modification which is $\mathbb{F}$ adapted $\mathbb{P}$ almost cadlag, if and only if $\Psi(X)=(\Psi(X_{t}),t\geq 0)$ can be chosen to form an $\mathbb{F}'$ adapted $\mathbb{P}'$ almost cadlag process. $X$ has a modification of finite variation, if and only if $\Psi(X)=(\Psi(X_{t}),t\geq 0)$ can be chosen to form a process with finite variation. $X$ is $\mathbb{F}$ (cadlag) predictable (resp. optional), if and only $\Psi(X)$ can be chosen to be $\mathbb{F}'$ (cadlag) predictable (resp. optional). $X$ is a martingale in $(\Omega,\mathcal{A},\mathbb{P}, \mathbb{F})$, if and only if $\Psi(X)$ (for its cadlag version) is a martingale in $(\Omega',\mathcal{A}',\mathbb{P}', \mathbb{F}')$.
\epro

Notice that, for random variable $U$ on $\Omega$, $\Psi(U)$ is \textit{à priori} an equivalent class in $\mathsf{L}^0(\Omega',\mathcal{A}',\mathbb{P}')$. But in computations, by abuse of notation, $\Psi(U)$ may be used to denote a particular member in the equivalent class $\Psi(U)$. 

\textbf{Proof.} We only consider the path regularity of $\Psi(X)$. Suppose then that $X$ is cadlag and bounded (which is not too restrictive). The restriction of $X$ on the rational set $(X_{t},t\in \mathbb{Q}_{+})$ has the same law as $(\Psi(X_{t}),t\in \mathbb{Q}_{+})$. By \cite[Chapitre III, $n^\circ$18]{DM}, to be the restriction on $\mathbb{Q}_{+}$ of a cadlag map defines a measurable set on $\mathbb{R}^{\mathbb{Q}_{+}}$. Hence, like $(X_{t},t\in \mathbb{Q}_{+})$, the random map $(\Psi(X_{t}),t\in \mathbb{Q}_{+})$ is the restriction on $\mathbb{Q}_{+}$ of a cadlag map, $\mathbb{P}'$ almost everywhere. Hence, $\liminf_{s\downarrow t}\Psi(X_{s}), t\geq 0,$ define a cadlag process (up to indistinguishability). Notice that, by the continuity of $\Psi$ under convergence in probability, we see that $\liminf_{s\downarrow t}\Psi(X_{s})$ is a member of the equivalent class $\Psi(X_{t})$. We prove thus that $(\liminf_{s\downarrow t}\Psi(X_{s}), t\geq 0)$ is a cadlag representative of $\Psi(X)$.

The reciprocal statement can be proved similarly. \ok

\brem
We introduce the stochastic basis isomorphism, because of Proposition \ref{sbi2} (cf. also Proposition \ref{sbi} below) about the martingale property under the isomorphism. Actually, the problem of filtration enlargement may find more easily a solution on an isomorphic probability space, than on the original one. 
\erem

\

\subsection{Filtration transformations under a map}

A common way to define isomorphisms of stochastic basis is to do it with a maps. In such a situation the isomorphic properties can be studied in a very transparent way, with or without completion of the filtrations. This subsection is based on \cite{Song-diff}. See section \ref{lsmfind} for an application.

\subsubsection{Inverse images of $\sigma$-algebras}\label{tribulimitinverse}

We begin with some facts on the $\sigma$-algebras. Let $E$ a set and $\mathcal{T}$ a $\sigma$-algebra on $E$. For $A\subset E$, we define$$
A\cap\mathcal{T}
=\{B\subset A: \exists C\in\mathcal{T}, B=A\cap C\}.
$$
The family $A\cap\mathcal{T}$ can be used as a family of subsets in $E$, or as a $\sigma$-algebra on $A$. Notice that, if $f$ and $g$ are two maps on $E$ such that $f=g$ on $A$, we have $A\cap\sigma(f)=A\cap\sigma(g)$. We need another fact on the $\sigma$-algebras. Let $\mathcal{T}_n, n\in\mathbb{N}^*$, be a decreasing sequence of $\sigma$-algebras on $E$. Let $F$ be another space and $\eta$ be a map from $F$ into $E$. We have$$
\eta^{-1}(\cap_{n\in\mathbb{N}^*}\mathcal{T}_n)
=\cap_{n\in\mathbb{N}^*}\eta^{-1}(\mathcal{T}_n).
$$ 
Obviously the right hand side term contains the left hand side term. Let $B$ be an element in $\cap_{n\in\mathbb{N}^*}\eta^{-1}(\mathcal{T}_n)$. For any $n\in\mathbb{N}^*$, there exists a $C_n\in \mathcal{T}_n$ such that $B=\eta^{-1}(C_n)$. Let $C=\limsup_{n\rightarrow \infty}C_n\in\cap_{n\in\mathbb{N}^*}\mathcal{T}_n$. We check that $$
\eta^{-1}(C)=\eta^{-1}(\cap_{n\geq 1}\cup_{m\geq n}C_m)
=\cap_{n\geq 1}\cup_{m\geq n}\eta^{-1}(C_m)=B.
$$
This proves the above identity. Consider a subset $D$ of $F$. Applying the previous result with the identity map from $D$ into $F$, we can state the above identity in a general form:

\bl\label{intersection}
We have the identity$$
D\cap\eta^{-1}(\cap_{n\in\mathbb{N}^*}\mathcal{T}_n)
=D\cap(\cap_{n\in\mathbb{N}^*}\eta^{-1}(\mathcal{T}_n))
=\cap_{n\in\mathbb{N}^*}(D\cap\eta^{-1}(\mathcal{T}_n)).
$$ 
\el

\subsubsection{$\sigma$-algebra completion}

We need the following lemma which describes the completion of a $\sigma$-algebra in term of the $\sigma$-algebra itself. This description will be useful when we compare the completion of a $\sigma$-algebra on the original space with the completion of a $\sigma$-algebra on an auxiliary space (cf. Proposition \ref{optionalinclusion}).

\bl\label{negtribu}
Let $\mathcal{T}_1, \mathcal{T}_2$ be two $\sigma$-algebras on some common space $\Omega$. Let $\nu$ be a probability measure defined on the two $\sigma$-algebras. Let $\mathcal{N}$ be the family of the $(\nu,\mathcal{T}_2)$ negligible sets. Then,$$
\mathcal{T}_1\vee\sigma(\mathcal{N})
=\{X\subset \Omega: \exists B\in\mathcal{T}_1, A\in\mathcal{T}_2, \nu[A]=1, X\cap A=B\cap A\}.
$$
\el

\textbf{Proof.} Denote the right hand side term of the above formula by $\mathcal{J}$. Then, $\Omega\in\mathcal{J}$. If $X\in\mathcal{J}$, let $B\in\mathcal{T}_1$ and $A\in\mathcal{T}_2$ such that $\nu[A]=1$ and $X\cap A=B\cap A$. Then, $X^c\cap A=B^c\cap A$, which means $X^c\in\mathcal{J}$. If $X_n\in\mathcal{J}$ for $n\in\mathbb{N}^*$, let $B_n\in\mathcal{T}_1$ and $A_n\in\mathcal{T}_2$ such that $\nu[A_n]=1$ and $X\cap A_n=B_n\cap A_n$. Set $A=\cap_{n\in\mathbb{N}^*}A_n$. Then, $\nu[A]=1-\nu[A^c]=1,$ while $$
(\cup_{n\in\mathbb{N}^*}X_n)\cap A
=\cup_{n\in\mathbb{N}^*}(X_n\cap A)
=\cup_{n\in\mathbb{N}^*}(B_n\cap A)
=(\cup_{n\in\mathbb{N}^*}B_n)\cap A),
$$
i.e. $\cup_{n\in\mathbb{N}^*}X_n\in\mathcal{J}$. The family $\mathcal{J}$ is a $\sigma$-algebra.

The $\sigma$-algebra $\mathcal{J}$ contains clearly $\mathcal{T}_1$. It also contains $\mathcal{N}$. Actually, for any $X\in\mathcal{N}$, there exists a $C\in\mathcal{T}_2$ such that $\nu[C]=0$ and $X\subset C$. Let $A=C^c$, we have $\nu[A]=1$ and $X\cap A=\emptyset=\emptyset\cap A$. This means that $X\in\mathcal{J}$. 

On the other hand, for any $X\in\mathcal{J}$, let $B\in\mathcal{T}_1$ and $A\in\mathcal{T}_2$ such that $\nu[A]=1$ and $X\cap A=B\cap A$. Then,$$
X=X\cap A+X\cap A^c
=B\cap A+X\cap A^c\in\mathcal{T}_1\vee\sigma(\mathcal{N}).
$$
This means that $\mathcal{J}\subset \mathcal{T}_1\vee\sigma(\mathcal{N})$. \ok

\subsubsection{Filtration completion and the inverse image}

In this subsection we consider two measurable spaces $(\Omega,\mathcal{A})$ and $(\Omega',\mathcal{A}')$ and a measurable map $\phi$ from $\Omega$ into $\Omega'$. Let $\mathbb{P}$ be a probability measure on $\mathcal{A}$ and $\check{\mathbb{G}}^0=(\check{\mathcal{G}}^0_{t}, t\geq 0)$ be a right continuous filtration in $\mathcal{A}'$ (no necessarily completed). We define $\check{\mathbb{P}}=\phi^{-1}(\mathbb{P})$ on $\mathcal{A}'$ and ${\mathcal{G}}^0_{t}=\phi^{-1}(\check{\mathcal{G}}^0_{t}), t\geq 0$ forming a filtration ${\mathbb{G}}^0$ in $\mathcal{A}$ which is right continuous, because of subsection \ref{tribulimitinverse}. We denote by $\mathbb{G}$ the completion of $\mathbb{G}^0$ under $\mathbb{P}$ with respect to $\mathcal{G}_{\infty}$, and by $\check{\mathbb{G}}$ the completion of $\check{\mathbb{G}}^0$ under $\check{\mathbb{P}}$ with respect to $\check{\mathcal{G}}^0_{\infty}$. In this setting, the map $\phi^{-1}$ defines an imbedding from $(\Omega,\mathcal{A},\mathbb{P}, \mathbb{G})$ to $(\Omega',\mathcal{A}',\check{\mathbb{P}}, \check{\mathbb{G}})$,

\bl\label{isoPsi}
Let $t\in[0,\infty]$. We have $\phi^{-1}(\check{\mathcal{G}}_{t})\subset {\mathcal{G}}_{t}$. For any $\mathcal{G}_{t}$ measurable function $f$, there exists a  $\check{\mathcal{G}}^0_{t}$ measurable function $g$ such that $f = g(\phi)$, $\mathbb{P}$ almost surely. We define a map $\Psi$ which maps the class of $f$ in $\mathsf{L}^0(\Omega,\mathcal{G}_{\infty},\mathbb{P})$ to the class of $g$ in $\mathsf{L}^0(\Omega',\check{\mathcal{G}}_{\infty},\check{\mathbb{P}})$. Then, $\Psi$ is an isomorphism from the stochastic basis $(\Omega,\mathbb{G},\mathbb{P})$ onto $(\Omega',\check{\mathbb{G}},\check{\mathbb{P}})$.
\el

\textbf{Proof.} We apply Lemma \ref{negtribu}. For $B'\in \check{\mathcal{G}}_{t}$, there exist $A'\in \check{\mathcal{G}}_{\infty}$ of probability 1 and $C'\in \check{\mathcal{G}}^0_{t}$ such that $B'\cap A' = C'\cap A'$, which implies $$
\phi^{-1}(B')\cap \phi^{-1}(A')=\phi^{-1}(B'\cap A')=\phi^{-1}(C'\cap A')=\phi^{-1}(C')\cap \phi^{-1}(A'),
$$ 
i.e., $\phi^{-1}(B')\in\mathcal{G}_{t}$. This proves the first assertion of the lemma.

For the second assertion, it is enough to prove that, for any element $B\in \mathcal{G}_{t}$, there exists an element $B'\in \check{\mathcal{G}}^0_{t}$ such that $\mathbb{P}[B\Delta \phi^{-1}(B')]=0$. Let then $B\in {\mathcal{G}}_{t}$. There exist $A\in {\mathcal{G}}_{\infty}$ of probability 1 and $C\in {\mathcal{G}}^0_{t}$ such that $B\cap A = C\cap A$. But $C=\phi^{-1}(C')$ for some $C'\in  \check{\mathcal{G}}^0_{t}$. We have $$
\mathbb{P}[B\Delta \phi^{-1}(C')] = \mathbb{P}[B\Delta C]\leq  \mathbb{P}[A^c]=0,
$$
proving the second assertion. For the last assertion, we note that $\Psi$ is well-defined, because, if $f = g(\phi)=h(\phi)$ $\mathbb{P}$, $g=h$ $\check{\mathbb{P}}$ almost surely. With this in mind, the last assertion can be checked with the second assertion. \ok

\brem
Note that we can not state an equality in the first assertion, because a $(\mathbb{P},\mathcal{G}_{\infty})$ negligible set may not be the inverse image of a set in $\Omega'$. 
\erem

Now we compare the optional and predictable $\sigma$-algebras on the two spaces.

\bl\label{optionalinclusion}
Let $\Phi$ to be the map from $\mathbb{R}_+\times\Omega$ into $\mathbb{R}_+\times(\Omega\times[0,\infty])$ with $\Phi(t,\omega)=(t,\phi(\omega))$. We have the following relationships:$$
\mathcal{O}(\mathbb{G}^0)
=
\Phi^{-1}(\mathcal{O}(\check{\mathbb{G}}^0))
\subset 
\Phi^{-1}(\mathcal{O}(\check{\mathbb{G}}))
\subset \mathcal{O}(\mathbb{G}).
$$
Similar results hold for the predictable $\sigma$-algebras.
\el

\textbf{Proof.}
Notice that $
\phi^{-1}(\check{\mathcal{G}}^0_t)
=
\mathcal{G}^0_t,
$
for every $t\in\mathbb{R}_+$, and by Lemma \ref{isoPsi},
$
\phi^{-1}(\check{\mathcal{G}}_t)
\subset
\mathcal{G}_t.
$
Then, for any (everywhere) càdlàg $\check{\mathbb{G}}^0$ (respectively $\check{\mathbb{G}}$) adapted process $X$, $X\circ\Phi$ is a càdlàg $\mathbb{G}^0$ (respectively $\mathbb{G}$) adapted process. We have therefore$$
\Phi^{-1}(\mathcal{O}(\check{\mathbb{G}}^0))
\subset \mathcal{O}(\mathbb{G}^0),\
\Phi^{-1}(\mathcal{O}(\check{\mathbb{G}}))
\subset \mathcal{O}(\mathbb{G}).
$$

For any $\mathbb{G}^0$ stopping time $S$, let $X=\ind_{[S,\infty)}$. For $t\in\mathbb{Q}_+$ let $f_t$ be a $\check{\mathcal{G}}^0_t$ measurable function bounded by 0 and 1, such that $X_t=f_t(\phi)$ (cf. \cite[Theorem 1.5]{HWY}). Let $$
\overline{S} (\omega,u) 
=\inf\{s\in\mathbb{Q}_+: f_s(\omega,u)=1\}.
$$
$\overline{S}$ is a $\check{\mathbb{G}}^0$ stopping time and $S=\overline{S}(\phi)$, i.e., $$
X=(\ind_{[\overline{S},\infty)})\circ\Phi\in \Phi^{-1}(\mathcal{O}(\check{\mathbb{G}}^0)).
$$
By \cite[Theorem 3.17]{HWY}, we conclude that $$
\mathcal{O}(\mathbb{G}^0)
\subset \Phi^{-1}(\mathcal{O}(\check{\mathbb{G}}^0))
\subset \Phi^{-1}(\mathcal{O}(\check{\mathbb{G}})). \  \ok
$$

\

\subsection{Generalized Girsanov formula}

Notice that the problem of filtration enlargement is a problem in law so that a filtration change without probability is meaningless. We look at therefore the implication of the probability in the problem, especially the consequence of a probability change. Clearly, the problem of filtration enlargement is invariant under equivalent probability change thanks to the Girsanov formula. But, in fact, we can do much more with the so-called generalized Girsanov formula. See section \ref{lsmintro}.

\brem
We know that a probability space change can be equivalent to a probability measure change. It is the case when one make an isomorphism between to probability spaces with the same "physical" basis $(\Omega,\mathcal{A})$ and with two different probability measures. We ask if such an equivalence exists between a filtration enlargement and a probability measure change. In fact, it is not really the case.
\erem

Consider a filtered measurable space $(\Omega,\mathcal{A},\mathbb{G})$. Let $\mu$ and $\nu$ be two probability measures on $\mathcal{A}$ (no usual condition supposed on the filtration). Let $\lambda$ denote the probability measure $\frac{\mu+\nu}{2}$ and let $h$ be the $(\lambda, \mathbb{G})$ uniformly integrable martingale representing the density processes of $\mathbb{P}$ with respect to $\lambda$ along the filtration $\mathbb{G}$. The process $h$ is supposed to be everywhere right continuous and $\lambda$ almost surely cadlag. We define$$
\dcb
\tau = \inf\{t\geq 0: h_{t} = 2\},\

\alpha = \frac{h}{2-h}\ind_{[0,\tau)}.
\dce
$$

\bethe\label{basicformulas}
Let $\rho$ to be a $\mathbb{G}$-stopping time with $\rho\leq \tau$. We have, for any
$\mathbb{G}$-stopping time
$\kappa$, for any $A \in
{\mathcal{G}}_{\kappa}$,$$
\mathbb{E}_{\mu}[\ind_A\ind_{\{\kappa < \rho\}}] = \mathbb{E}_{\nu}[\ind_A\ind_{\{\kappa < \rho\}}
\alpha_{\kappa}].
$$
Consequently, $\alpha\ind_{[0,\rho)}$ is a $(\nu,{\mathbb{G}})$ supermartingale. Moreover, $\alpha>0$ on $[0,\rho)$ under the probability $\mu$. For any non negative $\mathbb{G}$ predictable process $H$, $$
\mathbb{E}_{\mu}[H_\kappa\ind_{\{\kappa < \rho\}}] 
= \mathbb{E}_{\nu}[H_\kappa\ind_{\{\kappa < \rho\}}
\alpha_{\kappa}].
$$
Suppose in addition that $(\alpha\ind_{[0,\rho)})^\kappa$ is in class $(D)$ under $\nu$. Let $V$ be the non decreasing $\mathbb{G}$ predictable process associated with the supermartingale $\alpha\ind_{[0,\rho)}$. For any non negative $\mathbb{G}$ predictable process $H$, we have$$
\mathbb{E}_{\mu}[H_{\rho}\ind_{\{0<\rho\leq \kappa\}}] =
\mathbb{E}_{\nu}[\int_{0}^{\kappa}
H_s dV_s].
$$
Let $B$ be a $\mathbb{G}$ predictable process with bounded variation. We have$$
\mathbb{E}_{\mu}[B_{\kappa\wedge \rho}-B_0] =
\mathbb{E}_{\nu}[\int_{0}^{\kappa}
\alpha_{s-} dB_s].
$$
Consequently, $\mathbb{E}_{\mu}[\int_0^{\kappa}\ind_{\{0<s\leq \rho\}}\ind_{\{\alpha_{s-}=0\}}dB_s]=0$. Let $C$ be a $\mathbb{G}$ optional process having bounded variation on the open random interval $(0,\kappa\wedge \rho)$. For any bounded $\mathbb{G}$ predictable process $H$, We have$$
\mathbb{E}_{\mu}[\int_0^\infty H_s\ind_{(0,\kappa\wedge \rho)}(s) dC_s] =
\mathbb{E}_{\nu}[\int_{0}^\infty H_s\ind_{(0,\kappa\wedge \rho)}(s) \alpha_sdC_s].
$$
In particular, $$
(\ind_{(0,\kappa\wedge \rho)}\centerdot C)^{\mu\cdot\mathbb{G}-p}
=
\frac{1}{\alpha_-}\centerdot(\ind_{(0,\kappa\wedge \rho)}\alpha\centerdot C)^{\nu\cdot\mathbb{G}-p}.
$$
\ethe

\textbf{Proof.} For any positive $\mathcal{G}_{\kappa}$ measurable function $f$,$$
2\mathbb{E}_{\mu}[f \ind_{\{\kappa<\rho\}}]=\mathbb{E}_{\mu}[f \ind_{\{\kappa<\rho\}}h_{\kappa}]+\mathbb{E}_{\nu}[f \ind_{\{\kappa<\rho\}}h_{\kappa}],
$$
or equivalently
$$
\mathbb{E}_{\mu}[f (2-h_{\kappa})\ind_{\{\kappa<\rho\}}]=\mathbb{E}_{\nu}[f (2-h_{\kappa})\alpha_{\kappa}\ind_{\{\kappa<\rho\}}].
$$
This last identity is an equivalent form of the first formula of the lemma. To see the supermartingale property of $\alpha\ind_{[0,\rho)}$, it is enough to notice that, for $A\in\mathcal{G}_{s}, 0\leq s\leq t$, $$
\dcb
\mathbb{E}_{\nu}[\ind_A\alpha_t\ind_{\{t<\rho\}}]
=\mathbb{E}_{\mu}[\ind_A\ind_{\{t<\rho\}}]
\leq\mathbb{E}_{\mu}[\ind_A\ind_{\{s<\rho\}}]
=\mathbb{E}_{\nu}[\ind_A\alpha_s\ind_{\{s<\rho\}}].
\dce
$$
If $A=\{\upsilon\leq s\}$ with $\upsilon= \inf\{u\geq 0: h_{u} = 0\}$, the above expectation is null, which prove the positivity of $\alpha$ on $[0,\rho)$ under $\mu$.

The second formula of the lemma is a direct consequence of the first one. To prove the third formula of the theorem, we need only to check it on the processes $H$ of the form $\ind_A\ind_{(a,\infty)}$ with $0\leq a<\infty$ and $A\in\mathcal{G}_{a}$. We have 
$$
\dcb
&&\mathbb{E}_{\mu}[H_\rho\ind_{\{0<\rho\leq \kappa\}}]
=\mathbb{E}_{\mu}[\ind_A\ind_{\{a<\rho\}}\ind_{\{0<\rho\leq \kappa\}}]
=\mathbb{E}_{\mu}[\ind_A\ind_{\{a<\rho\}}\ind_{\{0<\rho\}}]
-\mathbb{E}_{\mu}[\ind_A\ind_{\{a<\rho\}}\ind_{\{\kappa<\rho\}}]\\
&=&\mathbb{E}_{\nu}[\ind_A\alpha_{a}\ind_{\{a<\rho\}}\ind_{\{0<\rho\}}]
-\mathbb{E}_{\nu}[\ind_A\alpha_{a\vee \kappa}\ind_{\{a<\rho\}}\ind_{\{\kappa<\rho\}}]

=\mathbb{E}_{\nu}[\ind_A\left(\alpha_{a}\ind_{\{a<\rho\}}-\alpha_{a\vee \kappa}\ind_{\{a\vee \kappa<\rho\}}\right)]
\\

&=&\mathbb{E}_{\nu}[\ind_A\left(V_ {a\vee \kappa}-V_ {a}\right)]
\ \mbox{ because $(\alpha\ind_{[0,\rho)})^\kappa$ is in class$(D)$}\\

&=&\mathbb{E}_{\nu}[\ind_A\int_0^\kappa \ind_{(a,\infty)}(s)dV_s]

=\mathbb{E}_{\nu}[\int_0^\kappa H_sdV_s].
\dce
$$
For the fourth formula, we write$$
\dcb
&&\mathbb{E}_{\mu}[B_{\kappa\wedge\rho}-B_0] 
=\mathbb{E}_{\mu}[(B_{\kappa}-B_0)\ind_{\{\kappa<\rho\}}]
+\mathbb{E}_{\mu}[(B_{\rho}-B_0)\ind_{\{\kappa\geq \rho\}}]\\
&=&\mathbb{E}_{\nu}[(B_{\kappa}-B_0)\alpha_\kappa\ind_{\{\kappa<\rho\}}]
+\mathbb{E}_{\mu}[(B_{\rho}-B_0)\ind_{\{0<\rho\leq \kappa\}}]\\
&=&\mathbb{E}_{\nu}[(B_{\kappa}-B_0)\alpha_\kappa\ind_{\{\kappa<\rho\}}]
+\mathbb{E}_{\nu}[\int_0^\kappa(B_{s}-B_u)dV_s]
=\mathbb{E}_{\nu}[\int_{0}^{\kappa}
\alpha_{s-} dB_s].
\dce
$$
Finally, for the last formulas, let $K=H\ind_{(0,\kappa\wedge \rho)}\centerdot C$. We note that $K_\kappa=K_{\kappa\wedge \rho-}$. Applying the preceding formulas, we can write $$
\dcb
&&\mathbb{E}_{\mu}[\int_0^\infty H_s\ind_{(0,\kappa\wedge \rho)}(s) dC_s] 
=\mathbb{E}_{\mu}[K_{\kappa\wedge \rho-}]
=\mathbb{E}_{\nu}[K_{\kappa-}\ind_{\{\kappa<\rho\}}\alpha_\kappa+\int_{0}^\kappa K_{s-}dV_s]\\

&=&\mathbb{E}_{\nu}[K_{\kappa}\ind_{\{\kappa<\rho\}}\alpha_\kappa+\int_{0}^\kappa K_{s-}dV_s]
=\mathbb{E}_{\nu}[\int_0^\kappa K_{s-}d(\alpha\ind_{[0,\rho)})_s+\int_0^\kappa (\alpha\ind_{[0,\rho)})_{s}dK_s+\int_{0}^\kappa K_{s-}dV_s]\\
&=&\mathbb{E}_{\nu}[\int_0^\kappa (\alpha\ind_{[0,\rho)})_{s}dK_s]
\ \mbox{ because $(\alpha\ind_{[0,\rho)})^\kappa$ is in class$(D)$,}\\
&=&\mathbb{E}_{\nu}[\int_0^\infty \alpha_s\ind_{(0,\kappa\wedge \rho)}(s)H_sdC_s].
\dce
$$
And also
$$
\dcb
&&\mathbb{E}_{\mu}[\int_0^\infty H_s d(\ind_{(0,\kappa\wedge \rho)}\centerdot C)^{\mu\cdot\mathbb{G}-p}_s] 

=\mathbb{E}_{\mu}[\int_0^\infty H_s \ind_{\{\alpha_{s-}>0\}}d(\ind_{(0,\kappa\wedge \rho)}\centerdot C)^{\mu\cdot\mathbb{G}-p}_s]\\

&=&\mathbb{E}_{\mu}[\int_0^\infty H_s\ind_{\{\alpha_{s-}>0\}}\ind_{(0,\kappa\wedge \rho)}(s) dC_s] 

=\mathbb{E}_{\nu}[\int_0^\infty H_s\ind_{\{\alpha_{s-}>0\}}\alpha_s\ind_{(0,\kappa\wedge \rho)}(s)dC_s]\\

&=&\mathbb{E}_{\nu}[\int_0^\infty H_s\ind_{\{\alpha_{s-}>0\}}d(\alpha\ind_{(0,\kappa\wedge \rho)}\centerdot C)^{\nu\cdot\mathbb{G}-p}_s]

=\mathbb{E}_{\nu}[\int_0^\infty \alpha_{s-}\frac{1}{\alpha_{s-}}H_s\ind_{\{\alpha_{s-}>0\}}d(\alpha\ind_{(0,\kappa\wedge \rho)}\centerdot C)^{\nu\cdot\mathbb{G}-p}_s]\\

&=&\mathbb{E}_{\mu}[\int_0^\infty H_s\frac{1}{\alpha_{s-}}d(\alpha\ind_{(0,\kappa\wedge \rho)}\centerdot C)^{\nu\cdot\mathbb{G}-p}_s]. \ \ok

\dce
$$

\

\subsection{Regular conditional distribution}\label{rcd}

The notion of the regular conditional distribution (cf. \cite{freedman, KS, RW}) will be necessary to present the next section. It is also a basic notion used to design the local solution method in section \ref{lsmfind}.

Consider a probability space $(\Omega,\mathcal{A},\mathbb{P})$. Let $\mathcal{B}$ be a sub-$\sigma$-algebra of $\mathcal{A}$. 

\bd
A function $\pi$ from $\Omega \times \mathcal{A}$ into $[0,1]$ is called a regular conditional probability on $\mathcal{A}$ given $\mathcal{B}$, if
\ebe
\item
for each $A\in\mathcal{A}$, the map $\omega \rightarrow \pi(\omega,A)$ is a version of $\mathbb{P}[A|\mathcal{B}]$;
\item
for each $\omega\in \Omega$, the map $A \rightarrow \pi(\omega,A)$ is a probability measure on $\mathcal{A}$.
\dbe
\ed

Let $X$ a measurable map from $(\Omega,\mathcal{A})$ into a measurable space $(\mathtt{C},\mathcal{C})$. 

\bd
A function $\overline{\pi}$ from $\Omega \times\mathcal{C}$ into $[0,1]$ is called a regular conditional distribution of $X$ given $\mathcal{B}$, if
\ebe
\item
for each $C\in\mathcal{C}$, the map $\omega \rightarrow \overline{\pi}(\omega,C)$ is a version of $\mathbb{P}[X\in C|\mathcal{B}]$;
\item
for each $\omega\in \Omega$, the map $C \rightarrow \overline{\pi}(\omega,C)$ is a probability measure on $\mathcal{C}$.
\dbe
\ed

Let $Y$ a measurable map from $(\Omega,\mathcal{A})$ into a measurable space $(\mathtt{D},\mathcal{D})$. 

\bd
A function $\tilde{\pi}$ from $\mathtt{D} \times\mathcal{C}$ into $[0,1]$ is called a regular conditional distribution of $X$ given $Y=y$, if
\ebe
\item[0.]
for each $C\in\mathcal{C}$, the map $y\rightarrow \tilde{\pi}(y,C)$ is  $\mathcal{D}$ measurable;
\item
for each $C\in\mathcal{C}$, the map $\omega \rightarrow \tilde{\pi}(Y(\omega),C)$ is a version of $\mathbb{P}[X\in C|\sigma(Y)]$;
\item
for each $y\in \mathtt{D}$, the map $C \rightarrow \tilde{\pi}(y,C)$ is a probability measure on $\mathcal{C}$.
\dbe
\ed

\brem
If $\tilde{\pi}$ exists, for non negative $\mathcal{C}\otimes\mathcal{D}$ measurable function $h(x,y)$, we have $$
\int h(x,Y) \tilde{\pi}(Y,dx)=\mathbb{E}[h(X,Y)| \sigma(Y)]
$$
almost surely.
\erem

\bethe
The regular conditional probability on $\mathcal{A}$ given $\mathcal{B}$ exists and is unique, if $(\Omega,\mathcal{A})$ is a Polish space with its Borel $\sigma$-algebra. The regular conditional distribution of $X$ given $\mathcal{B}$ and the regular conditional distribution of $X$ given $Y=y$ exist and are unique, if $(\mathtt{C},\mathcal{C})$ is a Polish space with its Borel $\sigma$-algebra. If in addition $\mathcal{B}$ is countably generated, a version of the regular conditional probability exists such that $$
\{\omega\in\Omega: \pi(\omega,B) = \ind_{\{\omega\in B\}}, \ \forall B\in \mathcal{B}\} \in\mathcal{B},
$$
and is of probability 1. Similarly, if $\sigma(Y) \subset \sigma(X)$ (i.e. $Y=\kappa(X)$ for some measurable $\kappa$), if $\mathcal{D}$ is countably generated, a version of $\tilde{\pi}$ exists such that $$
\{y\in\mathtt{D}: \tilde{\pi}(y, \kappa^{-1}(D)) = \ind_{\{y\in D\}}, \ \forall D\in \mathcal{D}\} \in\mathcal{D},
$$
and is of probability 1 for the law of $Y$.
\ethe

\textbf{Proof.} We only need to proceed along with the approach of \cite[Theorem 89.1]{RW}. To deal with the situation $\sigma(Y) \subset \sigma(X)$, we recall Doob's measurability theorem \cite[Theorem 1.5]{HWY} \ \ok

\brem
If the regular conditional probability on $\mathcal{A}$ given $\mathcal{B}$ exists, its image by $X$ will be a regular conditional distribution of $X$ given $\mathcal{B}$.
\erem

\brem
Consider another probability space $(\Omega',\mathcal{A}', \mathbb{P}')$ and a measurable map $\phi$ from $\Omega'$ into $\Omega$. Consider $\phi$ as a map into the space $(\Omega,\mathcal{B})$ and denoted by $\psi$. Suppose that $\mathbb{P}$ is the image of $\mathbb{P}'$ by $\phi$. Then, the regular conditional probability on $\mathcal{A}$ given $\mathcal{B}$ under $\mathbb{P}$ is a version of the regular conditional distribution of $\phi$ given $\psi$ under $\mathbb{P}'$. 
\erem

\
\subsection{Prediction process}

A ultimate goal beyong the study of the enlargement of filtration is to establish a "calculus" on the filtrations. People has thought that the prediction process is a suitable device to get such a calculus system. Here below is given a presentation of the prediction process based on {\cite[Theorem 13.1]{aldous}} and \cite{SongThesis}. (Notice that the local solution method, based on the change of probability measures, provides a system to compute on the filtrations, giving a partial answer to the question of "calculus".)

We consider a Polish space $\mathtt{S}$ and the associated Skorohod space $\mathsf{D}(\mathtt{S})$. Consider a cadlag adapted process $X$ defined on a stochastic basis $(\Omega,\mathbb{F}, \mathbb{P})$, which is also considered as a random variable taking value in $\mathsf{D}(\mathtt{S})$. Let $\mathsf{\Pi}$ denote the space of all probability measures on $\mathsf{D}(\mathtt{S})$, equipped with the topology of weak convergence (cf. \cite{billingsley} for the notion). 

\bethe\label{aldous}
	There exists a unique (up to indistinguishability) Skorohod process $\chi$, adapted to $\mathbb{F}$ taking values in $\mathsf{\Pi}$, such that, for any finite $\mathbb{F}$ stopping time $T$, $\chi_T$ is a regular conditional distribution for $X$ given $\mathcal{F}_T$ under $\mathbb{P}$.  
\ethe

(See \cite{KS, RW} for the notion of regular conditional distribution.) The process $\chi$ is called the prediction process associated with $X$. For $t\geq 0$, denote by $\pi_{t}$ the projection operator at $t$ from $\mathsf{D}(\mathtt{S})$ into $\mathtt{S}$. We have 

\bethe
For any finite $\mathbb{F}$ stopping time $T$, $\pi_{T}\chi_{T} = \delta_{X_T}$. The process $X$ is predictable, if and only if $\pi_{t}\chi_{t-}=\delta_{X_{t}}$ for $t\geq 0$. $X$ is quasi-left continuous, if and only if $\pi_{t}\chi_{t-}=\delta_{X_{t-}}$ for $t\geq 0$.
\ethe

\bethe
Let $H$ be a non negative measurable function on $\mathbb{R}_{+}\times \mathsf{D}(\mathtt{S})$. The process $(\chi_{t}(H_{t}), t\geq 0)$ is a version of the optional projection of the process $(H(t,X),t\geq 0)$ and The process $(\chi_{t-}(H_{t}), t\geq 0)$ is a version of the predictable projection of the process $(H(t,X),t\geq 0)$.
\ethe

\bethe
Suppose that $\mathbb{F}$ is the completion of the natural filtration of $X$. Then, $\mathbb{F}$ is also the completion of the natural filtration of $\chi$. Any optional process has an indistinguishable version of the form $(\varphi(t, \chi_{t}), t\geq 0)$. Any predictable process has an indistinguishable version of the form $(\varphi(t, \chi_{t-}), t\geq 0)$.
\ethe

\brem
The filtration acts on the stochastic calculus through the conditionings, especially through the optional and predictable projections. In this sense, the prediction process constitutes a representative of the filtration. As prediction processes are measure valued processes, one may expect to "compute" the filtrations via their prediction processes. But the things do not go so well, because the prediction processes are not so easy to handle as well.
\erem

\

\section{Enlargement of filtration in general}

We work in this section on a probability space $(\Omega, \mathcal{A},\mathbb{P})$ endowed with two filtrations $\mathbb{F}$ and $\mathbb{G}$ in an enlargement setting $\mathbb{F}\subset \mathbb{G}$.

\textbf{Problem} The fundamental question about the enlargement of filtration $\mathbb{F}\subset\mathbb{G}$, called \textit{the problem of filtration enlargement}, is whether a given $\mathbb{F}$ local martingale $X$ is a $\mathbb{G}$ semimartingale, and, if it is the case, what is the canonical decomposition of $X$ in $\mathbb{G}$ with its martingale part and its drift,  called \textit{enlargement of filtration formula}.

We have an immediate result.

\bethe\label{ssiG}
If an $\mathbb{F}$ local martingale $X$ is a $\mathbb{G}$ semimartingale, $X$ must be a $\mathbb{G}$ special semimartingale.
\ethe

\textbf{Proof.} We obtain this theorem with \cite[Theorem 8.6]{HWY} applied on the processs $X^*$. \ok

\subsection{Hypothesis $(H)$}

An immediate natural question about the problem of filtration enlargement is whether all $\mathbb{F}$ martingales can be $\mathbb{G}$ martingales. 

\bd
Consider a filtration enlargement setting $\mathbb{F}\subset \mathbb{G}$. We say that hypothesis $(H)$ holds for the filtration enlargement $\mathbb{F}\subset \mathbb{G}$, if any $\mathbb{F}$ local martingale is a $\mathbb{G}$ local martingale.
\ed

We have the following result (cf. \cite[Theorem 3]{BY} and on \cite[Chapter 3 section 2]{BJR}). This result constitutes one of the rare situations where a complete result can be obtained on an assumption formulated in a general setting of filtration enlargement (instead of a particular model).

\bethe
The following conditions are equivalent.
\ebe
\item
Hypothesis$(H)$ holds.

\item
For every $t\geq 0$, for every bounded $\mathcal{F}_{\infty}$ measurable random variable $\xi$, $\mathbb{E}[\xi | \mathcal{F}_{t}] = \mathbb{E}[\eta | \mathcal{G}_{t}]$.
\item
For every $t\geq 0$, for every bounded $\mathcal{G}_{t}$ measurable random variable $\eta$, $\mathbb{E}[\eta | \mathcal{F}_{t}] = \mathbb{E}[\eta | \mathcal{F}_{\infty}]$.

\item
For every $t\geq 0$, $\mathcal{F}_{\infty}$ is conditionally independent of $\mathcal{G}_{t}$, given $\mathcal{F}_{t}$.

\dbe
\ethe

\textbf{Proof.}
Hypothesis $(H)$ implies clearly the second condition. Let $\eta$ be a bounded $\mathcal{G}_{t}$ measurable random variable. For any bounded $\mathcal{F}_{\infty}$ measurable random variable $\xi$, applying the second property, we can write$$
\mathbb{E}[\mathbb{E}[\eta|\mathcal{F}_{t}]\ \xi]
=
\mathbb{E}[\eta \ \mathbb{E}[\xi|\mathcal{F}_{t}]]
=
\mathbb{E}[\eta \ \mathbb{E}[\xi|\mathcal{G}_{t}]]
=
\mathbb{E}[\eta \xi]
=
\mathbb{E}[\mathbb{E}[\eta|\mathcal{F}_{\infty}]\ \xi].
$$
This proves the third property. Let $\chi$ be a bounded $\mathcal{F}_{t}$ measurable random variable. We can write
$$
\mathbb{E}[\mathbb{E}[\eta\  \xi\ \chi]
=
\mathbb{E}[\mathbb{E}[\eta|\mathcal{F}_{\infty}]\ \xi \ \chi]
=
\mathbb{E}[\mathbb{E}[\eta|\mathcal{F}_{t}]\ \xi \ \chi]
=
\mathbb{E}[\mathbb{E}[\eta|\mathcal{F}_{t}]\ \mathbb{E}[\xi|\mathcal{F}_{t}] \ \chi],
$$
which implies the fourth property. Finally, if the fourth condition holds, 
$$
\mathbb{E}[\mathbb{E}[\eta\  \xi]
=
\mathbb{E}[\mathbb{E}[\eta\ \xi|\mathcal{F}_{t}]]
=
\mathbb{E}[\mathbb{E}[\eta|\mathcal{F}_{t}]\ \mathbb{E}[\xi|\mathcal{F}_{t}]]
=
\mathbb{E}[\eta\ \mathbb{E}[\xi|\mathcal{F}_{t}]],
$$
i.e., hypothesis $(H)$ holds.\ \ok

\textbf{Example.} When $L$ is an random variable independent of $\mathcal{F}_{\infty}$, the initial enlargement of $\mathbb{F}$ with $L$ satisfies hypothesis $(H)$.

\textbf{Example.}
Cox construction of random time gives rise of a progressive enlargement of filtration, which satisfies hypothesis $(H)$. See \cite[Chapter 3, section 3.2.2]{BJR}.

\brem
Hypothesis $(H)$ is a very strong condition on a filtration enlargement setting. However, later on, we will see by the local solution method that most of the known examples of filtration enlargement are connected with this hypothesis $(H)$ "locally" under  probability change. 
\erem

\

\subsection{Hypothesis $(H')$ and drift operator}

A natural extension of hypothesis $(H)$ is the following one.

\bd
Consider a filtration enlargement setting $\mathbb{F}\subset \mathbb{G}$. We say that hypothesis $(H')$ holds for the filtration enlargement $\mathbb{F}\subset \mathbb{G}$, if any $\mathbb{F}$ local martingale is a $\mathbb{G}$ semimartingale.
\ed

The fundamental result under hypothesis $(H')$ is about the drift operator. It is a second rare result that one can obtain in such a generality. According to Theorem \ref{ssiG}, the following notion is well-defined.

\bd
In the case of hypothesis $(H')$, for any $\mathbb{F}$ local martingale $X$, we denote by $\Gamma(X)$ the drift part of $X$ in the filtration $\mathbb{G}$. We call the map $X\rightarrow \Gamma(X)$ the drift operator of the enlargement setting $\mathbb{F}\subset \mathbb{G}$.
\ed

The operator $\Gamma$ takes values in the space $\mathsf{V}^0(\mathbb{G})$ of all $\mathbb{G}$ predictable processes with finite variations. For any element $A$ in the space $\mathsf{V}^0(\mathbb{G})$, we put$$
\varrho(A)=\mathbb{E}[1\wedge \left(\int_{0}^\infty |dA_{s}|\right)],
$$
which defines a tolopogy on the space $\mathsf{V}^0(\mathbb{G})$. Recall the normed space $\mathsf{H}^r(\mathbb{F})$ (for a $r\geq 1$) of $\mathbb{F}$ martingales.

\bethe
For any $t\geq 0$, the drift operator considered on the interval $[0,t]$ defines a continuous linear operator from $\mathsf{H}^r(\mathbb{F})$ into $\mathsf{V}^0(\mathbb{G})$.
\ethe

\textbf{Proof.} 
Let $t\geq 0$. Let $X$ be an $\mathbb{F}$ local martingale. Notice that, for any simple $\mathbb{G}$ predictable process $J$, the stochastic integral $J\centerdot X$ is well-defined in an elementary way. Using \cite[Proposition 2]{sss} (cf. also \cite[Chapter IV, Theorem 2]{protter}) we can check the identity$$
\varrho(\Gamma(X^t)) = \sup_{J}\mathbb{E}[1\wedge  |J\centerdot \Gamma(X^t)|],
$$
where $J$ runs over the family of all simple $\mathbb{G}$ predictable processes bounded by 1. Denoting the term on the right hand side by $\psi(\Gamma(X))$, denoting by $\overline{X}$ the martingale part of $X$ in $\mathbb{G}$, we have again $$
\psi(\Gamma(X^t))\leq \psi(X^t) + \psi(\overline{X}^t),
$$
with obvious interpretation of the right hand side terms. According to Proposition \ref{bracketcontrol}, there exists a universal constant $\gamma_{r}$ such that $$
\psi(\overline{X}^t) \leq (1+\gamma_{r})\|[X,X]^{1/2}\|_{r}.
$$
Hence, to prove the theorem, it is enough to prove that $\psi(X^t)$ tends to zero when $X$ tends to zero in $\mathsf{H}^r(\mathbb{F})$.

Suppose the opposite. Then, there exists a $\alpha>0$ such that, for all integer  $n\geq 1$, there exists a ${^n}\!X\in \mathsf{H}^r(\mathbb{F})$ and a simple $\mathbb{G}$ predictable process ${^n}\!H$ bounded by 1, such that $\|{^n}\!X\|_{\mathsf{H}^r}\leq \frac{1}{n}$ and $$
\mathbb{E}[1\wedge  \left|{^n}\!H \centerdot {^n}\!X^t\right|] \geq \alpha.
$$
On the other side, we can consider the maps $\phi_{n}:X\rightarrow {^n}\!H \centerdot {^n}\!X_{t}$, for each $n\geq 1$. It is a family of continuous linear map from $\mathsf{H}^r(\mathbb{F})$ into $\mathsf{L}^0(\Omega)$ (equipped with the usual distance $\mathfrak{t}(\xi)=\mathbb{E}[1\wedge |\xi|]$) with the property $\lim_{n\rightarrow\infty}\phi_{n}(X)=0$ for every $X\in \mathsf{H}^r(\mathbb{F})$. This implies, for any $0<\epsilon<\frac{\alpha}{4}$,$$
\cup_{n\geq 1}\{X \in \mathsf{H}^r(\mathbb{F}): \forall k\geq n, \mathfrak{t}(\phi_{k}(X))\leq \epsilon\} = \mathsf{H}^r(\mathbb{F}).
$$
By Baire property, there exists a $n_{0}$ for which $\{X: \forall k\geq n_{0}, \mathfrak{t}(\phi_{k}(X))\leq \epsilon\}$ contains a ball $\mathtt{B}(Y,a)$. We obtain then that, for any $X$ with $\|\!X\|_{\mathsf{H}^r}\leq a$, for any $k\geq n_{0}$,$$
\mathfrak{t}(\phi_{k}(X)) \leq \mathfrak{t}(\phi_{k}(Y)) + \mathfrak{t}(\phi_{k}(Y+X)) \leq 2\epsilon.
$$
In particular, $$
\mathbb{E}[1\wedge  \left|{^n}\!H \centerdot {^n}\!X^{t}\right|] =
\mathfrak{t}(\phi_{n}({^n}\!X))\leq 2\epsilon<\alpha,
$$
for $n>\frac{1}{a}\vee n_{0}$, a contradiction, which proves the theorem.
 \ \ok

\bethe
Suppose hypothesis $(H')$. For $0<v<2\wedge r$, there exists an equivalent probability measure $\mathbb{Q}$ such that, for any $\mathbb{F}$ local martingale $X$ in $\mathsf{H}^r(\mathbb{F},\mathbb{P})$, $X$ is a $\mathbb{G}$ semimartingale in $\mathsf{H}^v(\mathbb{G},\mathbb{Q})$.
\ethe

\

\subsection{The problem of "faux-amis"}

Suppose that a ${\mathbb{F}}$ local martingale $X$ remains a ${\mathbb{G}}$ (special) semimartingale. Let $X=M+V$ be the canonical decomposition in $\mathbb{G}$, where $M$ is a ${\mathbb{G}}$ local martingale and $V$ is a ${\mathbb{G}}$ predictable càdlàg process with finite variation. 

\bpro\label{jfg}
Let $J\in \mathcal{I}({\mathbb{P}},{\mathbb{F}},X)$ and denote $Y=J\stocint X$. If $Y$ is a $({\mathbb{P}},{\mathbb{G}})$ semimartingale, we have $J\in\mathcal{I}({\mathbb{P}},{\mathbb{G}},M)\cap\mathcal{I}({\mathbb{P}},{\mathbb{G}},V)$ and $Y=J\stocint M+J\stocint V$.
\epro

\proof 
The key point in this proposition is to prove the $M$-integrabiliy and the $V$-integrability of $J$. Without loss of the generality, suppose $Y\in\mathsf{H}^1_0({\mathbb{F}})$. Let $J_{\cdot,n}=\ind_{\{\sup_{1\leq i\leq d}|J_i|\leq n\}}J$ for $n\in\mathbb{N}$. By Proposition \ref{bracketcontrol}, for some constant $C$, $$
\mathbb{E}[\sqrt{[J_{\cdot,n}\stocint M,J_{\cdot,n}\stocint M]_\infty}]
\leq
C\mathbb{E}[\sqrt{[J_{\cdot,n}\stocint X,J_{\cdot,n}\stocint X]_\infty}]
\leq 
C\mathbb{E}[\sqrt{[Y,Y]_\infty}]<\infty, \ n\in\mathbb{N}.
$$
Taking $n\uparrow\infty$, we prove $J\in\mathcal{I}({\mathbb{P}},{\mathbb{G}},M)$. Because $Y$ is a semimartingale in $\mathbb{G}$, applying Proposition \ref{equalintegral},$$
\dcb
&&\ind_{\{\sup_{1\leq i\leq d}|J_i|\leq n\}}\stocint Y \ \mbox{ (in ${\mathbb{G}}$)}

=\ind_{\{\sup_{1\leq i\leq d}|J_i|\leq n\}}\stocint Y \ \mbox{ (in ${\mathbb{F}}$)}

=J_{\cdot,n}\stocint X \ \mbox{ (in ${\mathbb{F}}$)}\\

&=&J_{\cdot,n}\stocint X \ \mbox{ (in ${\mathbb{G}}$)}

=
J_{\cdot,n}\stocint M +J_{\cdot,n}\stocint V \ \mbox{ (in ${\mathbb{G}}$)}.
\dce
$$
As in \cite[Chapter III.6b]{JSh}, we represent the components $V_i$ of $V$ in the form $a_i\stocint F$, where the processes $a_i$ and $F$ are supposed to be ${\mathbb{G}}$ predictable and $F$ is càdlàg increasing. We have$$
\dcb
&&|\sum_{i=1}^dJ_ia_i|\ind_{\{\sup_{1\leq i\leq d}|J_i|\leq n\}}\stocint F\\
&=&
\mbox{sgn}(\sum_{i=1}^dJ_ia_i)\ind_{\{\sup_{1\leq i\leq d}|J_i|\leq n\}}J\stocint V\\
&=&\mbox{sgn}(\sum_{i=1}^dJ_ia_i)\ind_{\{\sup_{1\leq i\leq d}|J_i|\leq n\}}\stocint Y
-
\mbox{sgn}(\sum_{i=1}^dJ_ia_i)\ind_{\{\sup_{1\leq i\leq d}|J_i|\leq n\}}J\stocint M \ \mbox{ (in ${\mathbb{G}}$)}.
\dce
$$
As $Y$ and $J\stocint M$ are semimartingales in $\mathbb{G}$, by \cite[Remark(ii) of Definition 4.8]{CS} and \cite[Lemma 4.11 and Lemma 4.12]{CS}, the terms on the right hand side of this equality converge in probability when $n\uparrow \infty$. Consequently, $J\in\mathcal{I}({\mathbb{P}},{\mathbb{G}},V)$ (cf. \cite[Definition 3.7]{CS}) and the proposition follows. \ok

The next result is an extension of \cite[Proposition (2.1)]{Jeulin80}.

\bcor\label{corFA}
$Y=J\stocint X$ is a $({\mathbb{P}},{\mathbb{G}})$ semimartingale, if and only if $J\in\mathcal{I}({\mathbb{P}},{\mathbb{G}},M)\cap\mathcal{I}({\mathbb{P}},{\mathbb{G}},V)$.
\ecor

\textbf{The "Faux-Amis".}
Proposition \ref{jfg} is motivated by the "faux-amis" of \cite{JYfa}. (See also \cite{MY}).
\begin{quote}
	During 1977-1978, we have studied systematically Hypothesis $(H')$. At various occasions, we have "encountered" the following statements.

	\textbf{FA.1} 
	\textit{For $(H')$ to be realized, it is enough to have a generating system of $\mathbb{F}$ martingales to be $\mathbb{G}$ semimartingales.	
	}
	
	(For a long time, that statement have seemed very likely, especially in the case where the $\mathbb{F}$ martingales are generated by one martingale alone.)

	\textbf{FA.2} 
	\textit{
		If an $\mathbb{F}$ locally square integrable martingale $X$ is a  $\mathbb{G}$ semimartingale, the drift $A$ of $X$ in $\mathbb{G}$ satisfies $dA_{s}<\!\!< d\cro{X,X}^{\mathbb{F}}$.
	}
	
	(This second statement has originated, on the one hand, in the analogy which exists between the filtration enlargement and the change of probability measure, and on the other hand, in the famous Girsanov theorem.)
	
	Actually, as our notations want to say, the statements \textbf{FA.1} and \textbf{FA.2} are "faux-amis".
\end{quote}

\

Based on result such as Corollary \ref{corFA} or \cite[Proposition (2.1)]{Jeulin80}, a counter-example to \textbf{FA.1} is found in \cite{JYfa}.

\bl\label{processR}
Let $R$ be a positive measurable process with the following property: there exists a probability measure $\mu$ sur $\mathbb{R_{+}}$ such that $\mu[\{0\}]=0$, $\int_{0}^\infty x\mu(dx)=1$ and ${^{p}}\!f(R)=\int_{0}^\infty f(x)\mu(dx)$, for any bounded Borel function $f$. Then, for any predictable non decreasing process $A$, null at the origin, we have $$
\{A_{\infty}<\infty\} = \{\int_{0}^\infty R_{s} dA_{s}<\infty\}, \ \mbox{ almost surely.}
$$
\el

\textbf{Proof.} (It is a beautiful example of stochastic calculus. See \cite[Lemme(3.22)]{Jeulin80}.) For any non decreasing integrable process $B$ null at the origin, for any integer $n\geq 1$, let $T_{n}=\inf\{s>0: B^p_{s}\geq n\}$. Then,  $T_{n}$ is a predictable stopping time and therefore $$
\mathbb{E}[\ind_{\{T_{n} =\infty\}} B_{\infty}] 
\leq \mathbb{E}[B_{T_{n}-}] = \mathbb{E}[B^p_{T_{n}-}]\leq n.
$$ 
This shows that $B_{\infty}<\infty$ on $\cup_{n=1}^\infty\{T_{n} =\infty\}=\{B^p_{\infty}<\infty\}$. 
 Consider the two sets in the above equality. We note that the assumption of the lemma implies $(R\centerdot A)^p = A$, which yields the inclusion "$\subset$".

For an indicator random variable $J_{\infty}$, we compute the predictable projection of $J_{\infty}R$.  Let $J$ denote the martingale associated with $J_{\infty}$ and let $\jmath$ denote the minimum of the process $J$. For a predictable stopping time $T$, we have$$
\dcb
{^{p}}(J_{\infty}R)_{T}\ind_{\{T<\infty\}}
=
\mathbb{E}[J_{\infty}R_{T}|\mathcal{F}_{T-}]\ind_{\{T<\infty\}} = \int_{0}^\infty du \ \mathbb{E}[J_{\infty} \ind_{\{R_{T}>u\}}|\mathcal{F}_{T-}]\ind_{\{T<\infty\}}.
\dce
$$ 
Because $J_{\infty}$ is an indicator,$$
\dcb
\mathbb{E}[J_{\infty} \ind_{\{R_{T}>u\}}|\mathcal{F}_{T-}]
&=&
\mathbb{E}[(J_{\infty} - \ind_{\{R_{T}\leq u\}})^+|\mathcal{F}_{T-}]\\
&\geq&
(\mathbb{E}[J_{\infty} - \ind_{\{R_{T}\leq u\}}|\mathcal{F}_{T-}] )^+
=
(J_{T-} - \int_{0}^u \mu(dx) )^+.
\dce
$$
Integrating with respect to $u$, we get $
{^{p}}(J_{\infty}R)_{T}
\geq \phi(J_{T-}),
$
where $\phi(x)=\int_{0}^\infty(x - \int_{0}^u \mu(dx) )^+ du$. With this inequality in mind, we write the following inequalities concerning the processs $A$.
$$
\dcb
\mathbb{E}[\phi(\jmath)A_{\infty}]
\leq
\mathbb{E}[\int_{0}^\infty \phi(J_{s-})dA_{s}]
\leq
\mathbb{E}[ \int_{0}^\infty {^{p}}(J_{\infty} R)_{s}dA_{s}]
=
\mathbb{E}[ J_{\infty} \int_{0}^\infty R_{s}dA_{s}],
\dce
$$
i.e., $\phi(\jmath)A_{\infty}<\infty$ almost surely, whenever $J_{\infty} \int_{0}^\infty R_{s}dA_{s}$ is integrable. Note that $\phi$ is non decreasing, $\phi(0)=0, \phi(1)=1$, and strictly positive on $(0,1]$. The martingale $J$ being non negative, by \cite[Chapitre IV, $n^\circ$17]{DM2}, we obtain $\phi(\jmath)>0$ on $\{J_{\infty}=1\}$, or $$
\{J_{\infty}=1\}\subset \{A_{\infty}<\infty\}.
$$
Now, apply these inequalities with $J_{\infty}$ to be the indicator function of $\{(R\centerdot A)_{\infty} < n \}$ for every integer $n\geq 1$. The lemma results. \ \ok

Consider now a Brownian motion $X$ with its natural filtration $\mathbb{F}$. Consider the value $X_{1}$ of the Brownian motion $X$ at $t=1$. Let $R_{t}=\ind_{\{t<1\}}\frac{|X_{1}-X_{t}|}{\sqrt{t-1}}$. Then, with the independent increments, with the gaussian property, we check that the process $R$ satisfies the property of Lemma \ref{processR} with respect to $\mathbb{F}$ with the probability measure $\mu$ the law of $|X_{1}|$. 

Let $\mathbb{G}$ be the initial enlargement of $\mathbb{F}$ with the random variable $X_{1}$. Then, again by the gaussian property, the process $X$ is a $\mathbb{G}$ semimartingale whose $\mathbb{G}$ drift is given by $\int_{0}^t\ind_{\{s<1\}} \frac{X_{1}-X_{s}}{1-s}ds, t\geq 0$. Let $H$ be any $\mathbb{F}$ predictable process with $\int_{0}^t H_{s}^2 ds <\infty$ for all $t\geq 0$. According to Corollary \ref{corFA}, the stochastic integral $H\centerdot X$ is a $\mathbb{G}$ semimartingale, if and only if $$
\int_{0}^1 |H_{s}|\ind_{\{s<1\}} \frac{|X_{1}-X_{s}|}{1-s}ds = \int_{0}^1  \frac{|H_{s}|}{\sqrt{1-s}}  R_{s} ds <\infty.
$$
According to Lemma \ref{processR}, this integral is finite, if and only if $$
 \int_{0}^1  \frac{|H_{s}|}{\sqrt{1-s}}  ds <\infty.
$$
Now we introduce the (determinist) process $$
H_{s} = \frac{1}{\sqrt{1-s}} \frac{1}{(-\ln(1-s))^\alpha} \ind_{\{\frac{1}{2}<s<1\}},\ s\geq 0,
$$
for some $\frac{1}{2}< \alpha < 1$. We check directly that the last integral with this $H$ can not be finite, while for the same $H$, $\int_{0}^1 H_{s}^2 ds <\infty$. This proves that $H\centerdot X$ is not a $\mathbb{G}$ semimartingale.

\bethe
\textbf{FA.1} is wrong in general.
\ethe

\

\subsection{Invariance under isomorphism}

\bpro\label{sbi}
Let $\Psi$ be an isomorphism of the stochastic basis $(\Omega,\mathbb{G},\mathbb{P})$ onto a stochastic basis $(\Omega',\mathbb{G}',\mathbb{P}')$. Let $\mathbb{F}$ be a sub-filtration of $\mathbb{G}$. Let $\mathbb{F}'=\Psi(\mathbb{F})$, which is a sub-filtration of $\mathbb{G}'$. Then, a $(\mathbb{P},\mathbb{F})$ local martingale $M$ is a $(\mathbb{P},\mathbb{G})$ semimartingale, if and only if (a cadlag representative of) $\Psi(M)$ is a $(\mathbb{P}',\mathbb{G}')$ semimartingale. This means that the problem of filtration enlargement $\mathbb{F}\subset \mathbb{G}$ is equivalent to  the problem of filtration enlargement $\mathbb{F}'\subset \mathbb{G}'$.
\epro

\textbf{Proof.}
If $M$ is a $(\mathbb{P},\mathbb{G})$ semimartingale, it is special. Let $M=\overline{M}+V$ be its canonical decomposition in $(\mathbb{P},\mathbb{G})$. Then, $\Psi(M) = \Psi(\overline{M}) + \Psi(V)$, and, according to Proposition \ref{isota} and \ref{sbi2}, $\Psi(\overline{M})$ is a $(\mathbb{P}',\mathbb{G}')$ local martingale and $\Psi(V)$ is a $\mathcal{G}'$ predictable process with finite variation. This proves that $\Psi(M)$ is a $(\mathbb{P}',\mathbb{G}')$ semimartingale.

The inverse implication can be proved similarly. \ \ok

\

\subsection{A panorama of specific models and formulas}

To have an overview of the development of the theory of enlargement of filtration, we try to make a list of various formulas that we can find in the main literature. The list is not exhaustive. The applications are neglected.

\textbf{-Progressive enlargement of filtration}
\textit{Let $\mathbb{F}$ be a filtration and $\tau$ a random variable in $[0,\infty]$. The progressive enlargement of $\mathbb{F}$ with $\tau$ is the smallest filtration satisfying the usual conditions and making $\tau$ a stopping time.}

Here are some formulas on the progressive enlargement of filtration.
\ebe
\item
Barlow and Jeulin-Yor's formula for honnete times, cf. \cite{barlow, JY2}.
\item
Progressive enlargement and Jeulin-Yor's formula, cf. \cite{Jeulin80, DM3, JY2, yoeurp}.
\item
Yor's formula under martingale representation, cf. \cite{yor-n}.
\item
Cox construction of random time, cf. \cite{BJR}.
\item
Initial time or density approach, cf. \cite{CJZ, ElKJJ, ElKJJ2,  EJY, JC,  TXY}.
\item
Jeanblanc-Song and Li-Rutkowski formula, cf. \cite{JS1, JS2, LR}.
\item
Dynamic one-default model and Song's formula, cf. \cite{JS2, songmodel}.
\dbe

\

\textbf{-Initial enlargement of filtration}
\textit{Let $\mathbb{F}$ be a filtration and $L$ a random variable. The initial enlargement of $\mathbb{F}$ with $L$ is the smallest filtration satisfying the usual conditions and containing $L$ at origin. }

Here are some formulas on the initial enlargement of filtration.
\ebe
\item
Ito's idea to define the stochastic integral $\int_{0}^tB_{1}dB_{s}$, cf. \cite{AK, ito78, MY}.
\item
Initial enlargement with a partition, cf. \cite{yor-entro}.
\item
Initial enlargement and Jeulin's study, cf. \cite[Theorem (3.9)]{Jeulin80}
\item
Jeulin's formulas in various examples, cf. \cite{Jeulin80}
\item
Initial enlargement and Jacod's hypothesis, \cite{amen, jacod2, GP}.
\item
Formula with Malliavin calculus, cf. \cite{SongThesis, imkeller}.
\dbe

\

\textbf{-Enlargement by future infimum of a linear diffusion} This will be presented in detail below. See \cite{Jeulin80, Song-local}.

\

\textbf{-Multiple random times enlargement or successive enlargement} \textit{Let $\mathbb{F}$ be a filtration and $\tau_{1},\ldots,\tau_{k}$ (for a $k\in\mathbb{N}^*$) be $k$ random times. We call a successive enlargement of $\mathbb{F}$ by $\tau_{1},\ldots,\tau_{k}$ the smallest filtration satisfying the usual conditions making the random times $\tau_{1},\ldots,\tau_{k}$ stopping times.} As indicated by the name, a successive enlargement is the repetition of the progressive enlargement of filtration. However, the multiple presences of random times lead to combinational computation complexity. This problem has been involved in \cite{ElKJJ2, Jeulin80, pham, Song-splitting, Song-DC, TXY}.

\

\textbf{-Specification of the enlargement of filtration formula in applications} In applications, general enlargement of filtration formula should be specified under particular conditions. This specification is not always straightforward. The books \cite{JYC, MY} collect various examples where the enlargement of filtration formula is explicitly computed.

\

\textbf{-Random time specifications} \textit{Honest time, immersion time, initial time, pseudo-stopping time, pseudo-honest time, invariance time.} The importance of the progressive enlargement approach in credit risk modelling, the need to adapt progressive enlargements to specifical situations, leads to a great effort to reclassify the random times. \cite{ACDJ, CS, JC, NY1, JS1, JS2, LR}.

\

\textbf{-Unified methods} Besides the general study such as that about hypothesis $(H')$, the the theory of enlargement of filtration is composed mainly by the studied of concrete models and of the search of the enlargement of filtrations formulas. We notice a great diversity of methodologies to deal with different models. This can be explained, on the one hand, by the divers applications of the theory in different domains such as stochastic integral, stochastic differential equation, Brownian path decomposition, and nowaday, mathematical finance, etc. And on the other hand, it is because of our poor ability to operate on general filtrations. That said, it is nevertheless very natural to ask the question whether the diverse models can be worked on with a unified methodolog. This question has been the subject of  \cite{KP1, SongThesis, Song-local}. We will give a detailed presentation of the local solution method introduced in \cite{SongThesis}.

\

\textbf{-Martingale representation property}
The literature of the theory of enlargement of filtration is mainly marked by the search of enlargement of filtration formulas. But, in fact, any element in stochastic calculus is a potential subject in this theory. In \cite{BJR, BSJ,  CJZ, JC2, Kusuoka, TXY, JS3}, the martingale representation property has been studied in the context of progressive enlargement of filtration. As usual, the question has been tackled with different approaches. Once again, the local solution method used in \cite{JS3} provides the most general result.

\

\textbf{-No-arbitrage property in an enlarged filtration} The change of the information flow is a common problem in financial market modelling. The study of this problem leads to an intense effort to understand the non-arbitrage property under enlargement of filtration. See \cite{AFK, ACDJ, FJS, Song-drift, Song-NA}. Coming from the financial market modelling, the notion of non-arbitrage constitutes in fact an important characteristic of semimartingales and of the reference filtration.

\

\textbf{-Drift operator} We have seen that the drift operator is a well-defined notion in case of hypothesis $(H)$. Typically in a problem of filtration enlargement, one begins with the construction of an enlargement and deduce result in consequence. This mode of thinking, however, appears too rigid in many situations, because of the lack of knowledge on the reference filtration, and because the construction of the enlargement is usually not intrinsic for the issue. This point is not enough underlined in the literature. The paper \cite{Song-drift} gives a study where the non-arbitrage property in the enlarged filtration can be completely characterized in term of the drift operator alone, without explicit knowledge on the reference filtration. The works \cite{JS1, JS2} provide complementary information about this point. Another important aspect of the drift operator, as a characteristic of the enlargement, is its connexion with the notion of information and of the measurement of $\sigma$-algebras. See \cite{ankir, ankir2, ADI, Song-drift, yor-entro}.

\pagebreak

\begin{center}
\Large \textbf{Part II. 
Local solution method}
\end{center}

\pagebreak

\section{Local solution method : Introduction}\label{lsmintro}

As mentioned, We give below a full presentation of the methodology introduced in the thesis \cite{SongThesis} (cf. also \cite{Song-local}).

The thesis \cite{SongThesis} was motivated by two questions.
\ebe
\item
Why the  enlargement of filtration formulas look so much like Girsanov's formula ? What is the intrinsic relationship between the enlargement of filtration and the change of probability ?
\item
Is there a unified methodology which enable one to work on various different models ? Can one deduce the formulas in Jeulin's book with one and same general formula ? 
\dbe
We underline that these questions exist from the beginning of the literature of enlargement of filtration. And at the moment of \cite{SongThesis}, Jacod just published its famous criterion on the initial enlargement and Yoeurp just published its formal proof that the pregressive enlargement of filtration formula is a Girsanov formula. But, as much Jacod and Yoeurp undertake both with Girsanov theorem and their results are as general as possible in their respective fields, i.e., the initial enlargement or the progressive enlargement of filtration, as much the difference between their respective approaches is so apparent. This contrast  further reinforced the desire to understand the situation and to find a unified proof of both their formulas.

According to \cite{SongThesis}, we can make the following analysis. In the case, for example, of an initial enlargement of a filtration $\mathbb{F}$ with a random variable $L$, for an $\mathbb{F}$ martingale $M$, to discover if $M$ is a $\mathbb{G}$ local martingale, we are led to regard the conditioning $\mathbb{E}[M_{t}-M_{s}|\mathcal{F}_{s}\vee\sigma(L)]$ for $0\leq s<t$. This conditioning has an equivalent expression $\mathbb{E}_{\pi}[M_{t}-M_{s}|\mathcal{F}_{s}]$, where $\pi$ is the regular conditional probability given the random variable $L$. Note that, for the problem of filtration enlargement, this conditioning under $\pi$ should be considered in relation with the conditioning $\mathbb{E}[M_{t}|\mathcal{F}_{s}]=M_{s}$ under the initial probability measure $\mathbb{P}$. That being understood, a sufficient condition for the conditioning under $\pi$ to be computed with the conditioning under $\mathbb{P}$ is naturally to suppose that $\pi$ is absolutely continuous with respect to $\mathbb{P}$ on $\mathcal{F}_{\infty}$, or by Bayes formula that the conditional law of $L$ given $\mathcal{F}_{\infty}$ is absolutely continuous with respect to its own law.

More generally, consider an enlargement of a filtration $\mathbb{F}$ by a process $I$, i.e., the filtration $\mathbb{G}$ is generated by $\mathcal{F}_{t}\vee\sigma(I_{s},s\leq t)$ for $t\geq 0$. Suppose that under the probability measure $\mathbb{P}$, the process $I$ is independent of $\mathcal{F}_{\infty}$. Then, the problem of filtration enlargement has an immediate solution, because all $\mathbb{F}$ local martingales are $\mathcal{G}$ local martingale, thanks to the independence. Suppose that it is not the case, but instead, there exists a second probability measure $\mathbb{P}^0$ satisfying three conditions:
\ebe
\item
$\mathbb{P}$ coincides with $\mathbb{P}^0$ on $\mathcal{F}_{\infty}$.
\item
$\mathbb{P}$ is absolutely continuous with respect to $\mathbb{P}^0$ on $\mathcal{G}_{\infty}$.
\item
Under $\mathbb{P}^0$ the process $I$ is independent of $\mathcal{F}_{\infty}$.
\dbe
Then, the problem of filtration enlargement has a solution, because 
\ebe
\item[i.]
all $(\mathbb{P},\mathbb{F})$ local martingales are $(\mathbb{P}^0,\mathbb{F})$ local martingales, 
\item[ii.]
are $(\mathbb{P}^0,\mathbb{G})$ local martingales by independence, 
\item[iii.]
are $(\mathbb{P},\mathbb{G})$ semimartingales by Girsanov theorem. 
\dbe
By now, we have a prototype of a general problem solving method.

Clearly, this method is again primitive, whose implementation encounters several problems. First of all, a probability $\mathbb{P}^0$ satisfying the condition 3) will not exist in general. This is in fact a minor problem, because, thanks to Proposition \ref{sbi}, we can transfer the problem of filtration enlargement onto an enlarged probability space where such a probability $\mathbb{P}^0$ exists naturally.

A more serious problem is that the condition 2) is usually in conflict with the condition 3). A compromise is to remove the condition 2), and in consequence, in the step iii., we substitute the usual Girsanov theorem by the so-called generalized Girsanov theorem, i.e, the existence of a stopping time $T_{0}$ which marks the limit level of the validity of the Girsanov formula (i.e., the Girsanov formula is valid on $[0,T_{0})$, but invalid on $[T_{0},\infty)$). 

The above compromise raises a new problem, because the method is then limited at the level $T_{0}$. It is a real handicap. We notice that the paper \cite{yoeurp} has provided a Girsanov proof of the progressive enlargement formula, but has not involved the honest enlargement formula, because in particular of this limitation. An answer to this question is given in \cite{SongThesis}, where we propose to substitute the condition 3) by the existence of a probability measure $\mathbb{P}^t$ (for every given $t\geq 0$) satisfying
\ebe
\item[3$^\circ$]
Under $\mathbb{P}^t$ the process $I$ is independent of $\mathcal{F}_{\infty}$ on the time interval $[t,\infty)$.
\dbe
We then apply the generalized Girsanov theorem on $[t,\infty)$ (more precisely on $[t,T_{t})$). 

The consequence of this modification $3^\circ$ is that, instead of a solution to the problem of filtration enlargement, we get a collection of pieces of solutions on the random intervals $[t,T_{t}), t\geq 0$. It is this collection that is at the origin of the so-called \textit{local solution method}, a method which leads to a solution of the problem of filtration enlargement in two steps :
\ebe
\item[$_{\bullet}$]
To have a collection of pieces of solutions to the problem of filtration enlargement, called local solutions
\item[$_{\bullet\bullet}$]
To aggregate the local solutions into one global solution.
\dbe

\

\

\section{Local solution method : Pieces of semimartingales and their aggregration to form a global semimartingale}\label{noncovering}

We study the question of aggregating a collection of pieces of semimartingales to form a global semimartingale. The result is initiated in \cite{SongThesis} and achieved in \cite{Song99}. See also \cite{Song-local}. Before going into the details, we underline the similarity which exists between the semimartingales and the differentiable functions, because they are all integrators (cf. Bichteler-Dellacherie theorem). There are two aspects for a differentiable function $f$. It is differentiable because $\frac{df}{dt}$ exists, which is a local property. It is a function, because the integrals $\int_{0}^t \frac{df}{ds} ds$ exist, which is a global property. A semimartingale has the same features. In this section, we discuss the integrability issue.

\subsection{Assumption and result}\label{pieces}

We work on a probability space $(\Omega,\mathcal{B},\mathbb{P})$ equipped with a filtration $\mathbb{G}=(\mathcal{G}_t)_{t\geq 0}$ of sub-$\sigma$-algebras of $\mathcal{B}$, satisfying the usual condition (called also a stochastic basis). We consider in this section the processes satisfying the following conditions. Let $S$ be a real càdlàg $\mathbb{G}$ adapted process. Let $\mathtt{B}$ be a $\mathbb{G}$ random left interval (i.e. $\mathtt{B}=(T,U]$ with $T,U$ being two $\mathbb{G}$ stopping times). We say that $S$ is a $\mathbb{G}$ (special) semimartingale on $\mathtt{B}$, if $S^U - S^T$ is a $\mathbb{G}$ (special) semimartingale. We denote below $S^U - S^T$ by $\ind_\mathtt{B}\centerdot S$.

\bassumption\label{partialcovering}
We suppose
\ebe
\item[i.]
$S$ is a special semimartingale in its natural filtration.
\item[ii.]
There exists a sequence of random left intervals $(\mathtt{B}_i)_{i\in\mathbb{N}}$ on each of which $S$ is a special semimartingale. In any bounded open interval contained in $\cup_{i\in\mathbb{N}}\mathtt{B}_i$, there exist only a finite number of accumulation points of the right end points of third type in the sense of Lemma \ref{third-type}. 
\item[iii.]
There is a special semimartingale $\check{S}$ such that $\{S\neq \check{S}\}\subset \cup_{i\in\mathbb{N}}\mathtt{B}_i$.
\dbe
\eassumption

We will prove a theorem which characterizes, among the processes satisfying Assumption \ref{partialcovering}, those which are semimartingales. Before going into the details, however, let us consider a particular situation to illustrate the fundamental idea behind the theorem. Let $\pi_k, k\in\mathbb{N}$, be an increasing sequence of finite sets in $(0,1)$. For every $k$, the points of $\pi_k$ cuts $(0,1]$ into a number $n_k$ of left intervals, denoted by $B_{k,i}=(t^k_{i}, t^k_{i+1}]$, $i\leq n_k$. We suppose$$
\lim_{k\rightarrow\infty}\sup_{i\leq n_{k}}|t_{i+1}-t_{i}| = 0.
$$
Let $S$ to be an integrable process satisfying Assumption \ref{partialcovering} with respect to this family of left intervals $B_{k,i}$. Consider the quantity$$
\mbox{Var}(S) = \sup_{k\in\mathbb{N}}\sum_{i\leq n_{k}}\mathbb{E}[\left|\mathbb{E}[\ind_{B_{k,i}}\centerdot S |\mathcal{F}_{t^k_{i}} \right|].
$$
We know that if this quantity $\mbox{Var}(S)$ is finite, $S$ is a quasimartingale (cf. \cite[Definition 8.12]{HWY}), so that $S$ is a special semimartingale. We get therefore an answer to the question about the semimartingale issue of $S$, using information of the pieces of semimartingales $\ind_{B_{k,i}}\centerdot S$. But what the quantity $\mbox{Var}(S)$ measures on $S$ exactly ? Being a special semimartingale, $S$ has a canonical decomposition $S=M-A$ with $M$ martingale and $A$ predictable process with integrable variation. We can now rewrite the quantity $\mbox{Var}(S)$:$$
\mbox{Var}(S)
= \sup_{k\in\mathbb{N}}\sum_{i\leq n_{k}}\mathbb{E}[\left|\mathbb{E}[\ind_{B_{k,i}}\centerdot A |\mathcal{F}_{t^k_{i}} \right|]
=\mathbb{E}[\mbox{var}(A)] \ (\mbox{cf. \cite[Theorem 8.13]{HWY}}).
$$
We see that the quantity $\mbox{Var}(S)$, ignoring the martingale part $M$, measures exactly the drift part of $S$. In other words, the determining element making $S$ a semimartingale is the total variation of the drift process $A$.

\textbf{Associated random measure.} We now define an association between a process $S$ satisfying Assumption \ref{partialcovering} and a random measure $d\chi^{S,\cup}=d\chi^\cup$. Noting that $\ind_{\mathtt{B}_i}\centerdot S$ is a special semimartingale for any $i\geq 0$, we denote by $\chi^{\mathtt{B}_i}$ the drift process of $\ind_{\mathtt{B}_i}\centerdot S$. It is clear that the two random measures $\mathsf{d}\chi^{\mathtt{B}_i}$ and $\mathsf{d}\chi^{\mathtt{B}_j}$ coincides on $\mathtt{B}_i\cap\mathtt{B}_j$ for $i,j\geq 0$. We can therefore define with no ambiguity a $\sigma$-finite (signed) random measure $\mathsf{d}\chi^\cup$ on $\cup_{i\geq 0}\mathtt{B}_i$ such that $\ind_{\mathtt{B}_i}\mathsf{d}\chi^\cup=\mathsf{d}\chi^{\mathtt{B}_i}$.

\textbf{Distribution function.} We also introduce the notion of distribution function. We will say that a signed random measure $\lambda$ on $\mathbb{R}_+$ has a distribution function, if $\int_0^t|d\lambda|_s<\infty$ for any $t\geq 0$. In this case, the process $t\geq 0 \rightarrow \lambda([0,t])$ is called the distribution function of $\lambda$.

\textbf{Associated jump process.} For any $\mathbb{G}$ stopping time $R$, let $d_R=\inf\{s\geq R: s\notin \cup_{i\geq 0}\mathtt{B}_i\}$. Let$$
\mathtt{A}=\cup_{s\in\mathbb{Q}_+}(s, d_s],\ \mathtt{C}=\mathtt{A}\setminus \cup_{i\geq 0}\mathtt{B}_i.
$$
The sets $\mathtt{A}$ and $\mathtt{C}$ are $\mathbb{G}$ predictable sets. $\mathtt{C}$ is in addition a thin set contained in $\cup_{s\in\mathbb{Q}_+}[d_s]$. Denote by $C$ the following thin process $$
C_t=\ind_{\{t\in\mathtt{C}\}}\Delta_t(S-\check{S}), t\geq 0.
$$

We are ready to state the theorem of this section.

\bethe\label{semimartingaleproperty}
Suppose Assumption \ref{partialcovering}. For $S$ to be a semimartingale on the whole $\mathbb{R}_+$, it is necessary and it is sufficient that the random measures $\mathsf{d}\chi^\cup$ has a distribution function $\chi^\cup$ and the process $C$ is the jump process of a special semimartingale. 
\ethe

Note that the necessity of this theorem is clear (recalling that $S$ is locally in class$(D)$). As for the sufficiency, it will be the consequence of the different results we prove in the following subsections.

\

\subsection{Elementary integral}

We will need to do integral calculus on the sets $\mathtt{B}_{i }$ with respect to $S$ without semimartingale property. It will be the pathwise elementary integral, that we define below.

\subsubsection{Elementary integral on an interval}\label{eii}

We present here the notion of the integrals of elementary functions (which is applicable to elementary processes in an evident way).

Let $f$ be a real càdlàg function defined on an interval $[a,b]$ with $-\infty\leq a<b\leq \infty$ (where we make the convention that $(a,\infty]=(a,\infty)$). For any real number $c$, denote by $f^c$ the function 
$$
t\in[a,b]\rightarrow f(t\wedge c).
$$ 
A function $h$ on $(a,b]$ is called left elementary function, if $h$ can be written in such a form as : $h=\sum_{i=0}^nd_i\ind_{(x_i,x_{i+1}]}$ on $(a,b]$ for some real numbers $a=x_0<x_1<\ldots<x_{n}<x_{n+1}=b$ ($n\in\mathbb{N}$) and $d_0,d_1,\ldots,d_n$. 

\bd
We define $h\ind_{(a,b]}\centerdot f$ the elementary integral of $h$ with respect to $f$ to be the following function: $$
\left\{
\dcb
\mbox{for } t\leq a,& &(h\ind_{(a,b]}\centerdot f)(t) = 0,\\
\mbox{for } t\in(a,b],& &(h\ind_{(a,b]}\centerdot f)(t) = \sum_{i=0}^nd_i(f^{x_{i+1}}-f^{x_i})(t),\\
\mbox{for } t>b,& &(h\ind_{(a,b]}\centerdot f)(t) = (h\ind_{(a,b]}\centerdot f)(b).
\dce
\right.
$$
\ed

Then,

\bl\label{elem-integral} 
\ebe
\item
If $h$ has another representation $h=\sum_{i=0}^{m}e_i\ind_{(y_i,y_{i+1}]}$, we have$$
(h\ind_{(a,b]}\centerdot f) = \sum_{i=0}^{m}e_i(f^{y_{i+1}}-f^{y_i})\ \mbox{ on $(a,b]$.}
$$
\item
The elementary integral is bi-linear on the product set of the family of all left elementary functions and of the family of all real càdlàg functions on $(a,b]$.
\item
If $g$ is another left elementary function on $(a,b]$, we have$$
(g\ind_{(a,b]}\centerdot (h\ind_{(a,b]}\centerdot f))
=(gh\ind_{(a,b]}\centerdot f).
$$
\item
For any real number $c$, $(h\ind_{(a,b]}\centerdot f)^c=h\ind_{(a,b]}\centerdot f^c=(h\ind_{(a,c]})\ind_{(a,b]}\centerdot f
=h\ind_{(a,c\wedge b]}\centerdot f$. For any real number $c<d$ such that $(c,d]\subset(a,b]$, $(h\ind_{(c,d]})\ind_{(a,b]}\centerdot f=h\ind_{(c,d]}\centerdot f$.

\item
$\Delta_t(h\ind_{(a,b]}\centerdot f)=h(t)\ind_{a<t\leq b}\Delta_t f$ for $t\in\mathbb{R}$.
\dbe
\el

\textbf{Proof.}
\ebe
\item
Let $a=z_1<\ldots<z_k<z_{k+1}=b$ be a refinement of $a=x_0<x_1<\ldots<x_{n}<x_{n+1}=b$ and of $a=y_0<y_1<\ldots<y_{m}<y_{m+1}=b$. We note that $$
\dcb
\sum_{i=0}^nd_i(f^{x_{i+1}}-f^{x_i})&=&\sum_{i=0}^nh(x_{i+1})(f^{x_{i+1}}-f^{x_i}), \\
\sum_{i=0}^{m}e_i(f^{y_{i+1}}-f^{y_i})&=&\sum_{i=0}^{m}h(y_{i+1})(f^{y_{i+1}}-f^{y_i}).
\dce
$$ 
Denote by respectively $G$ and $F$ the above two expressions (as functions on $(a,b]$). Then, for any points $s<t$ in some $(z_i,z_{i+1}]$, supposing $(z_i,z_{i+1}]\subset (x_k,x_{k+1}]\cap (y_l,y_{l+1}]$, we have $$
\dcb
h(t)=h(z_{i+1})=h(x_{k+1})=h(y_{l+1}),\\
(f^{x_{j+1}}-f^{x_j})|_s^t=(f^{y_{h+1}}-f^{y_h})|_s^t=0, \ \mbox{for $j\neq k, h\neq l$}.

\dce
$$ 
Hence, $$
G|_s^t=h(x_{k+1})(f^{x_{k+1}}-f^{x_k})|_s^t
=h(t)(f(t)-f(s))
=h(y_{l+1})(f^{y_{l+1}}-f^{y_l})|_s^t
=F|_s^t.
$$
Since $\lim_{s\downarrow a}F(s)=\lim_{s\downarrow a}G(s)=0$, $F$ and $G$ coincides.
\item
The bi-linearity is the consequence of the first property.
\item
By the bi-linearity, we need only to check the third property for $h=\ind_{(x,x']}$ and $g=\ind_{(y,y']}$, where $x<x',y<y'$ are points in $(a,b]$. We have
$$
\dcb
&&(g\ind_{(a,b]}\centerdot (h\ind_{(a,b]}\centerdot f))
=(g\ind_{(a,b]}\centerdot (f^{x'}-f^{x}))\\
&=&(f^{x'}-f^{x})^{y'}-(f^{x'}-f^{x})^{y}
=(f^{x'}-f^{x})^{y'\wedge x'}-(f^{x'}-f^{x})^{y\wedge x'}\\
&=&(f^{x'}-f^{x})^{(y'\wedge x')\vee x}-(f^{x'}-f^{x})^{(y\wedge x')\vee x}
=f^{(y'\wedge x')\vee x}-f^{x}-f^{(y\wedge x')\vee x}+f^{x}\\
&=&f^{(y'\wedge x')\vee x}-f^{(y\wedge x')\vee x}.
\dce
$$
If $y>x'$ or $x>y'$, the above function is identically null just as $gh$ and $gh\ind_{(a,b]}\centerdot f$ do. Otherwise, it is equal to $$
f^{y'\wedge x'}-f^{y\vee x}
=gh\ind_{(a,b]}\centerdot f.
$$
\item
This fourth property is clear for $h=\ind_{(s,t]}$ for $a\leq s<t\leq b$. By the bi-linearity, it is valid in general.
\item
It is clear for $h=\ind_{(s,t]}$. By the bi-linearity, it is true in general. 
\dbe
The lemma is proved. \ok

\subsubsection{Limit of elementary integrals on intervals}

An interval, which is open at the left end and closed at the right end, will be called a left interval. Let $\mathtt{B}_i=(a_i,b_i], i\in\mathbb{N}$, be a sequence of non empty left intervals. 
Look at the union set $\cup_{i\geq 0}\mathtt{B}_i$. One of the following situations hold for a $x\in\mathbb{R}_+$:
\ebe
\item[-]
$x$ is in the interior of one of $\mathtt{B}_i$ for $i\geq 0$.
\item[-]
$x$ is in the interior of $\cup_{i\geq 0}\mathtt{B}_i$, but it is in the interior of no of the $\mathtt{B}_i$ for $i\geq 0$. In this case, $x$ is the right end of one of $\mathtt{B}_i$. There exists a $\epsilon>0$ such that $(x,x+\epsilon)\subset\cup_{i\geq 0}\mathtt{B}_i$, and for any $i\geq 0$, either $\mathtt{B}_i\subset (0,x]$, or $\mathtt{B}_i\subset (x,\infty)$.
\item[-]
$x$ is in $\cup_{i\geq 0}\mathtt{B}_i$, but it is not in its interior. Then, $x$ is the right end of one of $\mathtt{B}_i$ and there exists a sequence of points in $(\cup_{i\geq 0}\mathtt{B}_i)^c$ decreasing to $x$.
\item[-]
$x$ is in $(\cup_{i\geq 0}\mathtt{B}_i)^c$.
\dbe
Consider the right end points $b_i$. A point $b_i$ will be said of first type if $a_j<b_i<b_j$ for some $j$. It is of second type if it is not of first type, but $b_i=a_j$ for some $j$. It is of third type if it is not of first neither of second type.

\bl\label{third-type}
Let $f$ be a càdlàg function on $\mathbb{R}$. Let $a<b$. Suppose that $(a,b]\subset\cup_{i\geq 0}\mathtt{B}_i$. Suppose that the family of right end points of third type has only a finite number of accumulation points in $(a,b)$. Suppose that the limit of $\ind_{(a,b]}\ind_{\cup_{0\leq i\leq n}\mathtt{B}_i}\centerdot f$ exists when $n\uparrow\infty$ with respect to the uniform norm on compact intervals. Then, the limit is simply equal to $\ind_{(a,b]}\centerdot f$.
\el

\textbf{Proof.}
We denoted by $\ind_{(a,b]}\ind_{\cup_{0\leq i<\infty}\mathtt{B}_i}\centerdot f$ the limit of $\ind_{(a,b]}\ind_{\cup_{0\leq i\leq n}\mathtt{B}_i}\centerdot f$ (which is well-defined by subsection \ref{eii}).
\ebe
\item
Suppose firstly that there exists no right end point of third type in $(a,b)$. Let $a<s<t<b$. Then, $[s,t]$ is contained in $(\cup_{i\geq 0}\mathtt{B}_i)^\circ$ (where the exponent $^\circ$ denotes the interior of a set). Note that in the case we consider here, the right end points are interior points of some $(\cup_{0\leq i\leq n}\mathtt{B}_i)^\circ, n\geq 0$. So, $[s,t]\subset\cup_{n\geq 0}(\cup_{0\leq i\leq n}\mathtt{B}_i)^\circ$. There exists therefore a $N>0$ such that $[s,t]\subset\cup_{0\leq n\leq N}(\cup_{0\leq i\leq n}\mathtt{B}_i)^\circ$. We have$$
\dcb
&&f_t-f_s=(\ind_{(s,t]}\centerdot f)_t\\
&=&(\ind_{(s,t]}\ind_{\cup_{0\leq i\leq N}\mathtt{B}_i}\centerdot f)_t
=\lim_{n\uparrow\infty}(\ind_{(s,t]}\ind_{\cup_{0\leq i\leq n}\mathtt{B}_i}\centerdot f)_t,\ \mbox{a stationary limit,}\\
&=&\lim_{n\uparrow\infty}(\ind_{(s,t]}\ind_{(a,b]}\ind_{\cup_{0\leq i\leq n}\mathtt{B}_i}\centerdot f)_t\\
&=&\lim_{n\uparrow\infty}(\ind_{(s,t]}\centerdot (\ind_{(a,b]}\ind_{\cup_{0\leq i\leq n}\mathtt{B}_i}\centerdot f))_t,\ \mbox{cf. Lemma \ref{elem-integral},}\\

&=&(\ind_{(s,b]}\centerdot (\ind_{(a,b]}\ind_{\cup_{0\leq i<\infty}\mathtt{B}_i}\centerdot f))_t,\ \mbox{by convergence assumption,}\\
&=&(\ind_{(s,b]}\ind_{\cup_{0\leq i<\infty}\mathtt{B}_i}\centerdot f)_t-(\ind_{(s,b]}\ind_{\cup_{0\leq i<\infty}\mathtt{B}_i}\centerdot f)_s.
\dce
$$
As $(\ind_{(s,b]}\ind_{\cup_{0\leq i<\infty}\mathtt{B}_i}\centerdot f)_s\rightarrow 0$ when $s\downarrow a$, we obtain$$
f_t-f_a=(\ind_{(s,b]}\ind_{\cup_{0\leq i<\infty}\mathtt{B}_i}\centerdot f)_t, \forall a<t<b.
$$
Now to obtain the result stated in the lemma, we need only to check that $\Delta_b(\ind_{(s,b]}\ind_{\cup_{0\leq i<\infty}\mathtt{B}_i}\centerdot f)=\Delta_b f$. Notice that there exists a $N>0$ such that $b\in\cup_{0\leq i<n}\mathtt{B}_i$ for all $n\geq N$. We have$$
\dcb
\Delta_b(\ind_{(s,b]}\ind_{\cup_{0\leq i<\infty}\mathtt{B}_i}\centerdot f)
&=&\lim_{n\uparrow\infty}\Delta_b(\ind_{(s,b]}\ind_{\cup_{0\leq i\leq n}\mathtt{B}_i}\centerdot f)\\
&=&\lim_{n\uparrow\infty}\ind_{b\in\cup_{0\leq i\leq n}\mathtt{B}_i}\Delta_b f\ \mbox{(cf. Lemma \ref{elem-integral})}
=\Delta_b f.
\dce
$$
The lemma is proved when no right end point of third type exists.

\item
There exist a finite number of right end points of third type in $(a,b)$. Let $v_1<v_2<\ldots<v_k$ are the right end points of third type in $(a,b)$. Applying the preceding result, $$
\dcb
\lim_{n\uparrow\infty}\ind_{(a,b]}\ind_{\cup_{0\leq i\leq n}\mathtt{B}_i}\centerdot f
&=&\sum_{j=0}^k\lim_{n\uparrow\infty}\ind_{(v_j,v_{j+1}]}\ind_{\cup_{0\leq i\leq n}\mathtt{B}_i}\centerdot f \ (v_0=a,v_{k+1}=b)\\

&=&\sum_{j=0}^k\ind_{(v_j,v_{j+1}]}\centerdot f
=\ind_{(a,b]}\centerdot f.
\dce
$$
The lemma is true in this second case.
\item
There exist an infinite number of right end points of third type in $(a,b)$, but $b$ is the only accumulation point of these right end points of third type. We have, for $a<t<b$,
$$
\dcb
(\ind_{(a,b]}\ind_{\cup_{0\leq i< \infty}\mathtt{B}_i}\centerdot f)_t
&=&(\lim_{n\uparrow\infty}\ind_{(a,b]}\ind_{\cup_{0\leq i\leq n}\mathtt{B}_i}\centerdot f)_t\\
&=&(\lim_{n\uparrow\infty}\ind_{(a,t]}\ind_{\cup_{0\leq i\leq n}\mathtt{B}_i}\centerdot f)_t
=(\ind_{(a,t]}\centerdot f)_t
=(\ind_{(a,b]}\centerdot f)_t.
\dce
$$
As before, the two functions has the same jumps at $b$. The lemma is true in this third case.

\item
There exist an infinite number of right end points of third type in $(a,b)$, but $a$ is the only accumulation point of these right end points of third type. Let $a<s$.
$$
\dcb
&&(\ind_{(a,b]}\ind_{\cup_{0\leq i< \infty}\mathtt{B}_i}\centerdot f)-(\ind_{(a,b]}\ind_{\cup_{0\leq i< \infty}\mathtt{B}_i}\centerdot f)^s\\
&=&
\ind_{(s,b]}\centerdot(\ind_{(a,b]}\ind_{\cup_{0\leq i< \infty}\mathtt{B}_i}\centerdot f)
=\lim_{n\uparrow\infty}\ind_{(s,b]}\centerdot(\ind_{(a,b]}\ind_{\cup_{0\leq i\leq n}\mathtt{B}_i}\centerdot f)\\
&=&\lim_{n\uparrow\infty}\ind_{(s,b]}\ind_{(a,b]}\ind_{\cup_{0\leq i\leq n}\mathtt{B}_i}\centerdot f
=\lim_{n\uparrow\infty}\ind_{(s,b]}\ind_{\cup_{0\leq i\leq n}\mathtt{B}_i}\centerdot f

\ =\ \ind_{(s,b]}\centerdot f.

\dce
$$
Since $(\ind_{(a,b]}\ind_{\cup_{0\leq i< \infty}\mathtt{B}_i}\centerdot f)^s$ tends uniformly  to zero when $s\downarrow a$, the lemma is true in this fourth case.

\item
There exist an infinite number of right end points of third type in $(a,b)$, but $a,b$ are the only accumulation points of these right end points of third type. We have
$$
\dcb
\lim_{n\uparrow\infty}\ind_{(a,b]}\ind_{\cup_{0\leq i\leq n}\mathtt{B}_i}\centerdot f
&=&\lim_{n\uparrow\infty}\ind_{(a,a+\frac{b-a}{2}]}\ind_{\cup_{0\leq i\leq n}\mathtt{B}_i}\centerdot f
+\lim_{n\uparrow\infty}\ind_{(a+\frac{b-a}{2},b]}\ind_{\cup_{0\leq i\leq n}\mathtt{B}_i}\centerdot f\\

&=&\ind_{(a,a+\frac{b-a}{2}]}\centerdot f+\ind_{(a+\frac{b-a}{2},b]}\centerdot f
=\ind_{(a,b]}\centerdot f.
\dce
$$

\item
There exist an infinite number of right end points of third type in $(a,b)$, but there exist only a finite number of accumulation point of these right end points of third type in $(a,b)$. Let $v_1<v_2<\ldots<v_k$ be the accumulation points in $(a,b)$. Applying the preceding result, $$
\dcb
\lim_{n\uparrow\infty}\ind_{(a,b]}\ind_{\cup_{0\leq i\leq n}\mathtt{B}_i}\centerdot f
&=&\sum_{j=0}^k\lim_{n\uparrow\infty}\ind_{(v_j,v_{j+1}]}\ind_{\cup_{0\leq i\leq n}\mathtt{B}_i}\centerdot f \ (v_0=a,v_{k+1}=b)\\

&=&\sum_{j=0}^k\ind_{(v_j,v_{j+1}]}\centerdot f
=\ind_{(a,b]}\centerdot f.
\dce
$$
The lemma is proved. \ok
\dbe

\

\

\subsection{Properties of processes satisfying Assumption \ref{partialcovering}}

We consider a process $S$ satisfying Assumption \ref{partialcovering}.

\bl\label{UB}
(We assume only the two first conditions in Assumption \ref{partialcovering}.) Suppose that the random measure $\mathsf{d}\chi^\cup$ has a distribution function $\chi^\cup$. Then, $\ind_{\cup_{0\leq i\leq n}\mathtt{B}_i}\centerdot S, n\geq 0,$ converges in $\mathsf{H}^1_{\mathsf{loc}}$ to a semimartingale that we denote by $\ind_{\cup_{i\geq 0}\mathtt{B}_i}\centerdot S$. This semimartingale is special whose drift is $\chi^\cup$. 
\el

\textbf{Proof.} 
Remark that, under Assumption \ref{partialcovering}, denoting by $S^c$ the continuous martingale part of $S$ in its natural filtration, the bracket process $$
[S,S]_t=\cro{S^c,S^c}_t+\sum_{s\leq t}(\Delta_s S)^2, \ t\geq 0,
$$
is a well-defined increasing process in the natural filtration of $S$. The process $\sqrt{[S,S]}$ is locally integrable.

We notice that $\ind_{\cup_{0\leq i\leq n}\mathtt{B}_i}$ is a left elementary process. By Lemma \ref{third-type} and the inclusion-exclusion formula, $\ind_{\cup_{0\leq i\leq n}\mathtt{B}_i}\centerdot S$ is a special semimartingale with bracket process $\ind_{\cup_{0\leq i\leq n}\mathtt{B}_i}\centerdot [S,S]$ (computed pathwisely as quadratic variation). We use \cite[Corollaire(1.8)]{Jeulin80} to control the martingale part of $\ind_{\cup_{0\leq i\leq n}\mathtt{B}_i}\centerdot S$ by $\sqrt{\ind_{\cup_{0\leq i\leq n}\mathtt{B}_i}\centerdot [S,S]}$, and the drift part of $\ind_{\cup_{0\leq i\leq n}\mathtt{B}_i}\centerdot S$ by $\chi^\cup$. The two control processes are locally integrable. \ok

\bex
How to construct the sequence of stopping times $(T_{m})_{m\in\mathbb{N}}$ to realise the converges in $\mathcal{H}^1_{\mathsf{loc}}$ ?
\eex

\

We introduce the process
$g_{t} = \sup\{0\leq s < t; s \notin \cup_{n\geq 0}\mathtt{B}_n\},\ t > 0$ (when the set is empty, $g_{t}=0$). It is an increasing left continuous $\mathbb{G}$ adapted process, i.e. a predictable process. For any
$\epsilon > 0$, set
$\mathtt{A}_{\epsilon} = \{t\geq 0: t - g_{t} > \epsilon\}$. We can check that $\mathtt{A}_{\epsilon} \subseteq \mathtt{A}_{\delta}$ for $\delta<\epsilon$ and $\mathtt{A}=\cup_{\epsilon>0}\mathtt{A}_\epsilon$. Note that $0\notin \cup_{i\geq 0}\mathtt{B}_i$. 

Recall, for $\mathbb{G}$ stopping time $R$, $d_R=\inf\{s\geq R: s\notin \cup_{i\geq 0}\mathtt{B}_i\}$. Either $d_{R}$ belongs or not to $\cup_{i\geq 0}\mathtt{B}_i$. But always $S_{d_R}=\check{S}_{d_R}$ if $d_{R}<\infty$.  Define successively for every $n\geq 0$ the $\mathbb{G}$ stopping times: ($d_{R_{-1}}=0$)$$
\dcb
R_{n}=\inf\{s\geq d_{R_{n-1}}: s\in \mathtt{A}_\epsilon\}\ n\geq 0.
\dce
$$
For any $s\in\mathtt{A}_\epsilon$, there exists a $s'<s$ such that $[s',s]\subset\mathtt{A}_\epsilon$. Therefore, $R_n\notin \mathtt{A}_\epsilon$ and $R_n<d_{R_n}$ if $R_n<\infty$. Moreover, $R_n\in \cup_{i\geq 0}\mathtt{B}_i$ and $d_{R_n}-g_{d_{R_n}}> \epsilon$ and $(d_{R_n},d_{R_n}+\epsilon)\cap\mathtt{A}_\epsilon=\emptyset$. Consequently $d_{R_{n}}+\epsilon\leq R_{n+1}$ and $\lim_{k\rightarrow\infty}R_k=\infty$. We can write $\mathtt{A}_\epsilon=\cup_{n\geq 1} (R_{n},d_{R_{n}}]$ and hence $\ind_{\mathtt{A}_\epsilon}$ is a left elementary process on any finite interval. 

For a càdlàg process $X$, set $\mathtt{j}(X)=\{t>0:X_{t-}>0>X_{t} \mbox{ or }X_{t-}\leq 0<X_{t}\}$. We introduce the process (eventually taking infinite values)	
$$
\dcb
A_t&=&\sum_{0<s\leq t}\ind_{\{s\in\mathtt{A}\}}\left[\ind_{(S-\check{S})_{s-}>0}(S-\check{S})^-_{s}+\ind_{(S-\check{S})_{s-}\leq 0}(S-\check{S})^+_{s}\right]\\
&=&\sum_{0<s\leq t}\ind_{\{s\in\mathtt{A}\}} \ind_{\mathtt{j}(S-\check{S})}(s)\ |S-\check{S}|_{s}\ t\geq 0.
\dce
$$

\bl\label{BU+-}
Suppose that the distribution function $\chi^\cup$ exists and the process $C$ is the jump process of a special semimartingale. Then, the process $A$ is locally integrable, and the three $\mathbb{G}$ semimartingales $\ind_{\mathtt{A}_\epsilon} \centerdot (S-\check{S})$, $\ind_{\mathtt{A}_\epsilon}\centerdot (S-\check{S})^+$ and $\ind_{\mathtt{A}_\epsilon}\centerdot (S-\check{S})^-$, $n\geq 0,$ converge in $\mathcal{H}^1_{\mathsf{loc}}$ to semimartingales that we denote respectively by $\ind_{\mathtt{A}} \centerdot (S-\check{S})$, $\ind_{\mathtt{A}}\centerdot (S-\check{S})^+$ and $\ind_{\mathtt{A}}\centerdot (S-\check{S})^-$. They are special semimartingales. We have $$
\dcb
\ind_{\mathtt{A}}\centerdot (S-\check{S})^+
=\ind_{\{(S-\check{S})_{-}>0\}}\ind_{\mathtt{A}} \centerdot (S-\check{S}) + \ind_{\mathtt{A}}\centerdot A+\frac{1}{2}\ind_{\mathtt{A}} \centerdot l^\cup\\
\\
\ind_{\mathtt{A}}\centerdot (S-\check{S})^-
=-\ind_{\{(S-\check{S})_{-}\leq 0\}}\ind_{\mathtt{A}} \centerdot (S-\check{S}) + \ind_{\mathtt{A}}\centerdot A+\frac{1}{2}\ind_{\mathtt{A}} \centerdot l^\cup
\dce
$$
for some continuous increasing process $l^\cup$ null at the origin which increases only on the set $\{t\in\mathtt{A}: (S-\check{S})_{t-}=0\}$
\el

\textbf{Proof.} 
Let the process $C$ be the jump process of a special semimartingale $M+V$, where $M$ is a $\mathbb{G}$ local martingale and $V$ is a $\mathbb{G}$ predictable process with finite variation. By stopping $S,\check{S},M,V$ if necessary, we assume that $(S-\check{S})^*_\infty$ and $M^*$ are integrable, $(S-\check{S})_\infty$ exists, and the total variations of $\chi^\cup$ and of $V$ are integrable.

As the random measure $\mathsf{d}\chi^\cup$ has a distribution function $\chi^\cup$, the random measure $\ind_{(R,d_R]}\centerdot\mathsf{d}\chi^\cup$ also has a distribution function. The following limit computation is valid in $\mathsf{H}^1$: 
$$
\dcb
&&\ind_{(R,d_R]}\centerdot(\ind_{\cup_{i\geq 0}\mathtt{B}_i}\centerdot (S-\check{S}))\\
&=&\lim_{n\uparrow\infty}\ind_{(R,d_R]}\centerdot(\ind_{\cup_{0\leq i\leq n}\mathtt{B}_i}\centerdot (S-\check{S}))\\
&=&\lim_{n\uparrow\infty}(\ind_{(R,d_R]}\ind_{\cup_{0\leq i\leq n}\mathtt{B}_i})\centerdot (S-\check{S}).
\dce
$$
Notice that the limit in $\mathsf{H}^1$ implies the limit with respect to the uniform convergence in compact intervals in probability. Lemma \ref{third-type} is applicable. For any $R<t<d_R$,$$
(\ind_{(R,d_R]}\centerdot(\ind_{\cup_{i\geq 0}\mathtt{B}_i}\centerdot (S-\check{S})))_t
=(\ind_{(R,d_R]}\centerdot (S-\check{S}))_t
$$
With Lemma \ref{elem-integral}, we compute the jumps at $d_R$ : $$
\Delta_{d_R}(\ind_{(R,d_R]}\centerdot(\ind_{\cup_{i\geq 0}\mathtt{B}_i}\centerdot (S-\check{S})))
=
\ind_{\{R<d_R\}}\ind_{\{d_R\in \cup_{i\geq 0}\mathtt{B}_i\}}\Delta_{d_R}(S-\check{S}).
$$
These facts implies
$$
\dcb
\ind_{(R,d_R]}\centerdot (S-\check{S})=\ind_{(R,d_R]}\centerdot(\ind_{\cup_{i\geq 0}\mathtt{B}_i}\centerdot (S-\check{S}))+\ind_{\{R<d_R\}}\ind_{\{d_R\notin \cup_{i\geq 0}\mathtt{B}_i\}}\Delta_{d_R}(S-\check{S})\ind_{[d_R,\infty)}.
\dce
$$
Note that $
[(d_{R})_{\{R<d_R,d_R\notin \cup_{i\geq 0}\mathtt{B}_i\}}]
=(R,d_R]\setminus \cup_{i\geq 0}\mathtt{B}_i
$
is a $\mathbb{G}$ predictable set. This means that the $\mathbb{G}$ stopping time $(d_{R})_{\{R<d_R,d_R\notin \cup_{i\geq 0}\mathtt{B}_i\}}$ is a predictable stopping time. Let $\hat{d}_R$ denote this predictable stopping time. Consider the following jump process:
$$
\dcb
&&\ind_{\{R<d_R\}}\ind_{\{d_R\notin \cup_{i\geq 0}\mathtt{B}_i\}}\Delta_{d_R}(S-\check{S})\ind_{[d_R,\infty)}\\
&=&\ind_{\{R<d_R\}}\ind_{\{d_R\notin \cup_{i\geq 0}\mathtt{B}_i\}}C_{d_R}\ind_{[d_R,\infty)}\\

&=&\ind_{\{R<d_R\}}\ind_{\{d_R\notin \cup_{i\geq 0}\mathtt{B}_i\}}\Delta_{d_R}(M+V)\ind_{[d_R,\infty)}\\

&=&\Delta_{\hat{d}_{R}}(M+V)\ind_{[\hat{d}_{R},\infty)}\\

&=&\ind_{[\hat{d}_{R}]}\centerdot(M+V)\\
\dce
$$
Combining this equation with Lemma \ref{UB}, we see that $\ind_{(R,d_R]}\centerdot (S-\check{S})$ is a $\mathbb{G}$ special semimartingale whose drift is given by $\ind_{(R,d_R]}\centerdot\chi^\cup+\ind_{[\hat{d}_{R}]}\centerdot V$. 

Applying now the argument of the proof of Lemma \ref{UB}, we see that $\ind_{\mathtt{A}_\epsilon}\centerdot (S-\check{S})$ converges in $\mathsf{H}^1$ to a special semimartingale $\ind_{\mathtt{A}}\centerdot (S-\check{S})$.

Recall, for $t>0,$ $$
\ind_{(R,d_R]}\centerdot A_t = \sum_{R<s\leq d_R\wedge t} [1_{\{(S-\check{S})_{s-}>0\}}(S-\check{S})_{s}^{-}+1_{\{(S-\check{S})_{s-}\leq 0\}}(S-\check{S})_{s}^{+}]
= \sum_{R<s\leq d_R\wedge t,s\in\mathtt{j}(S-\check{S})} |S-\check{S}|_{s}.
$$
According to \cite[Chapter 9 \S 6]{HWY} applied to the semimartingale $(S-\check{S})^{d_R}\ind_{[R,\infty)}=\ind_{(R,d_R]}\centerdot (S-\check{S})+(S_R-\check{S}_R)\ind_{[R,\infty)}$,  
we know that the process 
$$
l^{(R,d_R]}_{t}(X) = 2[\ind_{(R,d_R]}\centerdot (S-\check{S})^{+}_t-\ind_{\{(S-\check{S})_{-}>0\}}\ind_{(R,d_R]} \centerdot (S-\check{S})_{t} - \ind_{(R,d_R]}\centerdot A_{t} ], \ t>0,
$$
is not decreasing, continuous and null at the origin, which increases only on the set $\{t\in(R,d_R]:(S-\check{S})_{t-}=0\}$. Note that, if we take another $\mathbb{G}$ stopping time $R'$, the random measure $\mathsf{d}l^{(R',d_{R'}]}$ coincide with the random measure $\mathsf{d}l^{(R,d_{R}]}$ on $(R,d_{R}]\cap(R',d_{R'}]$. Therefore, there exists a random measure $\mathsf{d}l^\cup$ on $\mathtt{A}$ such that $\ind_{(R,d_R]}\centerdot \mathsf{d}l^\cup=\mathsf{d}l^{(R,d_R]}$. (Note that $\mathsf{d}l^\cup$ is diffuse which does not charge $\mathtt{C}$.)

We have the following computation for any $\epsilon>0$:
$$
\dcb
&&\mathbb{E}[\int\ind_{\mathtt{A}_\epsilon}(s) (dl^\cup_s+2dA_s)]\\
&=&\mathbb{E}[\int\ind_{\cup_{n\geq 0}(R_n,d_{R_n}]}(s) (dl^\cup_s+2dA_s)]\\
&=&\lim_{k\uparrow\infty}\sum_{0\leq n\leq k}\mathbb{E}[l^{(R_n,d_{R_n}]}_{d_{R_n}}+2\ind_{(R_n,d_{R_n}]}\centerdot A_{d_{R_n}}]\\
&=&\lim_{k\uparrow\infty}\sum_{0\leq n\leq k}2\mathbb{E}[(S-\check{S})_{d_{R_n}}^{+}-(S-\check{S})_{R_n}^{+}-\ind_{\{(S-\check{S})_{-}>0\}}\ind_{(R_n,d_{R_n}]} \centerdot (S-\check{S})_{d_{R_n}} ]\\
&\leq&\lim_{k\uparrow\infty}\sum_{0\leq n\leq k}2\mathbb{E}[(S-\check{S})_{d_{R_n}}^{+}\ind_{\{R_n<\infty,d_{R_n}=\infty\}}+\ind_{\{(S-\check{S})_{-}>0\}}(\ind_{(R_n,d_{R_n}]} \centerdot |d\chi^\cup|_\infty +\ind_{[\hat{d}_{R}]}\centerdot |dV|_\infty)]\\
&&\mbox{because $S_{d_R}=\check{S}_{d_R}$ if $d_{R}<\infty$,}\\
&\leq&2\mathbb{E}[(S-\check{S})_{\infty}^{+}\sum_{n\geq 0}\ind_{\{R_n<\infty,d_{R_n}=\infty\}}]+2\mathbb{E}[\int_0^\infty (|d\chi^\cup_s|+|dV_s|) ]\\
&\leq&2\mathbb{E}[(S-\check{S})^*_{\infty}]+2\mathbb{E}[\int_0^\infty (|d\chi^\cup_s|+|dV_s|) ]<\infty\\
\dce
$$
Here the last inequality is because there exist only one $n$ such that $R_n<\infty,d_{R_n}=\infty$. Let $\epsilon\rightarrow 0$. $\ind_{\mathtt{A}_\epsilon}(s)$ tends to $\ind_{\mathtt{A}}(s)$. We conclude that $
\mathbb{E}[\int_0^\infty \ind_{\mathtt{A}}(s) (dl^\cup_s+dA_s)]<\infty.
$
That means that $A$ is finite valued and $\mathsf{d}l^\cup$ have a distribution functions $l^\cup$. 

It is now straightforward to see that $
\ind_{\mathtt{A}_\epsilon}\centerdot (S-\check{S})^+
$
converge in $\mathcal{H}^1$ to a limit that we denote by $\ind_{\mathtt{A}}\centerdot (S-\check{S})^+$, which is equal to
$$
\ind_{\mathtt{A}}\centerdot (S-\check{S})^+
=\ind_{\{(S-\check{S})_{-}>0\}}\ind_{\mathtt{A}} \centerdot (S-\check{S}) + \ind_{\mathtt{A}}\centerdot A+\frac{1}{2}\ind_{\mathtt{A}} \centerdot l^\cup
$$ 
The first part of the lemma is proved.

The other part of the lemma can be proved similarly. Notice that we obtain the same random measure $\mathsf{d}l^\cup$ in the decomposition of $(S-\check{S})^-$\ok

\bethe\label{Z-Ytheorem}
Suppose that the distribution function $\chi^\cup$ exists and the process $C$ is the jumps process of a $\mathbb{G}$ special semimartingale. Then, $S$ is a special semimartingale. More precisely, let \[
\begin{array}{lcl}
V^{+}&=&(S-\check{S})^{+}-(S-\check{S})_0^{+}-\ind_{\mathtt{A}}\centerdot(S-\check{S})^{+}\\
V^{-}&=&(S-\check{S})^{-}-(S-\check{S})_0^{-}-\ind_{\mathtt{A}}\centerdot(S-\check{S})^{-}
\end{array}
\]
$V^+, V^-$ are non decreasing locally integrable processes. They increase only outside of $\mathtt{A}$. We have :\[
\begin{array}{lll}
S&=&S_{0}+\ind_{\mathtt{A}}\centerdot (S-\check{S})+V^{+}-V^{-}+ \check{S}.
\end{array}
\]
\ethe

\textbf{Proof.} Let $X=(S-\check{S})^{+}$ and $X'=(S-\check{S})^{-}$. By stopping, we suppose that $X,X'$ are in class$(D)$, as well as $\ind_{\mathtt{A}}\centerdot X$. For a $\epsilon>0$, we compute the
difference :
\begin{eqnarray*}
	&&X - X_{0}-\ind_{\mathtt{A}_{\epsilon}}\centerdot X = X - X^{0} - \sum_{n\geq 0} (X^{d_{R_{n}}}-X^{R_{n}})\\
	&=& \sum_{n\geq 0} (X^{d_{R_{n}}}-X^{R_{n}})
	+\sum_{n\geq 0} (X^{R_{n}}-X^{d_{R_{n-1}}}) - \sum_{n\geq 1} (X^{d_{R_{n}}}-X^{R_{n}})\\
	&=& \sum_{n\geq 0} (X^{R_{n}}-X^{d_{R_{n-1}}})
\end{eqnarray*}
According to Lemma \ref{BU+-} $\ind_\mathtt{A}\centerdot X$ exists as the limit in $\mathcal{H}^1_{loc}$ of $\ind_{\mathtt{A}_{\epsilon}}\centerdot X$. Let $\epsilon$ tend to zero. The first term of the above identity tends to a process $V^+=X- X_{0} - \ind_{\mathtt{A}}\centerdot X$, uniformly on every compact interval in probability. In particular, $V^+$ is càdlàg and is in class$(D)$.

Consider the last term of the above identity. For $t>0$, let $N(t) = \sup\{k\geq 1; d_{R_{k}} \leq t\}$ ($N(t) = 0$ if the set is empty). Recall that, on the set $\{d_{R_{j}}<\infty\}$, $X_{d_{R_{j}}} = 0$. We have  
$$
\dcb
&&\sum_{n\geq 0} (X^{R_{n}}-X^{d_{R_{n-1}}})_t\\

&=&\sum_{n= 0}^{N(t)+1} X^{R_{n}}_{t}\\
&=&X^{R_{0}}_{t}+X^{R_{1}}_{t}+X^{R_{2}}_{t}+X^{R_{3}}_{t}+\ldots+X^{R_{N(t)}}_{t}+X^{R_{N(t)+1}}_{t}\\
&&\mbox{recall that $d_{R_{N(t)}}\leq t$ so that $R_{N(t)}\leq t$}\\
&=&X_{R_{0}}+X_{R_{1}}+X_{R_{2}}+X_{R_{3}}+\ldots+X_{R_{N(t)}}+X^{R_{N(t)+1}}_t\\
&=&\sum_{n=0}^{N(t)} X_{R_{n}} + X_{t}1_{\{d_{R_{N(t)}}\leq t < R_{N(t)+1}\}}+ X_{R_{N(t)+1}}1_{\{ R_{N(t)+1}\leq t\}}
\dce
$$ 
Notice that $X_{t}1_{\{d_{R_{N(t)}}\leq t < R_{N(t)+1}\}}=X_{t}1_{\{d_{R_{N(t)}}< t < R_{N(t)+1}\}}$ because $X_{d_{R_{N(t)}}}=0$. If $t \in \mathtt{A}$, for $\epsilon$ small enough, $t$
will belongs to
$\mathtt{A}_{\epsilon}$. As the interval $(d_{R_{N(t)}},R_{N(t)+1})$ is contained in the complementary of $\mathtt{A}_{\epsilon}$, we must have  $\ind_{\{d_{R_{N(t)}}< t
	\leq R_{N(t)+1}\}} = 0$. If $t \notin \mathtt{A}$, by Assumption \ref{partialcovering}, $X_{t} = 0$. In sum, for every $t>0$, there is a $\epsilon(t)$ such that, for $\epsilon < \epsilon(t)$, 
$$
\sum_{n\geq 0} (X^{R_{n}}-X^{d_{R_{n-1}}})_t
=\sum_{n=0}^{N(t)} X_{R_{n}} + X_{R_{N(t)+1}}\ind_{\{ R_{N(t)+1}\leq t\}}
=\sum_{n\geq 0} X_{R_{n}}\ind_{\{R_{n}\leq t\}}
$$
From this expression, we can write, for $0<s<t$,
$$
V^{+}_{t} - V^{+}_{s} 
=  \lim_{\epsilon \downarrow 0} (\sum_{n\geq 0} X_{R_{n}}\ind_{\{ R_{n}\leq t\}}-\sum_{n\geq 0} X_{R_{n}}\ind_{\{ R_{n}\leq s\}})
=  \lim_{\epsilon \downarrow 0} \sum_{n\geq 0} X_{R_{n}}\ind_{\{s<  R_{n}\leq t\}} \geq 0,
$$
i.e. $V^+$ is an increasing process. Moreover, for a fixed $a>0$, since $\mathtt{A}_a\subset \mathtt{A}_\epsilon$ for any $\epsilon<a$, since $R_n\notin\mathtt{A}_\epsilon$, the above expression implies $\ind_{\mathtt{A}_a}\centerdot V^+=0$. This argument shows that the random measure $\mathtt{d}V^+$ does not charge $\mathtt{A}$.  

According to Lemma \ref{BU+-} 
$$
\ind_{\mathtt{A}}\centerdot X =\ind_{\{(S-\check{S})_{-}>0\}}\ind_{\mathtt{A}} \centerdot (S-\check{S}) +\ind_{\mathtt{A}}\centerdot A+\frac{1}{2}\ind_{\mathtt{A}}\centerdot l^\cup
$$
and, therefore,
$$
X=X_0+\ind_{\{(S-\check{S})_{-}>0\}}\ind_{\mathtt{A}} \centerdot (S-\check{S}) +\ind_{\mathtt{A}}\centerdot A+\frac{1}{2}\ind_{\mathtt{A}}\centerdot l^\cup+V^+
$$
In the same way we prove$$
X'=X'_0-\ind_{\{(S-\check{S})_{-}\leq 0\}}\ind_{\mathtt{A}} \centerdot (S-\check{S}) +\ind_{\mathtt{A}}\centerdot A+\frac{1}{2}\ind_{\mathtt{A}}\centerdot l^\cup+V^-
$$
Writing finally $S=X-X'+\check{S}$, we prove the theorem.  \ok

\

\section{Local solution method : A method to find local solutions}\label{lsmfind}

In the previous section we have discussed the step $_{\bullet\bullet}$ of the local solution method. Especially, according to Theorem \ref{semimartingaleproperty}, for a given process $X$, if $X$ is $\mathbb{G}$ semimartingale on a rich collect of random intervals, $X$ will be a $\mathbb{G}$ semimartingale on the whole $\mathbb{R}_{+}$. In this section, we discuss the step $_{\bullet}$ of the local solution method. We establish a canonical method to find local solutions.

\subsection{The functions $h^u, u\geq 0$}\label{h-formula}

We equip now the probability space $(\Omega,\mathcal{A},\mathbb{P})$ with a second filtration $\mathbb{F}$ contained in the filtration $\mathbb{G}$. We need some topological properties behind the filtrations. We assume that $\Omega$ is a Polish space and $\mathcal{A}$ is its borel $\sigma$-algebra. We assume that there exists a filtration $\mathbb{F}^\circ=(\mathcal{F}^\circ_t)_{t\geq 0}$ of sub-$\sigma$-algebras of $\mathcal{A}$ such that $$
\mathcal{F}_t=\mathcal{F}^\circ_{t+}\vee\mathcal{N}, \ t\geq 0,
$$ 
where $\mathcal{N}$ is the family of $(\mathbb{P},\mathcal{A})$-negligible sets. We assume that there exists a measurable map $I$ from $\Omega$ into another polish space $\mathtt{E}$ equipped with its borel $\sigma$-algebra and a filtration $\mathbb{I}=(\mathcal{I}_t)_{t\geq 0}$ of countably generated sub-Borel-$\sigma$-algebras, such that $\mathcal{G}_{t}=\mathcal{G}^\circ_{t+}\vee \mathcal{N}, \ t\geq 0,$ where
$
\mathcal{G}^\circ_{t}=\mathcal{F}^\circ_{t} \vee I^{-1}(\mathcal{I}_t), t\geq 0.
$
We recall the following result (see \cite[Theorem 4.36]{HWY})

\bl\label{FF0}
For any $\mathbb{F}$ stopping time $T$, there exists a $\mathbb{F}^\circ_+$ stopping time $T^\circ$ such that $T=T^\circ$ almost surely. For any $A\in\mathcal{F}_T$, there exists a $A^\circ\in\mathcal{F}^\circ_{T^\circ+}$ such that $A\Delta A^\circ$ is negligible. For any $\mathbb{F}$ optional (resp. predictable) process $X$, there exists a $\mathbb{F}^\circ_+$ optional (resp. predictable) process $X^\circ$ such that $X$ and $X^\circ$ are indistinguishable. 

A similar result holds between the couple $\mathbb{G}$ and $\mathbb{G}^\circ_+=(\mathcal{G}^\circ_{t+})_{t\geq 0}$. 
\el

\brem
As a consequence of Lemma \ref{FF0}, we can speak about, for example, the predictable dual projection of an integrable increasing process with respect to the filtration $\mathbb{G}^\circ_+=(\mathcal{G}^\circ_{t+})_{t\geq 0}$ : this will mean that we compute the predictable dual projection in $\mathbb{G}$, then we take a version in $\mathbb{G}^\circ_+$.
\erem

Consider the product space $\Omega\times \mathtt{E}$ equipped with its product $\sigma$-algebra and its product filtration $\mathbb{J}$ composed of 
$$
\mathcal{J}_t=\sigma\{A\times B: A\in\mathcal{F}^\circ_t, \ B\in\mathcal{I}_t\}, t\geq 0.
$$
We introduce the map $\phi$ on $\Omega$ such that
$\phi(\omega) = (\omega, I(\omega)) \in \Omega \times \mathtt{E}$. Notice that, for $t\geq 0$, $\mathcal{G}^\circ_t=\phi^{-1}(\mathcal{J}_t)$. Therefore, for $C\in\mathcal{G}^\circ_{t+}$, there exist a sequence $(D_n)_{n\geq 1}$ of sets in respectively $\mathcal{J}_{t+\frac{1}{n}}$ such that $C=\phi^{-1}(D_n)$, which means $C=\phi^{-1}(\cap_{n\geq 1}\cup_{k\geq n}D_k)\in \phi^{-1}(\mathcal{J}_{t+})$. This observation yields the equality $\mathcal{G}^\circ_{t+}=\phi^{-1}(\mathcal{J}_{t+}), t\geq 0$ (cf. Lemma \ref{intersection}).

We equip the product space $\Omega \times \mathtt{E}$ with the the image probability $\mu$ of $\mathbb{P}$ by $\phi$. We introduce the identity map $\mathsf{i}$ on $\Omega$, the projection maps $\zeta(\omega,x)=x$ and $\iota(\omega,x)=\omega$ for $(\omega,x)\in\Omega\times \mathtt{E}$. For $t\geq 0$, let  
$$
\begin{array}{lll}
\pi_{t,F/F}(\omega, d\omega')
&=&
\mbox{regular conditional distribution of the map $\mathsf{i}$ under $\mathbb{P}$}
\\
&&\mbox{given the map $\mathsf{i}$ (itself) as a map valued in $(\Omega,\mathcal{F}^\circ_t)$}\\

\pi_{t,I/F}(\omega, dx')
&=&
\mbox{regular conditional distribution of $I$ under $\mathbb{P}$}
\\
&&\mbox{given the map $\mathsf{i}$ as a map valued in $(\Omega,\mathcal{F}^\circ_t)$}\\

\pi_{t,I/I}(x,dx')
&=&
\mbox{regular conditional distribution of $I$ under $\mathbb{P}$}
\\
&&\mbox{given $I$ (itself) as a map valued in $(\mathtt{E},\mathcal{I}_t)$}\\
\end{array}
$$
(See \cite{freedman, KS, RW} or subsection \ref{rcd} for the notion of regular conditional distribution.)

\brem
There exist situations where $\mathbb{F}^\circ$ is generated by a Borel map $Y$ from $\Omega$ into a polish space $\mathtt{F}$ equipped with a filtration $\mathbb{K}=(\mathcal{K}_t)_{t\geq 0}$ of sub-Borel-$\sigma$-algebras on $\mathtt{F}$, such that $\mathcal{F}^\circ_t=Y^{-1}(\mathcal{K}_t), t\geq 0$. Let $\tilde{\pi}_{t}(y,d\omega')$ be the regular conditional distribution of $\mathsf{i}$ under $\mathbb{P}$ given the map $Y$ considered as a map valued in the space $(\mathtt{F},\mathcal{K}_t)$. Then, $
\tilde{\pi}_t(Y(\omega),d\omega')
$
is a version of $\pi_{t,F/F}(\omega,d\omega')$, and its image measure by the map $I$ on $E$ is a version of $\pi_{t,I/F}(\omega,dx')$. See subsection \ref{rcd}.
\erem

For a fixed $u \geq 0$, for any $t \geq u,$ let $h_{t}^{u}(\omega, x)$
to be a function on $\Omega \times \mathtt{E}$ which is $[0,2]$ valued and ${\mathcal{J}}_{t}$-measurable, such that, for $\mathbb{P}$-almost all $\omega$,
$$
h_{t}^{u}(\omega, x) = \left. {\frac{2 \pi_{t,I/F}(\omega, dx)}
	{\pi_{t,I/F}(\omega, dx) +
		\int \pi_{u,I/F}(\omega, dx'') \pi_{u,I/I}(x'', dx)}}\right|
_{\mathcal{I}_{t}}.
$$
We introduce a family of probabilities $\nu^{u}$ on $\mathcal{J}_\infty$
indexed by
$u \geq 0$ determined by the equations :
\[
\int f(\omega, x) \nu^{u}(d\omega dx)
=
\mathbb{E}_{\mu}\left[ 
\int \pi_{u,F/F}(\iota, d\omega)
\int \pi_{u,I/I}(\zeta, dx) f(\omega, x)\right],
\]
where $f(\omega,x)$ represents a positive $\mathcal{J}_\infty$ measurable function. 

\brem
The measure $\nu^u$ coincides with $\mu$ on $\mathcal{J}_u$, on $\iota^{-1}(\mathcal{F}^\circ_\infty)$, and on $\zeta^{-1}(\mathcal{I}_\infty)$. The kernel $\int \pi_{u,F/F}(\iota, d\omega)
\int \pi_{u,I/I}(\zeta, dx)$ is a version of the regular conditional probability on $\mathcal{J}_{\infty}$, under the probability
$\nu^{u}$ given $\mathcal{J}_u$.  We notice that the regular conditonal probability makes $\zeta$ indenpendent of $\iota^{-1}(\mathcal{F}^\circ_{\infty})$. In this sense, the probability measure $\nu^u$ realizes the condition $3^\circ$ of section \ref{lsmintro}.
\erem

\bl\label{martingaletransfer} 
Let $M$ to be a
$\mathbb{F}^\circ_+$-adapted càdlàg $\mathbb{P}$ integrable process. Then, $M$ is a
$(\mathbb{P},\mathbb{F})$-martingale on $[u,\infty)$, if and only if $M(\iota)$ is a $(\nu^u,\mathbb{J}_+)$ martingale on $[u,\infty)$.
\el

\textbf{Proof.} (cf. Proposition \ref{sbi2}) Note that, because of the càdlàg path property, $M$ is a
$(\nu^u,\mathbb{J}_+)$ martingale, if and only if
\[
\dcb
\mathbf{E}_{\nu^{u}}[M_{t}(\iota)\ind_A(\iota)\ind_B(\zeta)]
=\mathbf{E}_{\nu^{u}}[M_{s}(\iota)\ind_A(\iota)\ind_B(\zeta)],
\dce 
\]
for any $t\geq s\geq u$ and for any $A\in\mathcal{F}^\circ_s$, $B\in\mathcal{I}_s$. This is equivalent to
\[
\dcb
&&\mathbf{E}_{\mu}\left[ 
\int \pi_{u,F/F}(\iota, d\omega)M_{t}(\omega)\ind_A(\omega)
\int \pi_{u,I/I}(\zeta, dx) \ind_B(x)\right]\\
&=&\mathbf{E}_{\mu}\left[ 
\int \pi_{u,F/F}(\iota, d\omega)M_{s}(\omega)\ind_A(\omega)
\int \pi_{u,I/I}(\zeta, dx) \ind_B(x)\right],
\dce 
\]
or equivalent to
\[
\dcb
&&\mathbf{E}_{\mathbb{P}}\left[ 
M_{t}\ind_A
\mathbb{E}[\int \pi_{u,I/I}(I, dx) \ind_B(x)|\mathcal{F}^\circ_u]\right]
=\mathbf{E}_{\mathbb{P}}\left[ 
M_{s}\ind_A
\mathbb{E}[\int \pi_{u,I/I}(I, dx) \ind_B(x)|\mathcal{F}^\circ_u]\right].
\dce 
\]
The last relation is true if and only if $M$ is a $(\mathbb{P},\mathbb{F})$-martingale on $[u,\infty)$. \ \ok

\bl\label{h+}
For $t
\geq u$, 
the function $h_{t}^{u}$ exists and it is a version of $\left.\frac{2\mu}{\mu+\nu^{u}}\right |_{{\mathcal{J}}_{t}}$ (respectively, $2-h_{t}^{u}=\left.\frac{2\nu^u}{\mu+\nu^{u}}\right |_{{\mathcal{J}}_{t}}$).
Consequently, the right limit version
$$
h_{t+}^{u} =
\lim_{\epsilon>0, \epsilon\rightarrow
	0}\sup_{s \in \mathbb{Q}_{+}, u\leq s\leq u+\epsilon}  h(\omega, x),\ t\geq u,
$$
is well-defined and the process $h^u_+=(h^u_{t+})_{t\geq u}$ is a càdlàg $(\mu+\nu^{u},\mathbb{J}_+)$ uniformly integrable martingale on $[u,\infty)$. We have $0\leq h^u_{+}\leq 2$ and, if
$$
\dcb
\tau^{u}(\omega, x)
&=&
\inf\{t \geq u: h_{t+}^{u} = 2\}\\
\upsilon^{u}(\omega, x)
&=&
\inf\{t \geq u: h_{t+}^{u} = 0\}\\
\dce
$$ 
$h^u_{t+}=2$ for $t\in[\tau^u,\infty)$ and $h^u_{t+}=0$ for $t\in[\upsilon^u,\infty)$, $(\mu+\nu^u)$ almost surely. We have also $$
\mu[\upsilon^u<\infty]=0,\ \ \nu^u[\tau^u<\infty]=0.
$$
\el

\textbf{Proof.} We prove first of all two identities : Let $t\geq u$. Let $H\geq 0$ be $\mathcal{J}_{t}$ measurable,
\[
\begin{array}{lcl}
\mathbb{E}_{\mu}[H]&=&\mathbb{E}_\mathbb{P}[\int H(\mathsf{i},x)\pi_{t,I/F}(\iota, dx)],\\
\mathbb{E}_{\nu^{u}}[H]&=&\mathbb{E}_{\mathbb{P}}[\int \pi_{u,I/F}(\mathsf{i}, dx'') \int H(\mathsf{i},x) \pi_{u,I/I}(x'', dx)].
\end{array}
\]
By monotone class theorem, it is enough to check them for $H=\ind_A(\iota)\ind_B(\zeta)$ where $A\in\mathcal{F}^\circ_t$ and $B\in\mathcal{I}_t$. But then, because$$
\ind_A(\iota)\int \ind_B(x)\pi_{t,I/F}(\iota, dx)
=\int \ind_A(\iota)\ind_B(x)\pi_{t,I/F}(\iota, dx),
$$ 
we have $$
\dcb
&&\mathbb{E}_{\mu}[\ind_A(\iota)\ind_B(\zeta)]
=\mathbb{E}_{\mathbb{P}}[\ind_A\ind_B(I)]
=\mathbb{E}_{\mathbb{P}}[\ind_A\mathbf{E}_{\mathbb{P}}[\ind_B(I)|\mathcal{F}_t]]\\
&=&\mathbb{E}_{\mathbb{P}}[\ind_A\int \ind_B(x)\pi_{t,I/F}(\mathsf{i}, dx)]
=\mathbb{E}_{\mathbb{P}}[\int \ind_A(\mathsf{i})\ind_B(x)\pi_{t,I/F}(\mathsf{i}, dx)],
\dce
$$
and
$$
\dcb
&&\mathbb{E}_{\nu^{u}}[\ind_A(\iota)\ind_B(\zeta)]
=\mathbb{E}_{\mu}\left[ 
\int \pi_{u,F/F}(\iota, d\omega)
\int \pi_{u,I/I}(\zeta, dx) \ind_A(\omega)\ind_B(x)\right]\\
&=&\mathbb{E}_{\mathbb{P}}\left[ 
\int \pi_{u,F/F}(\mathsf{i}, d\omega)\ind_A(\omega)
\int \pi_{u,I/I}(I, dx) \ind_B(x)\right]
=\mathbb{E}_{\mathbb{P}}\left[ 
\int \pi_{u,F/F}(\mathsf{i}, d\omega)\ind_A(\omega)
\mathbb{E}_{\mathbb{P}}[\int \pi_{u,I/I}(I, dx) \ind_B(x)|\mathcal{F}_u]\right]\\
&=&\mathbb{E}_{\mathbb{P}}\left[\ind_A
\mathbb{E}_{\mathbb{P}}[\int \pi_{u,I/I}(I, dx) \ind_B(x)|\mathcal{F}_u]\right]\\
&=&\mathbb{E}_{\mathbb{P}}[\ind_A\int \pi_{u,I/F}(\mathsf{i}, dx'') \int \pi_{u,I/I}(x'', dx)\ind_B(x)]
=\mathbb{E}_{\mathbb{P}}[\int \pi_{u,I/F}(\mathsf{i}, dx'') \int \pi_{u,I/I}(x'', dx)\ind_A(\mathsf{i})\ind_B(x)], 
\dce
$$
which concludes the identities.

Let $h_{t}^{u}$ be the density function $\left.\frac{2\mu}{\mu+\nu^{u}}\right |_{{\mathcal{J}}_{t}}$. For $A\in\mathcal{F}^\circ_t$, $B\in\mathcal{I}_t$, we have$$
2\mathbb{E}_{\mu}[\ind_A\ind_B]=\mathbb{E}_{\mu}[\ind_A\ind_Bh^u_t]+\mathbb{E}_{\nu^u}[\ind_A\ind_Bh^u_t]
$$
Applying the two identities for $\mathbb{E}_{\mu}$ and for $\mathbb{E}_{\nu^u}$, the above relation becomes$$
\dcb
&&2\mathbb{E}_\mathbb{P}[\int \pi_{t,I/F}(\mathsf{i}, dx)\ind_A(\mathsf{i})\ind_B(x)]\\
&=&\mathbb{E}_\mathbb{P}[\int \pi_{t,I/F}(\mathsf{i}, dx)\ind_A(\mathsf{i})\ind_B(x)h^u_t(\mathsf{i},x)]+\mathbb{E}_{\mathbb{P}}[\int \pi_{u,I/F}(\mathsf{i}, dx'') \int  \pi_{u,I/I}(x'', dx)\ind_A(\mathsf{i})\ind_B(x)h^u_t(\mathsf{i},x)]\\
\dce
$$
or
$$
\dcb
&&2\mathbb{E}_\mathbb{P}[\ind_A(\mathsf{i})\int \pi_{t,I/F}(\mathsf{i}, dx)\ind_B(x)]\\

&=&\mathbb{E}_\mathbb{P}[\ind_A(\mathsf{i})\left(\int \pi_{t,I/F}(\mathsf{i}, dx)+\int \pi_{u,I/F}(\mathsf{i}, dx'') \int  \pi_{u,I/I}(x'', dx)\right)\ind_B(x)h^u_t(\mathsf{i},x)]
\dce
$$
When $A$ runs over all $\mathcal{F}^\circ_t$, the above equation becomes$$
2\int \pi_{t,I/F}(\mathsf{i}, dx)\ind_B(x)
=
\left(\int \pi_{t,I/F}(\mathsf{i}, dx)+\int \pi_{u,I/F}(\mathsf{i}, dx'') \int  \pi_{u,I/I}(x'', dx)\right)\ind_B(x)h^u_t(\mathsf{i},x),
$$
$\mathbb{P}$-almost surely. Since $\mathcal{I}_t$ is countably generated, we conclude$$
\left.2\int \pi_{t,I/F}(\mathsf{i}, dx)\right|_{\mathcal{I}_t}
=
\left.\left(\int \pi_{t,I/F}(\mathsf{i}, dx)+\int \pi_{u,I/F}(\mathsf{i}, dx'') \int  \pi_{u,I/I}(x'', dx)\right)h^u_t(\mathsf{i},x)\right|_{\mathcal{I}_t},
$$ 
$\mathbb{P}$-almost surely. This is the first part of the lemma.

The second statements on $h^u_+$ are valid, because $h^u_+$ and $2-h^u_+$ are positive $(\mu+\nu^{u},\mathbb{J}_+)$ supermartingales. For the last assertion, we have$$
\dcb
\mu[\upsilon^u<\infty]
&=&\mathbb{E}_{\frac{1}{2}(\mu+\nu^u)}[\ind_{\{\upsilon^u\leq t\}}h^u_{t+}]=0\\ 
\nu^u[\tau^u<\infty]
&=&\mathbb{E}_{\frac{1}{2}(\mu+\nu^u)}[\ind_{\{\tau^u\leq t\}}(2-h^u_{t+})]=0.\ \ok
\dce
$$

Along with the function $h^u$, we introduce two processes on $[u,\infty)$:
\begin{eqnarray*}
	\alpha_{t}^{u}(\omega, x) 
	&=&
	\frac{h_{t+}^{u}(\omega, x)}
	{2-h_{t+}^{u}(\omega, x)} \ind_{\{t < \tau^{u}(\omega, x)\}},\\
	\beta_{t}^{u}(\omega, x) 
	&=&
	\frac{2-h_{t+}^{u}(\omega, x)}{h_{t+}^{u}(\omega, x)} \ind_{\{t < \upsilon^{u}(\omega, x)\}}.
\end{eqnarray*}

We have the following theorem, which is the adaptation of Theorem \ref{basicformulas} in the setting of couple of probability measures $(\mu, \nu^u)$.

\bethe\label{basicformulas2}
Fix $u\geq 0$. Let $\rho$ to be a $\mathbb{J}_{+}$-stopping time such that $u\leq \rho\leq \tau^u$. We have, for any
$\mathbb{J}_{+}$-stopping time
$\kappa$ with $u\leq \kappa$, for any $A \in
{\mathcal{J}}_{\kappa+}$,$$
\mathbb{E}_{\mu}[\ind_A\ind_{\{\kappa < \rho\}}] = \mathbb{E}_{\nu^{u}}[\ind_A\ind_{\{\kappa < \rho\}}
\alpha_{\kappa}^{u}].
$$
Consequently, $\ind_{[0,u)}+\alpha^{u}\ind_{[u,\rho)}$ is a $(\nu^u,{\mathbb{J}}_{+})$ supermartingale. Moreover, $\alpha^u>0$ on $[u,\rho)$ under the probability $\mu$. For any positive $\mathbb{J}_+$ predictable process $H$, $$
\mathbb{E}_{\mu}[H_\kappa\ind_{\{\kappa < \rho\}}] = \mathbb{E}_{\nu^{u}}[H_\kappa\ind_{\{\kappa < \rho\}}
\alpha_{\kappa}^{u}].
$$
Suppose in addition that $(\alpha^u\ind_{[u,\rho)})^\kappa$ is of class $(D)$ under $\nu^u$. Let $V$ be the increasing $\mathbb{J}_+$ predictable process associated with the supermartingale $\ind_{[0,u)}+\alpha^{u}\ind_{[u,\rho)}$ (see \cite[p.115 Theorem 13]{protter} and Lemma \ref{FF0}). For any positive $\mathbb{J}_+$ predictable process $H$, we have$$
\mathbb{E}_{\mu}[H_{\rho}\ind_{\{u<\rho\leq \kappa\}}] =
\mathbb{E}_{\nu^{u}}[\int_{u}^{\kappa}
H_s dV_s].
$$
Let $B$ be a $\mathbb{J}_+$ predictable process with bounded total variation. We have$$
\mathbb{E}_{\mu}[B_{\kappa\wedge \rho}-B_u] =
\mathbb{E}_{\nu^{u}}[\int_{u}^{\kappa}
\alpha^u_{s-} dB_s].
$$
Consequently, $\mathbb{E}_{\mu}[\int_u^{\kappa}\ind_{\{u<s\leq \rho\}}\ind_{\{\alpha^u_{s-}=0\}}dB_s]=0$. Let $C$ be a $\mathbb{J}_+$ optional process. Suppose that the random measure $\mathsf{d}C$ on the open random interval $(u,\kappa\wedge \rho)$ has bounded total variation. For any bounded $\mathbb{J}_+$ predictable process $H$, We have$$
\mathbb{E}_{\mu}[\int_0^\infty H_s\ind_{(u,\kappa\wedge \rho)}(s) dC_s] =
\mathbb{E}_{\nu^{u}}[\int_{0}^\infty H_s\ind_{(u,\kappa\wedge \rho)}(s) \alpha^u_sdC_s].
$$
In particular, $$
(\ind_{(u,\kappa\wedge \rho)}\centerdot C)^{\mu\cdot\mathbb{J}_+-p}
=
\frac{1}{\alpha^u_-}\centerdot(\ind_{(u,\kappa\wedge \rho)}\alpha^u\centerdot C)^{\nu^u\cdot\mathbb{J}_+-p}.
$$
\ethe

\textbf{Proof.} 
It is enough to adapt the proof of Theorem \ref{basicformulas} from the case $u=0$ to the case $u>0$. \ \ok

\

\subsection{Computation of local solutions}

We have now a theorem to find local solutions.

\bethe
\label{theorem:our-formula}
Let $M$ be a
bounded
$(\mathbb{P},\mathbb{F})$ martingale (assumed to be $\mathbb{F}^\circ_+$ adapted). Let $u\geq 0$. Let $\rho$ to be a $\mathbb{J}_{+}$-stopping time such that $u\leq \rho\leq \tau^u$. Let $\langle M(\iota),\alpha^{u}\ind_{[u,\rho)} \rangle^{\nu^u\cdot\mathbb{J}_+}$ denote the $\mathbb{J}_+$-predictable
bracket of the couple of $M(\iota)$ and $\ind_{[0,u)}+\alpha^{u}\ind_{[u,\rho)}$, computed under the probability $\nu^{u}$ in
the filtration ${\mathbb{J}}_{+}$. Then, there exists an increasing sequence $(\eta_n)_{n\geq 1}$ of $\mathbb{J}_+$ stopping times such that $\sup_{n\geq 1}\eta_n(\phi)\geq \rho(\phi)$ $\mathbb{P}$-almost surely, and for every $n\geq 1$, $M$ is a special
$\mathbb{G}$-semimartingale on the random left interval 
$(u,\eta_n\circ\phi]$ such that the process 
\[
(M^{\eta_n(\phi)}-M^{u}) 
- 
(\ind_{(u,\eta_n]}\frac{1}{\alpha^u_{-}}\centerdot \langle M(\iota),
\alpha^{u}\ind_{[u,\rho)}
\rangle^{\nu^u\cdot\mathbb{J}_+})\circ \phi 
- 
\ind_{\{u<\rho(\phi)\leq \eta_n(\phi)\}}\Delta_{\rho(\phi)}M \ind_{[\rho(\phi),\infty)}
\]
is a $(\mathbb{P},\mathbb{G})$ local martingale.

If $(\ind_{(u,\rho]}\frac{1}{\alpha^u_{-}}\centerdot \langle M(\iota),
\alpha^{u}\ind_{[u,\rho)}
\rangle^{\nu^u\cdot\mathbb{J}_+})\circ \phi$ has a distribution function under $\mathbb{P}$, the above decomposition formula remains valid if we replace $\eta_n$ by $\rho$.
\ethe

\brem
Notice that, according to \cite[Proposition II.2.4]{SongThesis}, different version of $\langle M(\iota),
\alpha^{u}\ind_{[u,\rho)}
\rangle^{\nu^u\cdot\mathbb{J}_+}$ under $\nu^u$ give rise of indistinguishable versions of $\langle M(\iota),
\alpha^{u}\ind_{[u,\rho)}
\rangle^{\nu^u\cdot\mathbb{J}_+}\circ \phi$ \ on $(u,\tau^u(\phi)]$ under $\mathbb{P}$. 
\erem

\textbf{Proof.} 
In this proof, instead of writing
$\nu^{u}$ (resp. $\tau^{u}, \upsilon^u, \alpha^{u}, \langle M(\iota),
\alpha^{u}\ind_{[u,\rho)}
\rangle^{\nu^u\cdot\mathbb{J}_+}$), we shall simply write
$\nu$ (resp.
$\tau, \upsilon,
\alpha, \langle M(\iota),
\alpha\ind_{[u,\rho)}
\rangle$).
For each integer $n
\geq 1$, set $$
\dcb
\gamma^+_{n} = \inf\{t \geq u:
\alpha_{t}>n\}\ \mbox{ and }\
\gamma^-_{n} = \inf\{t \geq u:
\alpha_{t}<1/n\},\ n\geq 1.
\dce
$$ 
Let $(\gamma^\circ_n)_{n\geq 1}$ be a sequence of $\mathbb{J}_{+}$-stopping times tending to $\infty$ under $\nu$ such that $\gamma^\circ_n\geq u$ for every $n\geq 1$, and $\cro{M(\iota),\alpha\ind_{[u,\rho)}}^{\gamma^\circ_n}$ has bounded total variation. 
Let $\eta_n=\gamma^+_{n}\wedge\gamma^-_{n}\wedge\gamma^\circ_{n}\wedge\rho$. Note that $\gamma^+_\infty=\tau$ on $\{\gamma^+_\infty<\infty\}$, $\gamma^-_\infty=\tau\wedge \upsilon$ under $\mu+\nu$. We have $$
\dcb
\mathbb{P}[\eta_\infty(\phi)<\rho(\phi)]
=\mu[\eta_\infty<\rho]
=\mathbb{E}_\nu[\ind_{\{\eta_\infty<\rho\}}\alpha_{\eta_\infty}]
=\mathbb{E}_\nu[\ind_{\{\tau\wedge\upsilon<\rho\}}\alpha_{\tau\wedge\upsilon}]
=0.
\dce
$$

Fix an $n\geq 1$. Let $\kappa$ be a $\mathbb{J}_{+}$-stopping time such that $u
\leq \kappa \leq \eta_{n}$. The stopping time $\kappa$ satisfies therefore the condition of Lemma \ref{basicformulas}. Recall also Lemma \ref{martingaletransfer} according to which $M(\iota)$ is a $(\nu^u,\mathbb{J}_+)$ martingale on $[u,\infty)$. 

Set $M^\flat=M-\Delta_{\rho(\phi)}M\ind_{[\rho(\phi),\infty)}$. Applying Lemma \ref{basicformulas}, we can write
\[
\begin{array}{lll}
&&\mathbb{E}_{\mathbb{P}}[(M^\flat_{\kappa(\phi)}- M^\flat_{u})]
\\
&=&\mathbb{E}_{\mu}[(M^\flat(\iota)_\kappa- M^\flat(\iota)_u)]
\\ 
&=& \mathbb{E}_{\mu}[(M(\iota)_\kappa- M(\iota)_u)\ind_{\{\kappa < \rho\}}] +
\mathbb{E}_{\mu}[(M^\flat(\iota)_\kappa- M^\flat(\iota)_u)1_{\{\kappa = \rho\}}] 
\\
&=&
\mathbb{E}_{\nu}[(M(\iota)_\kappa- M(\iota)_u)\alpha_{\kappa}\ind_{\{\kappa < \rho\}}] 
+
\mathbb{E}_{\mu}[(M(\iota)_{\rho-}- M(\iota)_u)\ind_{\{u<\kappa=\rho\}}]\\
&=&
\mathbb{E}_{\nu}[(M(\iota)_\kappa- M(\iota)_u)\alpha_{\kappa}\ind_{\{\kappa < \rho\}}] 
+
\mathbb{E}_{\nu}[\int_{u}^{\kappa}(M(\iota)_{-}- M(\iota)_{u})
dV]\\
&=&
\mathbb{E}_{\nu}[\int_{u}^{\kappa}
d\langle M(\iota), \alpha\ind_{[u,\rho)} \rangle_s]\\
&=&
\mathbb{E}_{\mu}[\int_{u}^{\kappa} \frac{1}{\alpha_{s-}}
d\langle M(\iota), \alpha\ind_{[u,\rho)} \rangle_s]\\
&&\mbox{ The above process is well-defined}\\
&&\mbox{ and it has bounded total variation, because $\frac{1}{\alpha_-}\leq n$ on $(u,\kappa]$}\\

&=&
\mathbb{E}_{\mathbb{P}}[(\int_{u}^{\kappa(\phi)} \frac{1}{\alpha(\phi)_{s-}}
d(\langle M(\iota), \alpha\ind_{[u,\rho)} \rangle\circ\phi)_s]
\dce
\]
This identity on $\mathbf{E}_{\mathbb{P}}[(M^\flat_{\kappa(\phi)}- M^\flat_{u})]$ implies that  (cf. \cite[Theorem 4.40]{HWY}) 
\[
(M^{\eta_n(\phi)}-M^{u}) 
- 
(\ind_{(u,\eta_n]}\frac{1}{\alpha_{-}}\centerdot \langle M(\iota),
\alpha\ind_{[u,\rho)}
\rangle)\circ \phi 
- 
\ind_{\{u<\rho(\phi)\leq \eta_n(\phi)\}}\Delta_{\rho(\phi)}M \ind_{[\rho(\phi),\infty)}
\]
is a $(\mathbb{P},\mathbb{G})$ local martingale. The first part of the theorem is proved. 

Suppose now that $(\ind_{(u,\rho]}\frac{1}{\alpha_{-}}\centerdot \langle M(\iota),
\alpha\ind_{[u,\rho)}
\rangle\circ \phi$ has a distribution function $C=(C_t)_{t\geq 0}$ under $\mathbb{P}$. It is clear that $\ind_{(u,\eta_n(\phi)]}\centerdot M^\flat$ converges in $\mathsf{H}^1_{loc}(\mathbb{P},\mathbb{G})$. Note that $\eta_\infty=\rho$ under $\mu$. Look at $\lim_{n\uparrow\infty}(M^\flat_{\eta_n}-M^\flat_u)$. If $\eta_n<\rho$ for all $n\geq 1$, this limit is equal to $(M^\flat_{\rho-}-M^\flat_u)=(M^\flat_{\rho}-M^\flat_u)$. If $\eta_n=\rho$ for some $n\geq 1$, this limit is equal to $(M^\flat_{\rho}-M^\flat_u)$ also. Look at the limit $$
\lim_{n\uparrow\infty}\ind_{(u,\eta_n(\phi)]}\frac{1}{\alpha(\phi)_{-}}\centerdot \langle M(\iota),
\alpha\ind_{[u,\rho)}
\rangle\circ \phi
=\ind_{\cup_{n\geq 1}(u,\eta_n(\phi)]}\frac{1}{\alpha(\phi)_{-}}\centerdot \langle M(\iota),
\alpha\ind_{[u,\rho)}
\rangle\circ \phi.
$$
We note that $(u,\rho]\setminus \cup_{n\geq 1}(u,\eta_n]=[\rho_{\{\forall n\geq 1, \eta_n<\rho\}}]$. The time $\rho_{\{\forall n\geq 1, \eta_n<\rho\}}$ is a $\mathbb{J}_+$ predictable stopping time. So, on the set $\{\rho_{\{\forall n\geq 1, \eta_n<\rho\}}<\infty\}$, $$
\dcb
&&\ind_{(u,\rho]\setminus \cup_{n\geq 1}(u,\eta_n]}\centerdot\langle M(\iota), \alpha\ind_{[u,\rho)} \rangle\\
&=&\Delta_{\rho_{\{\forall n\geq 1, \eta_n<\rho\}}}\langle M(\iota), \alpha\ind_{[u,\rho)} \rangle\\

&=&\mathbb{E}_{\nu}[\Delta_{\rho_{\{\forall n\geq 1, \eta_n<\rho\}}} M(\iota) \Delta_{\rho_{\{\forall n\geq 1, \eta_n<\rho\}}}(\alpha\ind_{[u,\rho)})|\mathcal{J}_{\rho_{\{\forall n\geq 1, \eta_n<\rho\}}-}]\\

&=&\mathbb{E}_{\nu}[\ind_{\{\forall n\geq 1, \eta_n<\rho<\infty\}}\Delta_{\rho} M(\iota) \Delta_{\rho}(\alpha\ind_{[u,\rho)})|\mathcal{J}_{\rho_{\{\forall n\geq 1, \eta_n<\rho\}}-}]\\

&=&\mathbb{E}_{\nu}[\ind_{\{\forall n\geq 1, \eta_n<\rho<\infty\}}\Delta_{\rho} M(\iota) (-\alpha_{\rho-})|\mathcal{J}_{\rho_{\{\forall n\geq 1, \eta_n<\rho\}}-}]\\

&=&\mathbb{E}_{\nu}[\Delta_{\rho_{\{\forall n\geq 1, \eta_n<\rho\}}} M(\iota) (-\alpha_{\rho_{\{\forall n\geq 1, \eta_n<\rho\}}-})|\mathcal{J}_{\rho_{\{\forall n\geq 1, \eta_n<\rho\}}-}]\ind_{\{\rho_{\{\forall n\geq 1, \eta_n<\rho\}}<\infty\}}\\

&=&\mathbb{E}_{\nu}[\Delta_{\rho_{\{\forall n\geq 1, \eta_n<\rho\}}} M(\iota) |\mathcal{J}_{\rho_{\{\forall n\geq 1, \eta_n<\rho\}}-}](-\alpha_{\rho_{\{\forall n\geq 1, \eta_n<\rho\}}-})\ind_{\{\rho_{\{\forall n\geq 1, \eta_n<\rho\}}<\infty\}}\\

&=&0\
\mbox{ because $\rho_{\{\forall n\geq 1, \eta_n<\rho\}}$ is predictable and $M$ is a bounded martingale}.

\dce
$$
We have actually
$$
\lim_{n\uparrow\infty}\ind_{(u,\eta_n(\phi)]}\frac{1}{\alpha(\phi)_{-}}\centerdot \langle M(\iota),
\alpha\ind_{[u,\rho)}
\rangle\circ \phi
=\ind_{(u,\rho(\phi)]}\frac{1}{\alpha(\phi)_{-}}\centerdot \langle M(\iota),
\alpha\ind_{[u,\rho)}
\rangle\circ \phi.
$$
The theorem is proved.
\ok

The following lemma will also be very useful.

\bl\label{alpha} Let $u \geq 0$. 
Let $\lambda$ be a ${\mathbb{J}_+}$-stopping time $\geq u$. Let $B$ be  a
$\mathcal{J}_{u}$-measurable set. If  for any $t\geq u$,  there is a
$\mathcal{J}_{t}$-measurable non negative function $l_{t}$ such that, $\mathbb{P}$-almost surely on
$\omega$
\[
\dcb
&&\int \pi_{t,I/F}(\omega, dx) f(x) \ind_{B}(\omega,x)\ind_{\{t \leq  \lambda(\omega,x)\}} \\
&=& \int \pi_{u,I/F}(\omega, dx'') \int\pi_{u,I/I}(x'', dx) f(x)\ind_{B}(\omega,x)\ind_{\{t \leq  \lambda(\omega,x)\}} l_{t}(\omega,x), 
\dce
\]
for all bounded $\mathcal{I}_{t}$-measurable function $f$,
then, \[\tau^{u}\ind_{B} \geq \lambda \ind_{B}, \ \  
\alpha_{t}^{u}\ind_{B}\ind_{\{t < \lambda\}} = \limsup_{t'\in\mathbb{Q}_+,t'\downarrow t} l_{t'}\ind_{B}\ind_{\{t < \lambda\}}, t\geq u.\]
\el

Note that, if $\lambda$ is a $\mathbb{J}$ stopping time, we can replace $\ind_{\{t\leq\lambda\}}$ par $\ind_{\{t<\lambda\}}$ in the assumption.

\textbf{Proof.} We have, for almost all $\omega$,
$$
\dcb
&&\left(\int \pi_{t,I/F}(\omega, dx)
+ \int \pi_{u,I/F}(\omega, dx'') \int\pi_{u,I/I}(x'', dx)\right) f(x)\ind_{B}(\omega,x)\ind_{\{t \leq \lambda(\omega,x)\}}\\
&=&\int \pi_{u,I/F}(\omega, dx'') \int\pi_{u,I/I}(x'', dx) f(x)\ind_{B}(\omega,x)\ind_{\{t \leq \lambda(\omega,x)\}} (l_{t}(\omega,x)+1).
\dce
$$
It follows that (noting that $\{t \leq \lambda\}\in\mathcal{J}_t$).
$$
\dcb
&&\left(\int \pi_{t,I/F}(\omega, dx)
+ \int \pi_{u,I/F}(\omega, dx'') \int\pi_{u,I/I}(x'', dx)\right) f(x)\ind_{B}(\omega,x)\ind_{\{t \leq \lambda(\omega,x)\}}h^u_t(\omega,x)\\
&=&2\int \pi_{t,I/F}(\omega, dx) f(x) \ind_{B}(\omega,x)\ind_{\{t \leq \lambda(\omega,x)\}}\\
&=& \int \pi_{u,I/F}(\omega, dx'') \int\pi_{u,I/I}(x'', dx) f(x)\ind_{B}(\omega,x)\ind_{\{t \leq \lambda(\omega,x)\}} 2l_{t}(\omega,x)\\
&=& \int \pi_{u,I/F}(\omega, dx'') \int\pi_{u,I/I}(x'', dx) f(x)\ind_{B}(\omega,x)\ind_{\{t \leq \lambda(\omega,x)\}} \frac{2l_{t}(\omega,x)}{l_{t}(\omega,x)+1}(l_{t}(\omega,x)+1)\\
&=&\left(\int \pi_{t,I/F}(\omega, dx)
+ \int \pi_{u,I/F}(\omega, dx'') \int\pi_{u,I/I}(x'', dx)\right) f(x)\ind_{B}(\omega,x)\ind_{\{t \leq \lambda(\omega,x)\}}\frac{2l_{t}(\omega,x)}{l_{t}(\omega,x)+1}.
\dce
$$
The above identity means, for almost all $\omega$, the (in)equalities
$$
h_{t}^{u}(\omega,x) \ind_{B}(\omega,x)\ind_{\{t \leq \lambda(\omega,x)\}} = 
\frac{2l_{t}(\omega,x)}{l_{t}(\omega,x)+1}\ind_{B}(\omega,x)\ind_{\{t \leq \lambda(\omega,x)\}}< 2,
$$ 
hold for $\left(\int \pi_{t,I/F}(\omega, dx) + \int \pi_{u,I/F}(\omega, dx'') \int\pi_{u,I/I}(x'', dx)\right)$ almost all $x$. By the proof of Lemma \ref{h+}, these (in)equalities also hold $\mu+\nu^u$ almost surely.  This is enough to conclude the lemma. \ok

\

\

\section{Classical results as consequence of the local solution method}

One of the motivation to introduce the local solution method was to provide unified proofs for the classical formulas. We will do it in this section with the three most important classical formulas. 

We employ the same notations as in subsection \ref{h-formula}. In particular,  $(\Omega,\mathcal{A},\mathbb{P})$ denotes a Polish probability space and $\mathbb{F}=(\mathcal{F}_{t})_{t\in \mathbb{R}_{+}}$ denotes a filtration on $\Omega$. We assume that there exists a filtration $\mathbb{F}^\circ=(\mathcal{F}^\circ_t)_{t\in \mathbb{R}_{+}}$ of sub-$\sigma$-algebras of $\mathcal{A}$ such that $$
\mathcal{F}_t=\mathcal{F}^\circ_{t+}\vee\mathcal{N}, \ t\geq 0,
$$ 
where $\mathcal{N}$ is the family of $(\mathbb{P},\mathcal{A})$-negligible sets.

\subsection{Jacod's criterion}\label{jcrit}

We consider a random variable $\xi$ taking values in a polish space $E$. We consider the enlargement of filtration given by $$
\mathcal{G}_t=\cap_{s>t}(\mathcal{F}_s\vee\sigma(\xi)),\ t\geq 0
$$
It is the so-called initial enlargement of filtration. In \cite{jacod2}, the following assumption is introduced

\bassumption($A'$)
For any $t\geq 0$, the conditional law of $\xi$ given the $\sigma$-algebra $\mathcal{F}_t$ is absolutely continuous with respect to the (unconditional) law $\lambda$ of $\xi$
\eassumption
It is proved then (cf. \cite[Lemme(1.8)]{jacod2}) that, if $q=q_t(x)=q_t(\omega,x)$ denotes the density function of the conditional law of $\xi$ given the $\sigma$-algebra $\mathcal{F}_t$ with respect to $\lambda(dx)$, there exists a jointly measurable version of $q$ such that $t\rightarrow q_t(\omega,x)$ is càdlàg, and for every $x$, the map $(\omega,t)\rightarrow q_t(\omega,x)$ is a $(\mathbb{P},\mathbb{F})$ martingale. 

\bethe\label{jacod-criterion} (\cite[Théorème(2.1)]{jacod2})
For every bounded $\mathbb{F}$ martingale $M$, it is a $\mathbb{G}$ special semimartingale whose drift is given by the formula:$$
\frac{1}{q_{-}(\xi)}\centerdot \left.(\cro{q(x),M}^{\mathbb{P}\cdot\mathbb{F}})\right|_{x=\xi}
$$ 
\ethe

Let us see how the above formula can be obtained by the local solution method. We have assumed the existence of $\mathbb{F}^\circ$ behind the filtration $\mathbb{F}$. Note that, if Assumption($A'$) holds, since $\mathcal{F}^\circ_t\subset\mathcal{F}_t, t\geq 0,$ the same assumption holds with respect to the filtration $\mathbb{F}^\circ$. Let $p=p_t(x)=p_t(\omega,x)$ denote the density function of the conditional law of $\xi$ given the $\mathcal{F}^\circ_t$ with respect to $\lambda(dx)$. We also assume that $q_t(\omega,x)$ is $\mathbb{J}_+$ adapted by Lemma \ref{FF0}.

Define the process $I$ by $I_t=\xi, \forall t\geq 0$ ($I$ is a constant process taking values in the space $C(\mathbb{R}_+,E)$ of all continuous function on $\mathbb{R}_+$ taking value in $E$, equipped with the natural filtration). Define the filtration $\mathbb{G}^\circ$ as in subsection \ref{h-formula}. Then, the filtrations $\mathbb{G}$ and $\mathbb{G}^\circ$ are related in the way : $\mathcal{G}_t=\mathcal{G}^\circ_{t+}\vee\mathcal{N}, t\geq 0,$ as in subsection \ref{h-formula}.

Let $u\geq 0$. We have clearly $\pi_{u,I/I}(y'',dy)=\delta_{y''}(dy)$. For $x\in E$, denote by $\mathfrak{c}_x$ the constant function at value $x$. The process $I$ is written as $I=\mathfrak{c}_\xi$. For $t\geq u$, for any bounded borel function $F$ on $I$, we have
$$
\dcb
\int \pi_{t,I/F}(\mathsf{i},dy)F(y)

=\mathbb{E}[F(I)|\mathcal{F}^\circ_t]

=\mathbb{E}[F(\mathfrak{c}_{\xi})|\mathcal{F}^\circ_t]
=\int F(\mathfrak{c}_{x})p_t(x)\lambda(dx)
\dce
$$
Therefore,
$$
\dcb
&&\int \pi_{t,I/F}(\mathsf{i},dy)F(y)\ind_{\{p_u(y_0)>0\}}
\\
&=&\int F(\mathfrak{c}_{x})\ind_{\{p_u((\mathfrak{c}_{x})_0)>0\}}p_t(x)\lambda(dx)

\\
&=&\int F(\mathfrak{c}_{x})\frac{p_t((\mathfrak{c}_{x})_0)}{p_u((\mathfrak{c}_{x})_0)}\ind_{\{p_u((\mathfrak{c}_{x})_0)>0\}}p_u(x)\lambda(dx)\\
\\
&=&\int \pi_{u,I/F}(\mathsf{i},dy)F(y)\frac{p_t(y_0)}{p_u(y_0)}\ind_{\{p_u(y_0)>0\}}\\
\\
&=&\int \pi_{u,I/F}(\mathsf{i},dy'')\int \pi_{u,I/I}(y'',dy)F(y)\frac{p_t(y_0)}{p_u(y_0)}\ind_{\{p_u(y_0)>0\}}
\dce
$$
By Lemma \ref{alpha}, we conclude that $\tau^u(\omega,,y)\ind_{\{p_u(\omega,y_0)>0\}}=\infty\ind_{\{p_u(\omega,y_0)>0\}}$ under $(\mu+\nu^u)$ and $$
\alpha^u_t(\omega,y)\ind_{\{p_u(\omega,y_0)>0\}}=\ind_{\{p_u(\omega,y_0)>0\}}\frac{p_{t+}(\omega,y_0)}{p_u(\omega,y_0)}, \ t\geq u
$$ 
From now on, let $u=0$. By the first formula in Lemma \ref{basicformulas}, the measure
$
\ind_{\{p_0(y_0)>0\}}\frac{p_{t+}(y_0)}{p_0(y_0)}d\nu^0
$
coincides with 
$\ind_{\{p_0(\iota,\zeta_0)>0\}}d\mu$ on $\mathcal{J}_{t+}$. By the second identity in the proof of Lemma \ref{h+}, we can write, for $A\in\mathcal{F}^\circ_t, B\in\mathcal{I}_t$ 
$$
\dcb
&&\int\ind_A(\omega)\ind_B(y)\ind_{\{p_0(\omega,y_0)>0\}}\frac{q_{t}(\omega,y_0)}{p_0(\omega,y_0)}d \nu^{0}(d\omega,dy)\\
&=&\mathbb{E}_{\mathbb{P}}[\ind_A(\mathsf{i})\int \pi_{0,I/F}(\mathsf{i}, dy'') \int \pi_{0,I/I}(y'', dy)\ind_B(y)\ind_{\{p_0(\mathsf{i},y_0)>0\}}\frac{q_{t}(\mathsf{i},y_0)}{p_0(\mathsf{i},y_0)}]\\

&=&\mathbb{E}_{\mathbb{P}}[\ind_A(\mathsf{i})\int \pi_{0,I/F}(\mathsf{i}, dy) \ind_B(y)\ind_{\{p_0(\mathsf{i},y_0)>0\}}\frac{q_{t}(\mathsf{i},y_0)}{p_0(\mathsf{i},y_0)}]\\

&=&\mathbb{E}_{\mathbb{P}}[\ind_A(\mathsf{i})\int  \ind_B(\mathfrak{c}_x)\ind_{\{p_0(\mathsf{i},(\mathfrak{c}_x)_0)>0\}}\frac{q_{t}(\mathsf{i},(\mathfrak{c}_x)_0)}{p_0(\mathsf{i},(\mathfrak{c}_x)_0)}p_0(\mathsf{i},x)\lambda(dx)]\\

&=&\mathbb{E}_{\mathbb{P}}[\ind_A(\mathsf{i})\int  \ind_B(\mathfrak{c}_x)\ind_{\{p_0(\mathsf{i},x)>0\}}q_t(\mathsf{i},x)\lambda(dx)]\\

&=&\mathbb{E}_{\mathbb{P}}[\ind_A(\mathsf{i})\mathbb{E}[ \ind_B(\mathfrak{c}_\xi)\ind_{\{p_0(\mathsf{i},\xi)>0\}}|\mathcal{F}_t]]\\

&=&\mathbb{E}_{\mathbb{P}}[\ind_A(\mathsf{i}) \ind_B(I)\ind_{\{p_0(\mathsf{i},\xi)>0\}}]\\

&=&\mathbb{E}_{\mu}[\ind_A(\iota) \ind_B(\zeta)\ind_{\{p_0(\iota,\zeta_0)>0\}}]\\
\dce
$$
On the other hand, with the martingale property of $q(x)$, we check that $\ind_{\{p_0(\omega,y_0)>0\}}\frac{q_{t}(\omega,y_0)}{p_0(\omega,y_0)}, t\geq 0$, is a $(\nu^u,\mathbb{J}_{+})$ martingale. The above identity remains valid, if we replace $t$ by $t+\epsilon, \epsilon>0,$ and if we replace $\ind_A(\iota) \ind_B(\zeta)$ by a $\ind_C, C\in\mathcal{J}_{t+}$. Let then $\epsilon\downarrow 0$. We prove thus that the two measures $\ind_{\{p_0(y_0)>0\}}\frac{p_{t+}(y_0)}{p_0(y_0)}d\nu^0$ and $\ind_{\{p_0(y_0)>0\}}\frac{q_{t}(y_0)}{p_0(y_0)}d\nu^0$ coincide all with 
$\ind_{\{p_0(\iota,\zeta_0)>0\}}d\mu$ on $\mathcal{J}_{t+}$. Consequently $\ind_{\{p_0(y_0)>0\}}\frac{p_{t+}(y_0)}{p_0(y_0)}=\ind_{\{p_0(y_0)>0\}}\frac{q_{t}(y_0)}{p_0(y_0)}$ $\nu^0$-almost surely.

Note that $\rho=\infty\ind_{\{p_0(y_0)>0\}}$ is a $\mathbb{J}_+$ stopping time with $0\leq \rho\leq \tau^0$. Using the product structure of the random measure $\pi_{0,F/F}(\omega,d\omega')\otimes\pi_{0,I/I}(y,dy')$, we check that $\ind_{\{p_0(y_0)>0\}}\frac{1}{p_0(y_0)}\left(\cro{M(\iota),q(x)}^{\mathbb{P}\cdot\mathbb{F}}\right)_{x=y_0}$ is a version of $\cro{M(\iota),\alpha^0\ind_{[0,\rho)}}^{\nu^0\cdot\mathbb{J}_+}$. As $$
\mathbb{E}[\ind_{\{p_0(\xi)=0\}}]
=\mathbb{E}[\mathbb{E}[\ind_{\{p_0(\xi)=0\}}|\mathcal{F}^\circ_0]]
=\mathbb{E}[\int \ind_{\{p_0(x)=0\}}p_0(x)\lambda(dx)]
=0
$$
we have $\rho(\phi)=\infty$ $\mathbb{P}$-almost surely.

Now Theorem \ref{theorem:our-formula} is applicable. We conclude that $M^{\infty}-M^0$ is a $(\mathbb{P},\mathbb{G})$ special semimartingale whose drift is given by
$$
\ind_{(0,\rho(\phi)]}\frac{1}{\alpha^0_-(\phi)}\centerdot \left(\cro{M(\iota),\alpha^0\ind_{[0,\rho)}}^{\nu^0\cdot\mathbb{J}_+}\right)\circ\phi
=
\frac{1}{q_{-}(\xi)}\centerdot \left.(\cro{q(x),M}^{\mathbb{P}\cdot\mathbb{F}})\right|_{x=\xi}
$$
This proves Theorem \ref{jacod-criterion}.

\subsection{Progressive enlargement of filtration}\label{progressive}

Consider a positive random variable $\mathfrak{t}$ (assumed not identically null). We introduce the process$$
\dcb
I= \ind_{(\mathfrak{t},\infty)}\in\mathtt{D}^g(\mathbb{R}_+,\mathbb{R})
\dce
$$
where $\mathtt{D}^g(\mathbb{R}_+,\mathbb{R})$ is the space of all càglàd functions. We equip this space with the distance $\mathfrak{d}(x,y), x,y\in \mathtt{D}^g(\mathbb{R}_+,\mathbb{R})$, as the Skorohod distance of $x_+,y_+$ in the space $\mathtt{D}(\mathbb{R}_+,\mathbb{R})$ of càdlàg functions. Let $\mathbb{G}$ be the filtration of the $\sigma$-algebras $$
\mathcal{G}_t=\cap_{s>t}(\mathcal{F}_{s}\vee\sigma(I_u:u\leq s)),\ t\geq 0
$$
It is the so-called progressive enlargement of filtration. In \cite[Proposition(4.12)]{Jeulin80}, it is proved :

\bethe\label{avant-tau}
For all bounded $(\mathbb{P},\mathbb{F})$ martingale $M$, the stopped process $M^\mathfrak{t}$ is a special $(\mathbb{P},\mathbb{G})$ semimartingale whose drift is given by $$
\ind_{(0,\mathfrak{t}]}\frac{1}{Z_-}\centerdot (\cro{N,M}^{\mathbb{P}\cdot \mathbb{F}}+\mathsf{B}^M),
$$
where $Z$ is the Azéma supermartingale $\mathbb{E}[t<\mathfrak{t}|\mathcal{F}_t], t\geq 0$, $N$ is the $\mathbb{F}$ martingale part of $Z$, and $\mathsf{B}^M$ is the $(\mathbb{P},\mathbb{F})$ predictable dual project of the jump process $\ind_{\{0<\mathfrak{t}\}}\Delta_\mathfrak{t}M \ind_{[\mathfrak{t},\infty)}$.
\ethe

We call the formula in this theorem the formula of progressive enlargement of filtration. Let us look at how the above formula can be established with the local solution method.

Let $\mathcal{I}$ be the natural filtration on $\mathtt{D}^g(\mathbb{R}_+,\mathbb{R})$ generated by the coordinates. For $x\in \mathtt{D}^g(\mathbb{R}_+,\mathbb{R})$, let $\rho(x)=\inf\{s\geq 0: x_s\neq x_0\}$. We consider $\rho$ as a function on $\Omega\times\mathtt{D}(\mathbb{R}_+,\mathbb{R})$ in identifying $\rho(x)$ with $\rho(\zeta(\omega,x))$.

Note that, for $t\geq 0$, $I^t=0$ on $\{t\leq \mathfrak{t}\}$ so that $$
\sigma(I_s:s\leq t)\cap\{\mathfrak{t}\geq t\}
=\{\emptyset,\Omega\}\cap\{\mathfrak{t}\geq  t\}
$$
We now compute $\alpha^u$ for $u=0$. For any bounded function $f\in\mathcal{I}_t$, $f(I)\ind_{\{t\leq \mathfrak{t}\}}=f(0)\ind_{\{t\leq \mathfrak{t}\}}$. Hence, $$
\dcb
&&\int\pi_{0,I/I}(I,dx)f(x)\ind_{\{t\leq\rho(x)\}}
=\mathbb{E}[f(I)\ind_{\{t\leq \mathfrak{t}\}}|\sigma(I_0)]\\
\\
&=&f(0)\mathbb{E}[\ind_{\{t\leq \mathfrak{t}\}}|\sigma(I_0)]
=f(0)\mathbb{P}[t\leq \mathfrak{t}]
=f(0)\mathbb{P}[t\leq \mathfrak{t}](1-I_0)
=\int f(x)\delta_0(dx)\ \mathbb{P}[t\leq \mathfrak{t}](1-I_0)
\dce
$$
We have then$$
\dcb
&&\int\pi_{0,I/F}(\mathsf{i},dx'')\int\pi_{0,I/I}(x'',dx)f(x)\ind_{\{t\leq \rho(x)\}}\\
&=&\int\pi_{0,I/F}(\mathsf{i},dx'')f(0)\mathbb{P}[t\leq \mathfrak{t}](1-x''_0)\\
&=&f(0)\mathbb{P}[t\leq\mathfrak{t}]\mathbb{E}[1-I_0|\mathcal{F}^\circ_0]\\
&=&f(0)\mathbb{P}[t\leq\mathfrak{t}]\\
\dce
$$
As for $\pi_{t,I/F}$, note that, by supermartingale property, if $Z_0=0$, $Z_t=0$. Notice also that, according to Lemma \ref{FF0}, $\widetilde{Z}_t=\mathbb{E}[\ind_{\{t\leq\mathfrak{t}\}}|\mathcal{F}_t], t\geq 0,$ may be chosen to be a Borel function, and to be assumed $\mathcal{F}^\circ_t$ adapted. Hence, we have
$$
\dcb
&&\int\pi_{t,I/F}(\mathsf{i},dx)f(x)\ind_{\{t\leq\rho(x)\}}\\
&=&\mathbb{E}[f(I)\ind_{\{t\leq\mathfrak{t}\}}|\mathcal{F}^\circ_t]\\
&=&f(0)\mathbb{E}[\ind_{\{t\leq\mathfrak{t}\}}|\mathcal{F}^\circ_t]\\
&=&f(0)\widetilde{Z}_t\\
&=&\frac{\widetilde{Z}_t}{\mathbb{P}[t\leq\mathfrak{t}]} f(0)\mathbb{P}[t<\mathfrak{t}]\\
&=&\int\pi_{0,I/F}(\mathsf{i},dx'')\int\pi_{0,I/I}(x'',dx)f(x)\ind_{\{t\leq\rho(x)\}}\frac{\widetilde{Z}_t}{\mathbb{P}[t\leq\mathfrak{t}]}\\
\dce
$$
By Lemma \ref{alpha}, we conclude that $\tau^0\geq \rho$ under $(\mu+\nu^0)$ and $$
\alpha^0_t(\omega,x)\ind_{\{t<\rho(x)\}}=\lim_{s\downarrow t}\frac{\widetilde{Z}_{s}}{\mathbb{P}[s\leq\mathfrak{t}]}\ind_{\{t<\rho(x)\}}
=\frac{Z_{t}}{\mathbb{P}[t<\mathfrak{t}]}\ind_{\{t<\rho(x)\}}, t\geq 0,
$$ 
according to \cite[p.63]{Jeulin80}. 

We now compute $\cro{M(\iota),\alpha^0\ind_{[0,\rho)}}^{\nu^0\cdot\mathbb{J}_+}$. We recall that $M$ and $Z$ are $(\nu^0,\mathbb{J}_+)$ semimartingales. Let $t_0=\inf\{s: \mathbb{P}[s<\mathfrak{t}]=0\}$. Then, $$
\int\pi_{0,I/I}(I,dx)\ind_{\{t_0<\rho(x)\}}
=\mathbb{E}[\ind_{\{t_0<\mathfrak{t}\}}|\sigma(I_0)]
=\mathbb{P}[t_0<\mathfrak{t}]=0
$$
We conclude that, under $\nu^0$, the process $v_t=\frac{1}{\mathbb{P}[t<\mathfrak{t}]}\ind_{\{t<\rho(x)\}}$ is well-defined $\mathbb{J}_+$ adapted with finite variation. By integration by parts formula, we write$$
d(Z_tv_t)=v_{t-}dZ_t+Z_{t-}dv_t+d[Z,v]_t
$$ 
Therefore, $$
d[M(\iota),Zv]_t
=v_{t-}d[M(\iota),Z]_t+Z_{t-}d[M(\iota),v]_t+d[M(\iota),[Z,v]]_t
$$
We deal successively the three terms at the right hand side of the above identity. Let $0\leq s<t$. Let $A\in\mathcal{F}^\circ_s$ and $B\in\mathcal{I}_s$. In the following computation we will consider the functions on $\Omega$ and respectively on $\mathtt{D}(\mathbb{R}_+,\mathbb{R})$ as functions on the product space $\Omega\times \mathtt{D}(\mathbb{R}_+,\mathbb{R})$. We give up the notations such as $\iota$ or $\zeta$. We will not denote all $\omega$ and $x$. Begin with the last term. $$
\dcb
&&\mathbb{E}_{\nu^0}[\ind_A\ind_B\left([M,[Z,v]]_t-[M,[Z,v]]_s\right)]\\

&=&\mathbb{E}[\int\pi_{0,I/I}(I,dx)\int\pi_{0,F/F}(\mathsf{i},d\omega)\ind_A\ind_B\left([M,[Z,v]]_t-[M,[Z,v]]_s\right)]\\

&=&\mathbb{E}[\int\pi_{0,I/I}(I,dx)\ind_B\int\pi_{0,F/F}(\mathsf{i},d\omega)\ind_A\sum_{s<u\leq t}\Delta_uM \Delta_uZ\Delta_uv]\\

&=&\mathbb{E}[\int\pi_{0,I/I}(I,dx)\ind_B\int\pi_{0,F/F}(\mathsf{i},d\omega)\ind_A\int_s^t\Delta_uvd[M,Z]_u]\\

&=&\mathbb{E}[\int\pi_{0,F/F}(\mathsf{i},d\omega)\ind_A\int_s^td[M,Z]_u\int\pi_{0,I/I}(I,dx)\ind_B\Delta_uv]\\

\dce
$$
Look at the integral $\int\pi_{0,I/I}(I,dx)\ind_B\Delta_uv$ (noting the $t\geq u>s\geq 0$). $$
\dcb
\int\pi_{0,I/I}(I,dx)\ind_B\Delta_uv
&=&\mathbb{E}[\ind_B(I)\left(\frac{1}{\mathbb{P}[u<\mathfrak{t}]}\ind_{u<\mathfrak{t}}-\frac{1}{\mathbb{P}[u\leq \mathfrak{t}]}\ind_{u\leq \mathfrak{t}}\right)|\sigma(I_0)]\\

&=&\mathbb{E}[\ind_B(0)\left(\frac{1}{\mathbb{P}[u<\mathfrak{t}]}\ind_{u<\mathfrak{t}}-\frac{1}{\mathbb{P}[u\leq \mathfrak{t}]}\ind_{u\leq \mathfrak{t}}\right)|\sigma(I_0)]\\

&=&\ind_B(0)\mathbb{E}[\frac{1}{\mathbb{P}[u<\mathfrak{t}]}\ind_{u<\mathfrak{t}}-\frac{1}{\mathbb{P}[u\leq \mathfrak{t}]}\ind_{u\leq \mathfrak{t}}|\sigma(I_0)]\\

&=&\ind_B(0)\mathbb{E}[\frac{1}{\mathbb{P}[u<\mathfrak{t}]}\ind_{u<\mathfrak{t}}-\frac{1}{\mathbb{P}[u\leq \mathfrak{t}]}\ind_{u\leq \mathfrak{t}}]\\
&=&0
\dce
$$
Consequently$$
\mathbb{E}_{\nu^0}[\ind_A\ind_B\left([M,[Z,v]]_t-[M(\iota),[Z,v]]_s\right)]
=0.
$$
Consider the second term.
$$
\dcb
&&\mathbb{E}_{\nu^0}[\ind_A\ind_B\int_s^tZ_{u-}d[M,v]_u]\\

&=&\mathbb{E}[\int\pi_{0,I/I}(I,dx)\int\pi_{0,F/F}(\mathsf{i},d\omega)\ind_A\ind_B\int_s^tZ_{u-}d[M,v]_u]\\

&=&\mathbb{E}[\int\pi_{0,I/I}(I,dx)\ind_B\int\pi_{0,F/F}(\mathsf{i},d\omega)\ind_A\sum_{s<u\leq t}Z_{u-}\Delta_uv\Delta_uM]\\

&=&\mathbb{E}[\int\pi_{0,I/I}(I,dx)\ind_B\int\pi_{0,F/F}(\mathsf{i},d\omega)\ind_A\int_s^tZ_{u-}\Delta_uvdM_u]\\
&&\mbox{ because under $\pi_{0,F/F}(\mathsf{i},d\omega)$ $v$ is a deterministic process}\\
&=&0
\ \mbox{ because $M$ is martingale under $\pi_{0,F/F}(\mathsf{i},d\omega)$}\\

\dce
$$
Consider the first term
$$
\dcb
&&\mathbb{E}_{\nu^0}[\ind_A\ind_B\int_s^tv_{u-}d[M,Z]_u]\\

&=&\mathbb{E}[\int\pi_{0,I/I}(I,dx)\int\pi_{0,F/F}(\mathsf{i},d\omega)\ind_A\ind_B\int_s^tv_{u-}d[M,Z]_u]\\

&=&\mathbb{E}[\int\pi_{0,I/I}(I,dx)\ind_B\int\pi_{0,F/F}(\mathsf{i},d\omega)\ind_A\int_s^tv_{u-}d[M,Z]_u]\\

&=&\mathbb{E}[\int\pi_{0,I/I}(I,dx)\ind_B\int\pi_{0,F/F}(\mathsf{i},d\omega)\ind_A\int_s^tv_{u-}d\cro{M,Z}^{\mathbb{P}\cdot\mathbb{F}}_u]\\
&&\mbox{ because under $\pi_{0,F/F}(\mathsf{i},d\omega)$ $v$ is a deterministic process}\\
&&\mbox{ while $M$ and $Z$ has a same behavior as under $(\mathbb{P},\mathbb{F})$}\\
\dce
$$
Combining these three terms, we obtain
$$
\dcb
&&\mathbb{E}_{\nu^0}[\ind_A\ind_B\int_s^td[M,Zv]_u]\\

&=&\mathbb{E}[\int\pi_{0,I/I}(I,dx)\ind_B\int\pi_{0,F/F}(\mathsf{i},d\omega)\ind_A\int_s^tv_{u-}d\cro{M,Z}^{\mathbb{P}\cdot\mathbb{F}}_u]\\

&=&\mathbb{E}_{\nu^0}[\ind_A\ind_B\int_s^t v_{u-} d\cro{M,N}^{\mathbb{P}\cdot\mathbb{F}}_u]\\

\dce
$$
By a usual limit argument, we conclude from this identity that $$
\cro{M,Zv}^{\nu^0\cdot\mathbb{J}_+}=v_-\centerdot\cro{M,N}^{\mathbb{P}\cdot\mathbb{F}}
$$
We can then compute 
$$
\dcb
d\cro{M,\alpha^0\ind_{[0,\rho)}}^{\nu^0\cdot\mathbb{J}_+}_t
&=&d\cro{M,Zv}^{\nu^0\cdot\mathbb{J}_+}_t\\

&=& v_{t-} d\cro{M,N}^{\mathbb{P}\cdot\mathbb{F}}_t\\

&=& \frac{1}{\mathbb{P}[t\leq \mathfrak{t}]}\ind_{\{t\leq \rho\}}d\cro{M,N}^{\mathbb{P}\cdot\mathbb{F}}_t\\

&=&\frac{Z_{t-}}{\mathbb{P}[t\leq \mathfrak{t}]} \ind_{\{t\leq \rho\}}\frac{1}{Z_{t-}}d\cro{M,N}^{\mathbb{P}\cdot\mathbb{F}}_u\\

&=&\alpha^0_{t-}\frac{1}{Z_{t-}}d\cro{M,N}^{\mathbb{P}\cdot\mathbb{F}}_u
\dce
$$
and we obtain
$$
\ind_{(0,\rho(\phi)]}\frac{1}{\alpha^0_-(\phi)}\centerdot \left(\cro{M(\iota),\alpha^0\ind_{[0,\rho)}}^{\nu^0\cdot\mathbb{J}_+}\right)\circ\phi
=\ind_{(0,\mathfrak{t}]}\frac{1}{Z_-}\centerdot \cro{M,N}^{\mathbb{P}\cdot\mathbb{F}}
$$

Theorem \ref{theorem:our-formula} is applicable. $\ind_{(0,\rho(\phi)]}\centerdot M$ is then the sum of a $(\mathbb{P},\mathbb{G})$ local martingale and the following process
$$
\dcb
&&\ind_{(0,\rho(\phi)]}\frac{1}{\alpha^0_-(\phi)}\centerdot \left(\cro{M(\iota),\alpha^0\ind_{[0,\rho)}}^{\nu^0\cdot\mathbb{J}_+	}\right)\circ\phi
+\ind_{\{0<\rho(\phi)\}}\Delta_{\rho(\phi)}M\ind_{[\rho(\phi),\infty)}\\

&=&
\ind_{(0,\mathfrak{t}]}\frac{1}{Z_-}\centerdot \cro{M,N}^{\mathbb{P}\cdot\mathbb{F}}
+\ind_{\{0<\mathfrak{t}\}}\Delta_{\mathfrak{t}}M\ind_{[\mathfrak{t},\infty)}
\dce
$$
This being done, to obtain the drift of $M$ in $\mathbb{G}$, we need only to compute the $(\mathbb{P},\mathbb{G})$ dual projection of the process$$
\ind_{\{0<\rho(\phi)\}}\Delta_{\rho(\phi)}M\ind_{[\rho(\phi),\infty)}
=\ind_{\{0<\mathfrak{t}\}}\Delta_{\mathfrak{t}}M\ind_{[\mathfrak{t},\infty)}
$$
which is given by $\ind_{(0,\mathfrak{t}]}\frac{1}{Z_-}\centerdot \mathsf{B}^M$, according to \cite[Lemme 4.]{JY2}. Theorem \ref{avant-tau} is proved.

\subsection{Honest time}\label{honest}

We consider the same situation as in the subsection \ref{progressive}. But we assume in addition that $\mathfrak{t}$ is a honest time (see \cite[Proposition (5.1)]{Jeulin80}). This means that there exists an increasing $\mathbb{F}$ adapted left continuous process $A$ such that $A_t\leq t\wedge \mathfrak{t}$ and $\mathfrak{t}=A_t$ on $\{\mathfrak{t}<t\}$. Let $\omega\in\{\mathfrak{t}\leq t\}$. Then, for any $\epsilon>0$, $\omega\in\{\mathfrak{t}< t+\epsilon\}$. So, $\mathfrak{t}(\omega)=A_{t+\epsilon}(\omega)$. Let $\epsilon\downarrow 0$. We see that $\mathfrak{t}(\omega)=A_{t+}(\omega)$. In fact, we can take$$
A_t=\sup\{0\leq s<t: \widetilde{Z}_s=1\},\ t>0,
$$
where $\widetilde{Z}={}^{\mathbb{P}\cdot\mathbb{F}-o}(\ind_{[0,\mathfrak{t}]})$.

We will now establish the formula of enlargement of filtration for honest time, using the local solution method. As in subsection \ref{progressive}, we assume that $A$ and $\widetilde{Z}$ are $\mathbb{F}^\circ$ adapted. The formula before $\mathfrak{t}$ has already been proved in Theorem \ref{avant-tau}. Let us consider the formula after $\mathfrak{t}$.

Fix $u>0$. For $t\geq u$, $$
\sigma(I_s:s\leq t)\cap\{\mathfrak{t}< t\}
=\sigma(\{\mathfrak{t}< s\}:s\leq t)\cap\{\mathfrak{t}< t\}
=\sigma(\mathfrak{t})\cap\{\mathfrak{t}< t\}
$$
So, for any bounded function $f\in\mathcal{I}_t$,$$
\dcb
&&\int\pi_{u,I/I}(I,dx)f(x)\ind_{\{\rho(x)< u\}}
=\mathbb{E}[f(I)\ind_{\{\mathfrak{t}< u\}}|\sigma(I_s:s\leq u)]\\
&=&\mathbb{E}[f(I^u)\ind_{\{\mathfrak{t}< u\}}|\sigma(I_s:s\leq  u)]\\
&=&f(I^u)\ind_{\{\mathfrak{t}< u\}}=f(I)\ind_{\{\rho(I)< u\}}\\
&=&\int f(x)\ind_{\{\rho(x)< u\}} \delta_{I}(dx)
\dce
$$
As for $\pi_{t,I/F}$, we have$$
\dcb
&&\int\pi_{t,I/F}(I,dx)f(x)\ind_{\{\rho(x)< u\}}
=\mathbb{E}[f(I)\ind_{\{\mathfrak{t}< u\}}|\mathcal{F}^\circ_t]\\
&=&\mathbb{E}[f(\ind_{[\mathfrak{t},\infty)})\ind_{\{\mathfrak{t}< u\}}|\mathcal{F}^\circ_t]\\
&=&\mathbb{E}[f(\ind_{[A_u,\infty)})\ind_{\{\mathfrak{t}< u\}}|\mathcal{F}^\circ_t]\\
&=&f(\ind_{[A_u,\infty)})\mathbb{E}[\ind_{\{\mathfrak{t}< u\}}|\mathcal{F}^\circ_t]\\
\dce
$$
Let $m^u_t=\mathbb{E}[\ind_{\{\mathfrak{t}< u\}}|\mathcal{F}^\circ_t], t\geq u$. We continue
$$
\dcb
&&\int\pi_{t,I/F}(I,dx)f(x)\ind_{\{\rho(x)< u\}}\ind_{\{m^u_u>0\}}\\
&=&f(\ind_{[A_u,\infty)})m^u_t\ind_{\{m^u_u>0\}}\\
&=&f(\ind_{[A_u,\infty)})\frac{m^u_t}{m^u_u}\ind_{\{m^u_u>0\}}m^u_u\\
&=&\int\pi_{u,I/F}(I,dx)f(x)\frac{m^u_t}{m^u_u}\ind_{\{m^u_u>0\}}\ind_{\{\rho(x)< u\}}\\
&=&\int\pi_{u,I/F}(I,dx'')\int\pi_{u,I/I}(x'',dx)f(x)\frac{m^u_t}{m^u_u}\ind_{\{m^u_u>0\}}\ind_{\{\rho(x)< u\}}\\
\dce
$$
By Lemma \ref{alpha}, we conclude that $\tau^u\ind_{\{\rho<u\}}\ind_{\{m^u_u>0\}}= \infty\ind_{\{\rho<u\}}\ind_{\{m^u_u>0\}}$ under $(\mu+\nu^u)$ and $$
\alpha^u_t(\omega,x)\ind_{\{m^u_u(\omega)>0\}}\ind_{\{\rho(x)<u\}}=\frac{m^u_{t+}(\omega)}{m^u_u(\omega)}\ind_{\{m^u_u(\omega)>0\}}\ind_{\{\rho(x)<u\}}, t\geq u
$$ 
It is straightforward to compute $\ind_{\{m^u_u(\iota)>0\}}\ind_{\{\rho<u\}}\ind_{(u,\infty)}\centerdot\cro{M(\iota),\alpha^u}^{\nu^u\cdot\mathbb{J}_+}$ which is given by$$
\ind_{\{m^u_u(\iota)>0\}}\ind_{\{\rho<u\}}\ind_{(u,\infty)}\centerdot\cro{M(\iota),\alpha^u}^{\nu^u\cdot\mathbb{J}_+}
=\ind_{\{\rho<u\}}\frac{1}{m^u_u(\iota)}\ind_{\{m^u_u(\iota)>0\}}\ind_{(u,\infty)}\centerdot\cro{M,m^u_+}^{\mathbb{P}\cdot\mathbb{F}}(\iota)
$$
Let us compute $\cro{M,m^u_+}^{\mathbb{P}\cdot\mathbb{F}}$. We begin with $m^u_+$. It is a bounded $(\mathbb{P},\mathbb{F})$ martingale. For $t\geq u$, $$
\dcb
m^u_{t+}&=&\mathbb{E}[\ind_{\{\mathfrak{t}< u\}}|\mathcal{F}_t]\\
&=&\mathbb{E}[\ind_{\{\mathfrak{t}< u\}}\ind_{\{\mathfrak{t}\leq  t\}}|\mathcal{F}_t]\\
&=&\mathbb{E}[\ind_{\{A_{t+}< u\}}\ind_{\{\mathfrak{t}\leq t\}}|\mathcal{F}_t]\\
&=&\ind_{\{A_{t+}< u\}}\mathbb{E}[\ind_{\{\mathfrak{t}\leq t\}}|\mathcal{F}_t]\\
&=&\ind_{\{A_{t+}< u\}}(1-Z_t)\\
&&\mbox{ or }\\
&=&\mathbb{E}[\ind_{\{\mathfrak{t}< u\}}\ind_{\{\mathfrak{t}<  t\}}|\mathcal{F}_t]\\
&=&\mathbb{E}[\ind_{\{A_{t}< u\}}\ind_{\{\mathfrak{t}< t\}}|\mathcal{F}_t]\\
&=&\ind_{\{A_{t}< u\}}\mathbb{E}[\ind_{\{\mathfrak{t}< t\}}|\mathcal{F}_t]\\
&=&\ind_{\{A_{t}< u\}}(1-\widetilde{Z}_t)\\

\dce
$$
Note that, if $t\in\{\widetilde{Z}=1\}$, $A_{t+\epsilon}\geq t$ for any $\epsilon>0$ so that $A_{t+}\geq t$. But $A_t\leq t$ so that $A_{t+}=t$. Let $T_u=\inf\{s\geq u: A_{s+}\geq u\}$. On the set $\{A_{u+}=u\}$, $T_u=u$. On the set $\{A_{u+}<u\}$, $T_u=\inf\{s\geq u: s\in\{\widetilde{Z}=1\}\}$, where we have $A_t=A_u$ for $u<t<T_u$, and $A_{T_u+}=T_u>u$ if $T_u<\infty$. Note that, according to [Jeulin Lemme(4.3)] $\{\widetilde{Z}=1\}$ is a closed set. It follows that $T_u\in\{\widetilde{Z}=1\}$ whenever $T_u<\infty$. Therefore, on the set $\{T_u<\infty\}$, $$
1-Z_{T_u}=\widetilde{Z}_{T_u}-Z_{T_u}={}^{\mathbb{P}\cdot\mathbb{F}-o}(\ind_{[\mathfrak{t}]})_{T_n}=\mathbb{E}[\ind_{\{\mathfrak{t}=T_u\}}|\mathcal{F}_{T_u}]
$$ 
Set $w^u_t=\ind_{\{A_{t+}< u\}}=\ind_{\{t<T_u\}}$ (For the last equality, it is evident if $t<u$. It is already explained for $t=u$. For $t>u$, the condition $A_{t+}<u$ is equivalent to say $A_{u+}<u$ and $A_{t+\epsilon}<u$ for some $\epsilon>0$, which implies $T_u\geq t+\epsilon$. Inversely, if $T_u> t$, we will have $A_{u+}<u$ so that $(u-\epsilon,u]\cap\{\widetilde{Z}=1\}=\emptyset$ for some $\epsilon>0$, because $\{\widetilde{Z}=1\}$ is a closed set. This means $A_{t+}\leq u-\epsilon$.). We compute for $t>u$$$
d(w^u(1-Z))_t=-w^u_{t-}dZ_t-(1-Z_{T_u-})\ind_{\{T_u\leq t\}}+\Delta_{T_u}Z\ind_{\{T_u\leq t\}}
=-w^u_{t-}dZ_t-(1-Z_{T_u})\ind_{\{T_u\leq t\}}.
$$
Consequently
$$
\dcb
d[M,m^u_+]_t
&=&-w^u_{t-}d[M,Z]_t-\Delta_{T_u}M(1-Z_{T_u-})\ind_{\{T_u\leq t\}}+\Delta_{T_u}M\Delta_{T_u}Z\ind_{\{T_u\leq t\}}\\

&=&-w^u_{t-}d[M,Z]_t-\Delta_{T_u}M(1-Z_{T_u})\ind_{\{T_u\leq t\}}
\dce
$$
Let $v$ denote the $(\mathbb{P},\mathbb{F})$ predictable dual projection of the increasing process $\Delta_{T_u}M(1-Z_{T_u})\ind_{\{u<T_u\leq t\}}, t>u$. Let $\widetilde{A}^{\mathfrak{t}}$ (resp. $A^{\mathfrak{t}}$) denotes the $(\mathbb{P},\mathbb{F})$ optional (resp. predictable) dual projection of $\ind_{\{0<\mathfrak{t}\}}\ind_{[\mathfrak{t},\infty)}$. For any bounded $\mathbb{F}$ predictable process $H$, we have
$$
\dcb
\mathbb{E}[\int_u^\infty H_udv_u]
&=&\mathbb{E}[H_{T_u}\Delta_{T_u}M(1-Z_{T_u})\ind_{\{u<T_u<\infty\}}]\\
&=&\mathbb{E}[H_{T_u}\Delta_{T_u}M\mathbb{E}[\ind_{\{\mathfrak{t}=T_u\}}|\mathcal{F}_{T_u}]\ind_{\{u<T_u<\infty\}}]\\
&=&\mathbb{E}[\mathbb{E}[H_{T_u}\Delta_{T_u}M\ind_{\{\mathfrak{t}=T_u\}}|\mathcal{F}_{T_u}]\ind_{\{u<T_u<\infty\}}]\\
&=&\mathbb{E}[\mathbb{E}[H_{\mathfrak{t}}\Delta_{\mathfrak{t}}M\ind_{\{\mathfrak{t}=T_u\}}|\mathcal{F}_{T_u}]\ind_{\{u<T_u<\infty\}}]\\

&=&\mathbb{E}[H_{\mathfrak{t}}\Delta_{\mathfrak{t}}M\ind_{\{\mathfrak{t}=T_u\}}\ind_{\{u<T_u<\infty\}}]\\

&=&\mathbb{E}[H_{\mathfrak{t}}\Delta_{\mathfrak{t}}M\ind_{\{\mathfrak{t}=T_u\}}\ind_{\{u<\mathfrak{t}<\infty\}}]\\

&=&\mathbb{E}[H_{\mathfrak{t}}\Delta_{\mathfrak{t}}M\ind_{\{\mathfrak{t}\leq T_u\}}\ind_{\{u<\mathfrak{t}<\infty\}}]\\
&&\mbox{ because $T_u\in\{\widetilde{Z}=1\}$ and $\mathfrak{t}$ is the last instant of the set $\{\widetilde{Z}=1\}$}\\
&=&\mathbb{E}[\int_0^\infty H_{s}\Delta_{s}M\ind_{(u,T_u]}(s)d(\ind_{[\mathfrak{t},\infty)})_s]\\

&=&\mathbb{E}[\int_0^\infty H_{s}\Delta_{s}M\ind_{(u,T_u]}(s)d\widetilde{A}^{\mathfrak{t}}_s]\\

&=&\mathbb{E}[\int_0^\infty H_{s}\Delta_{s}M\ind_{(u,T_u]}(s)d(\widetilde{A}^{\mathfrak{t}}-A^{\mathfrak{t}})_s]\\
&&\mbox{ because $M$ is a bounded martingale and $A^{\mathfrak{t}}$ is predictable}\\

&=&\mathbb{E}[\int_u^{T_u} H_{s}d[M,\widetilde{A}^{\mathfrak{t}}-A^{\mathfrak{t}}]_s]\\

&=&\mathbb{E}[\int_u^{T_u} H_{s}d\cro{M,\widetilde{A}^{\mathfrak{t}}-A^{\mathfrak{t}}}^{\mathbb{P}\cdot\mathbb{F}}_s]\\
\dce
$$
Putting these results together, we can write
$$
\dcb
&&\ind_{\{m^u_u(\iota)>0\}}\ind_{\{\rho(x)<u\}}\ind_{(u,\tau^u]}\frac{1}{\alpha^u_-}\centerdot \cro{M(\iota),\alpha^u}^{\nu^u\cdot\mathbb{J}_+}\\
&=&
\ind_{\{\rho(x)<u\}}\left(\ind_{(u,\infty)}\frac{m^u_u}{m^u_-} \frac{1}{m^u_u}\ind_{\{m^u_u>0\}}\centerdot\cro{M,m^u_+}^{\mathbb{P}\cdot\mathbb{F}}\right)(\iota)\\

&=&
\ind_{\{\rho(x)<u\}}\left(\ind_{(u,\infty)} \ind_{\{m^u_u>0\}}\frac{1}{m^u_-}\centerdot\cro{M,m^u_+}^{\mathbb{P}\cdot\mathbb{F}}\right)(\iota)\\

&=&
\ind_{\{\rho(x)<u\}}\left(\ind_{(u,\infty)} \ind_{\{m^u_u>0\}}\frac{1}{m^u_-}\centerdot(-\ind_{(u,T_u]}\centerdot\cro{M,N}^{\mathbb{P}\cdot\mathbb{F}}-\ind_{(u,T_u]}\centerdot\cro{M,\widetilde{A}^{\mathfrak{t}}-A^{\mathfrak{t}}}^{\mathbb{P}\cdot\mathbb{F}})\right)(\iota)\\

&=&
-\ind_{\{\rho(x)<u\}}\left(\ind_{(u,\infty)} \ind_{\{m^u_u>0\}}\ind_{(u,T_u]}\frac{1}{m^u_-}\centerdot(\cro{M,N}^{\mathbb{P}\cdot\mathbb{F}}+\cro{M,\widetilde{A}^{\mathfrak{t}}-A^{\mathfrak{t}}}^{\mathbb{P}\cdot\mathbb{F}})\right)(\iota)\\

&=&
-\ind_{\{\rho(x)<u\}}\left( \ind_{\{m^u_u>0\}}\ind_{(u,T_u]}\frac{1}{1-Z_-}\centerdot(\cro{M,N}^{\mathbb{P}\cdot\mathbb{F}}+\cro{M,\widetilde{A}^{\mathfrak{t}}-A^{\mathfrak{t}}}^{\mathbb{P}\cdot\mathbb{F}})\right)(\iota)
\dce
$$
Pull this identity back to $\Omega$, we have
$$
\dcb
&&\ind_{\{m^u_u>0\}}\ind_{\{\rho(\phi)<u\}}\ind_{(u,\tau^u(\phi)]}\frac{1}{\alpha^u_-(\phi)}\centerdot \cro{M(\iota),\alpha^u}^{\nu^u\cdot\mathbb{J}_+}(\phi)\\

&=&-\ind_{\{\mathfrak{t}<u\}} \ind_{\{m^u_u>0\}}\ind_{(u,T_u]}\frac{1}{1-Z_-}\centerdot(\cro{M,N}^{\mathbb{P}\cdot\mathbb{F}}+\cro{M,\widetilde{A}^{\mathfrak{t}}-A^{\mathfrak{t}}}^{\mathbb{P}\cdot\mathbb{F}})\\

&=&-\ind_{\{\mathfrak{t}<u\}} \ind_{\{1-\widetilde{Z}_u>0\}}\ind_{(u,T_u]}\frac{1}{1-Z_-}\centerdot(\cro{M,N}^{\mathbb{P}\cdot\mathbb{F}}+\cro{M,\widetilde{A}^{\mathfrak{t}}-A^{\mathfrak{t}}}^{\mathbb{P}\cdot\mathbb{F}})\\

&=&-\ind_{\{\mathfrak{t}<u\}} \ind_{(u,\infty)}\frac{1}{1-Z_-}\centerdot(\cro{M,N}^{\mathbb{P}\cdot\mathbb{F}}+\cro{M,\widetilde{A}^{\mathfrak{t}}-A^{\mathfrak{t}}}^{\mathbb{P}\cdot\mathbb{F}})
\dce
$$
because on the set $\{\mathfrak{t}<u\}$, $1-\widetilde{Z}_u>0$ and $T_u=\infty$.

Now we can apply Theorem \ref{theorem:our-formula}, according to which (noting that, according to \cite[Lemme(4.3) ii) and p.63]{Jeulin80}, $1-Z_->0$ on $(\mathfrak{t},\infty)$) $\ind_{\{\rho(\phi)<u\}}\ind_{\{m^u_u>0\}}\ind_{(u,\tau^u(\phi)]}\centerdot M$ is the sum of a $(\mathbb{P},\mathbb{G})$ local martingale and the following process with finite variation: 
$$
\ind_{\{\rho(\phi)<u\}}\ind_{\{m^u_u>0\}}\left(\ind_{(u,\tau^u(\phi)]}\frac{1}{\alpha^u_-(\phi)}\centerdot \left(\cro{M(\iota),\alpha^u}^{\nu^u\cdot\mathbb{J}_+}\right)\circ\phi
+\ind_{\{u<\tau^u(\phi)\}}\Delta_{\tau^u(\phi)}M\ind_{[\tau^u(\phi),\infty)}\right)
$$
In particular, noting that, on the set $\{m^u_u>0, \rho(\phi)<u\}=\{\mathfrak{t}<u\}$, $\tau(\phi)=\infty$ so that $$
\ind_{\{0<\tau^u(\phi)\}}\Delta_{\tau^u(\phi)}M\ind_{[\tau^u(\phi),\infty)}=0,
$$
the above formula implies that $\ind_{\{\mathfrak{t}<u\}}\ind_{(u,\infty)}\centerdot M$ is a $(\mathbb{P},\mathbb{G})$ special semimartingale whose drift is given by 
$$
\ind_{\{\mathfrak{t}<u\}}\ind_{(u,\infty)}\frac{1}{\alpha^u_-(\phi)}\centerdot \left(\cro{M(\iota),\alpha^u}^{\nu^u\cdot\mathbb{J}_+}\right)\circ\phi
=
-\ind_{\{\mathfrak{t}<u\}} \ind_{(u,\infty)}\frac{1}{1-Z_-}\centerdot(\cro{M,N}^{\mathbb{P}\cdot\mathbb{F}}+\cro{M,\widetilde{A}^{\mathfrak{t}}-A^{\mathfrak{t}}}^{\mathbb{P}\cdot\mathbb{F}})
$$

For $u\in\mathbb{Q}_+^*$, set $\mathtt{B}_u=\{\mathfrak{t}<u\}\cap(u,\infty)$. The $\mathtt{B}_u$ are $\mathbb{G}$ predictable left intervals. For every $u\in\mathbb{Q}_+^*$, $\ind_{\mathtt{B}_u}\centerdot M$ is a $(\mathbb{P,\mathbb{G}})$ special semimartingale with the above given drift. The first and the second conditions in Assumption \ref{partialcovering} are satisfied. 

The random measure generated by the above drifts on $\cup_{u\in\mathbb{Q}_+^*}\mathtt{B}_u=(\mathfrak{t},\infty)$ is$$
-\ind_{(\mathfrak{t},\infty)} \frac{1}{1-Z_-}\centerdot(\cro{M,N}^{\mathbb{P}\cdot\mathbb{F}}+\cro{M,\widetilde{A}^{\mathfrak{t}}-A^{\mathfrak{t}}}^{\mathbb{P}\cdot\mathbb{F}})
$$
Knowing the ${}^{\mathbb{P}\cdot\mathbb{F}-p}(\ind_{(\mathfrak{t},\infty)})=(1-Z_-)$, it can be checked that this random measure has a distribution function. Lemma \ref{UB} is applicable. It yields that $\ind_{(\mathfrak{t},\infty)}\centerdot M=\ind_{\cup_{u\in\mathbb{Q}_+^*}\mathtt{B}_u}\centerdot M$ is a $(\mathbb{P,\mathbb{G}})$ special semimartingale with the drift
$$
-\ind_{(\mathfrak{t},\infty)} \frac{1}{1-Z_-}\centerdot(\cro{M,N}^{\mathbb{P}\cdot\mathbb{F}}+\cro{M,\widetilde{A}^{\mathfrak{t}}-A^{\mathfrak{t}}}^{\mathbb{P}\cdot\mathbb{F}})
$$
We have thus proved the formula of enlargement of filtration for the honest time $\mathfrak{t}$ (see \cite[Théorème(5.10), Lemme(5.17)]{Jeulin80}).

\bethe\label{honest}
Let $\mathbb{G}$ be the progressive enlargement of $\mathbb{F}$ with a honest time $\mathfrak{t}$. For all bounded $(\mathbb{P},\mathbb{F})$ martingale $M$, it is a special $(\mathbb{P},\mathbb{G})$ semimartingale whose drift is given by $$
\ind_{(0,\mathfrak{t}]}\frac{1}{Z_-}\centerdot (\cro{N,M}^{\mathbb{P}\cdot \mathbb{F}}+\mathsf{B}^M)
-\ind_{(\mathfrak{t},\infty)} \frac{1}{1-Z_-}\centerdot(\cro{M,N}^{\mathbb{P}\cdot\mathbb{F}}+\cro{M,\widetilde{A}^{\mathfrak{t}}-A^{\mathfrak{t}}}^{\mathbb{P}\cdot\mathbb{F}})
$$
\ethe

\subsection{Girsanov theorem and formulas of enlargement of filtration}

We have proved the classical formulas with the local solution method. Compared with the initial proofs of the classical formulas, the local solution method is made on a product space which bear more structure and provides some interesting details, especially when we compute $\alpha^u$ and $\cro{M(\iota),\alpha^u\ind_{[u,\rho)}}^{\nu^u\cdot\mathbb{J}_+}$.

The local solution method shows another interesting point : the Girsanov's theorem in the enlargement of filtration. When one look at the classical examples of enlargement of filtration, one notes immediately some resemblance between the enlargement of filtration formula and the Girsanov theorem formula. It was a question at what extend these two themes are related. Clearly, Jacod's criterion is an assumption of type of Girsanov theorem. As for the progressive enlargement of filtration, \cite{yoeurp} gives a proof of the formula using Girsanov theorem in a framework of Föllmer's measure. However, the use of Föllmer's measure makes that interpretation in \cite{yoeurp} inadequate to explain the phenomenon in the other examples such as the case of honest time. In fact, not all examples of enlargement of filtration give a formula of type of Girsanov theorem (see Section \ref{futur-infermum}). However, locally, if we can compute the function $\alpha^u$ as in Section \ref{alpha}, we will have a formula on the interval $(u,\tau^u(\phi)]$ which is of type of Girsanov theorem, because Theorem \ref{theorem:our-formula} is a generalized Girsanov theorem (cf. \cite{kunita, yoeurp}). The local computation of $\alpha^u$ is possible in the case a honest time as it is shown in subsection \ref{honest}. That is also the case in most of the classical examples.

\

\

\section{A different proof of \cite[Théorème(3.26)]{Jeulin80}}

Théorème(3.26) in \cite{Jeulin80} is a result on initial enlargement of filtration. It is interesting to note that this result was considered as an example where general methodology did not apply. We will give a proof using the local solution method. 

Let $\Omega$ be the space $C_0(\mathbb{R}_+,\mathbb{R})$ of all continuous real functions on $\mathbb{R}_+$ null at zero. The space $\Omega$ is equipped with the Wiener measure $\mathbb{P}$. Let $X$ be the identity function on $\Omega$, which is therefore a linear brownian motion issued from the origin. Set $U_t=\sup_{s\leq t}X_s$. Let $\mathbb{F}^\circ$ be the natural filtration of $X$. Define the process $I$ by $I_t=U, \forall t\geq 0$ ($I$ is a constant process taking values in the space $\Omega$). Define the filtration $\mathbb{G}^\circ$ as in Subsection \ref{h-formula}, as well as the filtrations $\mathbb{F}$ and $\mathbb{G}$. 

For $z\in\Omega$, denote by $z_t$ its coordinate at $t\geq 0$ and by $\mathfrak{s}(z)$ the function $\sup_{v\leq t}z_v, t\geq 0$. Let $\rho_u(z)=\inf\{s\geq u: \mathfrak{s}_s(z)>\mathfrak{s}_u(z)\}$. Let $\mathfrak{T}_u=\rho_u(U)$ ($\rho_u$ is considered as a map without randomness, while $\mathfrak{T}_u$ is a random variable). For $x\in C(\mathbb{R}_+,\Omega)$, $x_t$ denotes the coordinate of $x$ at $t$. Let $\mathfrak{l}$ is the identity function on $\mathbb{R}_+$. 

\bethe\label{sup-adjonction}\cite[Théorème(3.26)]{Jeulin80}
Let $M$ be a $(\mathbb{P},\mathbb{F})$ martingale. Then, $M$ is a $\mathbb{G}$ semimartingale, if and only if the $\int_0^t\frac{1}{U_s-X_s}|d\cro{M,X}^{\mathbb{P}\cdot\mathbb{F}}_s|<\infty$ for all $t>0$. In this case, the drift of $M$ is given by $$
-\frac{1}{U-X}\left(1-\frac{(U-X)^2}{\mathfrak{T}-\mathfrak{l}}\right)\centerdot\cro{M,X}^{\mathbb{P}\cdot\mathbb{F}}
$$ 
\ethe

We use the notation $\mathfrak{c}_\xi$ introduced in subsection \ref{jcrit}. Note that $I=\mathfrak{c}_{U}$ and hence $
\sigma(I)_\infty=\sigma(U)=\sigma(I)_0.
$
We have $\pi_{u,I/I}(x'',dx)=\ind_{\{x''\in E\}}\delta_{x''}(dx)$. To compute $\pi_{t,I/F}$ we introduce $B$ to be the process $B_t=X_{\mathfrak{T}_u+t}-X_{\mathfrak{T}_u}$ for $t\geq 0$. We have the decomposition, $t\geq 0$,$$
U_t=U_t\ind_{\{t<\mathfrak{T}_u\}}+U_{t-\mathfrak{T}_u+\mathfrak{T}_u}\ind_{\{\mathfrak{T}_u\leq t\}}
=U_t^{u}\ind_{\{t<\mathfrak{T}_u\}}+(\mathfrak{s}_{t-\mathfrak{T}_u}(B)+U_u)\ind_{\{\mathfrak{T}_u\leq t\}}
=U_t^{u}+\mathfrak{s}_{t-\mathfrak{T}_u}(B)\ind_{\{\mathfrak{T}_u\leq t\}}
$$
This yields$$
\sigma(U)=\sigma(U^u, \mathfrak{T}_u, \mathfrak{s}(B))
$$
We note that $\mathfrak{s}(B)$ is independent of $\mathcal{F}_{\mathfrak{T}_u}$. Therefore, for any bounded borel function $f$ on $C(\mathbb{R}_+,\Omega)$, there exists a suitable borel function $G$ such that $$
\mathbb{E}_{\mathbb{P}}[f(\mathfrak{c}_U)|\mathcal{F}_{\mathfrak{T}_u}]
=G(U^u,\mathfrak{T}_u)
$$
Now we can compute the function $h^u_t$ for $t\geq u$. 
$$
\dcb
&&\int\pi_{t,I/F}(X,dx)f(x)\ind_{t\leq\rho_u(x_0)}\\ 
&=& \mathbb{E}_{\mathbb{P}}[f(I)\ind_{t\leq \mathfrak{T}_u}|\mathcal{F}_t^\circ]\\
&=&\mathbb{E}_{\mathbb{P}}[f(\mathfrak{c}_U)\ind_{t\leq \mathfrak{T}_u}|\mathcal{F}_t^\circ]\\

&=&\mathbb{E}_{\mathbb{P}}[\mathbb{E}_{\mathbb{P}}[f(\mathfrak{c}_U)|\mathcal{F}_{\mathfrak{T}_u}]\ind_{t\leq \mathfrak{T}_u}|\mathcal{F}_t^\circ]\\

&=&\mathbb{E}_{\mathbb{P}}[G(U^u, \mathfrak{T}_u)\ind_{t\leq \mathfrak{T}_u}|\mathcal{F}_t^\circ]\\

&=&\mathbb{E}_{\mathbb{P}}[G(U^u, T_{U_u}\circ\theta_t+t)\ind_{t\leq \mathfrak{T}_u}|\mathcal{F}_t^\circ]\\ &&\mbox{where $\theta$ is the usual translation operator and $T_a=\inf\{s\geq 0: X_s=a\}$}\\

&=&\int_0^\infty G(U^u, v+t)\frac{1}{\sqrt{2\pi}}\frac{(U_u-X_t)}{\sqrt{v}^3}e^{-\frac{(U_u-X_t)^2}{2v}}dv\ind_{t\leq \mathfrak{T}_u}\\

&&\mbox{ see \cite[p.80]{KS} for the density function of $T_a$}\\
&=&\int_0^\infty G(U^u, v+t)\frac{1}{\sqrt{2\pi}}\frac{(U_t-X_t)}{\sqrt{v}^3}e^{-\frac{(U_t-X_t)^2}{2v}}dv\ind_{t\leq \mathfrak{T}_u}\\

&=&\int_t^\infty G(U^u, s)\frac{1}{\sqrt{2\pi}}\frac{(U_t-X_t)}{\sqrt{s-t}^3}e^{-\frac{(U_t-X_t)^2}{2(s-t)}}ds\ind_{t\leq \mathfrak{T}_u}\\
\dce
$$
and, for a positive borel function $g$ on $\mathbb{R}_+$,
$$
\dcb
&&\int\pi_{u,I/F}(X,dx)f(x)g(\rho_u(x_0))\\ 
&=& \mathbb{E}_{\mathbb{P}}[f(I)g(\mathfrak{T}_u)|\mathcal{F}_u^\circ]\\
&=&\mathbb{E}_{\mathbb{P}}[f(\mathfrak{c}_U)g(\mathfrak{T}_u)|\mathcal{F}_u^\circ]\\

&=&\mathbb{E}_{\mathbb{P}}[G(U^u,\mathfrak{T}_u)g(\mathfrak{T}_u)|\mathcal{F}_u^\circ]\\

&=&\mathbb{E}_{\mathbb{P}}[G(U^u,T_{U_u}\circ\theta_u+u)g(T_{U_u}\circ\theta_u+u)|\mathcal{F}_u^\circ]\\

&=&\int_0^\infty G(U^u, v+u)g(v+u)\frac{1}{\sqrt{2\pi}}\frac{(U_u-X_u)}{\sqrt{v}^3}e^{-\frac{(U_u-X_u)^2}{2v}}dv\\

&=&\int_u^\infty G(U^u, s)g(s)\frac{1}{\sqrt{2\pi}}\frac{(U_u-X_u)}{\sqrt{s-u}^3}e^{-\frac{(U_u-X_u)^2}{2(s-u)}}ds\\
\dce
$$
Set $$
\Psi(y,a)=\frac{1}{\sqrt{2\pi}}\frac{y}{\sqrt{a}^3}e^{-\frac{y^2}{2a}},\ y>0, a>0
$$
The above two identities yields$$
\dcb
&&\int\pi_{t,I/F}(X,dx)f(x)\ind_{t\leq \rho_u(x_0)}\\
&=&\int_t^\infty G(U^u, s)\Psi(U_t-X_t,s-t)ds\ind_{t\leq \mathfrak{T}_u}\\

&=&\int_u^\infty G(U^u, s)\ind_{t\leq s}\Psi(U_t-X_t,s-t)ds\ind_{t\leq \mathfrak{T}_u}\\

&=&\int_u^\infty G(U^u, s)\ind_{t\leq s}\frac{\Psi(U_t-X_t,s-t)}{\Psi(U_u-X_u,s-u)}\frac{1}{\sqrt{2\pi}}\frac{(U_u-X_u)}{\sqrt{s-u}^3}e^{-\frac{(U_u-X_u)^2}{2(s-u)}}ds\ind_{t\leq \mathfrak{T}_u}\\

&=&
\int\pi_{u,I/F}(X,dx)f(x)\ind_{t\leq \rho_u(x_0)}\frac{\Psi(U_t-X_t,\rho_u(x_0)-t)}{\Psi(U_u-X_u,\rho_u(x_0)-u)}\ind_{t\leq \mathfrak{T}_u}\\

&=&
\int\pi_{u,I/F}(X,dx'')\int\pi_{u,I/I}(x'',dx)f(x)\ind_{t\leq \rho_u(x_0)}\frac{\Psi(U_t-X_t,\rho_u(x_0)-t)}{\Psi(U_u-X_u,\rho_u(x_0)-u)}\ind_{t\leq \mathfrak{T}_u}
\dce
$$
Let us identify $X(\iota)$ with $X$. According to Lemma \ref{alpha}, $\tau^u\geq \rho_u(x_0)$ and $$
\alpha^u_t\ind_{\{t<\rho_u(x_0)\}}
=\frac{\Psi(U_t-X_t,\rho_u(x_0)-t)}{\Psi(U_u-X_u,\rho_u(x_0)-u)}\ind_{\{t<\mathfrak{T}_u\}}\ind_{\{t<\rho_u(x_0)\}}
$$
Under the probability $\nu^u$, $X$ is a $\mathbb{J}_+$ brownian motion. Note that $\rho_u(x_0)$ is in $\mathcal{J}_0$ and $\mathfrak{T}_u$ is a stopping time in a brownian filtration. So, $\rho_u(x_0)\wedge \mathfrak{T}_u$ is a $\mathbb{J}_+$ predictable stopping time. The martingale part of $\alpha^u$ on $(u,\rho_u(x_0)\wedge \mathfrak{T}_u)$ is then given by$$
\dcb
&&-\frac{1}{\Psi(U_u-X_u,\rho_u(x_0)-u)}\frac{\partial}{\partial y}\Psi(U_t-X_t,\rho_u(x_0)-t) dX_t\\
&=&
-\frac{1}{\Psi(U_u-X_u,\rho_u(x_0)-u)}\Psi(U_t-X_t,\rho_u(x_0)-t)\frac{1}{U_t-X_t}(1-\frac{(U_t-X_t)^2}{\rho_u(x_0)-t}) dX_t\\
&=&
-\alpha^u_t\frac{1}{U_t-X_t}(1-\frac{(U_t-X_t)^2}{\rho_u(x_0)-t}) dX_t
\dce
$$
The continuity of the processes makes that the $\mathbb{J}_+$ predictable bracket between $M$ and $X$ under $\nu^u$ is the same as that in $(\mathbb{P},\mathbb{F})$. Theorem \ref{theorem:our-formula} is applicable. Hence, $\ind_{(u,\mathfrak{T}_u]}\centerdot M$ is a $(\mathbb{P},\mathbb{G})$ special semimartingale with drift given by $$
- \ind_{(u,\mathfrak{T}_u]}\frac{1}{U-X}(1-\frac{(U-X)^2}{\mathfrak{T}_u-\mathfrak{l}}) \centerdot \cro{M,X}
$$
Note that the above random measure has a distribution function, according to the argument in \cite[p.54]{Jeulin80}. Note also that we can replace $\frac{(U-X)^2}{\mathfrak{T}_u-\mathfrak{l}}$ by $\frac{(U-X)^2}{\mathfrak{T}_\cdot-\mathfrak{l}}$ in the above formula, because for $t\in(u,\mathfrak{T}_u]$, $\mathfrak{T}_u=\mathfrak{T}_t$.

Now we consider the random left intervals $\mathtt{B}_u=(u,\mathfrak{T}_u]$ for $u\in\mathbb{Q}_+$. The drift of $\ind_{(u,\mathfrak{T}_u]}\centerdot M$ generate a random measure $\mathsf{d}\chi^\cup$ on $\cup_{u\in\mathbb{Q}_+}(u,\mathfrak{T}_u]$, which is clearly diffuse. Note that $\cup_{u\in\mathbb{Q}_+}(u,\mathfrak{T}_u]\setminus \{U>x\}$ is a countable set. Hence $$
d\chi^\cup_t
=\ind_{\{U_t>X_t\}}\frac{1}{U_t-X_t}(1-\frac{(U_t-X_t)^2}{\mathfrak{T}_t-t}) d\cro{M,X}_t
$$ 
According to Lemma \ref{UB}, if $\mathsf{d}\chi^\cup$ has a distribution function, $\ind_{\cup_{u\in\mathbb{Q}_+}(u,\mathfrak{T}_u]}\centerdot M$ is a $(\mathbb{P},\mathbb{G})$ special semimartingale whose drift is the distribution function of $\mathsf{d}\chi^\cup$. We note that the set $\cup_{u\in\mathbb{Q}_+}(u,\mathfrak{T}_u]$ is a $\mathbb{F}$ predictable set so that the stochastic integral $\ind_{\cup_{u\in\mathbb{Q}_+}(u,\mathfrak{T}_u]}\centerdot M$ in $(\mathbb{P},\mathbb{G})$ is the same as that in $(\mathbb{P},\mathbb{F})$. But then, since $\ind_{(\cup_{u\in\mathbb{Q}_+}(u,\mathfrak{T}_u])^c}\centerdot\cro{X,X}=\ind_{\{U-X=0\}}\centerdot\cro{X,X}=0$ (recalling Levy's theorem that $U-X$ is the absolute value of a brownian motion), $\ind_{(\cup_{u\in\mathbb{Q}_+}(u,\mathfrak{T}_u])^c}\centerdot M=0$ according to the predictable representation theorem under Wiener measure. This means that $M\ (=\ind_{\cup_{u\in\mathbb{Q}_+}(u,\mathfrak{T}_u]}\centerdot M)$ itself is a $(\mathbb{P},\mathbb{G})$ special semimartingale with drift the distribution function of $\mathsf{d}\chi^\cup$. 

Inversely, if $M$ is a $(\mathbb{P},\mathbb{G})$ special semimartingale, clearly $\mathsf{d}\chi^\cup$ has a distribution function.

To achieve the proof of Theorem \ref{sup-adjonction}, we note that, according to \cite[Lemme(3.22)]{Jeulin80} $\mathsf{d}\chi^\cup$ has a distribution function if and only if $
\int_0^t\frac{1}{U_s-X_s}\left|d\cro{M,X}^{\mathbb{P}\cdot\mathbb{F}}_s\right|<\infty, \ t\geq 0.
$
We recall that $X$ does not satisfy that condition.

\

\section{A last passage time in \cite{JYC}}

The computation of $\alpha^u$ is based on absolute continuity of $\pi_{t,I/F}$ with respect to $\int \pi_{u,I/F}(\omega,dx'')\pi_{u,I/I}(x'',dx)$. The random time studied in this section is borrowed from \cite{JYC} proposed by Emery. We will compute its $\alpha^u$. This computation will show a good example how the reference measure can be random and singular.

\subsection{The problem setting}

Let $W$ be a brownian motion. Let $\xi=\sup\{0\leq t\leq 1: W_1-2W_t=0\}$. We consider the progressive enlargement of filtration with the random time $\xi$. According to \cite{JYC}, the Azéma supermartingale of $\xi$ is given by $$
\mathbb{P}[t<\xi|\mathcal{F}_t]
=1-h(\frac{|W_t|}{\sqrt{1-t}})
$$
where $
h(y)=\sqrt{\frac{2}{\pi}}\int_0^y s^2e^{-\frac{s^2}{2}}ds.
$
This will give the formula of enlargement of filtration on the random interval $[0,\xi]$ as it is shown in Theorem \ref{avant-tau}. So, we consider here only the enlargement of filtration formula on the time interval $(\xi,\infty)$. We consider this problem for $W$.

\subsection{Computing $\tau^u,\alpha^u$ for $0\leq u<1$}

We assume that $\mathbb{F}^\circ$ is the natural filtration of $W$. We introduce $
I= \ind_{[\xi,\infty)}\in\mathtt{D}(\mathbb{R}_+,\mathbb{R}).
$
Let $\mathcal{I}$ be the natural filtration on $\mathtt{D}(\mathbb{R}_+,\mathbb{R})$. For $x\in \mathtt{D}(\mathbb{R}_+,\mathbb{R})$, let $\rho(x)=\inf\{s\geq 0: x_s=1\}$. In the following computations, we will write explicitly variable $x\in \mathtt{D}(\mathbb{R_+},\mathbb{R})$, but omit writing $\omega$. We will not write the variables of the functions $\alpha^u, \tau^u$.

Fix $0\leq u<1$. For $u\leq t< 1$, for any bounded function $f\in\mathcal{I}_t$,$$
\dcb
&&\int\pi_{t,I/I}(I,dx)f(x)\ind_{\{\rho(x)\leq t\}}
=f(I)\ind_{\{\xi\leq t\}}
=\int f(x)\ind_{\{\rho(x)\leq t\}} \delta_{I}(dx)
\dce
$$
Also,
$$
\dcb
&&\int\pi_{t,I/F}(W,dx)f(x)\ind_{\{\rho(x)\leq u\}}
=\mathbb{E}[f(I)\ind_{\{\xi\leq u\}}|\mathcal{F}^\circ_t]
=\mathbb{E}[f(\ind_{[\xi,\infty)})\ind_{\{\xi\leq u\}}|\mathcal{F}^\circ_t]\\
\dce
$$
Let us compute the conditional law of $\xi$ given $\mathcal{F}^\circ_t$. It will be based on the following relation : for $0\leq a\leq u$, 
$$
\dcb
\{\xi< a\}
=\{\inf_{a\leq s\leq 1}W_s> \frac{W_1}{2}> 0\}
\cup \{\sup_{a\leq s\leq 1}W_s< \frac{W_1}{2}< 0\}
\dce
$$
We introduce some notations : $$
\dcb
\varkappa^{(t)}_a&=&\inf_{a\leq s\leq t}W_s\\ 
\theta^{(t)}_a&=&\sup_{a\leq s\leq t}W_s\\
g_t(c,b)&=&\ind_{\{c<b\}}\ind_{\{0<b\}}\frac{2(2b-c)}{\sqrt{2\pi t}^3}\exp\{-\frac{(2b-c)^2}{2t}\}\\
h_t(e,d,b)&=&\int_e^dg_t(c,b)dc =\int_e^{b\wedge d}\ind_{\{0<b\}}\frac{2(2b-c)}{\sqrt{2\pi t}^3}\exp\{-\frac{(2b-c)^2}{2t}\}dc\\
k_{t}(z,y)&=&\int_y^zg_{t}(y,b)db=\int_y^z\frac{2(2b-y)}{\sqrt{2\pi t}^3}\exp\{-\frac{(2b-y)^2}{2t}\}\ind_{\{0<b\}}db\\
\Gamma_{1-t}&=&\sup_{t\leq s<1}(-(W_s-W_t))\\
\Upsilon_{1-t}&=&\sup_{t\leq s<1}(W_s-W_t)
\dce
$$
Now we compute
$$
\dcb
\inf_{a\leq s\leq 1}W_s=\inf_{a\leq s\leq t}W_s\wedge \inf_{t\leq s\leq 1}W_s
=\inf_{a\leq s\leq t}W_s\wedge (W_t+\inf_{t\leq s\leq 1}(W_s-W_t))\\
=\inf_{a\leq s\leq t}W_s\wedge (W_t-\sup_{t\leq s\leq 1}(-(W_s-W_t)))
=\varkappa^{(t)}_a\wedge (W_t-\Gamma_{1-t})
\dce
$$
and similarly
$$
\sup_{a\leq s\leq 1}W_s=\sup_{a\leq s\leq t}W_s\vee \sup_{t\leq s\leq 1}W_s
=\sup_{a\leq s\leq t}W_s\vee (W_t+\sup_{t\leq s\leq 1}(W_s-W_t))
=\theta^{(t)}_a\vee (W_t+\Upsilon_{1-t})
$$
Hence,
$$
\dcb
&&\mathbb{E}[\ind_{\{\xi< a\}}|\mathcal{F}^\circ_t]\\
&=&\mathbb{P}[\varkappa^{(t)}_a\wedge (W_t-\Gamma_{1-t})> \frac{W_1-W_t+W_t}{2}> 0|\mathcal{F}^\circ_t]\\
&&+\mathbb{P}[\theta^{(t)}_a\vee (W_t+\Upsilon_{1-t})< \frac{W_1-W_t+W_t}{2}< 0|\mathcal{F}^\circ_t]\\

&=&\mathbb{P}[\varkappa^{(t)}_a\wedge (W_t-\Gamma_{1-t})> \frac{-(-(W_1-W_t))+W_t}{2}> 0|\mathcal{F}^\circ_t]\\
&&+\mathbb{P}[\theta^{(t)}_a\vee (W_t+\Upsilon_{1-t})< \frac{W_1-W_t+W_t}{2}< 0|\mathcal{F}^\circ_t]\\

&=&\int_{-\infty}^\infty dc\int_0^\infty db\ g_{1-t}(c,b)
\left[\ind_{\{\varkappa^{(t)}_a\wedge (W_t-b)> \frac{-c+W_t}{2}> 0\}}
+
\ind_{\{\theta^{(t)}_a\vee (W_t+b)< \frac{c+W_t}{2}< 0\}}\right]\\
&&\mbox{ cf. \cite[Proposition 8.1]{KS}}\\

&=&\int_{-\infty}^\infty dc\int_0^\infty db\ g_{1-t}(c,b)
\left[\ind_{\{W_t-2\varkappa^{(t)}_a\wedge (W_t-b)< c< W_t\}}
+
\ind_{\{-W_t+2\theta^{(t)}_a\vee (W_t+b)< c< -W_t\}}\right]\\
&=&\int_0^\infty db\ h_{1-t}(W_t-2\varkappa^{(t)}_a\wedge (W_t-b), W_t, b)
+\int_0^\infty db\ h_{1-t}(2\theta^{(t)}_a\vee (W_t+b)-W_t, -W_t, b)
\dce
$$
(Here the operation $\wedge$ is prior to the multiplication.) Note that the last expression is a continuous function in $a$. It is also equal to $\mathbb{E}[\ind_{\{\xi\leq a\}}|\mathcal{F}^\circ_t]$. $1-t>0$ being fixed, we can compute the differential with respect to $a\leq u$ under the integral. Note also $$
\frac{\partial}{\partial e}h_{1-t}(e,d,b)=0,\ \forall e>b\wedge d
$$
We have
$$
\dcb
&&d_a\mathbb{E}[\ind_{\{\xi\leq a\}}|\mathcal{F}^\circ_t]\\
&=&\int_0^\infty g_{1-t}(W_t-2\varkappa^{(t)}_a\wedge (W_t-b),b)
\ind_{\{W_t-2\varkappa^{(t)}_a\wedge (W_t-b)\leq b\wedge W_t\}}\ind_{\{\varkappa^{(t)}_a\leq W_t-b\}}db\ 2 d_a\varkappa^{(t)}_a\\
&&+
\int_0^\infty g_{1-t}(2\theta^{(t)}_a\vee (W_t+b)-W_t,b)
\ind_{\{2\theta^{(t)}_a\vee (W_t+b)-W_t\leq b\wedge(-W_t)\}}\ind_{\{b+W_t\leq \theta^{(t)}_a\}}db\ 2 (-d_a\theta^{(t)}_a)\\

&=&\ind_{\{\varkappa^{(t)}_a\geq 0\}}\int_{W_t-2\varkappa^{(t)}_a}^{W_t-\varkappa^{(t)}_a} g_{1-t}(W_t-2\varkappa^{(t)}_a,b)
db\ 2 d_a\varkappa^{(t)}_a\\
&&+
\ind_{\{\theta^{(t)}_a\leq 0\}}\int_{2\theta^{(t)}_a-W_t}^{\theta^{(t)}_a-W_t} g_{1-t}(2\theta^{(t)}_a-W_t,b)
db\ 2 (-d_a\theta^{(t)}_a)\\

&=&\ind_{\{\varkappa^{(t)}_a\geq 0\}}k_{1-t}(W_t-\varkappa^{(t)}_a,W_t-2\varkappa^{(t)}_a)\ 2 d_a\varkappa^{(t)}_a\\
&&+
\ind_{\{\theta^{(t)}_a\leq 0\}}k_{1-t}(\theta^{(t)}_a-W_t,2\theta^{(t)}_a-W_t)
db\ 2 (-d_a\theta^{(t)}_a)
\dce
$$
This yields
$$
\dcb
&&\mathbb{E}[f(\ind_{[\xi,\infty)})\ind_{\{\xi\leq u\}}|\mathcal{F}^\circ_t]\\

&=&\int_0^uf(\ind_{[a,\infty)})\ind_{\{\varkappa^{(t)}_a\geq 0\}}k_{1-t}(W_t-\varkappa^{(t)}_a,W_t-2\varkappa^{(t)}_a)\ 2 d_a\varkappa^{(t)}_a\\
&&+
\int_0^uf(\ind_{[a,\infty)}) \ind_{\{\theta^{(t)}_a\leq 0\}}k_{1-t}(\theta^{(t)}_a-W_t,2\theta^{(t)}_a-W_t)\ 2 (-d_a\theta^{(t)}_a)
\dce
$$
For $a< u\leq t<1$,$$
\dcb
\varkappa^{(t)}_a=\inf_{a\leq s\leq t}W_s=\varkappa^{(u)}_a\wedge \varkappa^{(t)}_u,\ \
\theta^{(t)}_a=\sup_{a\leq s\leq t}W_s=\theta^{(u)}_a\vee \theta^{(t)}_u\\
\dce
$$
So,
$$
\dcb
\ind_{\{a\leq u\}}d_a\varkappa^{(t)}_a
=\ind_{\{a\leq u\}}\ind_{\{\varkappa^{(u)}_a\leq \varkappa^{(t)}_u\}}d_a\varkappa^{(u)}_a,\ \

\ind_{\{a\leq u\}}d_a\theta^{(t)}_a
=\ind_{\{a\leq u\}}\ind_{\{\theta^{(u)}_a\geq \theta^{(t)}_u\}}d_a\theta^{(u)}_a
\dce
$$
Notice also $$
\dcb
d_a\varkappa^{(u)}_a=\ind_{\{W_a=\varkappa^{(u)}_a\}}d_a\varkappa^{(u)}_a,\ \
d_a\theta^{(u)}_a=\ind_{\{W_a=\theta^{(u)}_a\}}d_a\theta^{(u)}_a
\dce
$$
Let $v\in[u,1)$. We continue the computation for $u\leq t$ : 
$$
\dcb
&&\int\pi_{t,I/F}(W,dx)f(x)\ind_{\{\rho(x)\leq u\}}\ind_{\{t\leq v\}}=\mathbb{E}[f(\ind_{[\xi,\infty)})\ind_{\{\xi\leq u\}}\ind_{\{t\leq v\}}|\mathcal{F}^\circ_t]\\

&=&\int_0^uf(\ind_{[a,\infty)})\ind_{\{\varkappa^{(t)}_a\geq 0\}}k_{1-t}(W_t-\varkappa^{(t)}_a,W_t-2\varkappa^{(t)}_a)\ 2 d_a\varkappa^{(t)}_a\ind_{\{t\leq v\}}\\
&&+
\int_0^uf(\ind_{[a,\infty)})\ind_{\{\theta^{(t)}_a\leq 0\}} k_{1-t}(\theta^{(t)}_a-W_t,2\theta^{(t)}_a-W_t)\ 2 (-d_a\theta^{(t)}_a)\ind_{\{t\leq v\}}\\
\\
&=&\int_0^uf(\ind_{[a,\infty)})k_{1-t}(W_t-\varkappa^{(t)}_a,W_t-2\varkappa^{(t)}_a)\ind_{\{\varkappa^{(t)}_a\geq 0,\varkappa^{(u)}_a\leq \varkappa^{(t)}_u,W_a=\varkappa^{(u)}_a\}}\ 2 d_a\varkappa^{(u)}_a\ind_{\{t\leq v\}}\\
&&+
\int_0^uf(\ind_{[a,\infty)}) k_{1-t}(\theta^{(t)}_a-W_t,2\theta^{(t)}_a-W_t)\  \ind_{\{\theta^{(t)}_a\leq 0,\theta^{(u)}_a\geq \theta^{(t)}_u,W_a=\theta^{(u)}_a\}}2(-d_a\theta^{(u)}_a)\ind_{\{t\leq v\}}\\
\\
&=&\int_0^uf(\ind_{[a,\infty)})k_{1-t}(W_t-\varkappa^{(t)}_a,W_t-2\varkappa^{(t)}_a)\ind_{\{0\leq \varkappa^{(u)}_a\leq \varkappa^{(t)}_u,W_a=\varkappa^{(u)}_a\}}\ 2 d_a\varkappa^{(u)}_a\ind_{\{t\leq v\}}\\
&&+
\int_0^uf(\ind_{[a,\infty)})k_{1-t}(\theta^{(t)}_a-W_t,2\theta^{(t)}_a-W_t)\  \ind_{\{0\geq \theta^{(u)}_a\geq \theta^{(t)}_u,W_a=\theta^{(u)}_a\}}2(-d_a\theta^{(u)}_a)\ind_{\{t\leq v\}}\\
\\
&=&\int_0^uf(\ind_{[a,\infty)})\\
&&\left[k_{1-t}(W_t-\varkappa^{(t)}_a,W_t-2\varkappa^{(t)}_a)\ind_{\{0\leq \varkappa^{(u)}_a\leq \varkappa^{(t)}_u,W_a=\varkappa^{(u)}_a\}}
+k_{1-t}(\theta^{(t)}_a-W_t,2\theta^{(t)}_a-W_t)\  \ind_{\{0\geq \theta^{(u)}_a\geq \theta^{(t)}_u,W_a=\theta^{(u)}_a\}}\right]\\ 
&&(2 d_a\varkappa^{(u)}_a+2(-d_a\theta^{(u)}_a))\ind_{\{t\leq v\}}\\
\\
&=&\int_0^uf(\ind_{[a,\infty)})\\
&&\left[\frac{k_{1-t}(W_t-\varkappa^{(t)}_a,W_t-2\varkappa^{(t)}_a)}{k_{1-u}(W_u-\varkappa^{(u)}_a,W_u-2\varkappa^{(u)}_a)}\ind_{\{0\leq \varkappa^{(u)}_a\leq \varkappa^{(t)}_u,W_a=\varkappa^{(u)}_a\}}+
\frac{k_{1-t}(\theta^{(t)}_a-W_t,2\theta^{(t)}_a-W_t)}{k_{1-u}(\theta^{(u)}_a-W_u,2\theta^{(u)}_a-W_u)}  \ind_{\{0\geq \theta^{(u)}_a\geq \theta^{(t)}_u,W_a=\theta^{(u)}_a\}}\right]\\
&&\left[k_{1-u}(W_u-\varkappa^{(u)}_a,W_u-2\varkappa^{(u)}_a)\ind_{\{0\leq \varkappa^{(u)}_a,W_a=\varkappa^{(u)}_a\}}
+
k_{1-u}(\theta^{(u)}_a-W_u,2\theta^{(u)}_a-W_u) \ind_{\{0\geq \theta^{(u)}_a,W_a=\theta^{(u)}_a\}}\right]\\ 
&&(2 d_a\varkappa^{(u)}_a+2(-d_a\theta^{(u)}_a))\ind_{\{t\leq v\}}\\
\\
&=&\int\pi_{u,I/F}(W,dx)f(x)\ind_{\{\rho(x)\leq u\}}\Xi(u,t,W^t,\rho(x))\ind_{\{t\leq v\}}\\

&=&\int\pi_{u,I/F}(W,dx'')\int\pi_{t,I/I}(x'',dx)f(x)\ind_{\{\rho(x)\leq u\}}\ind_{\{t\leq v\}}\Xi(u,t,W^t,\rho(x))
\dce
$$
where $\Xi(u,t,W^t,a)$ is the function:$$
\dcb
\left[\frac{k_{1-t}(W_t-\varkappa^{(t)}_a,W_t-2\varkappa^{(t)}_a)}{k_{1-u}(W_u-\varkappa^{(u)}_a,W_u-2\varkappa^{(u)}_a)}\ind_{\{0\leq \varkappa^{(u)}_a\leq \varkappa^{(t)}_u,W_a=\varkappa^{(u)}_a\}}
+\frac{k_{1-t}(\theta^{(t)}_a-W_t,2\theta^{(t)}_a-W_t)}{k_{1-u}(\theta^{(u)}_a-W_u,2\theta^{(u)}_a-W_u)}  \ind_{\{0\geq \theta^{(u)}_a\geq \theta^{(t)}_u,W_a=\theta^{(u)}_a\}}\right]
\dce
$$
Applying Lemma \ref{alpha}, $\tau^u\ind_{\{\rho(x)\leq u\}}\geq v\ind_{\{\rho(x)\leq u\}}$ and, for $u\leq t<v$,
$$
\dcb
&&\alpha^u_t\ind_{\{\rho(x)\leq u\}}\ind_{\{0\leq \varkappa^{(u)}_{\rho(x)}\leq \varkappa^{(t)}_u,W_{\rho(x)}=\varkappa^{(u)}_{\rho(x)}\}}\\

&=&\frac{k_{1-t}(W_t-\varkappa^{(t)}_{\rho(x)},W_t-2\varkappa^{(t)}_{\rho(x)})}{k_{1-u}(W_u-\varkappa^{(u)}_{\rho(x)},W_u-2\varkappa^{(u)}_{\rho(x)})}\ind_{\{\rho(x)\leq u\}}\ind_{\{0\leq \varkappa^{(u)}_{\rho(x)}\leq \varkappa^{(t)}_u,W_{\rho(x)}=\varkappa^{(u)}_{\rho(x)}\}}
\dce
$$
and
$$
\dcb
&&\alpha^u_t\ind_{\{\rho(x)\leq u\}}\ind_{\{0\geq \theta^{(u)}_{\rho(x)}\geq \theta^{(t)}_u,W_{\rho(x)}=\theta^{(u)}_{\rho(x)}\}}\\

&=&\frac{k_{1-t}(\theta^{(t)}_{\rho(x)}-W_t,2\theta^{(t)}_{\rho(x)}-W_t)}{k_{1-u}(\theta^{(u)}_{\rho(x)}-W_u,2\theta^{(u)}_{\rho(x)}-W_u)}\ind_{\{\rho(x)\leq u\}}\ind_{\{0\geq \theta^{(u)}_{\rho(x)}\geq \theta^{(t)}_u,W_{\rho(x)}=\theta^{(u)}_{\rho(x)}\}}
\dce
$$

\subsection{The brackets computations}

Let us consider the semimartingale decomposition of the last expressions, which are composed of a part determined by the function $k_{1-t}$ and another part with indicator functions. Recall, for $z>0$,
$$
\dcb
k_{1-t}(z,y)
&=&\int_y^zg_{1-t}(y,b)db=\int_y^z\frac{2(2b-y)}{\sqrt{2\pi (1-t)}^3}\exp\{-\frac{(2b-y)^2}{2(1-t)}\}\ind_{\{0<b\}}db\\
&=&\frac{1}{\sqrt{2\pi (1-t)}^3}(1-t)\int_{y\vee 0}^z\frac{2(2b-y)}{(1-t)}\exp\{-\frac{(2b-y)^2}{2(1-t)}\}db\\

&=&\frac{1}{\sqrt{2\pi (1-t)}^3}(1-t)\int_{y\vee 0}^z\exp\{-\frac{(2b-y)^2}{2(1-t)}\}d\left(\frac{(2b-y)^2}{2(1-t)}\right)\\

&=&\frac{1}{\sqrt{2\pi (1-t)}^3}(1-t)\left.(-\exp\{-\frac{(2b-y)^2}{2(1-t)}\})\right|_{y\vee 0}^z\\
&=&\frac{1}{\sqrt{2\pi (1-t)}^3}(1-t)\left(-\exp\{-\frac{(2z-y)^2}{2(1-t)}\}+\exp\{-\frac{y^2}{2(1-t)}\}\right)
\dce
$$
Hence
$$
\dcb
&&k_{1-t}(W_t-\varkappa^{(t)}_{\rho(x)},W_t-2\varkappa^{(t)}_{\rho(x)})\\
&=&\frac{1}{\sqrt{2\pi (1-t)}^3}(1-t)\left(-\exp\{-\frac{W_t^2}{2(1-t)}\}+\exp\{-\frac{(W_t-2\varkappa^{(t)}_{\rho(x)})^2}{2(1-t)}\}\right)\\
\dce
$$
and
$$
\dcb
&&k_{1-t}(\theta^{(t)}_{\rho(x)}-W_t,2\theta^{(t)}_{\rho(x)}-W_t)\\
&=&\frac{1}{\sqrt{2\pi (1-t)}^3}(1-t)\left(-\exp\{-\frac{W_t^2}{2(1-t)}\}+\exp\{-\frac{(2\theta^{(t)}_{\rho(x)}-W_t)^2}{2(1-t)}\}\right)\\
\dce
$$
As $W=W(\iota)$ is $(\nu^u,\mathbb{J}_+)$ brownian motion, for the above processes, their $(\nu^u,\mathbb{J}_+)$ martingale part on $\{\rho\leq u\}\cap(u,1)$ are given by :
$$
\dcb
\frac{1}{\sqrt{2\pi (1-t)}^3}(1-t)\left(\frac{W_t}{1-t}\exp\{-\frac{W_t^2}{2(1-t)}\}-\frac{(W_t-2\varkappa^{(t)}_{\rho(x)})}{1-t}\exp\{-\frac{(W_t-2\varkappa^{(t)}_{\rho(x)})^2}{2(1-t)}\}\right)dW_t\\
\\
\frac{1}{\sqrt{2\pi (1-t)}^3}(1-t)\left(\frac{W_t}{1-t}\exp\{-\frac{W_t^2}{2(1-t)}\}-\frac{(W_t-2\theta^{(t)}_{\rho(x)})}{1-t}\exp\{-\frac{(W_t-2\theta^{(t)}_{\rho(x)})^2}{2(1-t)}\}\right)dW_t
\dce
$$

Consider next the indicator functions. Because $t\rightarrow \varkappa^{(t)}_u$ is decreasing, the random set $$
\{t\in[u,v]: \rho(x)\leq u, W_{\rho(x)}=\varkappa^{(u)}_{\rho(x)}\geq 0, \varkappa^{(u)}_{\rho(x)}\leq \varkappa^{(t)}_u\}
$$ 
is a random interval. Let $\lambda_\kappa$ be its supremum ($\sup\emptyset=-\infty$ by definition). Likely, the set$$
\{t\in[u,v]: \rho(x)\leq u, W_{\rho(x)}=\theta^{(u)}_{\rho(x)}\leq 0, \theta^{(u)}_{\rho(x)}\geq \theta^{(t)}_u\}
$$
is a random interval. Let $\lambda_\theta$ be its supremum. Since$$
\dcb
\{\rho(x)\leq u,W_{\rho(x)}=\varkappa^{(u)}_{\rho(x)}\geq 0\}\in\mathcal{J}_u\\
\{\rho(x)\leq u,W_{\rho(x)}=\theta^{(u)}_{\rho(x)}\leq 0\}\in\mathcal{J}_u
\dce
$$ 
$\lambda_\kappa\vee u$ and $\lambda_\theta\vee u$ are $\mathbb{J}_+$ stopping times inferior to $v\leq \tau^u$. 

Knowing the decomposition of $\alpha^u\ind_{[u,\lambda_\kappa\vee u)}$ and of $\alpha^u\ind_{[u,\lambda_\theta\vee u)}$, the continuity of the brownian motion gives a easy computation of the brackets :$$
\dcb
&&\ind_{\{u<t\leq \lambda_\kappa \}}\frac{1}{\alpha^u_{t-}}d\cro{W,\alpha^u\ind_{[u,\lambda_\kappa\vee u)}}_t
=\ind_{\{u<t\leq \lambda_\kappa\}}\frac{\left(\frac{W_t}{1-t}\exp\{-\frac{W_t^2}{2(1-t)}\}-\frac{(W_t-2\varkappa^{(t)}_{\rho(x)})}{1-t}\exp\{-\frac{(W_t-2\varkappa^{(t)}_{\rho(x)})^2}{2(1-t)}\}\right)}{\left(-\exp\{-\frac{W_t^2}{2(1-t)}\}+\exp\{-\frac{(W_t-2\varkappa^{(t)}_{\rho(x)})^2}{2(1-t)}\}\right)} dt\\
\\
&&\ind_{\{u<t\leq \lambda_\theta\}}\frac{1}{\alpha^u_{t-}}d\cro{W,\alpha^u\ind_{[u,\lambda_\theta\vee u)}}_t
=\ind_{\{u<t\leq \lambda_\theta\}}\frac{\left(\frac{W_t}{1-t}\exp\{-\frac{W_t^2}{2(1-t)}\}-\frac{(W_t-2\theta^{(t)}_{\rho(x)})}{1-t}\exp\{-\frac{(W_t-2\theta^{(t)}_{\rho(x)})^2}{2(1-t)}\}\right)}{\left(-\exp\{-\frac{W_t^2}{2(1-t)}\}+\exp\{-\frac{(2\theta^{(t)}_{\rho(x)}-W_t)^2}{2(1-t)}\}\right)} dt
\dce
$$

\subsection{The local solutions before time 1 with their drifts}

Pull the above quantities back to the initial space. We have
$$
\dcb
&&(\ind_{\{u<t\leq \lambda_\kappa\}})\circ\phi
=(\ind_{\{\rho(x)\leq u,W_{\rho(x)}=\varkappa^{(u)}_{\rho(x)}\geq 0\}}\ind_{\{u<t\leq v, \varkappa^{(u)}_{\xi}\leq \varkappa^{(t)}_u\}})\circ\phi\\
&=&\ind_{\{\xi\leq u,W_{\xi}=\varkappa^{(u)}_{\xi}\geq 0\}}\ind_{\{u<t\leq v, \varkappa^{(u)}_{\xi}\leq \varkappa^{(t)}_u\}}
=\ind_{\{\xi\leq u,W_{\xi}\geq 0\}}\ind_{\{u<t\leq v\}}
=\ind_{\{\xi\leq u,W_{1}> 0\}}\ind_{\{u<t\leq v\}}
\dce
$$
and
$$
\dcb
&&(\ind_{\{u<t\leq \lambda_\theta\}})\circ\phi
=(\ind_{\{\rho(x)\leq u,W_{\rho(x)}=\theta^{(u)}_{\rho(x)}\leq 0\}}\ind_{\{ u<t\leq v, \theta^{(u)}_{\xi}\geq \theta^{(t)}_u\}})\circ\phi\\
&=&\ind_{\{\xi\leq u,W_{\xi}=\theta^{(u)}_{\xi}\leq 0\}}\ind_{\{ u<t\leq v, \theta^{(u)}_{\xi}\geq \theta^{(t)}_u\}}
=\ind_{\{\xi\leq u,W_{\xi}\geq 0\}}\ind_{\{u<t\leq v\}}
=\ind_{\{\xi\leq u,W_{1}< 0\}}\ind_{\{u<t\leq v\}}
\dce
$$
and
$
\varkappa^{(t)}_\xi=\frac{W_1}{2},\ \theta^{(t)}_\xi=\frac{W_1}{2}.
$

$W$ is a continuous function, hence is $\mathbb{F}^\circ$-locally bounded. We can now apply Theorem \ref{theorem:our-formula} to conclude that the $\ind_{(u,\lambda_\kappa(\phi)]}\centerdot W$ and $\ind_{(u,\lambda_\theta(\phi)]}\centerdot W$ ($0\leq u\leq v<1$) are $(\mathbb{P},\mathbb{G})$ special semimartingales with the following decompositions $$
\dcb
&&\ind_{\{\xi\leq u\}}\ind_{\{W_1> 0\}}\ind_{\{u<t\leq v\}}dW_t\\
&=& d_t(\mbox{$\mathbb{G}$ local martingale})+
\ind_{\{\xi\leq u\}}\ind_{\{W_1> 0\}}\ind_{\{u<t\leq v\}}\frac{\left(\frac{W_t}{1-t}\exp\{-\frac{W_t^2}{2(1-t)}\}-\frac{(W_t-W_1)}{1-t}\exp\{-\frac{(W_t-W_1)^2}{2(1-t)}\}\right)}{\left(-\exp\{-\frac{W_t^2}{2(1-t)}\}+\exp\{-\frac{(W_t-W_1)^2}{2(1-t)}\}\right)} dt\\
\\&&\ind_{\{\xi\leq u\}}\ind_{\{W_1< 0\}}\ind_{\{u<t\leq v\}}dW_t\\
&=& d_t(\mbox{$\mathbb{G}$ local martingale})+\ind_{\{\xi\leq u\}}\ind_{\{W_1< 0\}}\ind_{\{u<t\leq v\}} \frac{\left(\frac{W_t}{1-t}\exp\{-\frac{W_t^2}{2(1-t)}\}-\frac{(W_t-W_1)}{1-t}\exp\{-\frac{(W_t-W_1)^2}{2(1-t)}\}\right)}{\left(-\exp\{-\frac{W_t^2}{2(1-t)}\}+\exp\{-\frac{(W_t-W_1)^2}{2(1-t)}\}\right)} dt
\dce
$$
From these decompositions, we obtain also the decomposition formula for $\ind_{\{\xi\leq u\}}\ind_{(u,v]}\centerdot W$

\subsection{The random measure generated by the drifts of the local solutions and its distribution function}

Let us consider $\alpha^u$ for $u\geq 1$. In this case, $I\in\mathcal{F}_u$. We obtain $\tau^u=\infty$ and $\alpha^u_t=1, t\geq u$. Clearly, $\ind_{(u,\infty)}\centerdot M$ is a $(\mathbb{P},\mathbb{G})$ martingale with a null drift.

Let $\mathtt{B}_i, i\geq 1,$ denote the left random intervals $\{\xi\leq u\}\cap(u,v], u,v\in\mathbb{Q}_+\cap[0,1)$. Let $\mathtt{B}_{0}$ denote the interval $(1,\infty)$. Then, $\cup_{i\geq 0}\mathtt{B}_i=(\xi,1)\cup(1,\infty)$. The family of $\mathtt{B}_i$ satisfies the second condition of Assumption \ref{partialcovering}. 

The drift of $W$ on the left intervals $\mathtt{B}_i$ yields a random measure concentrated on $(\xi,1)$:$$
d\chi^U_t
=\ind_{\{\xi<t<1\}}\frac{\left(\frac{W_t}{1-t}\exp\{-\frac{W_t^2}{2(1-t)}\}-\frac{(W_t-W_1)}{1-t}\exp\{-\frac{(W_t-W_1)^2}{2(1-t)}\}\right)}{\left(-\exp\{-\frac{W_t^2}{2(1-t)}\}+\exp\{-\frac{(W_t-W_1)^2}{2(1-t)}\}\right)} dt
$$ 
Since $W$ is continuous, $W-W_\xi=\ind_{(\xi,1)}\centerdot W+\ind_{(1,\infty)}\centerdot W$. According to Lemma \ref{UB}, $W-W_\xi$ is a $(\mathbb{P},\mathbb{G})$ semimartingale, whenever $\mathsf{d}\chi^\cup$ has a distribution function.

We consider then the question of distribution function. We begin with three inequalities : for $\xi\leq u\leq t<1$,
$$
\dcb
&&\left(-\exp\{-\frac{W_t^2}{2(1-t)}\}+\exp\{-\frac{(W_t-W_1)^2}{2(1-t)}\}\right)\\
&=&\exp\{-\frac{W_t^2}{2(1-t)}\}\left(\exp\{-\frac{(W_t-W_1)^2-W_t^2}{2(1-t)}\}-1\right)\\
&=&\exp\{-\frac{W_t^2}{2(1-t)}\}\left(\exp\{-\frac{-2W_tW_1+W_1^2}{2(1-t)}\}-1\right)\\
&=&\exp\{-\frac{W_t^2}{2(1-t)}\}\left(\exp\{\frac{2(W_t-W_\xi)W_1}{2(1-t)}\}-1\right)\\
&\geq&\exp\{-\frac{W_t^2}{2(1-t)}\}\frac{2(W_t-W_\xi)W_1}{2(1-t)}>0\\
\dce
$$
and
$$
\dcb
&&k_{1-t}(W_t-\varkappa^{(t)}_{a},W_t-2\varkappa^{(t)}_{a})\\
&=&\frac{1}{\sqrt{2\pi (1-t)}^3}(1-t)\left(-\exp\{-\frac{W_t^2}{2(1-t)}\}+\exp\{-\frac{(W_t-2\varkappa^{(t)}_{a})^2}{2(1-t)}\}\right)\\
&\leq &\frac{1}{\sqrt{2\pi (1-t)}^3}(1-t)\left(-\frac{(W_t-2\varkappa^{(t)}_{a})^2}{2(1-t)}+\frac{W_t^2}{2(1-t)}\right)\\
&= &\frac{1}{\sqrt{2\pi (1-t)}^3}(1-t)\frac{4W_t\varkappa^{(t)}_{a}-4(\varkappa^{(t)}_{a})^2}{2(1-t)}\\
&= &\frac{2}{\sqrt{2\pi (1-t)}^3}(W_t-\varkappa^{(t)}_{a})\varkappa^{(t)}_{a}\\

\dce
$$
and similarly
$$
\dcb
k_{1-t}(\theta^{(t)}_{a}-W_t,2\theta^{(t)}_{a}-W_t)
&\leq&\frac{2}{\sqrt{2\pi (1-t)}^3}(W_t-\theta^{(t)}_{a})\theta^{(t)}_{a}
\dce
$$
We now compute, for any $0<\epsilon<1$
$$
\dcb
&&\mathbb{E}[\int_\xi^{1-\epsilon}\frac{1}{|W_\xi-W_t|}dt]\\
&=&\mathbb{E}[\int_0^{1-\epsilon}\frac{1}{|W_\xi-W_t|}\ind_{\{\xi<t\}}dt]\\
&=&\mathbb{E}[\int_0^{1-\epsilon}\mathbb{E}[\frac{1}{|W_\xi-W_t|}\ind_{\{\xi<t\}}|\mathcal{F}^\circ_t] dt]\\

&=&\mathbb{E}[\int_0^{1-\epsilon}\left(\int_0^t\frac{1}{|W_a-W_t|}\ind_{\{\varkappa^{(t)}_a\geq 0\}}k_{1-t}(W_t-\varkappa^{(t)}_a,W_t-2\varkappa^{(t)}_a)\ 2 d_a\varkappa^{(t)}_a\right.\\
&&+
\left.\int_0^t\frac{1}{|W_a-W_t|}\ind_{\{\theta^{(t)}_a\leq 0\}} k_{1-t}(\theta^{(t)}_a-W_t,2\theta^{(t)}_a-W_t)\ 2 (-d_a\theta^{(t)}_a)\right) dt]\\

&=&\mathbb{E}[\int_0^{1-\epsilon}\left(\int_0^t\frac{1}{|\varkappa^{(t)}_a-W_t|}\ind_{\{\varkappa^{(t)}_a\geq 0\}}k_{1-t}(W_t-\varkappa^{(t)}_a,W_t-2\varkappa^{(t)}_a)\ 2 d_a\varkappa^{(t)}_a\right.\\
&&+
\left.\int_0^t\frac{1}{|\theta^{(t)}_a-W_t|}\ind_{\{\theta^{(t)}_a\leq 0\}} k_{1-t}(\theta^{(t)}_a-W_t,2\theta^{(t)}_a-W_t)\ 2 (-d_a\theta^{(t)}_a)\right) dt]\\

&\leq&\frac{2}{\sqrt{2\pi (1-t)}^3}\mathbb{E}[\int_0^{1-\epsilon}\left(\int_0^t\ind_{\{\varkappa^{(t)}_a\geq 0\}}|\varkappa^{(t)}_{a}|\ 2 d_a\varkappa^{(t)}_a
+
\int_0^t\ind_{\{\theta^{(t)}_a\leq 0\}} |\theta^{(t)}_{a}|\ 2 (-d_a\theta^{(t)}_a)\right) dt]\\

&\leq&\frac{2}{\sqrt{2\pi (1-t)}^3}\mathbb{E}[\int_0^{1-\epsilon}\left(\  (\varkappa^{(t)}_t)^2\ind_{\{\varkappa^{(t)}_t>0\}}
+
(\theta^{(t)}_t)^2\ind_{\{\theta^{(t)}_t<0\}}\right) dt]\\

&=&\frac{2}{\sqrt{2\pi (1-t)}^3} \int_0^{1-\epsilon} \mathbb{E}[W_t^2]\ dt\\
&<&\infty
\dce
$$
This implies that $$
\int_\xi^{1}\frac{1}{|W_\xi-W_t|}dt\ <\infty
$$
Applying the first of the above three inequalities, we see that the random measure $\mathsf{d}\chi^\cup$ has effectively a distribution function $\chi^\cup$. Consequently, $W-W_\xi$ is a $(\mathbb{P},\mathbb{G})$ special semimartingale with drift $\chi^\cup$. 

\subsection{Ending remark}
To end this example, we write the drift in another form : for $\xi<t\leq 1$,
$$
\dcb
&&\frac{\frac{W_t}{1-t}\exp\{-\frac{W_t^2}{2(1-t)}\}-\frac{(W_t-W_1)}{1-t}\exp\{-\frac{(W_t-W_1)^2}{2(1-t)}\}}{-\exp\{-\frac{W_t^2}{2(1-t)}\}+\exp\{-\frac{(W_t-W_1)^2}{2(1-t)}\}}\\

&=& \frac{\frac{W_t}{1-t}\exp\{-\frac{W_t^2}{2(1-t)}\}}{-\exp\{-\frac{W_t^2}{2(1-t)}\}+\exp\{-\frac{(W_t-W_1)^2}{2(1-t)}\}}
+\frac{-\frac{(W_t-W_1)}{1-t}\exp\{-\frac{(W_t-W_1)^2}{2(1-t)}\}}{-\exp\{-\frac{W_t^2}{2(1-t)}\}+\exp\{-\frac{(W_t-W_1)^2}{2(1-t)}\}}\\
\\
&=&\frac{\frac{W_t}{1-t}\exp\{-\frac{W_t^2}{2(1-t)}\}-\frac{(W_t-W_1)}{1-t}\exp\{-\frac{W_t^2}{2(1-t)}\}}{-\exp\{-\frac{W_t^2}{2(1-t)}\}+\exp\{-\frac{(W_t-W_1)^2}{2(1-t)}\}}
+\frac{\frac{W_t-W_1}{1-t}\exp\{-\frac{W_t^2}{2(1-t)}\}-\frac{(W_t-W_1)}{1-t}\exp\{-\frac{(W_t-W_1)^2}{2(1-t)}\}}{-\exp\{-\frac{W_t^2}{2(1-t)}\}+\exp\{-\frac{(W_t-W_1)^2}{2(1-t)}\}}\\
\\
&=&\frac{\frac{W_1}{1-t}\exp\{-\frac{W_t^2}{2(1-t)}\}}{-\exp\{-\frac{W_t^2}{2(1-t)}\}+\exp\{-\frac{(W_t-W_1)^2}{2(1-t)}\}}-\frac{W_t-W_1}{1-t}\\
\\
&=& 
\frac{\frac{W_1}{1-t}}{-1+\exp\{-\frac{(W_t-W_1)^2}{2(1-t)}+\frac{W_t^2}{2(1-t)}\}}-\frac{W_t-W_1}{1-t}\\
\\
&=&\frac{1}{-1+\exp\{\frac{2W_tW_1-W_1^2}{2(1-t)}\}}\frac{W_1}{1-t}-\frac{W_t-W_1}{1-t}
\dce
$$
Be ware of the component $-\frac{W_t-W_1}{1-t}\ind_{\{\xi<t\leq 1\}}dt$ in the above formula. This term is intrinsically related with the brownian motion bridge as it is indicated in \cite{AM}. In fact, using brownian motion bridge, we can have a different and shorter proof of the above enlargement of filtration formula.

\

\section{Expansion with the future infimum process}\label{futur-infermum}

This is an important result in the literature of the theory of enlargement of filtration. It is also a typical example which does not fit the classification of initial enlargement or progressive enlargement. The result is motivated by Pitman's theorem \cite{Pitman}. See \cite{Jeulin79, Jeulin80, MY, Rauscher, saisho, takaoka, Yor97}. The use of the enlargement of filtration to study Pitman's theorem was initiated in \cite{Jeulin79}, based ingeniously on a relationship between the enlargement by future infimum process and the enlargement by honest times (which can be dealed with within the framework of progressive enlargement). Such study continued in \cite{Jeulin80, Yor97,Rauscher} with various different approaches. We notice that all these approaches sit over properties specific to the problem. In this section we will show how the general local solution method works on this example.

\subsection{The setting}

Consider a regular conservative diffusion process $Z$ taking values in $(0,\infty)$ defined in a stochastic basis $(\Omega,\mathcal{A},\mathbb{F},\mathbb{P})$. See \cite{Revuz-Yor} for a presentation of general $\mathbb{R}$ valued diffusion. We recall only that the scale function of $Z$ is defined (modulo an affine transformation) as a continuous strictly increasing function $e$ on $(0,\infty)$ such that $$
\mathbb{P}_{x}[T_{b} < T_{a}] = \frac{e(x)-e(a)}{e(b)-e(a)},
$$ 
where $0<a<x<b<\infty$ and $T_{a}, T_{b}$ denote the corresponding hitting times.
We assume in this section that
\begin{enumerate}
\item
$\lim_{t\uparrow \infty} Z_{t} = \infty,$
\item
a scale function $e$ can be chosen such that $e(\infty)=0,e(0+)=-\infty$ (like Bessel(3)).
\end{enumerate}
Define the future infimum process $I = (I_{t})_{t>0}$ where $I_{t}=\inf_{s\geq t} Z_{s},\forall t>0$. Recall that (cf.
\cite[Proposition(6.29)]{Jeulin80}), for any finite $\mathbb{F}$ stopping time $T$,
\[ P[I_{T}\leq a|{\mathcal{F}}_{T}] = e(Z_{T}) \int_{0}^{Z_{T}}\  1_{(0,a]}(c) \ \upsilon(dc),
\]
where $\upsilon(c)=\frac{1}{e(c)}$. We consider the process $I$ as a map from $\Omega$ into the space $E=C_0(\mathbb{R}_+,\mathbb{R})$ of all continuous functions on $\mathbb{R}_+$. Let $\mathbb{I}$ be the natural filtration on $E$ generated by the coordinates functions. Let $\mathbb{F}^\circ$ be the natural filtration of $Z$. Define the filtration $\mathbb{G}^\circ$ as in Subsection \ref{h-formula}, as well as the filtrations $\mathbb{F}$ and $\mathbb{G}$.   

\bl\label{lemma:Z=I}
For $0<s<t$, let  $U(s,t)=\inf_{s\leq v\leq t}Z_{v}$. Then, $\mathbb{P}[I_{t} = U(s,t)]=0$. Consequently, \[
\mathbb{P}[\mbox{\textup{there are at least two distinct points $v,v'\geq s$ such that $Z_{v}=I_{s},Z_{v'}=I_{s}$}}]=0.\]
\el 

\textbf{Proof.} To prove the first assertion, we write
\[\begin{array}{lcl}
\mathbb{P}[I_{t} = U(s,t)]&=& P[P[I_{t} = a|{\mathcal{F}}_{t}]_{a=U(s,t)}]
=\mathbb{E}[e(Z_{t}) \int_{0}^{Z_{t}} \upsilon(dc) 1_{\{c=U(s,t)\}}]=0.
\end{array}\]
The second assertion is true because it is overestimated by $\mathbb{P}[\cup_{t\in \mathbb{Q}_{+},t>s}\{U(s,t)=I_{t}\}].$ \ok

As usual, for $x\in E$ and $s\geq 0$, the coordinate of $x$ at $s$ is denoted by $x_s$. We define $D_{u}(x) = \inf\{s>u: x_s 
\neq x_u\}$. For $0\leq u<t$, $I_u=U(u,t)\wedge I_t$. Therefore, $$
\{t\leq D_{u}(I)\}=\{I_u=I_t\}=\{U(u,t)\geq I_{t}\}
$$
and on the set $\{t\leq D_{u}(I)\}$, the stopped process $I^{t}$ coincides with the process $U^{u}\wedge I_{t}=U^{u}\wedge I_{u}=I^u$, where $U^{u}$ denotes the process $(U(s\wedge u,u))_{s\geq 0}$.

\subsection{Computation of $\alpha^u$}

To apply the results in section \ref{lsmfind}, assume that $\Omega$ is a Polish space, $\mathcal{A}$ is its Borel $\sigma$-algebra and $Z_{t}, t\geq 0,$ are Borel functions. Let $0\leq u<t$. Let $f$ be a bounded $\mathcal{I}_t$ measurable function ($f$ depending only on $I^t$). Consider the random measure $\pi_{u,I/I}$.
\[\begin{array}{lcl}
&&\int \pi_{u,I/I}(I,dx) f(x) \ind_{\{t\leq D_{u}(x)\}} 
=\mathbb{E}[f(I)\ind_{\{t\leq D_{u}(I)\}}|\sigma(I_s: 0\leq s\leq u)]\\
&=&\mathbb{E}[f(I^t)\ind_{\{t\leq D_{u}(I)\}}|\sigma(I_s: 0\leq s\leq u)]\ \mbox{ because $f\in\mathcal{I}_t$,}\\
&=&\mathbb{E}[f(I^u)\ind_{\{t\leq D_{u}(I)\}}|\sigma(I_s: 0\leq s\leq u)]\ \mbox{ because of $t\leq D_{u}(I)$,}\\
&=&f(I^u)\mathbb{E}[\ind_{\{t\leq D_{u}(I)\}}|\sigma(I_s: 0\leq s\leq u)]\\
&=& f(I^{u}) q(u,t,I^{u}),
\end{array}\]
where $q(u,t,x)$ is a $\mathcal{I}_{u}$-measurable positive function such that $$
q(u,t,I)
=
\mathbb{P}[t\leq D_{u}(I)|\sigma(I_s: 0\leq s\leq u)].
$$
The function $q(u,t,x)$ can be chosen l\`adc\`ag decreasing in $t$. 
Set $\eta(u,t,x) = \ind_{\{q(u,t,x)>0\}}, x\in E$ ($q(u,t)$ and $\eta$ depending only on $I^u$). In the following computations, the parameters $u,t$ being fixed, we omit them from writing. We have 
$$
\int\pi_{u,I/I}(I,dx)(1-\eta(x)) \ind_{\{t\leq D_{u}(x)\}}
=\mathbb{E}[\ind_{\{t\leq D_{u}(I)\}}|\sigma(I_s: 0\leq s\leq u)](1-\eta(I)) 
= 0.
$$
Put $\pi_{u,I/I}(x'',dx)$ together with $\pi_{u,I/F}(\mathsf{i},dx'')$. 
\[\begin{array}{lcl}
&&\int \pi_{u,I/F}(\mathsf{i},dx'') \int \pi_{u,I/I}(x'',dx) f(x)\eta(x) \ind_{\{t\leq D_{u}(x)\}}
\\
&=&\int \pi_{u,I/F}(\mathsf{i},dx'')  f(x''^{u})\eta(x''^u) q(x''^{u})\\
&=&\mathbb{E}[f(I^{u})\eta(I^u) q(I^{u})|\mathcal{F}^\circ_u]\\
&=&\mathbb{E}[f(U^{u}\wedge I_u)\eta(U^{u}\wedge I_u) q(U^{u}\wedge I_u)|\mathcal{F}^\circ_u]\\

&=&e(Z_{u}) \int_{0}^{Z_{u}}\ d(e(c)^{-1})\ f(U^{u}\wedge c)\eta(U^{u}\wedge c) q(U^{u}\wedge c).
\end{array}\]
Consider next the random measure $\pi_{t,I/F}$.
$$
\dcb
&&\int \pi_{t,I/F}(\mathsf{i},dx) f(x)\eta(x) \ind_{\{t\leq D_{u}(x)\}}\\
&=&\mathbb{E}[f(I)\eta(I) \ind_{\{t\leq D_{u}(I)\}}|\mathcal{F}^\circ_t]\\
&=&\mathbb{E}[f(U^u\wedge I_t)\eta(U^u\wedge I_t) \ind_{\{U(u,t)\geq I_t\}}|\mathcal{F}^\circ_t]\\
&=&\mathbb{E}[e(Z_{t})  \ \int_{0}^{Z_{t}}  \upsilon(dc) \ f(U^{u}\wedge c)\eta(U^{u}\wedge c)\ind_{\{U(u,t)\geq c\}}|\mathcal{F}^\circ_t]\\
&=&e(Z_{t})  \ \int_{0}^{Z_{t}}  \upsilon(dc) \ f(U^{u}\wedge c)\eta(U^{u}\wedge c)\ind_{\{U(u,t)\geq c\}}\\
&&\mbox{ noting that $c\leq U(u,t)\leq Z_u\wedge Z_t$}\\
&&\mbox{ and the coordinate of $(U^u\wedge c)$ at $u$ is $Z_u\wedge c$}\\
&=&\frac{e(Z_{t})}{e(Z_{u})}\ e(Z_{u}) \int_{0}^{Z_{u}}  \upsilon(dc) \ f(U^{u}\wedge c)\eta(U^{u}\wedge c) \ind_{\{U(u,t)\geq (U^u\wedge c)_u\}}q(U^{u}\wedge c) \frac{1}{q(U^{u}\wedge c)}\\
&=&\frac{e(Z_{t})}{e(Z_{u})}\ \int \pi_{u,I/F}(\mathsf{i},dx'') \int \pi_{u,I/I}(x'',dx)f(x)
\eta(x)\ind_{\{t\leq D_{u}(x)\}}\frac{1}{q(x)}
\ind_{\{U(u,t)\geq x_u\}}.
\dce
$$ 
Lemma \ref{alpha} is applicable. For $x\in E$ and $\omega\in\Omega$, let \[\begin{array}{lcl}
W_{t}(x) &=& \frac{1}{q(u,t+,x)},\\
\psi(\omega,x) &=& \inf\{s\geq u: Z_{s}(\omega) < x_u\},\\
\lambda(x) &=& \inf\{s\geq u:   q(u,s,x) = 0\},\\
\rho(\omega,x) &=&
\psi(\omega,x)\wedge \lambda(x) \wedge D_{u}(x).
\end{array}\] 
We have
$$
\tau^{u}(\omega,x)\geq \lambda(x)\wedge D_{u}(x),\
\alpha_{t}^{u}\ind_{\{\rho>t\}} = \frac{e(Z_{t})}{e(Z_{u})} W_{t}\ind_{\{\rho>t\}}, \ u\leq t.
$$

\subsection{Local solutions}

The integration par parts formula gives the following identity under the probability $\nu^{u}$ on the interval $(u,\infty)$:
$$\begin{array}{lcl}
	\alpha^{u}_{t}\ind_{\{t<\rho\}}& =& \frac{1}{e(Z_{u})}e(Z_{t\wedge \rho})W_{t}^{\rho-}
	-\frac{1}{e(Z_{u})}e(Z_{\rho}) W_{\rho-}\ind_{\{\rho\leq t\}} \\
	&=&\frac{1}{e(Z_{u})}e(Z_{u})W_{u}^{\rho-} + \frac{1}{e(Z_{u})}\int_{u}^{t\wedge \rho}W_{s-}de(Z)_s + \frac{1}{e(Z_{u})}\int_{u}^{t\wedge \rho}e(Z)_s dW^{\rho-}_s 
	-\frac{1}{e(Z_{u})}e(Z_{\rho}) W_{\rho-}\ind_{\{\rho\leq t\}}.
\end{array}
$$
Here $e(Z)$ is viewed as a $(\nu^{u},{\mathbb{J}}_{+})$ local martingale (cf. Lemma \ref{martingaletransfer}).

Let $\langle e(Z) \rangle$ be the quadratic variation of $e(Z)$ under $\mathbb{P}$ in the filtration $\mathbb{F}$. We notice
that  it is $\mathbb{F}$-predictable. It is straightforward to check that $\langle e(Z) \rangle$ is also a version of the
predictable quadratic variation of $e(Z)$ with respect to
$(\nu^{u},{\mathbb{J}}_{+})$ (cf. Lemma \ref{martingaletransfer}). It results that
$$
\cro{e(Z),\alpha^{u}\ind_{[u,\rho)}}^{\nu^{u}\cdot{\mathbb{J}}_{+}}
=
e(Z_{u})^{-1}W_{-}\ind_{(u,\rho]}\centerdot \langle e(Z)
\rangle
$$ on the interval $[u,\infty)$. Pull these quantities back the the original space, we obtain$$
\dcb
\psi(\phi)&=&\inf\{s>u: Z_{s}<I_{u}\} = \infty\\
\lambda(I)&=&\inf\{s>u:   q(u,s,I) = 0\}\geq D_{u}(I)\ \mbox{ $\mathbb{P}$-a.s.}\\
\rho(\phi)&=&D_u(I)
\dce
$$
where the second inequality comes from the fact $
\mathbb{P}[q(u,t,I)=0,t\leq D_u(I)]=0.
$
Applying Theorem \ref{theorem:our-formula}, we conclude that $\ind_{(u,D_u(I)]}\centerdot e(Z)$ is a $(\mathbb{P},\mathbb{G})$ special semimartingale with drift given by:
\[
\mathsf{d}\chi^{(u,D_u(I)]} = \ind_{(u,D_u(I)]} \frac{1}{e(Z)} d\langle e(Z) \rangle.
\]

\subsection{Aggregation of local solutions}

Consider $\mathtt{B}_u=(u,D_u(I)]$ for $u\in\mathbb{Q}_+$. Note that $I$ increases only when $Z=I$. Therefore, according to Lemma \ref{lemma:Z=I}, $D_{u}(I) = \inf\{s>u; Z_{s}=I_{s}\}$. Hence, the interior of $\cup_{u\in\mathbb{Q}_+}(u, D_u(I)]$ is $\{Z>I\}$. The family of $\mathtt{B}_u$ satisfies Assumption \ref{partialcovering}. These random measures $\mathsf{d}\chi^{(u,D_u(I)]}$ yields a diffuse random measure $$
\mathsf{d}\chi^\cup=
\ind_{\{Z>I\}} \frac{1}{e(Z)} \mathsf{d}\langle e(Z) \rangle
$$ 
on $\cup_{u\in\mathbb{Q}_+}(u, D_u(I)]$ which has clearly a distribution function. 

\bl\label{g-decomposition}
Let $\mathtt{G}=\{Z>I\}$. $e(Z)$ is a $\mathcal{G}$-semimartingale. We have the semimartingale decomposition : \[
e(Z)=e(Z_{0})-e(I_{0})+ M + \frac{1}{e(Z)}\cdot \langle e(Z) \rangle+V^++e(I),
\]
where $M$ is a continuous 
local $\mathbb{G}$-martingale whose quadratic variation is given by $\ind_{\mathtt{G}}\cdot \langle e(Z) \rangle$, and $V^+$ is the increasing process defined in Lemma \ref{Z-Ytheorem}. 
\el

\textbf{Proof.} The difference between $\mathtt{G}$ and $\cup_{u\in\mathbb{Q}_+}(u, D_u(I)]$ is contained in $\cup_{u\in\mathbb{Q}_+}[D_u(I)]$ which is a countable set. Let $\mathtt{A}$ be the set used in Theorem \ref{Z-Ytheorem}. The difference between $\mathtt{A}$ and $\cup_{u\in\mathbb{Q}_+}(u, D_u(I)]$ is equally countable. 

Now we can check that the process $e(Z)$ satisfies Assumption \ref{pieces} with respect to the family $(\mathtt{B}_u, u\in\mathbb{Q}_+)$ and with respect to the $(\mathbb{P},\mathbb{G})$ special semimartingale $e(I)$. Theorem \ref{Z-Ytheorem} is applicable to conclude the lemma. We make only a remark that $d\cro{e(Z)}$ does not charge $\mathtt{A}\setminus \mathtt{G}$ which is a countable set while $\cro{e(Z)}$ is continuous.
\ok

\subsection{Computation of $V^+$}

Lemma \ref{g-decomposition} solves the enlargement of filtration problem. Now we want to compute the process $V^+$. 

\bl\label{lemma:non-charge}
$\int 1_{\mathtt{G}^{c}} d\langle e(Z) \rangle = 0$. Consequently, if $X=e(Z)-e(I)$ and $l^{0}(X)$ denotes the local time of $X$ at 0, we have$$
V^+=\frac{1}{2}l^{0}(X)
$$
\el

\textbf{Proof.} By Lemme \ref{Z-Ytheorem}, $X=X_0+\ind_\mathtt{G}\centerdot X+V^+$. We can calculate the quadratic variation $[X]$ in $\mathcal{G}$ in two ways : 
\[
\dcb
d[X] = d[\ind_\mathtt{G}\centerdot X + V^+] = d[\ind_\mathtt{G}\centerdot X] = \ind_\mathtt{G}d\langle e(Z) \rangle\\
d[X] = d[e(Z)-e(I)] = d\langle e(Z) \rangle.
\dce
\]
This yields $\int 1_{\mathtt{G}^{c}} d\langle e(Z) \rangle = 0$. Using the formula in Lemma \ref{g-decomposition}, applying \cite[Chapter VI, Theorem 1.7]{Revuz-Yor}, $$
\frac{1}{2}l^0_t
=\ind_{\{X=0\}}\centerdot\left( \frac{1}{e(Z)}\cdot \langle e(Z) \rangle+V^+\right)
=\ind_{\{X=0\}}\centerdot V^+ = V^+
$$
\ok

\bl\label{lemma:tribu}
The processes of the form $Hf(I_{t})1_{(t,u]}$, where $0<t<u<\infty$, $H \in {\mathcal{F}}_{t}$, $f$ is a bounded continuous function on $]0,\infty[$, generate the $\mathbb{G}$-predictable $\sigma$-algebra.
\el

\textbf{Proof.} It is because $I_s=U(s,t)\wedge I_t$ for any $0\leq s\leq t$ and $U(s,t)\in\mathcal{F}_t$. \ok

For any $a>0,b>0$, set $\delta_b(a) = \inf\{s>b: Z_{s}\leq a\}$. Then, on the set $\{a< Z_{b},\delta_b(a)<\infty\}$, $I_{s}=I_{b}\leq a$ for $b\leq s\leq \delta_b(a)$,
and consequently, $e(I_{u\wedge\delta_b(a)})- e(I_{t\wedge\delta_b(a)}) = 0$ for $b\leq t\leq u$. On the other hand, on the set $\{a< Z_{b}, \delta_b(a) = \infty\}=\{a< I_{b}\}$. In sum, 
\[
\ind_{\{a< Z_{b}\}}(e(I_{u\wedge\delta_b(a)})- e(I_{t\wedge\delta_b(a)}))
=\ind_{\{a< I_{b}\}}(e(I_{u})- e(I_{t})), \ b\leq t\leq u.
\]
Note that $l^0$ is a continuous process increasing only when $Z=I$. The above argument holds also for $l^0$:
\[
\ind_{\{a< Z_{b}\}}(l^0_{u\wedge\delta_b(a)}- l^0_{t\wedge\delta_b(a)})
=\ind_{\{a< I_{b}\}}(l^0_{u}- l^0_{t}), \ b\leq t\leq u.
\]

\label{lemma:projection} 

\bl \label{corollary:projection}
For any $a>0$, for any finite $\mathbb{F}$ stopping time $\beta$, the $\mathbb{F}$-predictable dual projection of the increasing process $\ind_{\{a< I_{\beta}\}}\ind_{(\beta,\infty)}\centerdot e(I)$
is given by
$$
-\ind_{(\beta,\infty)}\ind_{\{U(\beta,\cdot)> a\}}\frac{1}{2e(Z)} \centerdot \langle e(Z) \rangle
$$
\el

\textbf{Proof.} Let $a>0$ to be a real number. Let $V$ be a bounded $\mathbb{F}$ stopping time such that, stopped at $V$, the two processes $(\ln(-e(a))-\ln(-e(Z_{s}))-1) \centerdot e(Z_{s})$ and $\frac{1}{2e(Z_{s})} \centerdot \langle e(Z) \rangle_{s}$ are uniformly integrable. Let $R,T$ be two others $\mathbb{F}$ stopping times such that $R\leq T\leq V$. Notice that, on the set $\{a\leq I_{R}\}$, $e(a)\leq e(I_T)\leq e(I_V)< 0$. Under these condition, we have the following computations:
\[\begin{array}{lcl}
&&\mathbb{E}[e(I_{V})\ind_{\{a< I_{R}\}}|{\mathcal{F}}_{T}]\\
&=& \mathbb{E}[\ind_{\{a< U(R,V)\}} \mathbb{E}[e(I_{V})\ind_{\{a< I_{V}\}}|{\mathcal{F}}_{V}] |{\mathcal{F}}_{T}]\\
&=& \ind_{\{a< Z_{R}\}}\mathbb{E}[\ind_{\{V< \delta_R(a)\}}e(Z_{V}) \int_{a}^{ Z_{V}}  e(c) \ d(e(c)^{-1})
|{\mathcal{F}}_{T}]\\ 
&=& \ind_{\{a< Z_{R}\}}\mathbb{E}[\ind_{\{V< \delta_R(a)\}} e(Z_{V})(\ln(-e(a))-\ln(-e(Z_{V})))
|{\mathcal{F}}_{T}]\\ 
&=& \ind_{\{a< Z_{R}\}}\mathbb{E}[e(Z_{V\wedge \delta_R(a)})(\ln(-e(a))-\ln(-e(Z_{V\wedge \delta_R(a)})))
|{\mathcal{F}}_{T}].
\end{array}\]
We have the same computation when $V$ is replaced by $T$. By Ito's formula (in the filtration $\mathbb{F}$), 
\[\begin{array}{lcl}
&&d\left( e(Z_{s})(\ln(-e(a))-\ln(-e(Z_{s})))\right)\\
&=&(\ln(-e(a))-\ln(-e(Z_{s}))-1) de(Z)_{s} - \frac{1}{2e(Z_{s})} d\langle e(Z) \rangle_{s}, \ s\geq 0.
\end{array}\]
Taking the difference of the above computations, we obtain
$$
\dcb
&&\mathbb{E}[\ind_{\{a< I_{R}\}}\int_T^Vde(I)_s|{\mathcal{F}}_{T}]
=\mathbb{E}[(e(I_{V})-e(I_T))\ind_{\{a< I_{R}\}}|{\mathcal{F}}_{T}]\\
&=& \ind_{\{a< Z_{R}\}}\mathbb{E}[\int_{T\wedge \delta_R(a)}^{V\wedge \delta_R(a)}(\ln(-e(a))-\ln(-e(Z_{s}))-1) de(Z)_{s} - \int_{T\wedge \delta_R(a)}^{V\wedge \delta_R(a)}\frac{1}{2e(Z_{s})} d\langle e(Z) \rangle_{s}
|{\mathcal{F}}_{T}]\\
&=& -\ind_{\{a< Z_{R}\}}\mathbb{E}[\int_{T}^{V}\ind_{\{s<  \delta_R(a)\}}\frac{1}{2e(Z_{s})} d\langle e(Z) \rangle_{s}
|{\mathcal{F}}_{T}]\\
&=& -\mathbb{E}[\int_{T}^{V}\ind_{\{U(R,s)> a\}}\frac{1}{2e(Z_{s})} d\langle e(Z) \rangle_{s}|{\mathcal{F}}_{T}]
\dce
$$
This identity yields that, for any positive bounded elementary $\mathbb{F}$ predictable process $H$,
$$
\dcb
\mathbb{E}[\ind_{\{a< I_{R}\}}\int_R^VH_sde(I)_s]

&=& -\mathbb{E}[\int_{R}^{V}H_s\ind_{\{U(R,s)> a\}}\frac{1}{2e(Z_{s})} d\langle e(Z) \rangle_{s}]
\dce
$$
It is to notice that we can find an increasing sequence $(V_n)_{n\geq 1}$ of $\mathbb{F}$ stopping times tending to the infinity, and each of the $V_n$ satisfies that assumption on $V$. Replace $V$ by $V_n$ and let $n$ tend to the infinity, monotone convergence theorem then implies
$$
\dcb
\mathbb{E}[\ind_{\{a< I_{R}\}}\int_R^\infty H_sde(I)_s]

&=& -\mathbb{E}[\int_{R}^{\infty}H_s\ind_{\{U(R,s)> a\}}\frac{1}{2e(Z_{s})} d\langle e(Z) \rangle_{s}]
\dce
$$
Notice that the left hand side is finite. Therefore, the right hand side is concerned by an integrable random variable. Now, replace $R$ by $\beta\wedge V_n$ and let $n$ tend to the infinity. The dominated convergence theorem implies 
$$
\dcb
\mathbb{E}[\ind_{\{a< I_{\beta}\}}\int_\beta^\infty H_sde(I)_s]
&=& -\mathbb{E}[\int_{\beta}^{\infty}H_s\ind_{\{U(\beta,s)> a\}}\frac{1}{2e(Z_{s})} d\langle e(Z) \rangle_{s}]
\dce
$$
This identity proves the lemma. \ok

\bethe
$l^{0}(X) = 2e(I)-2e(I_0)$.
\ethe

\textbf{Proof.} Using the occupation formula (cf. \cite{Revuz-Yor}), for
$0<t<u$, 
\[ l^{0}_{u}(X) - l^{0}_{t}(X)=\lim_{\epsilon \downarrow 0} \frac{1}{\epsilon} \int_{0}^{\epsilon} (l^a_u-l^a_t) da 
=\lim_{\epsilon \downarrow 0} \frac{1}{\epsilon} \int_{t}^{u} 1_{\{X_{v}\leq
	\epsilon\}}d\langle X
\rangle_{v}
=
\lim_{\epsilon \downarrow 0} \frac{1}{\epsilon} \int_{t}^{u} 1_{\{e(Z_{v})-e(I_{v})\leq \epsilon\}}d\langle e(Z)
\rangle_{v}  ,
\]
almost surely. Let $a>0$. Let $0\leq R\leq T\leq V$ be three $\mathbb{F}$ stopping times such that, stopped at $V$, the processes $Z$ and $\cro{e(Z)}$ are bounded. Note that this condition implies that the family of random variables $\{l^a_{s\wedge V}: s\geq 0, 0\leq a\leq 1\}$ is uniformly integrable. By convex combination, the family $\{\frac{1}{\epsilon}\int_0^\epsilon (l^a_{V}-l^a_{T})da : 0< \epsilon\leq 1\}$ also is uniformly integrable. Consequently, we can write the identity:$$
E[Hf(e(I_{T})) \ind_{\{a< Z_{R}\}}(l^0_{V\wedge \delta_R(a)}-l^0_{T\wedge \delta_R(a)})]
=
\lim_{\epsilon\downarrow 0}E[Hf(e(I_{T})) \ind_{\{a< Z_{R}\}}\frac{1}{\epsilon} \int_{T\wedge \delta_R(a)}^{V\wedge \delta_R(a)} 1_{\{e(Z_{v})-e(I_{v})\leq
	\epsilon\}}d\langle e(Z)
\rangle_{v} ]
$$
for any positive continuously differentiable function $f$ with compact support and any positive bounded random variable $H$ in ${\mathcal{F}}_{t}$. We now compute the limit at right hand side. 
\[\begin{array}{ll}
&E[Hf(e(I_{T})) \ind_{\{a< Z_{R}\}}\frac{1}{\epsilon} \int_{T\wedge \delta_R(a)}^{V\wedge \delta_R(a)} 1_{\{e(Z_{v})-e(I_{v})\leq
	\epsilon\}}d\langle e(Z)
\rangle_{v} ]\\
=&E[H \ind_{\{a< Z_{R}\}} \frac{1}{\epsilon} \int_{T\wedge \delta_R(a)}^{V\wedge \delta_R(a)} d\langle e(Z)
\rangle_{v}\ f(e(U(T,v))\wedge
e(I_{v}))1_{\{e(Z_{v})-\epsilon\leq
	e(I_{v})\}} ]\\
=&E[H \ind_{\{a< Z_{R}\}} \frac{1}{\epsilon} \int_{T\wedge \delta_R(a)}^{V\wedge \delta_R(a)} d\langle e(Z)
\rangle_{v} \ e(Z_{v}) \int_{e^{-1}(e(Z_{v})-\epsilon)}^{Z_{v}}  d(e(c)^{-1})f(e(U(T,v))\wedge
e(c)) ]\\
=&E[H \ind_{\{a< Z_{R}\}} \frac{1}{\epsilon} \int_{T}^{V} \ind_{\{v<\delta_R(a)\}}d\langle e(Z)
\rangle_{v} \ e(Z_{v}) \int_{e^{-1}(e(Z_{v})-\epsilon)}^{Z_{v}}  d(e(c)^{-1})f(e(U(T,v))\wedge
e(c)) ]\\ 
=&-2E[H\ind_{\{a< I_{R}\}} \frac{1}{\epsilon} \int_{T}^{V} de(I_{v}) e(Z_{v})^{2}
\int_{e^{-1}(e(Z_{v})-\epsilon)}^{Z_{v}} d(e(c)^{-1})\  f(e(U(T,v))\wedge
e(c)) ]\\
&\mbox{ according to Lemma \ref{corollary:projection}}\\
=&-2E[H\ind_{\{a< I_{R}\}} \frac{1}{\epsilon} \int_{T}^{V} de(I_{v}) e(I_{v})^{2}
\int_{e^{-1}(e(I_{v})-\epsilon)}^{I_{v}} d(e(c)^{-1})\  f(e(U(T,v))\wedge
e(c)) ]\\
\end{array}\]
The above last term is divided into two parts
$$
\dcb
-2E[H\ind_{\{a< I_{R}\}} \frac{1}{\epsilon} \int_{T}^{V} de(I_{v}) e(I_{v})^{2}
\int_{e^{-1}(e(I_{v})-\epsilon)}^{I_{v}} d(e(c)^{-1})\  f(e(U(T,v))\wedge
e(I_v)) ]\\
\\
-2E[H\ind_{\{a< I_{R}\}} \frac{1}{\epsilon} \int_{T}^{V} de(I_{v}) e(I_{v})^{2}
\int_{e^{-1}(e(I_{v})-\epsilon)}^{I_{v}} d(e(c)^{-1})\  (f(e(U(T,v))\wedge
e(c))-f(e(U(T,v))\wedge
e(I_v))) ]
\dce
$$
The first part is computed as follows :
$$
\dcb
&&-2E[H\ind_{\{a< I_{R}\}} \frac{1}{\epsilon} \int_{T}^{V} de(I_{v}) e(I_{v})^{2}
\int_{e^{-1}(e(I_{v})-\epsilon)}^{I_{v}} d(e(c)^{-1})\  f(e(I_T)) ]\\
&=&2E[H\ind_{\{a< I_{R}\}} f(e(I_T)) \frac{1}{\epsilon} \int_{T}^{V} de(I_{v}) e(I_{v})^{2}
(-\frac{1}{e(I_v)}+\frac{1}{e(I_v)-\epsilon})]\\
&=&2E[H\ind_{\{a< I_{R}\}} f(e(I_T)) \frac{1}{\epsilon} \int_{T}^{V} de(I_{v}) e(I_{v})^{2}
\frac{\epsilon}{e(I_v)(e(I_v)-\epsilon)}]\\
&=&2E[H\ind_{\{a< I_{R}\}} f(e(I_T)) \int_{T}^{V} de(I_{v}) e(I_{v})^{2}\frac{1}{e(I_v)(e(I_v)-\epsilon)}]\\
\dce
$$
When $\epsilon$ decreases down to zero, this quantity increases to $$
2E[H\ind_{\{a\leq I_{R}\}} f(e(I_T))(e(I_{V})-e(I_T))]
$$
The second part is overestimated by
$$
\dcb
&&\left|-2E[H\ind_{\{a< I_{R}\}} \frac{1}{\epsilon} \int_{T}^{V} de(I_{v}) e(I_{v})^{2}
\int_{e^{-1}(e(I_{v})-\epsilon)}^{I_{v}} d(e(c)^{-1})\  (f(e(U(T,v))\wedge
e(c))-f(e(U(T,v))\wedge
e(I_v))) ]\right|\\

&\leq&-2E[H\ind_{\{a< I_{R}\}} \frac{1}{\epsilon} \int_{T}^{V} de(I_{v}) e(I_{v})^{2}
\int_{e^{-1}(e(I_{v})-\epsilon)}^{I_{v}} d(e(c)^{-1})\|f'\|_\infty  |e(c)-e(I_v)| ]\\

&\leq&-2E[H\ind_{\{a< I_{R}\}} \frac{1}{\epsilon} \int_{T}^{V} de(I_{v}) e(I_{v})^{2}
\int_{e^{-1}(e(I_{v})-\epsilon)}^{I_{v}} d(e(c)^{-1})\|f'\|_\infty  \epsilon ]\\

&\leq&-2E[H\ind_{\{a< I_{R}\}} \int_{T}^{V} de(I_{v}) e(I_{v})^{2}
\int_{e^{-1}(e(I_{v})-\epsilon)}^{I_{v}} d(e(c)^{-1})\|f'\|_\infty]\\

&=&2E[H\ind_{\{a< I_{R}\}} \int_{T}^{V} de(I_{v}) e(I_{v})^{2}
(-\frac{1}{e(I_v)}+\frac{1}{e(I_v)-\epsilon})\|f'\|_\infty]\\

&=&2E[H\ind_{\{a< I_{R}\}} \int_{T}^{V} de(I_{v}) e(I_{v})^{2}
\frac{\epsilon}{e(I_v)(e(I_v)-\epsilon)}\|f'\|_\infty]\\

&\leq&\epsilon 2E[H\ind_{\{a< I_{R}\}} \int_{T}^{V} de(I_{v}) e(I_{v})^{2}
\frac{1}{e(I_v)^2}\|f'\|_\infty]\\

&\leq&\epsilon 2E[H\ind_{\{a< I_{R}\}} (e(I_{V})-e(I_T))\|f'\|_\infty]\\

\dce
$$
which tends to zero when $\epsilon$ goes down to zero. In sum, we have
$$
E[Hf(e(I_{T})) \ind_{\{a< Z_{R}\}}(l^0_{V\wedge \delta_R(a)}-l^0_{T\wedge \delta_R(a)})]
=
2E[H\ind_{\{a< I_{R}\}} f(e(I_T))(e(I_{V})-e(I_T))]
$$
Using the formula preceding Lemma \ref{corollary:projection}, we write$$
E[Hf(e(I_{T})) \ind_{\{a< I_{R}\}}(l^0_{V}-l^0_{T})]
=
2E[H\ind_{\{a< I_{R}\}} f(e(I_T))(e(I_{V})-e(I_T))]
$$
or equivalently
$$
E[H\ind_{\{a< I_{R}\}}f(e(I_{T}))\int_T^V dl^0_{s}]
=
2E[H\ind_{\{a\leq I_{R}\}} f(e(I_T))\int_T^V de(I)_s]
$$
Recall that, according to Lemma \ref{lemma:tribu}, the family of all processes of the form $\ind_{\{a\leq I_R\}}Hf(I_T)\ind_{(T,V]}$ generates all $\mathbb{G}$ predictable processes. The above identity means then that, on $\mathcal{P}(\mathbb{G})$, the Dolean-Dade measure of $l^0$ and of $2e(I)$ coincide. That is the theorem. \ok

\bethe
The process $e(Z)$ is a $\mathbb{G}$-semimartingale. Its canonical decomposition is given by :\[
e(Z) = e(Z_0) +M + 2e(I)-2e(I_0) + \frac{1}{e(Z)}\cdot \langle e(Z) \rangle,
\]
where $M$ is a $\mathbb{G}$ local martingale. 
\ethe

Applying Ito's formula, we obtain also an equivalent result

\bcor
The process
$\frac{1}{e(Z)}-\frac{2}{e(I)}$ is a $\mathbb{G}$ local  martingale.
\ecor

\pagebreak

\

\begin{center}
\Large \textbf{Part III. 
Drift operator and the viability}
\end{center}

\pagebreak

\section{Introduction}

A financial market is modeled by a triplet $(\mathbb{P},\mathbb{F},S)$ of a probability measure $\mathbb{P}$, of an information flow $\mathbb{F}=(\mathcal{F}_t)_{t\in\mathbb{R}_+}$, and of an $\mathbb{F}$ adapted asset process $S$. The basic requirement about such a model is its viability (in the sense of no-arbitrage). There are situations where one should consider the asset process $S$ in an enlarged information flow $\mathbb{G}=(\mathcal{G}_t)_{t\geq 0}$ with $\mathcal{G}_t\supset\mathcal{F}_t$. The viability of the new market $(\mathbb{P},\mathbb{G},S)$ is not guaranteed. The purpose of this paper is to find such conditions that the viability will be maintained despite the expansion of the information flow.

Concretely, we introduce the notion of the full viability for information expansions (cf. subsection \ref{viabilitydef} Assumption \ref{fullviabilitydef}). This means that, for any $\mathbb{F}$ special semimartingale asset process $S$, if $(\mathbb{P},\mathbb{F},S)$ is viable, the expansion market $(\mathbb{P},\mathbb{G},S)$ also is viable. Under the assumption of martingale representation property in $(\mathbb{P},\mathbb{F})$, we prove that (cf. Theorem \ref{fullviability}) the full viability is equivalent to the following fact: there exist a (multi-dimensional) $\mathbb{G}$ predictable process $\varphi$ and a (multi-dimensional) $\mathbb{F}$ local martingale $N$, such that (1) for any $\mathbb{F}$ local martingale $X$, the expression $\transp\varphi \centerdot[N,X]^{\mathbb{F}\cdot p}$ is well-defined and $X-\transp\varphi \centerdot[N,X]^{\mathbb{F}\cdot p}$ is a $\mathbb{G}$ local martingale; (2) the continuous increasing process $\transp{\varphi}(\centerdot[N^c,\transp N^c]){\varphi}$ is finite; (3) the jump increasing process 
$
(
\sum_{0<s\leq t}\left(\frac{\transp {\varphi}_{s}\Delta_{s} N}{1+\transp{\varphi}_{s}\Delta_{s} N}
\right)^2
)^{1/2}
$,
$t\in\mathbb{R}_+$,
is $(\mathbb{P},\mathbb{G})$ locally integrable.

It is to note that, if no jumps occurs in $\mathbb{F}$, continuous semimartingale calculus gives a quick solution to the viability problem of the information expansion. The situation becomes radically different when jumps occur, especially because we need to compute and to compare the different projections in $\mathbb{F}$ and in $\mathbb{G}$. To have an idea about the implication of jumps in the study of the market viability, the articles \cite{KC2010, K2012} give a good illustration. In this paper we come to a satisfactory result in a general jump situation, thanks to a particular property derived from the martingale representation. In fact, when a process $W$ has a the martingale representation property, the jump $\Delta W$ of this process can only take a finite number of \texttt{"}predictable\texttt{"} values. We refer to \cite{song-mrp-drift} for a detailed account, where it is called the finite predictable constraint condition (which has a closed link with the notion of multiplicity introduced in \cite{BEKSY}).  

Usually the martingale representation property is mentioned to characterize a specific process (a Brownian motion, for example). But, in this paper, what is relevant is the stochastic basis $(\mathbb{P},\mathbb{F})$ having a martingale representation property, whatever representation process is. One of the fundamental consequences of the finite predictable constraint condition is the possibility to find a finite family of very simply locally bounded mutually \texttt{"}avoiding\texttt{"} processes which have again the martingale representation property. This possibility reduces considerably the computation complexity and gives much clarity to delicate situations.

The viability property is fundamental for financial market modeling. There exists a huge literature (cf. for example, \cite{choulli, Fon, imkeller, ing, kabanov, KC2010, K2012,  Kar, loew, Sch1, Sch2, takaoka}). Recently, there is a particular attention on the viability problem related to expansions of information flow (cf. \cite{AFK,ACDJ,FJS,Song-NA}). It is to notice that, however, the most of the works on expansions of information flow follow two specific ideas : the initial enlargement of filtration or the progressive enlargement of filtration (cf. \cite{JYC, Jeulin80, JYexample, protter,MY} for definition and properties). In this paper, we take the problem in a very different perspective. We obtain general result, without assumption on the way that $\mathbb{G}$ is constructed from $\mathbb{F}$. It is known (cf. \cite{SongThesis, Song-local}) that the initial or progressive enlargement of filtration are particular situations covered by the so-called local solution method. The methodology of this paper does not take part in this category, adding a new element in the arsenal of filtration analysis.

The concept of information is a fascinating, but also a difficult notion, especially when we want to quantify it. The framework of enlargement of filtrations $\mathbb{F}\subset\mathbb{G}$ offers since long a nice laboratory to test the ideas. In general, no common consensus exists how to quantify the difference between two information flows $\mathbb{F}$ and $\mathbb{G}$. The notion of entropy has been used there (see for example \cite{ADI, yor-entro}). But a more convincing measurement of information should be the drift operator $\Gamma(X)$, i.e. the operator which gives the drift part of the $\mathbb{F}$ local martingale $X$ in $\mathbb{G}$ (cf. Lemma \ref{linearG}).  This observation is strengthened by the result of the present paper. We have seen that, in the case of our paper, the drift operator takes the form $\Gamma(X)=\transp\varphi \centerdot[N,X]^{\mathbb{F}\cdot p}$ for two processes $\varphi$ and $N$ (the drift multiplier assumption, cf. Assumption \ref{assump1}), and the full viability of the information expansion is completely determined by the size of the positive quantities $\transp\varphi\centerdot[N^c,\transp N^c]\varphi$ and $\frac{1}{1+\transp\varphi\Delta N}$, which have all the appearance of a measure. See \cite{ankir2} for complementary discussion. See also \cite{JS3} for a use of $\Gamma$ in the study of the martingale representation property in $\mathbb{G}$.

\

\section{Notations and vocabulary}\label{vocabulary}

We employ the vocabulary of stochastic calculus as defined in \cite{HWY, jacod} with the specifications below.

\

\textbf{Probability space and random variables}

A stochastic basis $(\Omega, \mathcal{A},\mathbb{P},\mathbb{F})$ is a quadruplet, where $(\Omega, \mathcal{A},\mathbb{P})$ is a probability space and $\mathbb{F}$ is a filtration of sub-$\sigma$-algebras of $\mathcal{A}$, satisfying the usual conditions. 

The relationships involving random elements are always in the almost sure sense. For a random variable $X$ and a $\sigma$-algebra $\mathcal{F}$, the expression $X\in\mathcal{F}$ means that $X$ is $\mathcal{F}$-measurable. The notation $\mathbf{L}^p(\mathbb{P},\mathcal{F})$ denotes the space of $p$-times $\mathbb{P}$-integrable $\mathcal{F}$-measurable random variables.

\

\textbf{The processes}

The jump process of a c\` adl\`ag process $X$ is denoted by $\Delta X$, whilst the jump at time $t\geq 0$ is denoted by $\Delta_tX$. By definition, $\Delta_0X=0$ for any c\` adl\`ag process $X$. When we call a process $A$ a process having finite variation, we assume automatically that $A$ is c\`adl\`ag. We denote then by $\mathsf{d}A$ the (signed) random measure that $A$ generates. 

An element $v$ in an Euclidean space $\mathbb{R}^d$ ($d\in\mathbb{N}^*$) is considered as a vertical vector. We denote its transposition by $\transp v$. The components of $v$ will be denoted by $v_h, 1\leq h\leq d$. 

We deal with finite family of real processes $X=(X_h)_{1\leq h\leq d}$. It will be considered as $d$-dimensional vertical vector valued process. The value of a component $X_h$ at time $t\geq 0$ will be denoted by $X_{h,t}$. When $X$ is a semimartingale, we denote by $[X,\transp X]$ the $d\times d$-dimensional matrix valued process whose components are $[X_i,X_j]$ for $1\leq i,j\leq d$.

\

\textbf{The projections}

With respect to a filtration $\mathbb{F}$, the notation ${^{\mathbb{F}\cdot p}}\bullet$ denotes the predictable projection, and the notation $\bullet^{\mathbb{F}\cdot p}$ denotes the predictable dual projection. 

\

\textbf{The martingales and the semimartingales}

Fix a probability $\mathbb{P}$ and a filtration $\mathbb{F}$. For any $(\mathbb{P},\mathbb{F})$ special semimartingale $X$, we can decompose $X$ in the form (see \cite[Theorem 7.25]{HWY}) :$$
\dcb
X=X_0+X^m+X^v,\
X^m=X^c+X^{da}+X^{di},
\dce
$$
where $X^m$ is the martingale part of $X$ and $X^v$ is the drift part of $X$, $X^c$ is the continuous martingale part, $X^{da}$ is the part of compensated sum of accessible jumps, $X^{di}$ is the part of compensated sum of totally inaccessible jumps. We recall that this decomposition of $X$ depends on the reference probability and the reference filtration. We recall that every part of the decomposition of $X$, except $X_0$, is assumed null at $t=0$.

\

\textbf{The stochastic integrals}

In this paper we employ the notion of stochastic integral only about the predictable processes. The stochastic integral are defined as 0 at $t=0$. We use a point \texttt{"}$\centerdot$\texttt{"} to indicate the integrator process in a stochastic integral. For example, the stochastic integral of a real predictable process ${H}$ with respect to a real semimartingale $Y$ is denoted by ${H}\centerdot Y$, while the expression $\transp{K}(\centerdot[X,\transp X]){K}$ denotes the process$$
\int_0^t \sum_{i=1}^k\sum_{j=1}^k{K}_{i,s}{K}_{j,s} \mathsf{d}[X_i,X_j]_s,\ t\geq 0,
$$
where ${K}$ is a $k$-dimensional predictable process and $X$ is a $k$-dimensional semimartingale. The expression $\transp{K}(\centerdot[X,\transp X]){K}$ respects the matrix product rule. The value at $t\geq 0$ of a stochastic integral will be denoted, for example, by $\transp{K}(\centerdot[X,\transp X]){K}_t$.

The notion of the stochastic integral with respect to a multi-dimensional local martingale $X$ follows \cite{Jacodlivre}. We say that a (multi-dimensional) $\mathbb{F}$ predictable process is integrable with respect to $X$ under the probability $\mathbb{P}$ in the filtration $\mathbb{F}$, if the non decreasing process $\sqrt{\transp{H}(\centerdot[X,\transp X]){H}}$ is $(\mathbb{P},\mathbb{F})$ locally integrable. For such an integrable process ${H}$, the stochastic integral $\transp{H}\centerdot X$ is well-defined and the bracket process of $\transp{H}\centerdot X$ can be computed using \cite[Remarque(4.36) and Proposition(4.68)]{Jacodlivre}. Note that two different predictable processes may produce the same stochastic integral with respect to $X$. In this case, we say that they are in the same equivalent class (related to $X$).

The notion of multi-dimensional stochastic integral is extended to semimartingales. We refer to \cite{JacShi} for details.

\textbf{Caution.}
Note that some same notations are used in different parts of the paper for different meaning. Also some definitions will be repeated in different parts of the paper.

\

\section{Three fundamental concepts}

Three notions play particular roles in this paper.

\subsection{Enlargements of filtrations and Hypothesis$(H')$ }

Let $\mathbb{F}=(\mathcal{F}_t)_{t\geq 0}$ and $\mathbb{G}=(\mathcal{G}_t)_{t\geq 0}$ be two filtrations on a same probability space such that $\mathcal{F}_t\subset\mathcal{G}_t$. We say then that $\mathbb{G}$ is an expansion (or an enlargement) of the filtration $\mathbb{F}$. Let $T$ be a $\mathbb{G}$ stopping time. We introduce the Hypothesis$(H')$ (cf. \cite{JYC, Jeulin80, JYexample, protter, MY}):   

\bassumption
(\textbf{Hypothesis$(H')$} on the time horizon $[0,T]$) All $(\mathbb{P},\mathbb{F})$ local martingale is a $(\mathbb{P},\mathbb{G})$ semimartingale on $[0,T]$.
\eassumption

Whenever Hypothesis$(H')$ holds, the associated drift operator can be defined (cf. \cite{song-mrp-drift}).

\bl\label{linearG}
Suppose hypothesis$(H')$ on $[0,T]$. Then there exists a linear map $\Gamma$ from the space of all $(\mathbb{P},\mathbb{F})$ local martingales into the space of c\`adl\`ag $\mathbb{G}$-predictable processes on $[0,T]$, with finite variation and null at the origin, such that, for any $(\mathbb{P},\mathbb{F})$ local martingale $X$, $\widetilde{X}:=X-\Gamma(X)$ is a $(\mathbb{P},\mathbb{G})$ local martingale on $[0,T]$. Moreover, if $X$ is a $(\mathbb{P},\mathbb{F})$ local martingale and $H$ is an $\mathbb{F}$ predictable $X$-integrable process, then $H$ is $\Gamma(X)$-integrable and $\Gamma(H\centerdot X)=H\centerdot \Gamma(X)$ on $[0,T]$. The operator $\Gamma$ will be called the drift operator.
\el

\

\subsection{The martingale representation property}\label{MRPsection}

Let us fix a stochastic basis $(\Omega, \mathcal{A},\mathbb{P},\mathbb{F})$. We consider a multi-dimensional stochastic process $W$. We say that $W$ has the martingale representation property in the filtration $\mathbb{F}$ under the probability $\mathbb{P}$, if $W$ is a $(\mathbb{P},\mathbb{F})$ local martingale, and if all $(\mathbb{P},\mathbb{F})$ local martingale is a stochastic integral with respect to $W$. We say that the martingale representation property holds in the filtration $\mathbb{F}$ under the probability $\mathbb{P}$, if there exists a local martingale $W$ which possesses the martingale representation property. In this case we call $W$ the representation process.

\subsubsection{The choice of representation process}\label{www}

Recall the result in \cite{song-mrp-drift}. Suppose the martingale representation property in $(\mathbb{P},\mathbb{F})$. Reconstituting the original representation process if necessary, we can find a particular representation process in the form $W=(W',W'',W''')$, where $W', W''$ denote respectively the processes defined in \cite[Formulas (4) and (5)]{song-mrp-drift} and $W'''$ denote the process $X^\circ$ in \cite[Section 4.5]{song-mrp-drift}. The processes $W',W'',W'''$ are locally bounded $(\mathbb{P},\mathbb{F})$ local martingales; $W'$ is continuous; $W''$ is purely discontinuous with only accessible jump times; $W'''$ is purely discontinuous with only totally inaccessible jump times; the three components $W',W'',W'''$ are mutually pathwisely orthogonal; the components of $W'$ are mutually pathwisely orthogonal; the components of $W'''$ are mutually pathwisely orthogonal. Let $\mathsf{n}',\mathsf{n}'',\mathsf{n}'''$ denote respectively the dimensions of the three components $W',W'',W'''$. We know that, if $d$ denotes the dimension of the original representation process, $\mathsf{n}'=d$ and $\mathsf{n}''=1+d$. (Notice that some components may by null.) Let $H$ be a $\mathbb{F}$ predictable $W$-integrable process. The process $H$ is naturally cut into three components $(H',H'',H''')$ corresponding to $(W',W'',W''')$. The pathwise orthogonality implies that $H'_h$ is $W'_h$-integrable for $1\leq h\leq d$, $H''$ is $W''$-integrable, and $H'''_h$ is $W'''_h$-integrable for $1\leq h\leq \mathsf{n}'''$.

The finite predictable constraint condition holds (cf. \cite{song-mrp-drift}). There exists a $\mathsf{n}'''$-dimensional $\mathbb{F}$ predictable process $\alpha'''$ such that$$
\Delta W'''_h = \alpha'''_h \ind_{\{\Delta W'''_h\neq 0\}},\
1\leq h\leq \mathsf{n}'''.
$$
Let $(T_n)_{1\leq n<\mathsf{N}^a}$ ($\mathsf{N}^a$ being a finite or infinite integer) be a sequence of strictly positive $(\mathbb{P},\mathbb{F})$ predictable stopping times such that $[T_n]\cap [T_{n'}]=\emptyset$ for $n\neq n'$ and $\{s\geq 0:\Delta_sW''\neq 0\}\subset\cup_{n\geq 1}[T_n]$. For every $1\leq n<\mathsf{N}^a$, there exists a partition $(A_{n,0},A_{n,1},\ldots,A_{n,d})$ such that $
\mathcal{F}_{T_n}=\mathcal{F}_{T_n-}\vee\sigma(A_{n,0},A_{n,1},A_{n,2},\ldots,A_{n,d})
$
(a finite multiplicity according to \cite{BEKSY}). Denote $p_{n,k}=\mathbb{P}[A_{n,k}|\mathcal{F}_{T_n-}], 0\leq k\leq d$. We have  
$$
W''_k
=
\sum_{n=1}^{\mathsf{N}^a-}\frac{1}{2^n}(\ind_{A_{n,k}}-p_{n,k})\ind_{[T_n,\infty)}.
$$
Let us denote by $\mathsf{a}_{n}$ the vector $(\ind_{A_{n,h}})_{0\leq h\leq d}$ and by $p_n$ the vector $(p_{n,h})_{0\leq h\leq d}$, so that $\Delta_{T_n}W''=\frac{1}{2^n}(\mathsf{a}_{n}-p_n)\ind_{[T_n,\infty)}$.

\subsubsection{Coefficient in martingale representation}

If the martingale representation property holds, the $(\mathbb{P},\mathbb{F})$ local martingale $X$ takes all the form $\transp H\centerdot W$ for some $W$-integrable predictable process. We call (any version of) the process $H$ the coefficient of $X$ in its martingale representation with respect to the process $W$. This appellation extends naturally to vector valued local martingales.

When we make computation with the martingale representation, we often need to extract information about a particular stopping time from an entire stochastic integral. The following lemma is proved in \cite[Lemma 3.1]{song-mrp-drift}.

\bl\label{single-jump}
Let $R$ be any $\mathbb{F}$ stopping time. Let $\xi\in\mathbf{L}^1(\mathbb{P},\mathcal{F}_{R})$. Let  $H$ denote any coefficient of the $(\mathbb{P},\mathbb{F})$ martingale  $\xi\ind_{[R,\infty)}-(\xi\ind_{[R,\infty)})^{\mathbb{F}\cdot p}$ in its martingale representation with respect to $W$.
\ebe
\item
If $R$ is predictable, the two predictable processes $H$ and $H\ind_{[R]}$ are in the same equivalent class related to $W$, whose value is determined by the equation on $\{R<\infty\}$
$$
\transp H_{R} \Delta_{R}W=
\xi-\mathbb{E}[\xi|\mathcal{F}_{R-}].
$$
\item
If $R$ is totally inaccessible, the process $H$ satisfies the equations  on $\{R<\infty\}$
$$
\transp H_{R} \Delta_{R}W=\xi,\ \mbox{ and } \
\transp H_{S} \Delta_{S}W =0 \mbox{ on $\{S\neq R\}$},
$$
for any $\mathbb{F}$ stopping time $S$.
\dbe
\el

\

\subsection{Structure condition}\label{SC-def}

Let given a stochastic basis $(\Omega, \mathcal{A},\mathbb{P},\mathbb{F})$.

\bd
Let $T>0$ be an $\mathbb{F}$ stopping time. We say that a multi-dimensional $(\mathbb{P},\mathbb{F})$ special semimartingale $S$ satisfies the structure condition in the filtration $\mathbb{F}$ under the probability $\mathbb{P}$ on the time horizon $[0,T]$, if there exists a real $(\mathbb{P},\mathbb{F})$ local martingale $D$ such that, on the time interval $[0,T]$, $D_0=0, \Delta D<1$, $[S^m_{i}, D]^{\mathbb{F}\cdot p}$ exists, and $S^{v}_{i} = [S^m_{i}, D]^{\mathbb{F}\cdot p}$ for all components $S_i$. We will call $D$ a structure connector. 
\ed

\bd\label{scdef}
Let $T>0$ be a $\mathbb{F}$ stopping time. We call a strictly positive $\mathbb{F}$ adapted real process $X$ with $X_0=1$, a local martingale deflator on the time horizon $[0,T]$ for a (multi-dimensional) $(\mathbb{P},\mathbb{F})$ special semimartingale $S$, if the processes $X$ and $X S$ are $(\mathbb{P},\mathbb{F})$ local martingales on $[0,T]$. 
\ed

\brem
The notion of \texttt{"}structure condition\texttt{"} exists in the literature, particularly in \cite{choulli, Sch2}. The above Definition \ref{scdef} employs the same name for something slightly different from the original one, Our definition adapts especially to the study of this paper. We use the same appellation because Definition \ref{scdef} can be considered as an extension of the original one. \ok
\erem

We recall that the existence of local martingale deflators and the structure condition are conditions equivalent to the no-arbitrage conditions \texttt{NUPBR} and \texttt{NA1} (cf. \cite{KC2010, takaoka}). We know that, when the no-arbitrage condition \texttt{NUPBR} is satisfied, the market is viable, and vice versa.

\bethe\label{deflator-connector}
Let $T>0$ be an $\mathbb{F}$ stopping time. A (multi-dimensional) special semimartingale $S$ possesses a local martingale deflator $X$ in $(\mathbb{P},\mathbb{F})$ on the time horizon $[0,T]$, if and only if $S$ satisfies the structure condition on the time horizon $[0,T]$ with a structure connector $D$. In this case, $X=\mathcal{E}(-D)$.
\ethe

\textbf{Proof.} We know that a strictly positive local martingale is always a Dolean-Dade exponential (cf.\cite[Theorem 9.41]{HWY} or \cite{jacod,choulli-stricker}).  The lemma is a consequence of the integration by parts formula$$
XS = X_0S_0+ S_-\centerdot X + X_-\centerdot S - X_-\centerdot [S,D].
$$ 
In particular, if $X$ and $XS$ are local martingales, $[S,D]$ is locally integrable. \ok

\

\section{Main results}\label{DMA}

Together with a given stochastic basis $(\Omega,\mathcal{A}, \mathbb{P},\mathbb{F})$ let $\mathbb{G}$ be an enlargement of $\mathbb{F}$. 

\subsection{Drift multiplier assumption and full viability}\label{viabilitydef}

Our study is based on the following two notions.

\bassumption\label{fullviabilitydef}
(\textbf{Full viability} on $[0,T]$) Let $T$ be a $\mathbb{G}$ stopping time. The expanded information flow $\mathbb{G}$ is fully viable on $[0,T]$. This means that, for any $\mathbb{F}$ asset process $S$ (multi-dimensional special semimartingale with strictly positive components) satisfying the structure condition in $\mathbb{F}$, the process $S$ satisfies the structure condition in the expanded market environment $(\mathbb{P},\mathbb{G})$ on the time horizon $[0,T]$. 
\eassumption

\brem\label{locallyboundedM}
As indicated in \cite{song-mrp-drift}, the full viability implies that, for any $(\mathbb{P},\mathbb{F})$ locally bounded local martingale $M$, $M$ satisfies the structure condition in $(\mathbb{P},\mathbb{G})$ on $[0,T]$.
\erem

\bassumption\label{assump1} (\textbf{Drift multiplier assumption}) Let $T$ be a $\mathbb{G}$ stopping time.
\
\ebe 
\item
The Hypothesis$(H')$ is satisfied on the time horizon $[0,T]$ with a drift operator $\Gamma$.

\item
There exist $N=(N_1,\ldots,N_\mathsf{n})$ an $\mathsf{n}$-dimensional $(\mathbb{P},\mathbb{F})$ local martingale, and ${\varphi}$ an $\mathsf{n}$ dimensional $\mathbb{G}$ predictable process such that, for any $(\mathbb{P},\mathbb{F})$ local martingale $X$, $[N,X]^{\mathbb{F}\cdot p}$ exists, ${\varphi}$ is $[N,X]^{\mathbb{F}\cdot p}$-integrable, and $$
\Gamma(X)=\transp{\varphi}\centerdot [N,X]^{\mathbb{F}\cdot p}
$$
on the time horizon $[0,T]$. 
\dbe
\eassumption

We will need frequently the following consequence of the drift multiplier assumptions \ref{assump1}.

\bl\label{A-G-p}
For any $\mathbb{F}$ adapted c\`adl\`ag process $A$ with $(\mathbb{P},\mathbb{F})$ locally integrable variation, we have $$
A^{\mathbb{G}\cdot p}=A^{\mathbb{F}\cdot p}+\Gamma(A-A^{\mathbb{F}\cdot p})=A^{\mathbb{F}\cdot p}+\transp{\varphi}\centerdot[N,A]^{\mathbb{F}\cdot p}
$$
on $[0,T]$. In particular, for $R$ an $\mathbb{F}$ stopping time, for $\xi\in\mathbf{L}^1(\mathbb{P},\mathcal{F}_{R})$, $$
(\xi\ind_{[R,\infty)})^{\mathbb{G}\cdot p}
=
(\xi\ind_{[R,\infty)})^{\mathbb{F}\cdot p}+\transp{\varphi}\centerdot(\Delta_{R}N\xi\ind_{[R,\infty)})^{\mathbb{F}\cdot p}
$$
on $[0,T]$. If $R$ is $\mathbb{F}$ totally inaccessible, $R$ also is $\mathbb{G}$ totally inaccessible on $[0,T]$.
\el

\bassumption\label{1+fin}
For any $\mathbb{F}$ predictable stopping time $R$, for any positive random variable $\xi\in\mathcal{F}_R$, we have $\{\mathbb{E}[\xi|\mathcal{G}_{R-}]>0, R\leq T, R<\infty\}=\{\mathbb{E}[\xi|\mathcal{F}_{R-}]>0, R\leq T, R<\infty\}$. 
\eassumption

\brem
Clearly, if the random variable $\xi$ is already in $\mathcal{F}_{R-}$ (or if $\mathcal{F}_{R-}=\mathcal{F}_{R}$), the above set equality holds. Hence, a sufficient condition for Assumption \ref{1+fin} to be satisfied is that the filtration $\mathbb{F}$ is quasi-left-continuous.
\erem

\

\subsection{The theorems}

The two notions of full viability and the drift multiplier assumption are closely linked. According to \cite{song-mrp-drift}, under the martingale representation property, the full viability on $[0,T]$ implies the drift multiplier assumption. The aim of this paper is to refine the above result to have an exact relationship between the drift multiplier assumption and the viability of the expanded information flow. We will prove the two theorems. Let $T$ be a $\mathbb{G}$ stopping time.

\bethe\label{main}
Suppose that $(\mathbb{P},\mathbb{F})$ satisfies the martingale representation property. Suppose the drift multiplier assumption. Let $S$ be any $(\mathbb{P},\mathbb{F})$ special semimartingale satisfying the structure condition in $(\mathbb{P},\mathbb{F})$ with a structure connector $D$. If the process $
\transp{\varphi}\centerdot[N^c,\transp N^c]{\varphi}
$
is finite on $[0,T]$ and if the process $$
\sqrt{\sum_{0<s\leq t\wedge T}\frac{1}{(1+\transp{\varphi}_{s}\Delta_{s} N)^2}\left(
\Delta_{s}D
+
\transp {\varphi}_{s}\Delta_{s} N  
\right)^2},\ t\in\mathbb{R}_+,
$$
is $(\mathbb{P},\mathbb{G})$ locally integrable, then, $S$ satisfies the structure condition in $(\mathbb{P},\mathbb{G})$. 
\ethe

\bethe\label{fullviability}
Suppose that $(\mathbb{P},\mathbb{F})$ satisfies the martingale representation property. Then, $\mathbb{G}$ is fully viable on $[0,T]$, if and only if the drift multiplier assumption is satisfied (with the processes $\varphi$ and $N$) such that
\begin{equation}\label{fn-sur-fn}
\dcb
\transp{\varphi}(\centerdot[N^c,\transp N^c]){\varphi}\ \mbox{ is a finite process on $[0,T]$ and }\\
\\
\sqrt{\sum_{0<s\leq t\wedge T}\left(\frac{\transp {\varphi}_{s}\Delta_{s} N}{1+\transp{\varphi}_{s}\Delta_{s} N}
\right)^2},\ t\in\mathbb{R}_+,\
\mbox{ is $(\mathbb{P},\mathbb{G})$ locally integrable.}
\dce
\end{equation}

\ethe

\bcor\label{commondeflator}
Suppose the $(\mathbb{P},\mathbb{F})$ martingale representation property. If condition (\ref{fn-sur-fn}) is satisfied, there exists a common $(\mathbb{P},\mathbb{G})$ deflator for all $(\mathbb{P},\mathbb{F})$ local martingales.
\ecor

\

\section{Structure condition decomposed under the martingale representation property}

We now begin the proof of Theorem \ref{main} and Theorem \ref{fullviability}. Recall that, when $(\mathbb{P},\mathbb{F})$ possesses the martingale representation property, we can choose the representation process to ease the computations. We suppose from now on the drift multiplier assumption \ref{assump1} and the following one.

\bassumption\label{assump-mrt}
$(\mathbb{P},\mathbb{F})$ is supposed to satisfy the martingale representation property, with a representation process $W$ of the form $W=(W',W'',W''')$ satisfying the conditions in subsection \ref{www} with respectively the dimensions $d, 1+d, \mathsf{n}'''$. 
\eassumption

Let $S$ be a multi-dimensional $\mathbb{F}$ asset process satisfying the structure condition in $\mathbb{F}$ with an $\mathbb{F}$ structure connector $D$. Set $M:=S^m$ (in $\mathbb{F}$). The $(\mathbb{P},\mathbb{G})$ canonical decomposition of $S$ on $[0,T]$ is given by $$
S=\widetilde{M}+[D,M]^{\mathbb{F}\cdot p}+\transp{\varphi}\centerdot [N,M]^{\mathbb{F}\cdot p}.
$$
The structure condition for $S$ in the expanded market environment $(\mathbb{P},\mathbb{G})$ takes the following form : there exists a $\mathbb{G}$ local martingale $Y$ such that $Y_0=0, \Delta Y<1$, $[Y,\widetilde{M}]^{\mathbb{G}\cdot p}$ exists, and 
\begin{equation}\label{structure-condition}
[Y,\widetilde{M}]^{\mathbb{G}\cdot p}=[D,M]^{\mathbb{F}\cdot p}+\transp{\varphi}\centerdot [N,M]^{\mathbb{F}\cdot p}
\end{equation}
on the time horizon $[0,T]$.

We consider the following specific structure conditions. (recall $\widetilde{X}=X-\Gamma(X)$.) 
\ebe
\item[. ]\textbf{Continuous structure condition related to $D$.}
For $1\leq h\leq d$, there exists a $\mathbb{G}$ predictable $\widetilde{W}'$-integrable process $K'_h$ such that, on the time horizon $[0,T]$,
\begin{equation}\label{structure-condition-c}
K'_h\centerdot[\widetilde{W}'_h, \widetilde{W}'_h]^{\mathbb{G}\cdot p}
=[D, W'_h]^{\mathbb{F}\cdot p}+\transp{\varphi}\centerdot [N, W'_h]^{\mathbb{F}\cdot p}.
\end{equation}
\vspace{13pt}

\item[. ]\textbf{Accessible structure condition related to $D$.}
There exists a $\mathbb{G}$ predictable $\widetilde{W}''$-integrable process $K''$  such that $\transp K''\Delta \widetilde{W}''<1$, and, on the time horizon $[0,T]$,
\begin{equation}\label{structure-condition-da}
\transp K''\centerdot[\widetilde{W}'',\transp \widetilde{W}'']^{\mathbb{G}\cdot p}
=[D,\transp {W}'']^{\mathbb{F}\cdot p}+\transp{\varphi}\centerdot [N,\transp {W}'']^{\mathbb{F}\cdot p}.
\end{equation}
\vspace{13pt}

\item[. ]\textbf{Totally inaccessible structure condition related to $D$.}
For $1\leq h\leq \mathsf{n}'''$, there exists a $\mathbb{G}$ predictable $\widetilde{W}'''_h$-integrable process $K'''_h$  such that $K'''_h\Delta \widetilde{W}'''_h<1$, and, on the time horizon $[0,T]$,
\begin{equation}\label{structure-condition-di}
 K'''_h\centerdot[\widetilde{W}'''_h,\widetilde{W}'''_h]^{\mathbb{G}\cdot p}
=[D,{W}'''_h]^{\mathbb{F}\cdot p}+\transp{\varphi}\centerdot [N,{W}'''_h]^{\mathbb{F}\cdot p}.
\end{equation}\vspace{17pt}
\dbe
Note that the above conditions assume in particular that all the stochastic integrals exist.

\bl\label{piece-ensemble}
Let $S$ be a multi-dimensional $\mathbb{F}$ asset process satisfying the structure condition in $\mathbb{F}$ with an $\mathbb{F}$ structure connector $D$. If the group of the conditions (\ref{structure-condition-c}), (\ref{structure-condition-da}), (\ref{structure-condition-di}) related to $D$ are satisfied, the structure condition (\ref{structure-condition}) for $S$ in $\mathbb{G}$ is satisfied.
\el

\textbf{Proof.} 
Write the martingale representation of $M:=S^m$ (in $\mathbb{F}$): $$
\dcb
M&=&H'\centerdot W'+H''\centerdot W''+H'''\centerdot W'''.
\dce
$$ 
for some $W$-integrable $\mathbb{F}$ predictable processes $(H', H'', H''')$. Let $K'_h,K'',K'''_h$ be the solutions of respectively (\ref{structure-condition-c}), (\ref{structure-condition-da}), (\ref{structure-condition-di}). Set $K':=(K'_h)_{1\leq h\leq d}$, $K''':=(K'''_h)_{1\leq h\leq \mathsf{n}'''}$ and define $$
Y:=\transp K'\centerdot\widetilde{W}'
+\transp K''\centerdot\widetilde{W}''
+\transp K'''\centerdot\widetilde{W}'''.
$$
Note that, with the drift multiplier assumption, $\Gamma(W'')$ has only $\mathbb{F}$ predictable jumps and $\Gamma(W'),\Gamma(W''')$ are continuous so that $\Delta\widetilde{W}'''=\Delta{W}'''$. With the pathwise orthogonality of the processes ${W}'',{W}'''_h, 1\leq h\leq \mathsf{n}''',$ (cf. subsection \ref{www}), we see that $\Delta Y<1$. With the integrability of $H',H'',H'''$ with respect to separately $W',W'',W'''$ (cf. subsection \ref{www}) and \cite[Lemma 2.2]{JS3},$$
\widetilde{M}=H'\centerdot \widetilde{W}'+H''\centerdot \widetilde{W}''+H'''\centerdot \widetilde{W}'''.
$$
Therefore,
$$
\dcb
[Y,\transp \widetilde{M}]
&=&
[\transp K'\centerdot\widetilde{W}'
+\transp K''\centerdot\widetilde{W}''
+\transp K'''\centerdot\widetilde{W}''', \ \ H'\centerdot \widetilde{W}'+H''\centerdot \widetilde{W}''+H'''\centerdot \widetilde{W}''']\\

&=&
\sum_{h=1}^{d} K'_h\centerdot[\widetilde{W}'_h,\widetilde{W}'_h]H'_h
+\transp K''\centerdot[\widetilde{W}'',\transp \widetilde{W}'']H''
+\sum_{h=1}^{\mathsf{n}'''} K'''_h\centerdot[\widetilde{W}'''_h, \widetilde{W}'''_h]H'''_h.
\dce
$$
Let $L>0$ and define $\mathtt{B}=\{|H|\leq L\}$. We can write
$$
\dcb
&&(\ind_{\mathtt{B}}\centerdot[Y,\transp \widetilde{M}])^{\mathbb{G}\cdot p}\\

&=&
\left(
\sum_{h=1}^{d} K'_h\centerdot[\widetilde{W}'_h,\widetilde{W}'_h]H'_h\ind_{\mathtt{B}}
+\transp K''\centerdot[\widetilde{W}'',\transp \widetilde{W}'']H''\ind_{\mathtt{B}}
+\sum_{h=1}^{\mathsf{n}'''} K'''_h\centerdot[\widetilde{W}'''_h, \widetilde{W}'''_h]H'''_h\ind_{\mathtt{B}}
\right)^{\mathbb{G}\cdot p}
\\

&=&

\sum_{h=1}^{d} K'_h\centerdot[\widetilde{W}'_h,\widetilde{W}'_h]^{\mathbb{G}\cdot p}H'_h\ind_{\mathtt{B}}
+\transp K''\centerdot[\widetilde{W}'',\transp \widetilde{W}'']^{\mathbb{G}\cdot p}H''\ind_{\mathtt{B}}
+\sum_{h=1}^{\mathsf{n}'''} K'''_h\centerdot[\widetilde{W}'''_h, \widetilde{W}'''_h]^{\mathbb{G}\cdot p}H'''_h\ind_{\mathtt{B}}
\\

&=&
+\sum_{h=1}^{d}(\centerdot[D,\transp {W}'_h]^{\mathbb{F}\cdot p}H'_h\ind_{\mathtt{B}}+\transp{\varphi}\centerdot [N,\transp {W}'_h]^{\mathbb{F}\cdot p}H'_h\ind_{\mathtt{B}})\\
&&
+\centerdot[D,\transp {W}'']^{\mathbb{F}\cdot p}H''\ind_{\mathtt{B}}+\transp{\varphi}\centerdot [N,\transp {W}'']^{\mathbb{F}\cdot p}H''\ind_{\mathtt{B}}\\
&&
+
\sum_{h=1}^{\mathsf{n}'''}(\centerdot[D,\transp {W}'''_h]^{\mathbb{F}\cdot p}H'''_h\ind_{\mathtt{B}}+\transp{\varphi}\centerdot [N,\transp {W}'''_h]^{\mathbb{F}\cdot p}H'''_h\ind_{\mathtt{B}})\\

&=&
+\ind_{\mathtt{B}}\centerdot[D,\transp M]^{\mathbb{F}\cdot p}+\ind_{\mathtt{B}}\transp{\varphi}\centerdot [N,\transp M]^{\mathbb{F}\cdot p}.
\dce
$$
This formula for any $L>0$ shows firstly that $[Y,\transp \widetilde{M}]^{\mathbb{G}\cdot p}$ exists and then $$
[Y,\transp \widetilde{M}]^{\mathbb{G}\cdot p}
=
[D,\transp M]^{\mathbb{F}\cdot p}+\transp{\varphi}\centerdot [N,\transp M]^{\mathbb{F}\cdot p}.\ \ok
$$

Below we will solve separately the three structure conditions (\ref{structure-condition-c}), (\ref{structure-condition-da}), (\ref{structure-condition-di}) related to $D$. 

\

\section{Solution of the continuous structure condition}

The continuous structure condition (\ref{structure-condition-c}) has a quick answer. Let$$
\dcb
D = J'\centerdot W'+J''\centerdot W''+J'''\centerdot W''',\\
N = \zeta'\centerdot W'+\zeta''\centerdot W''+\zeta'''\centerdot W'''.
\dce
$$
be the martingale representation of $D,N$ in $\mathbb{F}$.

\bethe\label{answer-c}
The continuous structure condition (\ref{structure-condition-c}) related to $D$ has a solution if and only if the process $
\transp{\varphi}\centerdot[N^c,\transp N^c]\varphi
$
is a finite process on $[0,T]$. In this case, $K'_h=J'_h+(\transp{\varphi}\zeta')_h$ for $1\leq h\leq d$ are particular solutions. 
\ethe

\textbf{Proof.} 
Recall equation (\ref{structure-condition-c}), for $1\leq h\leq d$,	$$
\transp K'_h\centerdot[\widetilde{W}'_h, \widetilde{W}'_h]^{\mathbb{G}-p}
=[D,W'_h]^{\mathbb{F}-p}+\transp{\varphi}\centerdot [N, W'_h]^{\mathbb{F}-p}.
$$
With the continuity, the equation takes another form
$$
\dcb
\transp K'_h\centerdot[{W}'_h, {W}'_h]
&=&J'_h\centerdot [W'_h, W'_h]+(\transp{\varphi}\zeta')_h\centerdot [W'_h,W'_h].
\dce
$$
Hence, if the continuous structure condition related to $D$ has a solution $K'$, necessarily $K'_h=J'_h+(\transp{\varphi}\zeta')_h$ almost surely under the random measure $\mathsf{d}[W'_h,W'_h]$ for $1\leq h\leq d$, and $(\transp{\varphi}\zeta')_h$ is $\widetilde{W}'_h$-integrable ($J'_h$ being by assumption $W'_h$-integrable), i.e., $
\transp{\varphi}\centerdot[N^c,\transp N^c]\varphi
$
is a finite process on $[0,T]$.

Conversely, if the process $
\transp{\varphi}\centerdot[N^c,\transp N^c]\varphi
$
is finite on $[0,T]$, define $K'_h=J'_h+(\transp{\varphi}\zeta')_h, 1\leq h\leq d,$ on $[0,T]$. It forms a solution of the continuous structure condition related to $D$.
\ok

\

\section{Solution of the accessible structure condition}

\subsection{Equations at the stopping times $T_n$}

Recall (cf. subsection \ref{www}) $(T_n)_{1\leq n<\mathsf{N}^a}$ ($\mathsf{N}^a\leq \infty$) a sequence of strictly positive $(\mathbb{P},\mathbb{F})$ predictable stopping times such that $[T_n]\cap [T_{n'}]=\emptyset$ for $n\neq n'$ and $$
\{s\geq 0:\Delta_sW''\neq 0\}\subset\cup_{n\geq 1}[T_n].
$$  
A $\mathbb{G}$ predictable process $K''$ satisfies the equation (\ref{structure-condition-da}) if and only if 
$$
\dcb
&&\transp K''\centerdot[\widetilde{W}'',\transp \widetilde{W}'']^{\mathbb{G}\cdot p}
=
[D,\transp {W}'']^{\mathbb{F}\cdot p}+\transp{\varphi}\centerdot [N,\transp {W}'']^{\mathbb{F}\cdot p}\\

&=&
\transp J''\centerdot[W'',\transp {W}'']^{\mathbb{F}\cdot p}+\transp{\varphi}\zeta''\centerdot [W'',\transp {W}'']^{\mathbb{F}\cdot p}
=
(\transp J''+\transp{\varphi}\zeta'')\centerdot [W'',\transp {W}'']^{\mathbb{F}\cdot p}
\dce
$$
on $[0,T]$. Computing the jumps at $\mathbb{F}$ stopping times $T_n$, we can also say that $K''$ satisfies the equation (\ref{structure-condition-da}) if and only if, for every $1\leq n<\mathsf{N}^a$, on $\{T_n\leq T, T_n<\infty\}$, $K''_{T_n}$ satisfies the equation 
\begin{equation}\label{at-Tn}
\transp K''_{T_n} \mathbb{E}[\Delta_{T_n}\widetilde{W}'' \transp \Delta_{T_n}\widetilde{W}''|\mathcal{G}_{T_n-}]
=(\transp J''+\transp\overline{\varphi}\zeta'')_{T_n}\mathbb{E}[\Delta_{T_n}W'' \transp \Delta_{T_n}W''|\mathcal{F}_{T_n-}] 
\end{equation}
on $\{T_n\leq T, T_n<\infty\}$. (Recall that $\widetilde{W}''$ has no other jumps than that at the $T_n$'s.)

\

\subsection{Conditional expectation at predictable stopping times $T_n$}

For a fixed $1\leq n<\mathsf{N}^a$, let $(A_{n,0},A_{n,1},\ldots,A_{n,d})$ be the partition which satisfies the relation $
\mathcal{F}_{T_n}=\mathcal{F}_{T_n-}\vee\sigma(A_{n,0},A_{n,1},\ldots,A_{n,d}).
$
(cf. subsection \ref{www}) Denote $p_{n,h}=\mathbb{P}[A_{n,h}|\mathcal{F}_{T_n-}]$ and $\overline{p}_{n,h}=\mathbb{P}[A_{n,h}|\mathcal{G}_{T_n-}]$ for $0\leq h\leq d$.  

\bl\label{F-Tn}
We have
\ebe
\item[. ]
For any finite random variable $\xi\in\mathcal{F}_{T_n}$, the conditional expectation $\mathbb{E}[\xi|\mathcal{F}_{T_n-}]$ is well-defined. Let 
$$
\mathfrak{a}_n(\xi)_h=\ind_{\{p_{n,h}>0\}}\frac{1}{p_{n,h}}\mathbb{E}[\ind_{A_{n,h}}\xi|\mathcal{F}_{T_n-}], \ 0\leq h\leq d.
$$
We have $\xi=\sum_{h=0}^d\mathfrak{a}_n(\xi)_h\ind_{A_{n,h}}$.

\item[. ]
Denote the vector valued random variable $n_{n,h}:=\mathfrak{a}_n(\Delta_{T_n}N)_h, 0\leq h\leq d$. We have $$
(1+\transp{\varphi}_{T_n}n_{n,h})p_{n,h}
=
\mathbb{E}[\ind_{A_h}|\mathcal{G}_{T_n-}]
=\overline{p}_h,
$$
for $0\leq h\leq d$, on $\{T_n\leq T, T_n<\infty\}$.
  
\dbe
\el

\textbf{Proof.}
The first assertion of the lemma is the consequence of the relation $
\mathcal{F}_{T_n}=\mathcal{F}_{T_n-}\vee\sigma(A_{n,0},A_{n,1},\ldots,A_{n,d}).
$
The second assertion follows from a direct computation of $(\ind_{A_{n,h}}\ind_{[T_n,\infty)})^{\mathbb{G}\cdot p}$ using Lemma \ref{A-G-p}.  \ok

\bl\label{1+fn>0}
For $ 1\leq n<\mathsf{N}^a$ we have $
1+\transp{\varphi}_{T_n}\Delta_{T_n} N> 0.
$
\el

\textbf{Proof.} We compute, for $0\leq h\leq d$,$$
\dcb
0\leq\mathbb{E}[\ind_{\{1+\transp{\varphi}_{T_n}\Delta_{T_n} N \leq 0\}}\ind_{A_{n,h}}|\mathcal{G}_{T_n-}]

&=&\ind_{\{1+\transp{\varphi}_{T_n}n_{n,h} \leq 0\}}(1+\transp{\varphi}_{T_n}n_{n,h})p_h\leq 0
\dce
$$
It follows that $\ind_{\{1+\transp{\varphi}_{T_n}\Delta_{T_n} N \leq 0\}}\ind_{A_{n,h}}=0$ for $0\leq h\leq d$, i.e., $1+\transp{\varphi}_{T_n}\Delta_{T_n} N>0$. \ok

\bl\label{samekernels}
For $ 1\leq n<\mathsf{N}^a$, on $\{T_n\leq T, T_{n}<\infty\}$, suppose the set equality 
\begin{equation}\label{ptop}
\{0\leq h\leq d: p_{n,h}>0\}=\{0\leq h\leq d: \overline{p}_{n,h}>0\}.
\end{equation}
Then, on $\{T_n\leq T, T_{n}<\infty\}$,  the two matrix $\mathbb{E}[\Delta_{T_n}{W}'' \transp \Delta_{T_n}{W}''|\mathcal{F}_{T_n-}]$ and $\mathbb{E}[\Delta_{T_n}\widetilde{W}'' \transp \Delta_{T_n}\widetilde{W}''|\mathcal{G}_{T_n-}]$ have the same kernel space, which is the space $\mathscr{K}_n$ of $a\in\mathbb{R}\times \mathbb{R}^d$ such that $a$ (as function of its components) is constant on the set $\{0\leq h\leq d: p_{n,h}>0\}$, and the same image space $\mathscr{K}^\perp_n$. There exists a matrix valued $\mathcal{G}_{T_n-}$ measurable random variable $\mathsf{G}_n$ such that $\mathbb{E}[\Delta_{T_n}\widetilde{W}\transp \Delta_{T_n}\widetilde{W}|\mathcal{G}_{T_n-}]\mathsf{G}_n$ is the orthogonal projection $\mathsf{P}_n$ onto $\mathscr{K}^\perp_n$. 
\el

\brem
The martingale representation property in $\mathbb{F}$ implies that, for any $\mathbb{F}$ predictable stopping time $R$, $\mathcal{F}_{R-}=\mathcal{F}_{R}$ on $\{R<\infty, R\neq T_{n}, 1\leq n<\mathsf{N}^a\}$. Therefore, because of Lemma \ref{F-Tn}, the validity of the set equalities (\ref{ptop}) in Lemma \ref{samekernels} is equivalent to Assumption \ref{1+fin}.
\erem

\textbf{Proof.}
Write $\Delta_{T_n}W''_h=\Delta_{T_n}\widetilde{W}''_h+\Delta_{T_n}\Gamma(W''_h)$ and take the conditioning with respect to $\mathcal{G}_{T_n-}$ on $\{T_n\leq T, T_{n}<\infty\}$. We obtain$$
\Delta_{T_n}\Gamma(W''_h)
=
\mathbb{E}[\Delta_{T_n}W''_h|\mathcal{G}_{T_n-}]
=
\mathbb{E}[\frac{1}{2^n}(\ind_{A_{n,h}}-p_{n,h})|\mathcal{G}_{T_n-}]
=
\frac{1}{2^n}(\overline{p}_{n,h} - p_{n,h}),
$$
so that $\Delta_{T_n}\widetilde{W}''_h=\frac{1}{2^n}(\ind_{A_{n,h}}-\overline{p}_{n,h})$. With this in mind, as well as the set equality $\{0\leq h\leq d: p_{n,h}>0\}=\{0\leq h\leq d: \overline{p}_{n,h}>0\}$, we conclude the first assertion of the lemma. The second assertion can be concluded with \cite[Lemma 5.14]{song-mrp-drift}. \ok

\

\subsection{An intermediate result}

Note that $\widetilde{W}''$ is a $\mathbb{G}$ purely discontinous local martingale. In fact, $\widetilde{W}''=\ind_{\cup_n[T_n]}\centerdot\widetilde{W}''$.

\bl\label{thmprimeprime}
$K''$ solves the accessible structure condition (\ref{structure-condition-da}) related to $D$, if and only if Assumption \ref{1+fin} holds and, for every $1\leq n<\mathsf{N}^a$, on $\{T_n\leq T, T_n<\infty\}$,  
$$
\transp K''_{T_n} \mathsf{P}_n
=(\transp J''+\transp{\varphi}\zeta'')_{T_n}\mathbb{E}[\Delta_{T_n}W''\transp \Delta_{T_n}W''|\mathcal{F}_{T_n-}]\mathsf{G}_n	
$$
and the process $
\sum_{n=1}^{\mathsf{N}^a-}K''_{T_n}\ind_{[T_n]}
$
is $\widetilde{W}''$ integrable on $[0,T]$, i.e.
\begin{equation}\label{VKI}
\sqrt{
\sum_{n=1}^{\mathsf{N}^a-}\ind_{\{T_n\leq t\wedge T\}}
\left(
(\transp J''+\transp{\varphi}\zeta'')_{T_n}\mathbb{E}[\Delta_{T_n}W\transp \Delta_{T_n}W|\mathcal{F}_{T_n-}]\mathsf{G}_n
\Delta_{T_n}\widetilde{W}''
\right)^2
}
\end{equation}
is $(\mathbb{P},\mathbb{G})$ locally integrable.
\el

\textbf{Proof.}
\textbf{\texttt{"}If\texttt{"} part.} Suppose Assumption \ref{1+fin}. We note then that  the set equality (\ref{ptop}) in Lemma \ref{samekernels} holds on $\{T_n\leq T,  T_{n}<\infty\}$ for every $1\leq n<\mathsf{N}^a$. With Lemma \ref{samekernels}, the first formula of the lemma implies formula (\ref{at-Tn}), and hence equation (\ref{structure-condition-da}). Note that, on every $\{T_n\leq T,  T_{n}<\infty\}$, $\Delta_{T_n}\widetilde{W}''_h=\frac{1}{2^n}(\ind_{A_{n,h}}-\overline{p}_{n,h}), 0\leq h\leq d$. This implies that, for any $a\in\mathscr{K}_n$, $\transp a \Delta_{T_n}\widetilde{W}''=0$ (noting that $\ind_{A_{n,h}}=0$ if $\overline{p}_{n,h}=0$), i.e., $\Delta_{T_n}\widetilde{W}''\in\mathscr{K}_n^\perp$, which implies$$
\transp K''_{T_n}\Delta_{T_n}\widetilde{W}''=\transp K''_{T_n}\mathsf{P}_n\Delta_{T_n}\widetilde{W}''.
$$
Together with the first formula, it implies that $K''$ is $\widetilde{W}''$ integrable on $[0,T]$ if and only if expression (\ref{VKI}) is $(\mathbb{P},\mathbb{G})$ locally integrable. 

It remains to prove the inequality $\transp K''_{T_n}\Delta_{T_n}\widetilde{W}''<1$ on $\{T_n\leq T,  T_{n}<\infty\}$, if $K''_{{T_n}}$ is given by the first formula of the lemma. The first formula implies formula (\ref{at-Tn}) which implies
$$
\transp K''_{T_n} \mathbb{E}[\Delta_{T_n}\widetilde{W}''(\ind_{A_{n,h}}-\overline{p}_{n,h})|\mathcal{G}_{T_n-}]
=(\transp J''+\transp{\varphi}\zeta'')_{T_n}\mathbb{E}[\Delta_{T_n}W''(\ind_{A_{n,h}}-{p}_{n,h})|\mathcal{F}_{T_n-}], 
$$
for all $0\leq h\leq d$, on $\{T_n\leq T, T_n<\infty\}$, or equivalently,
$$
\dcb
&&\transp K''_{T_n} \mathbb{E}[\Delta_{T_n}\widetilde{W}''\ind_{A_{n,h}}|\mathcal{G}_{T_n-}]
=
(\transp J''+\transp{\varphi}\zeta'')_{T_n}\mathbb{E}[\Delta_{T_n}W''\ind_{A_{n,h}}|\mathcal{F}_{T_n-}],
\dce
$$
because $$
\dcb
\mathbb{E}[\Delta_{T_n}\widetilde{W}''|\mathcal{G}_{T_n-}]=0,\ \
\mathbb{E}[\Delta_{T_n}W''|\mathcal{F}_{T_n-}]=0.
\dce
$$
We denote by $\mathsf{a}_{n}$ the vector $(\ind_{A_{n,h}})_{0\leq h\leq d}$, by $p_n$ the vector $(p_{n,h})_{0\leq h\leq d}$, by $\overline{p}_n$ the vector $(\overline{p}_{n,h})_{0\leq h\leq d}$, to write $$
\Delta_{T_n}W''=\frac{1}{2^n}(\mathsf{a}_{n}-p_n)\ind_{[T_n,\infty)},\
\Delta_{T_n}\widetilde{W}''=\frac{1}{2^n}(\mathsf{a}_{n}-\overline{p}_n)\ind_{[T_n,\infty)}.
$$
For $0\leq h\leq d$, let $d_{n,h}:=\mathfrak{a}_n(\Delta_{T_n}D)_h$, $n_{n,h}:=\mathfrak{a}_n(\Delta_{T_n}N)_h$ to write
$$
\dcb
&&(\transp J''+\transp{\varphi}\zeta'')_{T_n}\mathbb{E}[\Delta_{T_n}W''\ind_{A_{n,h}}|\mathcal{F}_{T_n-}]\\
&=&
\mathbb{E}[\Delta_{T_n}D\ind_{A_{n,h}}|\mathcal{F}_{T_n-}]
+
\transp{\varphi}_{T_n}\mathbb{E}[\Delta_{T_n}N\ind_{A_{n,h}}|\mathcal{F}_{T_n-}]
=
(d_{n,h}+\transp{\varphi}_{T_n}n_{n,h})
p_{n,h}.
\dce
$$
Let $(\epsilon_0,\ldots,\epsilon_d)$ be the canonical basis in $\mathbb{R}\times\mathbb{R}^d$. By Lemma \ref{F-Tn}, 
$$
\transp K''_{T_n} \mathbb{E}[\Delta_{T_n}\widetilde{W}\ind_{A_{n,h}}|\mathcal{G}_{T_n-}]
=
\frac{1}{2^n}\transp K''_{T_n}(\epsilon_{h}-\overline{p}_n) (1+\transp \varphi_{T_n} n_{n,h}){p}_{n,h}.
$$
Putting them together we obtain the equality 
$$
\frac{1}{2^n}\transp K''_{T_n}(\epsilon_{h}-\overline{p}_n) (1+\transp \varphi_{T_n} n_{n,h}){p}_{n,h}
=
(d_{n,h}+\transp{\varphi}_{T_n}n_{n,h})
p_{n,h},
$$
on $\{T_n\leq T,  T_{n}<\infty\}$, so that, because $\Delta_{T_n}D<1$ and $1+\transp{\varphi}_{T_n}\Delta_{T_n} N> 0$ (cf. Lemma \ref{1+fn>0}),$$
\frac{1}{2^n}\transp K''_{T_n}(\epsilon_{h}-\overline{p}_n)
=
\frac{d_{n,h}+\transp{\varphi}_{T_n}n_{n,h}}{1+\transp{\varphi}_{T_n}n_{n,h}}<1,
$$
on $\{p_{n,h}>0\}\cap A_{n,h}=A_{n,h}$. This proves $\transp K''_{T_n}\Delta_{T_n}\widetilde{W}''<1$ on ${A_{n,h}}$ (for every $0\leq h\leq d$). 

\textbf{\texttt{"}Only if\texttt{"} part.} 
The accessible structure condition (\ref{structure-condition-da}) related to $D$ implies formula (\ref{at-Tn}). Hence, as the above computation shows,
$$
\frac{1}{2^n}\transp K''_{T_n}(\epsilon_{h}-\overline{p}_n) (1+\transp \varphi_{T_n} n_{n,h}){p}_{n,h}
-
(1+\transp{\varphi}_{T_n}n_{n,h})
p_{n,h},
=
(d_{n,h}-1)
p_{n,h},
$$
on $\{T_n\leq T,  T_{n}<\infty\}$. The fact $\Delta_{T_n}D-1<0$ implies $(d_{n,h}-1)p_{n,h}\leq 0$. Note equally that $$
0=\ind_{\{d_{n,h}-1=0\}}(d_{n,h}-1)\ind_{A_{n,h}}=\ind_{\{d_{n,h}-1=0\}}(\Delta_{T_n}D-1)\ind_{A_{n,h}}.
$$
This means that $\ind_{\{d_{n,h}-1=0\}}\ind_{A_{n,h}}=0$. Taking conditional expectation with respect to $\mathcal{F}_{T_n-}$ we have also $\ind_{\{d_{n,h}-1=0\}}p_{n,h}=0$, i.e., on $\{p_{n,h}>0\}$, $d_{n,h}-1<0$. This, combining with the above identity of $(d_{n,h}-1)p_{n,h}$, implies that $(1+\transp{\varphi}_{T_n}n_{n,h})>0$ whenever $p_{n,h}>0$. Hence, as a consequence of Lemma \ref{F-Tn}, the set equality (\ref{ptop}) in Lemma \ref{samekernels} holds, which implies, on the one hand, Assumption \ref{1+fin}, and on the other hand, the conclusions in Lemma \ref{samekernels}. Now, we can repeat the reasoning in the \texttt{"}If\texttt{"} part to achieve the proof of the lemma.  
\ok

\brem
It is interesting to note that the accessible structure condition (\ref{structure-condition-da}) related to $D$ implies that $d_{n,h}-1<0$ and $(1+\transp{\varphi}_{T_n}n_{n,h})>0$, whenever $p_{n,h}>0$.
\erem

\

\subsection{The integrability and conclusion}

The integrability condition (\ref{VKI}) in Lemma \ref{thmprimeprime} looks awful. Using Lemma \ref{F-Tn} we now give a pleasant interpretation of the formula (\ref{VKI})	. Recall that $\mathsf{G}_n$ denotes the $\mathcal{G}_{T_n-}$ measurable random matrix which inverses the matrix $\mathbb{E}[\Delta_{T_n}\widetilde{W}''\transp \Delta_{T_n}\widetilde{W}''|\mathcal{G}_{T_n-}]$ on the space $\mathscr{K}^\perp_n$.

\bl
Under the condition of Lemma \ref{samekernels}, we have
$$
\dcb
&&\sum_{n=1}^{\mathsf{N}^a-}\ind_{\{T_n\leq t\wedge T\}}
\left(
\left(\transp J''_{T_n}
+
\transp{\varphi}_{T_n} \zeta''_{T_n}\right)\mathbb{E}[\Delta_{T_n}W'' \transp\Delta_{T_n}W''|\mathcal{F}_{T_n-}]
\mathsf{G}_n\Delta_{T_n}\widetilde{W}''\right)^2\\

&=&\sum_{n=1}^{\mathsf{N}^a-}\ind_{\{T_n\leq t\wedge T\}}\frac{1}{(1+\transp{\varphi}_{T_n}\Delta_{T_n} N)^2}\left(
\Delta_{T_n}D
+
\transp {\varphi}_{T_n}\Delta_{T_n} N 
\right)^2.
\dce
$$
\el

\textbf{Proof.}
Note that, for $0\leq h\leq d$, $W''_h$ is a bounded process with finite variation. $W''_h$ is always a $\mathbb{G}$ special semimartingale whatever hypothesis$(H')$ is valid or not. We denote always by $\widetilde{W}''$ the $\mathbb{G}$ martingale part of $W''$.

Consider the space $\mathtt{E}=\{0,1,2,\ldots,d\}$. We fix $1\leq n<\mathsf{N}^a$ and endow $\mathtt{E}$ with two (random) probability measures $$
\dcb
\mathsf{m}[\{h\}]:=p_{n,h}\\
\overline{\mathsf{m}}[\{h\}]:=(1+\transp{\varphi}_{T_n}n_{n,h})p_{n,h},\ 0\leq h\leq d.
\dce
$$
Let $\boldsymbol{\epsilon}=(\epsilon_h)_{0\leq h\leq d}$ denote the canonical basis in $\mathbb{R}\times\mathbb{R}^{d}$. Let $d_{n,h}:=\mathfrak{a}_n(\Delta_{T_n}D)_h$, $n_{n,h}:=\mathfrak{a}_n(\Delta_{T_n}N)_h$, and $$
\dcb
&&w_{n,h}:=\mathfrak{a}_n(\Delta_{T_n}W'')_h
=
\frac{1}{2^n}\ind_{\{p_{n,h}>0\}}\left({\epsilon}_h - p_n\right).
\dce
$$ 
(Recall the notations $p_{n}, \overline{p}_{n}, \mathsf{a}_{n}$ in the proof of Lemma \ref{thmprimeprime}.) We define then the function $\mathsf{d}:=\sum_{h=0}^{d}d_{n,h}\ind_{\{h\}}$, $\mathsf{n}:=\sum_{h=0}^{d}n_{n,h}\ind_{\{h\}}$, and $$
\mathsf{w}:
=
\frac{1}{2^n}\sum_{h=0}^{d}{\epsilon}_h\ind_{\{h\}}
-
\frac{1}{2^n}\sum_{h=0}^{d}\ind_{\{h\}}p_n
=
\frac{1}{2^n}\sum_{h=0}^{d}{\epsilon}_h\ind_{\{h\}}
-
\frac{1}{2^n}p_n,
$$ 
on $\mathtt{E}$. As $\mathbb{E}_{{\mathsf{m}}}[\ind_{\{p_{n}=0\}}\ind_{\{h\}}]=0$, $\mathsf{w}$ is $\mathsf{m}-a.s.$ and $\overline{\mathsf{m}}-a.s.$ equal to
$$
\sum_{h=0}^{d}w_{n,h}\ind_{\{h\}}
=
\frac{1}{2^n}\sum_{h=0}^{d}\ind_{\{p_{n,h}>0\}}{\epsilon}_h\ind_{\{h\}}
-
\frac{1}{2^n}\sum_{h=0}^{d}\ind_{\{p_{n,h}>0\}}\ind_{\{h\}}p_n.
$$
$\mathsf{w}$ is a $(1+d)$-dimensional vector valued function. We denote by $\mathsf{w}_k$ its $k$th component, which is the real function$$
\mathsf{w}_k
=
\frac{1}{2^n}\ind_{\{k\}}
-
\frac{1}{2^n}p_{n,k},
$$ 
on $\mathtt{E}$. Be careful: do not confound it with $\mathsf{w}(h)$ which is a vector. We have $$
\mathbb{E}_{{\mathsf{m}}}[\mathsf{w}_k]=0,\
\mathbb{E}_{\overline{\mathsf{m}}}[\mathsf{w}_k]=\frac{1}{2^n}(\overline{p}_{n,k}-p_{n,k}),
$$
so that, on $\{T_{n}\leq T, T_n<\infty\}$,$$
\Delta_{T_n}\widetilde{W}''=\frac{1}{2^n}(\mathsf{a}_n-\overline{p}_n)
= \frac{1}{2^n}(\mathsf{a}_n-p_n - (\overline{p}_n - p_n)) 
=
\Delta_{T_n}{W}'' - \mathbb{E}_{\overline{\mathsf{m}}}[\mathsf{w}].
$$
For a function $F$ we compute.$$
\dcb
&&\mathbb{E}[F(\Delta_{T_n}{W}'')|\mathcal{G}_{T_n-}]
=
\sum_{h=0}^{d}F(w_{n,h})\mathbb{E}[\ind_{A_{n,h}}|\mathcal{G}_{T_n-}]

=
\mathbb{E}_{\overline{\mathsf{m}}}[F(\mathsf{w})].
\dce
$$
Similarly, 
$\mathbb{E}[F(\Delta_{T_n}{W}'')|\mathcal{F}_{T_n-}]
=\mathbb{E}_{{\mathsf{m}}}[F(\mathsf{w})].
$
Let $$
\dcb
x:&=&\left(J''_{T_n}
+
\transp\zeta''_{T_n}{\varphi}_{T_n} \right),\\
y:&=&\mathsf{G}_n \mathbb{E}[\Delta_{T_n}W'' \transp\Delta_{T_n}W''|\mathcal{F}_{T_n-}] x.
\dce
$$
Then, for all $z\in\mathbb{R}\times\mathbb{R}^{d}$, we write
\begin{equation}\label{zWWyx}
\transp z \mathbb{E}[\Delta_{T_n}\widetilde{W}''\transp \Delta_{T_n}\widetilde{W}''|\mathcal{G}_{T_n-}] y
=
\transp z \mathbb{E}[\Delta_{T_n}{W}''\transp \Delta_{T_n}{W}''|\mathcal{F}_{T_n-}] x.
\end{equation} 
As$$
\mathbb{E}[\Delta_{T_n}\widetilde{W}''\transp \Delta_{T_n}\widetilde{W}''|\mathcal{G}_{T_n-}]
=
\mathbb{E}[\Delta_{T_n}{W}''\transp \Delta_{T_n}\widetilde{W}''|\mathcal{G}_{T_n-}],
$$
the equation (\ref{zWWyx}) becomes
$$
\dcb
\mathbb{E}[(\transp z\Delta_{T_n}{W}'')(\transp \Delta_{T_n}\widetilde{W}'' y)|\mathcal{G}_{T_n-}]
=
\mathbb{E}[(\transp z \Delta_{T_n}{W}'')(\transp \Delta_{T_n}{W}''x)|\mathcal{F}_{T_n-}],
\dce
$$
or equivalently
$$
\dcb
\mathbb{E}_{\overline{\mathsf{m}}}[(\transp z\mathsf{w})\transp(\mathsf{w}-\mathbb{E}_{\overline{\mathsf{m}}}[\mathsf{w}]) y]
=
\mathbb{E}_{{\mathsf{m}}}[(\transp z \mathsf{w})(\transp \mathsf{w}x)].
\dce
$$
Set the function $\mathsf{q}:=\sum_{h=0}^{\mathsf{n}''}(1+\transp{\varphi}_{T_n}n_k)\ind_{\{h\}}$ on $\mathtt{E}$. We note that $\mathsf{q}=\frac{\mathsf{d}\overline{\mathsf{m}}}{\mathsf{d}{\mathsf{m}}}$. We have therefore
$$
\dcb
\mathbb{E}_{{\mathsf{m}}}[(\transp z\mathsf{w})\mathsf{q}(\transp\mathsf{w}-\mathbb{E}_{\overline{\mathsf{m}}}[\transp\mathsf{w}]) y]
=
\mathbb{E}_{{\mathsf{m}}}[(\transp z \mathsf{w})(\transp \mathsf{w}x)].
\dce
$$
For any vector $a=(a_0,a_1,\ldots,a_d)\in\mathbb{R}\times \mathbb{R}^d$ such that $\transp a p_n=0$, we have$$
\transp a \mathsf{w}
=
\sum_{k=0}^da_k\mathsf{w}_k
=
\sum_{k=0}^da_k(\frac{1}{2^n}\ind_{\{k\}}
-
\frac{1}{2^n}p_{n,k})
=
\frac{1}{2^n}\sum_{k=0}^da_k\ind_{\{k\}}.
$$
This means that the functions of the form $(\transp a \mathsf{w})$ generate the space of all functions on $\mathtt{E}$ with null ${\mathsf{m}}$-expectation. But,$$
\mathbb{E}_{{\mathsf{m}}}[\mathsf{q}(\transp\mathsf{w}-\mathbb{E}{\overline{\mathsf{m}}}[\transp\mathsf{w}]) y]
=
\mathbb{E}_{\overline{\mathsf{m}}}[(\transp\mathsf{w}-\mathbb{E}_{\overline{\mathsf{m}}}[\transp\mathsf{w}]) y]
=0.
$$
Hence,$$
\mathsf{q}(\transp\mathsf{w}-\mathbb{E}_{\overline{\mathsf{m}}}[\transp\mathsf{w}]) y
=
\transp \mathsf{w}x\  \mbox{ or  }\
(\transp\mathsf{w}-\mathbb{E}_{\overline{\mathsf{m}}}[\transp\mathsf{w}]) y
=
\frac{1}{\mathsf{q}}\transp \mathsf{w}x,
\
\mathsf{m}-a.s..
$$
(cf. Lemma \ref{samekernels}.) Regarding the values at every for $0\leq h\leq d$ with $p_{n,h}>0$,
$$
(\transp w_{n,h}-\mathbb{E}_{\overline{\mathsf{m}}}[\transp\mathsf{w}]) y
=
\frac{1}{(1+\transp{\varphi}_{T_n}n_k)}\transp w_{n,h} x,
$$
Consider now the process
$$
\dcb
\sum_{n=1}^{\mathsf{N}^a-}\ind_{\{T_n\leq t\wedge T\}}
\left(
\left(\transp J''_{T_n}
+
\transp{\varphi}_{T_n} \zeta''_{T_n}\right)\mathbb{E}[\Delta_{T_n}W'' \transp\Delta_{T_n}W''|\mathcal{F}_{T_n-}]
\mathsf{G}_n\Delta_{T_n}\widetilde{W}''\right)^2.
\dce
$$
We have$$
\dcb
&&
\left(\transp J''_{T_n}
+
\transp{\varphi}_{T_n} \zeta''_{T_n}\right)\mathbb{E}[\Delta_{T_n}W'' \transp\Delta_{T_n}W''|\mathcal{F}_{T_n-}]
\mathsf{G}_n\Delta_{T_n}\widetilde{W}''\\

&=&
\transp y \Delta_{T_n}\widetilde{W}''
=
\sum_{h=0}^{d}\transp\Delta_{T_n}\widetilde{W}'' y\ind_{A_{n,h}}
=
\sum_{h=0}^{d}\transp (w_{n,h}-\mathbb{E}_{\overline{\mathsf{m}}}[\mathsf{w}])y\ind_{A_{n,h}}\\

&=&
\sum_{h=0}^{d}\frac{1}{(1+\transp{\varphi}_{T_n}n_k)}\transp w_{n,h} x\ind_{A_{n,h}}.
\dce
$$
It implies
$$
\dcb
&&
\left(
\left(\transp J''_{T_n}
+
\transp{\varphi}_{T_n} \zeta''_{T_n}\right)\mathbb{E}[\Delta_{T_n}W'' \transp\Delta_{T_n}W''|\mathcal{F}_{T_n-}]
\mathsf{G}_n\Delta_{T_n}\widetilde{W}''
\right)^2\\

&=&
\left(
\sum_{h=0}^{d}\frac{1}{(1+\transp{\varphi}_{T_n}n_k)}\transp w_{n,h} x\ind_{A_{n,h}}
\right)^2

=
\sum_{h=0}^{d}\frac{1}{(1+\transp{\varphi}_{T_n}n_k)^2}\left(
\transp w_{n,h} x 
\right)^2\ind_{A_{n,h}}\\

&=&
\sum_{h=0}^{d}\frac{1}{(1+\transp{\varphi}_{T_n}\Delta_{T_n} N)^2}\left(
\transp x\ \Delta_{T_n}W'' 
\right)^2\ind_{A_{n,h}}\\

&=&
\frac{1}{(1+\transp{\varphi}_{T_n}\Delta_{T_n} N)^2}\left(
\transp \left(J''_{T_n}
+
\transp\zeta''_{T_n}{\varphi}_{T_n} \right) \Delta_{T_n}W'' 
\right)^2\\

&=&
\frac{1}{(1+\transp{\varphi}_{T_n}\Delta_{T_n} N)^2}\left(
\Delta_{T_n}D
+
\transp {\varphi}_{T_n}\Delta_{T_n} N 
\right)^2.
\dce
$$
The lemma is proved.
\ok

As a corollary of Lemma \ref{thmprimeprime} and the above one, we state

\bethe\label{Rthmprimeprime}
The accessible structure condition (\ref{structure-condition-da}) related to $D$ is satisfied, if and only if Assumption \ref{1+fin} holds and the process $$
\sqrt{\sum_{n=1}^{\mathsf{N}^a-}\ind_{\{T_n\leq t\wedge T\}}\frac{1}{(1+\transp{\varphi}_{T_n}\Delta_{T_n} N)^2}\left(
\Delta_{T_n}D
+
\transp {\varphi}_{T_n}\Delta_{T_n} N  
\right)^2}
$$
is $(\mathbb{P},\mathbb{G})$ locally integrable.
\ethe

\

\section{Solution of the totally inaccessible structure condition}

\subsection{Equations at the stopping times $S_{n}$}

Let $(S_n)_{1\leq n<\mathsf{N}^i}$ ($\mathsf{N}^i\leq \infty$) be a sequence of $(\mathbb{P},\mathbb{F})$ totally inaccessible stopping times such that $[S_n]\cap [S_{n'}]=\emptyset$ for $n\neq n'$ and $\{s\geq 0:\Delta_sW'''\neq 0\}\subset\cup_{n\geq 1}[S_n]$. A $\mathbb{G}$ predictable process $K'''$ satisfies the equation (\ref{structure-condition-di}) if and only if, for $1\leq h<\mathsf{n}'''$, 
\begin{equation}\label{di-GF}
\dcb
K'''_h\centerdot[\widetilde{W}'''_h,\widetilde{W}'''_h]^{\mathbb{G}\cdot p}

=
K'''_h\centerdot[{W}'''_h,{W}'''_h]^{\mathbb{G}\cdot p}
=
(J'''_h+\transp{\varphi}\zeta'''_h)\centerdot [W'''_h, {W}'''_h]^{\mathbb{F}\cdot p}
\dce
\end{equation}
on $[0,T]$.

\bl\label{KequSn}
A $\mathbb{G}$ predictable process $K'''$ satisfies the equation (\ref{structure-condition-di}) if and only if, for $1\leq n<\mathsf{N}^i$, for $1\leq h\leq \mathsf{n}'''$, $K'''_{h,S_n}$ satisfies the equation 
$$
\dcb
&&
(1+\transp{\varphi}_{S_n}\mathtt{R}_n)K'''_{h,S_n}\mathbb{E}[\ind_{\{\Delta_{S_n}W'''_h\neq 0\}}|\mathcal{G}_{S_n-}]
=
(J'''_h+\transp{\varphi}\zeta'''_h)_{S_n}\mathbb{E}[\ind_{\{\Delta_{S_n}W'''_h\neq 0\}}|\mathcal{F}_{S_n-}]
\dce
$$
on $\{S_n\leq T, S_n<\infty\}$, where $\mathsf{R}_n=\mathbb{E}[\Delta_{S_n}N|\mathcal{F}_{S_n-}]$.
\el

\textbf{Proof.}
Let $1\leq n<\mathsf{N}^i, 0\leq h\leq \mathsf{n}'''$. We define $g_{n,h}$ to be an $\mathbb{F}$ (resp. $\overline{g}_{n,h}$ a $\mathbb{G}$) predictable process such that $$
\dcb
g_{n,h,S_n}:=\mathbb{E}[(\Delta_{S_n} W'''_h)^2|\mathcal{F}_{S_n-}] \
\mbox{ resp. }
\overline{g}_{n,h,S_n}:=\mathbb{E}[(\Delta_{S_n} W'''_h)^2|\mathcal{G}_{S_n-}].
\dce
$$
Let $f$ denote the coefficient of the $(\mathbb{P},\mathbb{F})$ martingale $\Delta_{S_n} W'''_h\ind_{[S_n,\infty)}-(\Delta_{S_n} W'''_h\ind_{[S_n,\infty)})^{\mathbb{F}\cdot p}$ in its martingale representation with respect to $W$. By pathwise orthogonality, $f$ has all components null but $f_{h}$ and by Lemma \ref{single-jump} $f_h$ can be modified to be bounded. We have, on the one hand,$$
\dcb
&&f_{h} K'''_h\centerdot[{W}'''_h,{W}'''_h]^{\mathbb{G}\cdot p}
=
K'''_h\centerdot[f_{h} \centerdot{W}'''_h,{W}'''_h]^{\mathbb{G}\cdot p}
=
K'''_h\centerdot[\transp f \centerdot{W},{W}'''_h]^{\mathbb{G}\cdot p}\\
&=&
K'''_h\centerdot[\Delta_{S_n}{W}'''_h\ind_{[S_n,\infty)}-(\Delta_{S_n}{W}'''_h\ind_{[S_n,\infty)})^{\mathbb{F}\cdot p},{W}'''_h]^{\mathbb{G}\cdot p}\\

&=&
K'''_h\centerdot((\Delta_{S_n}{W}'''_h)^2\ind_{[S_h,\infty)})^{\mathbb{G}\cdot p}
=
K'''_h\overline{g}_{n,h}\centerdot\left(\ind_{[S_h,\infty)}\right)^{\mathbb{G}\cdot p}
\dce
$$
on $[0,T]$. On the other hand, 
$$
\dcb
&&
f_{h}(J'''_h+\transp{\varphi}\zeta'''_h)\centerdot [W'''_h, {W}'''_h]^{\mathbb{F}\cdot p}
=
(J'''_h+\transp{\varphi}\zeta'''_h)g_{n,h}\centerdot\left(\ind_{[S_n,\infty)}\right)^{\mathbb{F}\cdot p}
\dce
$$
on $[0,T]$. All put together, equation (\ref{di-GF}) implies
\begin{equation}\label{di-GF2}
K'''_h\overline{g}_{n,h}\centerdot\left(\ind_{[S_n,\infty)}\right)^{\mathbb{G}\cdot p}
=
(J'''_h+\transp{\varphi}\zeta'''_h)g_{n,h}\centerdot\left(\ind_{[S_n,\infty)}\right)^{\mathbb{F}\cdot p}
\end{equation}
on $[0,T]$ for any $1\leq n<\mathsf{N}^i, 0\leq h\leq \mathsf{n}'''$. For the converse conclusion, we note that$$
\dcb
&&[{W}'''_h,{W}'''_h]^{\mathbb{F}\cdot p}
=
\sum_{n=1}^{\mathsf{N}^a-}\left((\Delta_{S_n}{W}'''_h)^2\ind_{[S_n,\infty)}\right)^{\mathbb{F}\cdot p}
=
\sum_{n=1}^{\mathsf{N}^a-}g_{n,h}\centerdot\left(\ind_{[S_n,\infty)}\right)^{\mathbb{F}\cdot p},
\dce
$$
and
$$
\dcb
&&[{W}'''_h,{W}'''_h]^{\mathbb{G}\cdot p}
=
\sum_{n=1}^{\mathsf{N}^a-}\left((\Delta_{S_n}{W}'''_h)^2\ind_{[S_n,\infty)}\right)^{\mathbb{G}\cdot p}
=
\sum_{n=1}^{\mathsf{N}^a-}\overline{g}_{n,h}\centerdot\left(\ind_{[S_n,\infty)}\right)^{\mathbb{G}\cdot p}.
\dce
$$
Consequently, if the process $K'''$ satisfied all equation (\ref{di-GF2}) for $1\leq k<\mathsf{N}^i, 0\leq h\leq \mathsf{n}'''$, the process $K'''$ is $[{W}'''_h,{W}'''_h]^{\mathbb{G}\cdot p}$-integrable and equation (\ref{di-GF}) is satisfied.

Consider the equations (\ref{di-GF2}). Following Lemma \ref{A-G-p}, we compute on $[0,T]$ :
\begin{equation}\label{1+fNGFS}
\dcb
&&
(\ind_{[S_n,\infty)})^{\mathbb{G}\cdot p}
=
(\ind_{[S_n,\infty)})^{\mathbb{F}\cdot p}+\transp{\varphi}_{S_n}(\Delta_{S_n}N\ind_{[S_n,\infty)})^{\mathbb{F}\cdot p}
=
\left(1+\transp{\varphi}\mathtt{r}_n\right)\centerdot(\ind_{[S_n,\infty)})^{\mathbb{F}\cdot p},
\dce
\end{equation}
where $\mathtt{r}_n$ is a $\mathbb{F}$ predictable process such that $(\mathsf{r}_n)_{S_n}=\mathsf{R}_n$. Hence, on $[0,T]$, for any $\mathbb{G}$ predictable set $\mathtt{A}$ such that $\ind_{\mathtt{A}}(1+\transp{\varphi}\mathtt{r}_n)$ is bounded, equation (\ref{di-GF2}) implies
$$
\dcb
&&\ind_{\mathtt{A}}(1+\transp{\varphi}\mathtt{r}_n)
K'''_h\overline{g}_{n,h}\centerdot\left(\ind_{[S_k,\infty)}\right)^{\mathbb{G}\cdot p}
=
\ind_{\mathtt{A}}(1+\transp{\varphi}\mathtt{r}_n)(J'''_h+\transp{\varphi}\zeta'''_h)g_{n,h}\centerdot (\ind_{[S_n,\infty)})^{\mathbb{F}\cdot p}\\

&=&
\ind_{\mathtt{A}}(J'''_h+\transp{\varphi}\zeta'''_h)g_{n,h}\centerdot (\ind_{[S_n,\infty)})^{\mathbb{G}\cdot p}.
\dce
$$
This is equivalent to
$$
\dcb
&&(1+\transp{\varphi}\mathtt{R}_n)K'''_{h,S_n}\overline{g}_{h,n,S_n}
=
(J'''_h+\transp{\varphi}\zeta'''_h)_{S_n}g_{n,h,S_n},\ \mbox{ on $\{S_n\leq T, S_n<\infty\}$}.
\dce
$$
Let $\alpha'''$ be the $\mathbb{F}$ predictable process in subsection \ref{www} such that, for $1\leq h\leq \mathsf{n}'''$, $$
\Delta W'''_h = \alpha'''_h\ind_{\{\Delta W'''_h\neq 0\}}.
$$ 
We compute
$$
\dcb
&&(\alpha'''_{h,S_n})^2
(1+\transp{\varphi}_{S_n}\mathtt{R}_n)K'''_{h,S_n}\mathbb{E}[\ind_{\{\Delta_{S_n}W'''_h\neq 0\}}|\mathcal{G}_{S_n-}]
=
(1+\transp{\varphi}_{S_n}\mathtt{R}_n)K'''_{h,S_n}\mathbb{E}[(\Delta_{S_n}W'''_h)^2|\mathcal{G}_{S_n-}]\\
&=&
( J'''_h+\transp{\varphi}\zeta'''_h)_{S_n}\mathbb{E}[(\Delta_{S_n}W'''_h)^2|\mathcal{F}_{S_n-}]
=
(\alpha'''_{h,S_n})^2( J'''_h+\transp{\varphi}\zeta'''_h)_{S_n}\mathbb{E}[\ind_{\{\Delta_{S_n}W'''_h\neq 0\}}|\mathcal{F}_{S_n-}].  
\dce
$$
The lemma is proved, because on $\{\alpha'''_{h,S_n}=0\}$, $$
\mathbb{E}[\ind_{\{\Delta_{S_n}W'''_h\neq 0\}}|\mathcal{G}_{S_n-}]
=
\mathbb{E}[\ind_{\{\Delta_{S_n}W'''_h\neq 0\}}|\mathcal{F}_{S_n-}]
=0.\ \ok
$$

\

\subsection{Conditional expectations at stopping times $S_n$}\label{exp-at-jump}

For a fixed $1\leq n<\mathsf{N}^i$, applying the martingale representation property, applying \cite[Lemme(4.48)]{jacod} with the finite $\mathbb{F}$ predictable constraint condition in subsection \ref{www}, we see that, on $\{S_n<\infty\}$,$$
\mathcal{F}_{S_n}=\mathcal{F}_{S_n-}\vee\sigma(\Delta_{S_n}W''')
=\mathcal{F}_{S_n-}\vee\sigma(\{\Delta_{S_n}W'''_1\neq 0\},\ldots,\{\Delta_{S_n}W'''_{\mathsf{n}'''}\neq 0\}).
$$
(Note that $\{S_n<\infty\}\subset\{\Delta_{S_n}W'''\neq 0\}$.) We set$$
B_{n,k}:=\{\Delta_{S_n}W'''_k\neq 0\},\ q_{n,k}=\mathbb{P}[B_{n,k}|\mathcal{F}_{S_n-}],\ 
\overline{q}_{n,k}=\mathbb{P}[B_{n,k}|\mathcal{G}_{S_n-}],\ 1\leq k\leq {\mathsf{n}'''}.
$$
Note that, by our choice of $W'''$ in subsection \ref{www}, the $B_{n,k}$ form a partition on $\{S_n<\infty\}$ (cf. \cite{song-mrp-drift}).

\bl\label{F-Sn}
Let $1\leq n<\mathsf{N}^i$.   
\ebe
\item[. ]
For any finite random variable $\xi\in\mathcal{F}_{S_n}$, the conditional expectation $\mathbb{E}[\xi|\mathcal{F}_{S_n-}]$ is well-defined. Let $$
\mathfrak{i}_n(\xi)_k=\ind_{\{q_{n,k}>0\}}\frac{1}{q_{n,k}}\mathbb{E}[\ind_{B_{n,k}}\xi|\mathcal{F}_{S_n-}], \ 1\leq k\leq \mathsf{n}'''.
$$
We have $\xi=\sum_{h=1}^{\mathsf{n}'''}\mathfrak{i}_n(\xi)_h\ind_{B_h}$. 
 
\item[. ]
Denote the vector valued random variable $l_{n,k}:=\mathfrak{i}_n(\Delta_{S_n}N)_k, 1\leq k\leq {\mathsf{n}'''}$. We have $$
(1+\transp{\varphi}_{S_n}l_{n,k})q_{n,k}
=
(1+\transp{\varphi}_{S_n}\mathtt{R}_{n})\mathbb{E}[\ind_{B_{n,k}}|\mathcal{G}_{S_n-}]
=
(1+\transp{\varphi}_{S_n}\mathtt{R}_{n})\overline{q}_{n,k}
$$
for $1\leq k\leq {\mathsf{n}'''}$, on $\{S_n\leq T, S_n<\infty\}$, where $\mathtt{R}_n$ is the vector valued process introduced in Lemma \ref{KequSn}.

\item[. ]
We have $(1+\transp{\varphi}_{S_n}\mathtt{R}_{n})>0$ almost surely on $\{S_n\leq T, S_n<\infty\}$.

\dbe

\el

\textbf{Proof.}
The proof of the first assertion is straightforward. To prove the second assertion, we introduce $\mathbb{F}$ predictable processes $H,G$ such that $H_{S_n}=l_{n,k}, G_{S_n}=q_{n,k}$. We apply then Lemma \ref{A-G-p} to write, for any $1\leq k\leq d$, $$
\dcb
(\ind_{B_{n,k}}\ind_{[S_n,\infty)})^{\mathbb{G}\cdot p}

&=&(1+\transp{\varphi}H)G\centerdot(\ind_{[S_n,\infty)}
)^{\mathbb{F}\cdot p}
\dce
$$
on $[0,T]$. Apply again formula (\ref{1+fNGFS}) to obtain
$$
\dcb
((1+\transp{\varphi}_{S_n}\mathtt{R}_{n})\ind_{B_k}\ind_{[S_n,\infty)})^{\mathbb{G}\cdot p}

&=&((1+\transp{\varphi}_{S_n}l_{n,k})q_{n,k}\ind_{[S_n,\infty)}
)^{\mathbb{G}\cdot p}\\

\dce
$$
on $[0,T]$, which proves the second formula. Consider $(1+\transp{\varphi}_{S_n}\mathtt{R}_{n})$. We compute on $[0,T]$ $$
\dcb
0\leq (\ind_{\{1+\transp{\varphi}_{S_n}\mathtt{R}_{n}\leq 0\}}\ind_{[S_n,\infty)})^{\mathbb{G}\cdot p}

&=&\ind_{\{1+\transp{\varphi}\mathtt{r}_{n}\leq 0\}}(1+\transp{\varphi}\mathtt{r}_{n})\centerdot(\ind_{[S_n,\infty)})^{\mathbb{F}-p}\leq 0.
\dce
$$
This yields $
\mathbb{E}[\ind_{\{1+\transp{\varphi}_{S_n}\mathtt{R}_{n}\leq 0\}}\ind_{\{S_n\leq T,S_n<\infty\}}]
=0,
$
proving the third assertion.  \ok

\bl\label{1+fn2>0}
For $ 1\leq n<\mathsf{N}^i$ we have $
1+\transp{\varphi}_{S_n}\Delta_{S_n} N> 0
$
on $\{S_n\leq T, S_n<\infty\}$.
\el

\textbf{Proof.} We compute, for $1\leq h\leq \mathsf{n}'''$,$$
\dcb
&&0\leq\mathbb{E}[\ind_{\{1+\transp{\varphi}_{S_n}\Delta_{S_n} N \leq 0\}}\ind_{B_{n,h}}|\mathcal{G}_{S_n-}]

=\ind_{\{1+\transp{\varphi}_{S_n}l_{n,h} \leq 0\}}\overline{q}_{n,h}
\\

&=&\ind_{\{1+\transp{\varphi}_{S_n}l_{n,h} \leq 0\}}\frac{1+\transp{\varphi}_{S_n}l_{n,h}}{1+\transp{\varphi}_{S_n}\mathtt{R}_n}q_{n,h}\leq 0
\dce
$$
 on $\{S_n\leq T, S_n<\infty\}$. It follows that $\ind_{\{1+\transp{\varphi}_{S_n}\Delta_{S_n} N \leq 0\}}\ind_{B_{n,h}}=0$ for $1\leq h\leq \mathsf{n}'''$, i.e., $1+\transp{\varphi}_{S_n}\Delta_{S_n} N>0$. \ok

\

\subsection{Consequences on the totally inaccessible structure condition}

Note that $\widetilde{W}'''$ is a $\mathbb{G}$ purely discontinous local martingale. This is because $W'''$ is the limit in martingale space of $(\mathbb{P},\mathbb{F})$ local martingales with finite variation (cf. \cite[Theorem 6.22]{HWY}), and therefore the same is true for $\widetilde{W}'''$ by \cite[Proposition (2,2)]{Jeulin80}.

\bethe\label{KWG2}
The totally inaccessible structure condition (\ref{structure-condition-di}) related to $D$ is satisfied for all $1\leq h\leq \mathsf{n}'''$, if and only if the process
$$
\sqrt{\sum_{n=1}^{\mathsf{N}^i-}\ind_{\{S_n\leq t\wedge T\}}\frac{1}{(1+\transp{\varphi}_{S_n}\Delta_{S_n} N)^2}\left(
\Delta_{S_n}D
+
\transp {\varphi}_{S_n}\Delta_{S_n} N  
\right)^2}
$$
is $(\mathbb{P},\mathbb{G})$ locally integrable. In this case, a solution process $K'''$ is given by
\begin{equation}\label{sol-i}
\dcb
&&
K'''_{h}
=
\frac{J'''_h+\transp{\varphi}\zeta'''_h}{1+\transp{\varphi}\zeta'''_h\alpha'''_h}\ind_{\{1+\transp{\varphi}\zeta'''_h\alpha'''_h\neq 0\}},\ \mbox{ $\mathsf{d}[W'''_h,W'''_h]-a.s.$ on $[0,T]$, $1\leq h\leq \mathsf{n}'''$.} 
\dce
\end{equation}
\ethe

\textbf{Proof.} 
Suppose the integrability condition in the theorem and define $K'''_h$, $1\leq h\leq \mathsf{n}'''$, by (\ref{sol-i}). As$$
l_{n,h}\ind_{B_{n,h}}
=
\Delta_{S_n}N\ind_{B_{n,h}}
=
\zeta_{S_n}\Delta_{S_n}W\ind_{B_{n,h}}
=
\zeta'''_{h,S_n}\Delta_{S_n}W'''_h\ind_{B_{n,h}}
=
(\zeta'''_h\alpha'''_h)_{S_n}\ind_{B_{n,h}},
$$
the formula (\ref{sol-i}) implies, on $\{S_n\leq T, S_n<\infty\}$ for any $1\leq n<\mathsf{N}^i$, $$
\dcb
&&
K'''_{h,S_n}\ind_{B_{n,h}}
=
\frac{(J'''_h+\transp{\varphi}\zeta'''_h)_{S_n}}{(1+\transp{\varphi}_{S_n}l_{n,h})}\ind_{B_{n,h}}, 
\dce
$$
(Noting that the random measure $\mathsf{d}[W'''_h,W'''_h]$ charges the set $B_{n,h}\cap[S_n]$). Take the conditioning with respect to $\mathcal{G}_{S_n-}$ with help of Lemma \ref{F-Sn}.
\begin{equation}\label{lrqK}
\dcb
&&
(1+\transp{\varphi}_{S_n}l_{n,k})q_{n,k} K'''_{h,S_n}
=
(J'''_h+\transp{\varphi}\zeta'''_h)_{S_n}q_{n,k}, 
\dce
\end{equation}
i.e., the equations in Lemma \ref{KequSn} are satisfied. We prove hence that $K'''_h$ is a solution of equation (\ref{structure-condition-di}). 

We now prove that $K'''_h$ is $\widetilde{W}'''_h$-integrable on $[0,T]$. For any $1\leq n<\mathsf{N}^i$, on the set $B_{n,h}\cap\{S_n\leq T, S_n<\infty\}$, $J'''_h\Delta_{S_n}W'''_h=\Delta_{S_n}D$, $(\transp{\varphi}\zeta'''_h)_{S_n}\Delta_{S_n}W'''_h=\transp{\varphi}_{S_n}\Delta_{S_n} N$ so that 
\begin{equation}\label{KjumpW}
\dcb
K'''_{h,S_n}\Delta_{S_n}W'''_h
&=&\frac{(J'''_h+\transp{\varphi}\zeta'''_h)_{S_n}}{(1+\transp{\varphi}\zeta'''_h\alpha'''_h)_{S_n}}\Delta_{S_n}W'''_h\\
&=&
\frac{J'''_h\Delta_{S_n}W'''_h+(\transp{\varphi}\zeta'''_{h})_{S_n}\Delta_{S_n}W'''_h}{1+(\transp{\varphi}\zeta'''_h)_{S_n}\Delta_{S_n}W'''_h}
=
\frac{\Delta_{S_n}D+\transp{\varphi}\Delta_{S_n} N}{1+\transp{\varphi}\Delta_{S_n} N}\ind_{B_{n,h}}.
\dce
\end{equation}
This proves the $\widetilde{W}'''_h$-integrability of $K'''_h$ on $[0,T]$. (Recall $\Delta\widetilde{W}'''=\Delta{W}'''$.) 

We finally check if $K'''_h\Delta \widetilde{W}'''_h=K'''_h\Delta W'''_h<1$ on $[0,T]$. But, on $B_{n,h}\cap\{S_n\leq T, S_n<\infty\}$,
$$
\dcb
&&
K'''_{h,S_n}\Delta_{S_n}W'''_h
=
\frac{\Delta_{S_n}D+\transp{\varphi}\Delta_{S_n} N}{1+\transp{\varphi}\Delta_{S_n} N}
\

<
\frac{1+\transp{\varphi}\Delta_{S_n} N}{1+\transp{\varphi}\Delta_{S_n} N} =1,
\dce
$$
because $1+\transp{\varphi}\Delta_{S_n} N>0$ (cf. Lemma \ref{1+fn2>0}). The totally inaccessible structure condition related to $D$ is satisfied by $K'''_h$, $1\leq h\leq \mathsf{n}'''$.

Conversely, suppose that $K'''_h$ is a solution of the totally inaccessible structure condition (\ref{structure-condition-di}) related to $D$. The formula in Lemma \ref{KequSn} is satisfied so as formula (\ref{lrqK}) (with help of Lemma \ref{F-Sn}). Multiply formula (\ref{lrqK}) by $\ind_{B_{n,h}}$ on $\{S_n\leq T, S_n<\infty\}$, we obtain
$$
\dcb
&&
K'''_{h,S_n}\ind_{B_{n,k}}
=
\frac{(J'''_h+\transp{\varphi}\zeta'''_h)_{S_n}}{1+\transp{\varphi}_{S_n}l_{n,h}}\ind_{B_{n,k}}
=
\frac{(J'''_h+\transp{\varphi}\zeta'''_h)_{S_n}}{(1+\transp{\varphi}\zeta'''_h\alpha'''_h)_{S_n}}\ind_{B_{n,k}},
\dce
$$
for $1\leq n<\mathsf{N}^i$, $1\leq h\leq \mathsf{n}'''$. This implies in turn the formula (\ref{KjumpW}). The $\widetilde{W}'''_h$-integrability of $K'''_h$ on $[0,T]$ for $1\leq h\leq \mathsf{n}'''$ implies finally the integrability conditon of the theorem. \ok

\

\section{Final conclusions}

\textbf{Proof of Theorem \ref{main}}
It is direct consequence of Lemma \ref{piece-ensemble} together with Theorem \ref{answer-c}, Theorem \ref{Rthmprimeprime} and Theorem \ref{KWG2}. \ok

\

\textbf{Proof of Theorem \ref{fullviability}}
We have the martingale representation property with representation process $W=(W',W'',W''')$.

\textbf{Necessary part}

Suppose the full viability on $[0,T]$. Then, the drift operator satisfied the drift multiplier assumption on $[0,T]$, as it is proved in \cite[Theorem 5.5]{song-mrp-drift}. Let us prove that the group of the conditions (\ref{structure-condition-c}), (\ref{structure-condition-da}), (\ref{structure-condition-di}) related to $D$ in $\mathbb{G}$ are satisfied.

First of all, there exists a $\mathbb{G}$ structure connector $Y$ for $W'$ on $[0,T]$ (cf. Remark \ref{locallyboundedM}), i.e. (cf. formula (\ref{structure-condition}))
$$
[Y,\widetilde{W}']^{\mathbb{G}\cdot p}=\transp{\varphi}\centerdot [N,W']^{\mathbb{F}\cdot p}.
$$
By the continuity, we can replace $[Y,\transp\widetilde{W}']^{\mathbb{G}\cdot p}$ by $\transp K'\centerdot[\widetilde{W}',\transp \widetilde{W}']^{\mathbb{G}\cdot p}$ for some $\mathbb{G}$ predictable $\widetilde{W}'$-integrable process $K'$. We have proved the continuous structure condition (\ref{structure-condition-c}) related to $D=0$ in $\mathbb{G}$.

If $\mu$ denotes the jump measure of $\widetilde{W}''$, let us check that the conditions in \cite[Theorem 3.8]{song-mrp-drift} is satisfied. Recall $\boldsymbol{\epsilon}=(\epsilon_h)_{0\leq h\leq d}$ the canonical basis in $\mathbb{R}\times\mathbb{R}^{d}$, $\mathsf{a}_{n}$ the vector $(\ind_{A_{n,h}})_{0\leq h\leq d}$, $1\leq n<\mathsf{N}^a$. We have the identity on $\{T_n<\infty\}$
$$
\Delta_{T_n}\widetilde{W}''=\frac{1}{2^n}(\mathsf{a}_{n}-\overline{p}_n)
=
\frac{1}{2^n}\sum_{h=0}^d(\epsilon_h-\overline{p}_n)\ind_{A_{n,h}}.
$$
Conforming to the notation of \cite[Theorem 3.8]{song-mrp-drift}, we have $$
\alpha_{n,h}\frac{1}{2^n}(\epsilon_h-\overline{p}_n), \
\gamma_{n,i}= \frac{1}{2^n}(\epsilon_i-\overline{p}_{n,i}\mathbf{1}), \
0\leq i\leq d.
$$
where $\transp \mathbf{1}=(1,1,\ldots,1)\in \mathbb{R}\times\mathbb{R}^{d}$. Let $v$ be a vector in $\mathbb{R}\times\mathbb{R}^{d}$ orthogonal to the $\gamma_{n,i}$. Then$$
v_i=\overline{p}_{n,i}\transp \mathbf{1}v,\ 0\leq i\leq d,
$$
i.e., $v$ is proportional to $\overline{p}_n$. The vectors $\gamma_i$ together with $\overline{p}_n$ span whole $\mathbb{R}\times\mathbb{R}^{d}$, which is the condition of \cite[Theorem 3.8]{song-mrp-drift}.

By the full viability, there exists a $\mathbb{G}$ structure connector $Y$ for $W''$ on $[0,T]$, i.e. 
$$
[Y,\widetilde{W}'']^{\mathbb{G}\cdot p}=\transp{\varphi}\centerdot [N,W'']^{\mathbb{F}\cdot p}.
$$
Applying \cite[Lemma 3.1 and Theorem 3.8]{song-mrp-drift}, we can replace $[Y,\transp\widetilde{W}'']^{\mathbb{G}\cdot p}$ by $\transp K''\centerdot[\widetilde{W}'',\transp \widetilde{W}'']^{\mathbb{G}\cdot p}$ for some $\mathbb{G}$ predictable $\widetilde{W}''$-integrable process $K''$, proving the accessible structure condition (\ref{structure-condition-da}) related to $D=0$ in $\mathbb{G}$.

Notice that $\Delta W'''_h$, $1\leq h\leq \mathsf{n}'''$, satisfies clearly the condition in \cite[Theorem 3.9]{song-mrp-drift} (with $\mathsf{n}=1$). We can repeat the above reasoning to prove the totally inaccessible structure condition (\ref{structure-condition-di}) related to $D=0$ in $\mathbb{G}$.

Now apply Theorem \ref{answer-c}, Theorem \ref{Rthmprimeprime} and Theorem \ref{KWG2}. We prove the condition (\ref{fn-sur-fn}) in Theorem \ref{fullviability}.

\textbf{Sufficient part}

Conversely, suppose the drift multiplier assumption and the condition (\ref{fn-sur-fn}). Apply Lemma \ref{piece-ensemble} together with Theorem \ref{answer-c}, Theorem \ref{Rthmprimeprime} and Theorem \ref{KWG2}, then translate the conclusion with Theorem \ref{deflator-connector} in term of deflators. We conclude that any $\mathbb{F}$ local martingale has a $\mathbb{G}$ deflator. Let now $S$ be an $\mathbb{F}$ special semimartingale with a $\mathbb{F}$ deflator $D$. Then, $(D,DS)$ has a $\mathbb{G}$ deflator $Y$, i.e. $S$ has a $\mathbb{G}$ deflator $DY$. Apply again Theorem \ref{deflator-connector}. We conclude the proof.\ok

\

\textbf{Proof of Corollary \ref{commondeflator}}
We note that, in the case $D\equiv 0$, the proof of Lemma \ref{piece-ensemble} gives a common $(\mathbb{P},\mathbb{G})$ structure connector to all $(\mathbb{P},\mathbb{F})$ local martingales. Corollary \ref{commondeflator} is therefore the consequence of Theorem \ref{deflator-connector}. \ok

\

\

\end{document}